\documentclass[10pt,reqno]{amsart}
\usepackage{amssymb,mathtools,empheq}
\usepackage{tikz-cd}
\usepackage[margin=1in]{geometry}
\usepackage{float}
\usepackage{xcolor}
\usepackage{hyperref}
\colorlet{energyblue}{blue!75!black}
\colorlet{singularpink}{magenta!85!black}
\colorlet{remainderorange}{orange!85!black}

\makeatletter
\def\l@subsection{\@tocline{2}{0pt}{2pc}{5pc}{}}
\makeatother

\theoremstyle{definition}
\newtheorem{definition}{Definition}
\newtheorem{lemma}[definition]{Lemma}
\newtheorem{proposition}[definition]{Proposition}
\newtheorem{corollary}[definition]{Corollary}

\theoremstyle{remark}
\newtheorem{remark}[definition]{Remark}

\theoremstyle{plain}
\newtheorem{theorem}[definition]{Theorem}

\numberwithin{equation}{section}
\numberwithin{definition}{section}

\title{An FBSDE Approach to the Two-Dimensional Abelian--Higgs Model}
\author[Zexi He]{Zexi He$^{*}$}
\thanks{$^{*}$School of Mathematical Sciences, Fudan University, Shanghai 200433, China.
  \href{mailto:26110180014@m.fudan.edu.cn}{\texttt{26110180014@m.fudan.edu.cn}}.}

\begin{document}
	
	\begin{abstract}
We develop an iterated forward–backward stochastic differential equation (FBSDE) approach to the two-dimensional infinite-volume Abelian–Higgs model in Coulomb gauge. We prove tightness of the joint renormalized lattice gauge–Higgs laws, after continuum extension, as the lattice spacing tends to zero and the volume tends to infinity. We represent the lattice measure as the joint terminal law of the FBSDEs and establish uniform moment bounds for the corresponding terminal shifts. For the Higgs field, we combine a low–high frequency decomposition with covariant energy estimates to obtain bounds with sufficiently mild dependence on the gauge shift. For the gauge field, we establish high-frequency estimates for the effective force and combine them with quantitative diamagnetic control in the infrared.
	\end{abstract}
	
	\maketitle
	
	\begingroup
	\setlength{\emergencystretch}{2em}
	\tableofcontents
	\endgroup
	
\section{Introduction}
\label{sec:introduction}

	The Abelian–Higgs model, also known as the $U(1)$ Higgs  model, is formally described in Euclidean quantum field theory by the  measure  on
the space of distributions
$\mathcal{D}'(\mathbb{R}^2;\mathbb{R}^2\times\mathbb{C})$:
\begin{equation}
	\label{eq:intro-formal-abelian-higgs-measure}
	\mu^{\mathrm{YMH}}(dA,d\phi)
	\overset{\mathrm{formally}}{\propto}
	\exp\bigl(-S(A,\phi)\bigr)
	\mathcal{D}A\,\mathcal{D}\phi,
\end{equation}
where
\begin{equation}
	\label{eq:intro-continuum-abelian-higgs-action}
	S(A,\phi)
	=
	\int_{\mathbb{R}^2}
	\left(\frac12\lvert F_A(x)\rvert^2
	+\frac12\lvert D_A\phi(x)\rvert^2
	+\frac14\lvert\phi(x)\rvert^4\right)\,dx.
\end{equation}
Here, for a smooth gauge field
$A\colon\mathbb{R}^2\to\mathbb{R}^2$ and a smooth matter field
$\phi\colon\mathbb{R}^2\to\mathbb{C}$, the continuum covariant derivative
is defined by
\begin{equation}
	\label{eq:continuum-covariant-derivative}
	D_{A,j}\phi
	:=
	\partial_j\phi+iA_j\phi,
	\qquad
	D_A\phi
	:=
	(D_{A,1}\phi,D_{A,2}\phi),
\end{equation}
with $F_{A,jk}:=\partial_j A_k-\partial_k A_j$ and
$|F_A|^2:=\frac12\sum_{j,k=1}^2|F_{A,jk}|^2$.
For any smooth real-valued function $g$, the gauge transformation
\begin{equation}
	\label{eq:intro-continuum-gauge-transformation}
	(A,\phi)\longmapsto(A^g,\phi^g)
	:=\bigl(A+\nabla g,e^{-ig}\phi\bigr)
\end{equation}
preserves the action. More precisely,
\begin{equation}
	\label{eq:intro-continuum-gauge-invariance}
	D_{A^g}\phi^g=(D_A\phi)^g,\qquad
	F_{A^g}=F_A,\qquad S(A^g,\phi^g)=S(A,\phi).
\end{equation}
Several difficulties arise in giving a rigorous meaning to the formal expression~\eqref{eq:intro-formal-abelian-higgs-measure}. The formal Lebesgue measure $\mathcal{D}A\,\mathcal{D}\phi$ has no rigorous
infinite-dimensional counterpart. Under the natural Gaussian reference
measures, $A$ and $\phi$ are distribution-valued, so the nonlinear products
in the action are not classically defined, leading to ultraviolet difficulties.
In infinite volume, the spatial integral in
\eqref{eq:intro-continuum-abelian-higgs-action} may also diverge,
introducing infrared difficulties.
Gauge invariance also makes the gauge-field quadratic form degenerate, so integration over gauge-equivalent configurations introduces an infinite gauge-volume factor; see ~\eqref{eq:gauge-volume-factorization} and \cite[\S9.4]{Pes95}.

To regularize the model and remove gauge redundancy, we work on the
periodic lattice
$\Lambda_{M,\varepsilon}=\varepsilon\mathbb Z^2/(2\pi M\mathbb Z^2)$.
We regard $A_j(x)$ as the value of the gauge field on the oriented edge from
$x$ to $x+\varepsilon e_j$, while $\phi(x)$ is attached to the vertex $x$.
After imposing the Coulomb gauge $\operatorname{div}^\varepsilon A=0$ and
introducing the renormalization constants $a_{M,\varepsilon}$ and
$b_{M,\varepsilon}$, the resulting rigorous lattice probability measure is
given by
\begingroup
\renewcommand{\theHequation}{intro-renormalized-lattice-ymh}
\begin{equation}
	\tag{$\mathrm{YMH}_{M,\varepsilon}$}
	\label{eq:intro-renormalized-lattice-abelian-higgs-measure}
	\begin{aligned}
		&\mu_{M,\varepsilon}^{\mathrm{YMH}}(dA,d\phi)\propto\\
		&\quad\exp\Bigg(-\int_{\Lambda_{M,\varepsilon}}
		\Bigg(
		\frac12\lvert\nabla^\varepsilon A(x)\rvert^2
		+\frac12\lvert D_A^\varepsilon\phi(x)\rvert^2
		+\frac14\lvert\phi(x)\rvert^4
		-\left(\frac12a_{M,\varepsilon}+b_{M,\varepsilon}\right)
		\lvert\phi(x)\rvert^2
		\Bigg)\,dx\Bigg)
		\lambda_{M,\varepsilon}^{\mathrm{GF}}(dA)\,d\phi.
	\end{aligned}
\end{equation}
\endgroup
Here \(\lambda_{M,\varepsilon}^{\mathrm{GF}}\) is the induced Lebesgue measure on
the gauge-fixed slice. The measure \(\mu_{M,\varepsilon}^{\mathrm{YMH}}\) is a
Coulomb-gauge realization of the massless non-compact Gaussian formulation (G)
of \cite[Section~4]{BFS79b}. This gauge-fixing construction is inspired by
\cite{BCCKM26} and admits an equivalent Faddeev--Popov formulation (see
Remark~\ref{rem:helmholtz-fp-gaussian-averaging}). For the precise construction,
see Section~\ref{sec:lattice-abelian-higgs-and-variational-problems}, in
particular~\eqref{eq:rigorous-lattice-abelian-higgs-measure}
and~\eqref{eq:renormalized-lattice-abelian-higgs-measure}.

Fix $0<\nu\ll1$ and put
$\rho(x):=(1+|\nu x|^2)^{-2}$.
Let $\mathcal E^\varepsilon$ be the lattice-to-continuum extension
operator defined in \eqref{eq:lattice-sharp-extension-definition}.
Our main result, proved in Section~\ref{sec:main-result}, is that the lattice Abelian--Higgs measures form a tight family.
\begin{theorem}
	\label{thm:intro-subsequential-limits-lattice-ymh}
	The family
	$\bigl\{\mathcal E^\varepsilon_\#
	\mu_{M,\varepsilon}^{\mathrm{YMH}}\bigr\}_{
		\substack{
			M\geq1,\;0<\varepsilon<1\\
			2\pi M/\varepsilon\in2\mathbb N+1
	}}$
	is tight  on
	$
	\mathcal C^{-\nu}(\rho^{1+\nu};\mathbb R^2)
	\times\mathcal C^{-\nu}(\rho^{2+\nu};\mathbb C)$.
\end{theorem}
This theorem establishes subsequential compactness of the joint gauge–Higgs laws as both the ultraviolet and infrared cutoffs are removed.  The proof relies on the following FBSDE characterization and uniform estimates for the associated terminal shifts.
\subsection{Iterated FBSDE representation}
\label{subsec:introduction-iterated-fbsde}
The iterated FBSDE representation arises from the conditional
structure of the lattice measure. We write
$\mu_{M,\varepsilon}^{\mathrm{YMH}}$ as a Gibbs perturbation
of the product of the pure Yang--Mills measure and a Gaussian
free field (GFF).
Integrating out the Higgs field  defines
the  gauge marginal.
Applying the Bou\'e--Dupuis
variational formula \cite{BD98} to the conditional Higgs
law and the gauge-field marginal leads to two nested
stochastic control problems. Their optimality conditions
yield an inner FBSDE for the Higgs field and an outer FBSDE
for the gauge field.

We use the notation introduced in
Sections~\ref{sec:preliminaries}--\ref{sec:lattice-abelian-higgs-and-variational-problems}, in particular $0<\alpha\ll\nu\ll1$.
Here $P_{\chi_\varepsilon}^{M}$ is the pinning operator,
$\Pi_{\chi_\varepsilon}^{M}$ is the projection operator, and
$K^{M,\varepsilon}:=\nabla_\perp^\varepsilon(-\Delta^\varepsilon)^{-1}$.
The following theorem summarizes the FBSDE representation
and our uniform estimates.
\begin{theorem}[Lattice FBSDE representation and uniform  bounds]
	\label{thm:intro-lattice-fbsde-uniform-bounds}
	The measure
	$\mu_{M,\varepsilon}^{\mathrm{YMH}}$ admits the following FBSDE representation:
	\begin{empheq}[left=\empheqlbrace]{align}
		Z_t^{M,\varepsilon}
		&=-\bigl(\mathcal{L}_{B_1^{M,\varepsilon}}^\varepsilon+1\bigr)^{-1}
		\int_0^t\mathbb E_s^\phi\Bigl[
		\bigl(
		:\!\mathcal{L}_{C^{M,\varepsilon}}^\varepsilon\!:_{M,\varepsilon}
		-\mathcal{L}_{B_1^{M,\varepsilon}}^\varepsilon-1
		\bigr)\Phi^{M,\varepsilon}
		+:\!\lvert\Phi^{M,\varepsilon}\rvert^2
		\Phi^{M,\varepsilon}\!:_{M,\varepsilon}
		\Bigr]\,ds,
		\label{eq:intro-lattice-scalar-fbsde}\\[0.5ex]
		B_t^{M,\varepsilon}
		&=-P_{\chi_\varepsilon}^{M}K^{M,\varepsilon}\bigl(P_{\chi_\varepsilon}^{M}K^{M,\varepsilon}\bigr)^*
		\int_0^t\mathbb E_s^A\Bigl[
		\mathbb E^\phi\bigl[
		\Pi_{\chi_\varepsilon}^{M}
		\operatorname{Im}\bigl(
		\overline{\Phi^{M,\varepsilon}}D_{C^{M,\varepsilon}}^\varepsilon\Phi^{M,\varepsilon}
		\bigr)\bigr]
		\Bigr]\,ds,
		\label{eq:intro-lattice-gauge-fbsde}
	\end{empheq}
	for $t\in[0,1]$, with
	\begin{equation}
		\label{eq:intro-lattice-terminal-fields}
		\begin{aligned}
			C^{M,\varepsilon}
			:=
			A^{M,\varepsilon,0}+B_1^{M,\varepsilon}+H^{M,\varepsilon},
			\qquad
			\Phi^{M,\varepsilon}
			:=
			\phi_{B_1^{M,\varepsilon}}^{M,\varepsilon}+Z_1^{M,\varepsilon}.
		\end{aligned}
	\end{equation}
	Here $A^{M,\varepsilon,0}$ and $\phi_{B_1^{M,\varepsilon}}^{M,\varepsilon}$
	are the reference Gaussian and conditionally Gaussian fields, respectively,
	and $H^{M,\varepsilon}$ is a random vector taking values in $[0,M^{-1})^2$.
	The joint terminal law satisfies
	\begin{equation}
		\label{eq:intro-lattice-fbsde-law-representation}
		\mu_{M,\varepsilon}^{\mathrm{YMH}}
		=
		\operatorname{Law}_{\mathbb{P}}
		\bigl(C^{M,\varepsilon},\Phi^{M,\varepsilon}\bigr).
	\end{equation}
	For every $p\geq1$, the terminal shifts satisfy the uniform bounds
	\begin{equation}
		\label{eq:intro-lattice-terminal-shift-uniform-bounds}
			\sup_{M,\varepsilon}
			\mathbb{E}\left[
			\lVert B_1^{M,\varepsilon}\rVert_{H^{1,\varepsilon}(\rho)}^p
			\right]
			<\infty,
			\qquad
			\sup_{M,\varepsilon}
			\mathbb{E}\left[
			\lVert Z_1^{M,\varepsilon}\rVert_{
				H_{B_1^{M,\varepsilon}}^{1-\alpha,\varepsilon}(\rho^2)}^2
			\right]
			<\infty.
	\end{equation}
\end{theorem}
The FBSDE representation is proved in
Section~\ref{subsec:variational-representation-and-fbsdes}, while the uniform
estimates are established in
Sections~\ref{sec:estimate-higgs-field} and~\ref{sec:estimate-gauge-field}. An advantage of the iterated FBSDE formulation is that it allows us to estimate the two fields successively. Since the gauge and Higgs Wiener noises are independent, we  first estimate the Higgs FBSDE conditionally on  the gauge variables and then use these estimates in the gauge-field FBSDE. This contrasts with the stochastic quantization approach of \cite{BC24}, where the gauge and Higgs fields are treated as a coupled system.

\subsection{Main difficulties and proof strategy}
\label{subsec:introduction-main-difficulties-proof-strategy}
\mbox{}\par\nobreak

\indent\textbf{Estimates for the Higgs field.}
The main difficulty is to obtain conditional estimates for the Higgs
terminal shift with sufficiently mild dependence on the gauge shift.
Following \cite{GH19}, we use a low--high frequency decomposition,
controlling the high-frequency part by covariant Schauder estimates
and the low-frequency part by a covariant energy estimate exploiting
the coercivity of the cubic interaction.
A key stochastic input is the uniform control of the Wick powers
of the covariant Gaussian reference field, with only a small-power
dependence on the gauge shift.

\indent\textbf{Estimates for the gauge field.}
The infrared singularity of the inverse Laplacian obstructs direct volume-uniform weighted bounds for the low-frequency parts of both the unpinned Gaussian reference field and the gauge-field terminal shift. We address this by pinning the reference field and controlling the low-frequency shift through quantitative diamagnetic estimates adapted from \cite{CC24}.

For the high-frequency gauge shift, we use the FBSDE Schauder estimate. A key observation is its spatial stationarity, which allows us to use rapidly decaying spatial weights when estimating the projected effective force. The regular contributions to the effective force
are controlled by the Higgs-field estimates, while the main
irregular mixed term requires further use of the high-frequency
Higgs equation. Comparing the resulting covariant expression with
its free Gaussian counterpart produces more regular errors and reveals an exact cancellation of divergent counterterms.

A further key estimate is the uniform logarithmic bound for the
expectation of the covariant current:
\begin{equation}
	\label{eq:intro-mean-covariant-current-log-bound}
	\left\|
	\mathbb E^\phi\left[
	\operatorname{Im}\bigl(
	\overline{\phi_B^\varepsilon}
	D_B^\varepsilon\phi_B^\varepsilon
	\bigr)\right]
	\right\|_{H^{-1,\varepsilon}(\rho)}
	\lesssim
	\log\left(2+\|B\|_{H^{1,\varepsilon}(\rho)}\right).
\end{equation}
To establish this bound, we express this expectation
through a directional derivative of the trace of the covariant heat kernel.
The discrete Feynman--Kac formula then reduces the problem to
estimates for line integrals along random-walk bridges.

\subsection{Relation to earlier work}
\label{subsec:introduction-related-literature}

The rigorous study of the Abelian--Higgs model by lattice methods goes back to \cite{BFS79b,BFS80,BFS81}, followed by the renormalization-group work in \cite{Bal82a,Bal82b,Bal83} and \cite{Kin86a,Kin86b}. These works established ultraviolet stability and continuum and infinite-volume results for gauge-invariant observables in related lattice formulations. In contrast to these earlier results on gauge-invariant observables, we establish tightness of the joint gauge-fixed fields.

More recent approaches to gauge theories include stochastic quantization. \cite{S21} studied the Langevin dynamics of a two-dimensional Abelian gauge theory, while \cite{BC24} proved global well-posedness for the two-dimensional stochastic Abelian--Higgs equations. Related work on non-Abelian gauge theories includes the two-dimensional stochastic Yang--Mills dynamics \cite{CCHS22} and the three-dimensional stochastic Yang--Mills--Higgs equations \cite{CCHS24}.
More recently, \cite{BCHZ26} established global well-posedness of the two-dimensional non-Abelian stochastic Yang--Mills--Higgs dynamics on the space of gauge orbits and constructed its unique invariant measure.

For a compact lattice Abelian--Higgs model, \cite{CC24} obtained finite-volume ultraviolet estimates for the gauge-field marginal using a loop expansion and a quantitative diamagnetic inequality.

The variational approach to Euclidean field measures was developed for the $\Phi_3^4$ model in \cite{BG20}. The associated FBSDE method was subsequently used to study the infinite-volume $\Phi_2^4$ model in \cite{BG22} and the sine--Gordon model in \cite{GM24,Mey26}. The work most closely related to ours is \cite{BCCKM26}, which constructed the Abelian--Higgs measure on the two-dimensional torus using an iterated Bou\'e--Dupuis representation. This provides the starting point for our lattice FBSDE formulation.

\subsection{Discussion and future directions}
\label{subsec:introduction-discussion-future-directions}

The present work establishes tightness but not uniqueness of the
subsequential limits. A natural next step is to formulate a continuum
version of the iterated FBSDE and identify the limit points through
its terminal law. Establishing uniqueness at that level would be a
further step toward convergence of the full lattice family. Other directions include using the FBSDE approach to verify the Osterwalder–Schrader axioms for the limit points and extending the method to non-Abelian gauge theories.

	\medskip
\noindent\textbf{Use of AI-assisted tools.}
During the preparation of this manuscript, the author used OpenAI's ChatGPT (GPT-5.5, GPT-5.6, and GPT-6 models) and Codex to discuss mathematical arguments,
check formulas and proof drafts, and improve the exposition and \LaTeX{} typesetting.
The finite-volume tightness argument was developed by the author with AI assistance, including an AI-generated initial proof of the finite-volume version of Proposition~\ref{lem:mean-covariant-current-logarithmic-bound}.
AI played a major role in extending the results to infinite volume, building on the author's analytical strategy, particularly the low–high frequency decomposition and covariant energy estimates.
Specifically, AI proposed the use of a pinning operator and supplied initial proofs of Proposition~\ref{prop:gauge-low-shift-uniform-moments}, its supporting results in Appendix~\ref{app:positive-expansion-gaussian-domination}, and Lemma~\ref{lem:gauge-stationary-weight-transfer}.
AI also assisted extensively in extending finite-volume estimates to weighted settings in Sections~\ref{sec:estimate-higgs-field}, \ref{sec:estimate-gauge-field}, and~\ref{sec:stochastic-estimate} and Appendix~\ref{app:proof-gauge-irregular-Zhi1-unified-free-comparison}, and with related technical results on covariant Sobolev spaces in Section~\ref{sec:covariant-sobolev-spaces} and interpolation in Appendix~\ref{app:interpolation-theorems}.
All AI-assisted mathematical arguments were independently verified by the author, who takes full responsibility for the manuscript.

\medskip
\noindent\textbf{Acknowledgements.}
The author is grateful to Professor Shanjian Tang for his guidance and to Professor Weijun Xu and Da Li  for helpful discussions on variational methods and FBSDEs. He is also grateful to Dr. Abdulwahab Mohamed for sharing slides on the variational approach to the Abelian–Higgs model, which helped inspire the iterated FBSDE approach developed in this paper. The author also thanks Professors Ajay Chandra and Ilya Chevyrev, as well as Drs. Nikolay Barashkov, Andreas Koller, and Abdulwahab Mohamed, for inviting him to present this work at an online meeting and for their valuable comments and discussions.
This work was supported by the National Natural Science Foundation of China under Grant No. 11631004.

	\section{Preliminaries}
	\label{sec:preliminaries}
	\subsection{Parameters and weight}
	\label{subsec:parameters}
	Fix $0<\nu\ll1$, and set
	\begin{equation}
		\label{eq:regularity-parameters}
		\alpha_1:=\nu^{100},
		\qquad
		\alpha:=\nu^{1000},
		\qquad
		\kappa_2:=\alpha^4,
		\qquad
		\kappa_1:=\alpha^{10},
		\qquad
		\kappa:=\alpha^{20}.
	\end{equation}
	We further set
	\begin{equation}
		\label{eq:large-constant-K0}
		l_0:=100,
		\qquad
		K_0:=\alpha^{-2}.
	\end{equation}
	Then
	\begin{equation}
		\label{eq:parameter-hierarchy}
		0<\kappa\ll\kappa_1\ll\kappa_2
		\ll\frac{1}{K_0}\ll\alpha\ll\alpha_1\ll1.
	\end{equation}
	We define the weight
	\begin{equation}
		\label{eq:polynomial-weight-definition}
		\rho(x):=\langle \nu x\rangle^{-4}
		=\bigl(1+|\nu x|^2\bigr)^{-2}.
	\end{equation}
	Throughout the paper, unless stated otherwise, $X\lesssim Y$ for nonnegative quantities $X$ and $Y$
	means that $X\leq CY$ for some constant $C>0$ depending only on the
	parameters fixed in
	\eqref{eq:regularity-parameters}--\eqref{eq:parameter-hierarchy}.
	Any additional dependence is displayed as a subscript.
	We have
	\begin{equation}
		\label{eq:weight-basic-estimates}
		\bigl|(D^\varepsilon)^\alpha\rho^a(x)\bigr|
		\lesssim_\alpha\nu^{|\alpha|}\rho^a(x),
		\qquad
		\frac{\rho(x)}{\rho(y)}
		\lesssim\rho^{-1}(x-y),
	\end{equation}
	for all $x,y\in\Lambda_\varepsilon$, $|a|\le l_0$, and multi-indices
	$\alpha\in\mathbb{N}_0^2$.
	\begin{remark}
		The parameter choices above imply the relations used below:
		\begin{equation}
			\label{eq:kappa-over-alpha-relation}
			\frac{\kappa}{\alpha}=\alpha^9\kappa_1\ll\kappa_1,
			\qquad
			\frac{\kappa_1}{\alpha^2}=\alpha^4\kappa_2\ll\kappa_2,
			\qquad
			\frac{1}{\alpha}\ll K_0,
			\qquad
			\kappa K_0=\alpha^{18}\ll1.
		\end{equation}
		The dependence on $\nu$ in~\eqref{eq:polynomial-weight-definition} ensures
		\eqref{eq:weighted-covariant-coercivity} for sufficiently small $\nu$,
		and hence the energy estimate in Lemma~\ref{lem:higgs-low-energy-estimate}.
	
	\end{remark}

	\subsection{Discrete Fourier analysis}
	\label{subsec:discrete-fourier-helmholtz}
	Let $0<\varepsilon<1$ and $M\geq1$ satisfy
	$2\pi M/\varepsilon\in2\mathbb{N}+1$. We set
	\begin{equation}
	\begin{aligned}
		\Lambda_\varepsilon&=\varepsilon\mathbb{Z}^2,
		&\widehat\Lambda_\varepsilon
		&=\mathbb{R}^2\big/\bigl(\tfrac{2\pi}{\varepsilon}\mathbb{Z}^2\bigr),\\
		\Lambda_{M,\varepsilon}
		&=\varepsilon\mathbb{Z}^2/(2\pi M\mathbb{Z}^2),
		&\widehat\Lambda_{M,\varepsilon}
		&=M^{-1}\mathbb{Z}^2\big/\bigl(\tfrac{2\pi}{\varepsilon}\mathbb{Z}^2\bigr).
	\end{aligned}
	\end{equation}
	We set
	$Q_{M,\varepsilon}=M^{-1}\mathbb{Z}^2\cap
	[-\pi/\varepsilon,\pi/\varepsilon)^2$.
	Every function in $\mathbb{C}^{\Lambda_{M,\varepsilon}}$ is identified
	with its $2\pi M\mathbb{Z}^2$-periodic extension to $\Lambda_\varepsilon$.
	For lattice functions and $j=1,2$, define
	\begin{equation}
		\label{eq:lattice-integration-and-forward-difference}
		\begin{aligned}
			&\int_{\Lambda_\varepsilon}f(x)\,dx
			:=\varepsilon^2
			\sum_{x\in\Lambda_\varepsilon}f(x),\\
			&(\tau_jf)(x):=
			f(x+\varepsilon e_j),
			\qquad
			(\tau_j^*f)(x):=
			f(x-\varepsilon e_j),\\
			&D_j^\varepsilon=\partial_j^\varepsilon
			:=
			\varepsilon^{-1}(\tau_j-\operatorname{Id}),
			\qquad
			D_j^{\varepsilon *}
			=(\partial_j^\varepsilon)^*
			:=
			\varepsilon^{-1}(\tau_j^*-\operatorname{Id}).
		\end{aligned}
	\end{equation}
	The integral is understood whenever the sum is absolutely convergent;
	on $\Lambda_{M,\varepsilon}$, the same convention denotes the finite
	sum over one period. We write
	$\lvert\Lambda_{M,\varepsilon}\rvert:=
	\varepsilon^2\#\Lambda_{M,\varepsilon}$ for its volume.
 For square-summable functions, define
	the following inner products; the difference operators act on arbitrary
	lattice functions:
	\begin{equation}
		\label{eq:lattice-inner-product-and-adjoints}
		\begin{aligned}
			&\langle f,g\rangle_\varepsilon
			:=
			\operatorname{Re}\int_{\Lambda_\varepsilon}
			\overline{f(x)}g(x)\,dx,
			\qquad
			\langle f,g\rangle_{\mathbb{C},\varepsilon}
			:=
			\int_{\Lambda_\varepsilon}\overline{f(x)}g(x)\,dx,\\
			&\nabla^\varepsilon f
			:=D^\varepsilon f
			:=
			(D_1^\varepsilon f,D_2^\varepsilon f),
			\qquad
			(D^\varepsilon)^*A
			:=
			\sum_{j=1}^2D_j^{\varepsilon *}A_j,\\
			&\nabla_\perp^\varepsilon\eta
			:=
			\bigl(D_2^{\varepsilon *}\eta,
			-D_1^{\varepsilon *}\eta\bigr),\\
			&\operatorname{div}^\varepsilon A
			:=
			-(D^\varepsilon)^*A,
			\qquad
			-\Delta^\varepsilon
			:=
			\sum_{j=1}^2D_j^{\varepsilon *}D_j^\varepsilon
			=(D^\varepsilon)^*D^\varepsilon.
		\end{aligned}
	\end{equation}
	The finite-volume inner products
	$\langle\cdot,\cdot\rangle_{M,\varepsilon}$ and
	$\langle\cdot,\cdot\rangle_{\mathbb{C},M,\varepsilon}$ are defined
	by replacing $\Lambda_\varepsilon$ with $\Lambda_{M,\varepsilon}$.
	For $p\in[1,\infty)$, set
	$\lVert f\rVert_{L^{p,\varepsilon}(\Lambda_\varepsilon)}
	:=
	\bigl(\int_{\Lambda_\varepsilon}|f(x)|^p\,dx\bigr)^{1/p}$;
	for $p=\infty$, use the supremum norm. The norms on
	$\Lambda_{M,\varepsilon}$ are defined by summing over one period.
	We write
$
	\mathbb{R}_0^{\Lambda_{M,\varepsilon}}
	:=
	\left\{f\in\mathbb{R}^{\Lambda_{M,\varepsilon}}:
	\int_{\Lambda_{M,\varepsilon}}f(x)\,dx=0\right\}.
$
	
	Fix a nonnegative $\chi\in C_c^\infty((-1,1)^2)$ which is strictly
	positive on $[-1/2,1/2]^2$. For $x\in\Lambda_\varepsilon$, set
	\begin{equation}
	\chi_\varepsilon(x):=
	\frac{\chi(x)}{\varepsilon^2\sum_{y\in\Lambda_\varepsilon}\chi(y)},
	\qquad \langle1,\chi_\varepsilon\rangle_\varepsilon=1.
	\end{equation}
 For a real lattice field, set
	\begin{equation}
	\label{eq:lattice-pinning-projection}
	\begin{aligned}
	\overline A^{\chi_\varepsilon}
	:=
	\langle A,\chi_\varepsilon\rangle_\varepsilon
	=\varepsilon^2\sum_{x\in\Lambda_\varepsilon}A(x)\chi_\varepsilon(x),\quad
	P_{\chi_\varepsilon}A
	:=A-\overline A^{\chi_\varepsilon}.
	\end{aligned}
	\end{equation}
	 Whenever $\chi_\varepsilon$ is used as a function
	on $\Lambda_{M,\varepsilon}$, it denotes the periodization
	$\sum_{n\in\mathbb Z^2}\chi_\varepsilon(x+2\pi Mn)$; thus
	$\int_{\Lambda_{M,\varepsilon}}\chi_\varepsilon(x)\,dx=1$.
 Define
	\begin{equation}
	\begin{aligned}
	\bigl(\mathbb{R}^2\bigr)_{0,\mathrm{div}}^{\Lambda_{M,\varepsilon}}
	&:=
	\left\{A\in(\mathbb{R}^2)^{\Lambda_{M,\varepsilon}}:
	\operatorname{div}^\varepsilon A=0,\
	\int_{\Lambda_{M,\varepsilon}}A(x)\,dx=0\right\},\\
	\bigl(\mathbb{R}^2\bigr)_{\chi_\varepsilon,\mathrm{div}}^{\Lambda_{M,\varepsilon}}
	&:=
	\left\{A\in(\mathbb{R}^2)^{\Lambda_{M,\varepsilon}}:
	\operatorname{div}^\varepsilon A=0,\ \overline A^{\chi_\varepsilon}=0\right\},\\
	\mathcal A_{\mathrm{gf}}^{M,\varepsilon}
	&:=
	\bigl(\mathbb{R}^2\bigr)_{\chi_\varepsilon,\mathrm{div}}^{\Lambda_{M,\varepsilon}}
	+[0,M^{-1})^2.
	\end{aligned}
	\end{equation}
	
	For fixed $\varepsilon$, the discrete Schwartz space is
	\begin{equation}
	\mathcal S(\Lambda_\varepsilon)
	:=
	\left\{f:\Lambda_\varepsilon\to\mathbb C:
\sup_{x\in\Lambda_\varepsilon}
	\langle x\rangle^N|f(x)|<\infty
	\text{ for every }N\in\mathbb N_0\right\}.
	\end{equation}
Its continuous complex-linear dual $\mathcal S'(\Lambda_\varepsilon)$
	is identified with the polynomially growing sequences via the bilinear
	pairing $f(\varphi):=\int_{\Lambda_\varepsilon}f(x)\varphi(x)\,dx$:
	\begin{equation}
	\mathcal S'(\Lambda_\varepsilon)
	=
	\left\{f:\Lambda_\varepsilon\to\mathbb C:
	|f(x)|\leq C\langle x\rangle^N
	\text{ for some }C>0,\ N\in\mathbb N_0\right\}.
	\end{equation}
	For $\xi\in\widehat\Lambda_\varepsilon$ and
	$f\in\mathcal S(\Lambda_\varepsilon)$, set
	\begin{equation}
		\label{eq:lattice-fourier-transform}
		\mathrm e_{\xi}(x):=
		\exp(i\xi\cdot x),
		\qquad
		(\mathcal F_\varepsilon f)(\xi)
		:=
		\int_{\Lambda_\varepsilon}
		f(x)\overline{\mathrm e_{\xi}(x)}\,dx.
	\end{equation}
	For $f\in\mathbb C^{\Lambda_{M,\varepsilon}}$ and
	$k\in\widehat\Lambda_{M,\varepsilon}$, define
	\begin{equation}
		\label{eq:finite-volume-fourier-transform}
		\widehat f(k)\equiv(\mathcal F_{M,\varepsilon}f)(k)
		:=
		\int_{\Lambda_{M,\varepsilon}}
		f(x)\overline{\mathrm e_k(x)}\,dx.
	\end{equation}
	Throughout the paper, $\widehat f(k)$ denotes
	$(\mathcal F_{M,\varepsilon}f)(k)$ whenever $f$ is a
	$2\pi M$-periodic lattice field.
	The inverse Fourier transform of a smooth function on the dual torus is
	\begin{equation}
		\label{eq:lattice-inverse-fourier-transform}
		(\mathcal F_\varepsilon^{-1}g)(x)
		:=
		\frac1{(2\pi)^2}
		\int_{\widehat\Lambda_\varepsilon}
		g(\xi)\,\mathrm e_{\xi}(x)\,d\xi.
	\end{equation}
	For $f,g\in L^{2,\varepsilon}(\Lambda_\varepsilon)$,
	Parseval's identity reads
	\begin{equation}
	\langle f,g\rangle_{\mathbb C,\varepsilon}
	=\frac1{(2\pi)^2}\int_{\widehat\Lambda_\varepsilon}
	\overline{\widehat f(\xi)}\widehat g(\xi)
	\,d\xi.
	\end{equation}
	For $f\in\mathbb C^{\Lambda_{M,\varepsilon}}$,
	we have, as a distribution on the dual torus,
	\begin{equation}
	\mathcal F_\varepsilon f
	=\frac1{M^2}
	\sum_{k\in\widehat\Lambda_{M,\varepsilon}}
	\widehat f(k)\,\delta_k,
	\qquad
	f(x)=\frac1{\lvert\Lambda_{M,\varepsilon}\rvert}
	\sum_{k\in\widehat\Lambda_{M,\varepsilon}}
	\widehat f(k)e^{ik\cdot x}.
	\end{equation}
	In particular, its Fourier support is contained in
	$\widehat\Lambda_{M,\varepsilon}$, and the finite-volume Parseval identity is
	\begin{equation}
	\langle f,g\rangle_{\mathbb C,M,\varepsilon}
	=\frac1{\lvert\Lambda_{M,\varepsilon}\rvert}
	\sum_{k\in\widehat\Lambda_{M,\varepsilon}}
	\overline{\widehat f(k)}\widehat g(k).
	\end{equation}

	Define the smooth periodic symbols
	\begin{equation}
		\label{eq:lattice-frequency-multiplier}
		\xi_j^\varepsilon
		:=
		\frac{\exp(i\varepsilon\xi_j)-1}{i\varepsilon},
		\qquad
		(\mathcal F_\varepsilon\partial_j^\varepsilon f)(\xi)
		=i\xi_j^\varepsilon(\mathcal F_\varepsilon f)(\xi).
	\end{equation}
	We write $\xi^\varepsilon
	=(\xi_1^\varepsilon,\xi_2^\varepsilon)$.
	Since $\xi_j^\varepsilon$ is  complex, set
	\begin{equation}
		\label{eq:perpendicular-gradient-symbol}
		(\xi^\varepsilon)^\perp
		:=
		\left(-\overline{\xi_2^\varepsilon},
		\overline{\xi_1^\varepsilon}\right),
		\qquad
		(\mathcal F_\varepsilon\nabla_\perp^\varepsilon f)(\xi)
		=i(\xi^\varepsilon)^\perp(\mathcal F_\varepsilon f)(\xi).
	\end{equation}
 The same multiplier formulas hold for finite-volume
	Fourier coefficients at $k\in\widehat\Lambda_{M,\varepsilon}$.
	
	The finite-volume Fourier transform yields the following result.
	\begin{lemma}[Discrete Helmholtz decomposition]
		\label{lem:discrete-helmholtz-decomposition}
		For every $A:\Lambda_{M,\varepsilon}\to\mathbb R^2$, there exist unique
		$f,\eta\in\mathbb R_0^{\Lambda_{M,\varepsilon}}$ and a unique constant
		vector $h\in\mathbb R^2$ such that
		\begin{equation}
			\label{eq:lattice-helmholtz-decomposition}
			A=\nabla^\varepsilon f+\nabla_\perp^\varepsilon\eta+h.
		\end{equation}
	\end{lemma}
	
	\subsection{Covariant derivatives and Laplacians}
	\label{sec:covariant-derivatives-and-laplacians}
	\leavevmode\par\nopagebreak

	\smallskip
	\noindent\textit{Continuum operators.}
	For the continuum covariant derivative $D_A$ defined
	in~\eqref{eq:continuum-covariant-derivative}, its $L^2$-adjoint and the
	associated covariant Laplacian are given by
	\begin{equation}
		\label{eq:continuum-covariant-laplacian}
		D_{A,j}^*=-\partial_j-iA_j,
		\qquad
		\mathcal{L}_A
		:=
		D_A^*D_A
		=
		\sum_{j=1}^2D_{A,j}^*D_{A,j}.
	\end{equation}
	
	\smallskip
	\noindent\textit{Lattice operators.}
	Let $A,B\colon\Lambda_\varepsilon\to\mathbb{R}^2$ be lattice gauge fields.
	For $j=1,2$, define
	\begin{equation}
		\label{eq:prelim-covariant-derivative}
		\begin{aligned}
			D_{A,j}^\varepsilon\phi(x)
			&:=
			\frac{
				\exp\bigl(i\varepsilon A_j(x)\bigr)
				\phi(x+\varepsilon e_j)-\phi(x)
			}{\varepsilon},\\
			D_{A,j}^{\varepsilon *}\phi(x)
			&:=
			\frac{
				\exp\bigl(-i\varepsilon A_j(x-\varepsilon e_j)\bigr)
				\phi(x-\varepsilon e_j)-\phi(x)
			}{\varepsilon}.
		\end{aligned}
	\end{equation}
	For lattice functions $f,g\colon\Lambda_\varepsilon\to\mathbb{C}$ and
	$j=1,2$, the discrete covariant Leibniz rules read
	\begin{equation}
		\label{eq:covariant-leibniz-rules}
		\begin{alignedat}{2}
			D_{B,j}^\varepsilon(fg)
			&=(\tau_jf)D_{B,j}^\varepsilon g
			+(D_j^\varepsilon f)g
			&&=fD_{B,j}^\varepsilon g
			+(D_j^\varepsilon f)U_{B,j}^\varepsilon\tau_jg,\\
			D_{B,j}^{\varepsilon *}(fg)
			&=(\tau_j^*f)D_{B,j}^{\varepsilon *}g
			+(D_j^{\varepsilon *}f)g
			&&=fD_{B,j}^{\varepsilon *}g
			+(D_j^{\varepsilon *}f)\tau_j^*(U_{B,j}^{\varepsilon *}g).
		\end{alignedat}
	\end{equation}
	Here $U_{B,j}^\varepsilon$ and $U_{B,j}^{\varepsilon *}$ are defined
	in~\eqref{eq:rho-A-definition}.
	Set
	\begin{equation}
		\label{eq:lattice-covariant-laplacian}
		D_A^\varepsilon
		:=
		(D_{A,1}^\varepsilon,D_{A,2}^\varepsilon),
		\qquad
		\mathcal{L}_A^\varepsilon
		:=
		(D_A^\varepsilon)^*D_A^\varepsilon
		=
		\sum_{j=1}^2
		D_{A,j}^{\varepsilon *}D_{A,j}^\varepsilon.
	\end{equation}
	We have $\mathcal{L}_0^\varepsilon=-\Delta^\varepsilon$ and
	$\langle\mathcal{L}_A^\varepsilon\phi,\phi\rangle_\varepsilon
	=\lVert D_A^\varepsilon\phi\rVert_{L^{2,\varepsilon}}^2$.
	
	Define $r_A^\varepsilon$ componentwise by
	\begin{equation}
		\label{eq:r-A-definition}
		r_{A,j}^\varepsilon(x)
		:=
		\frac{\exp\bigl(i\varepsilon A_j(x)\bigr)-1}{i\varepsilon},
		\qquad
		r_A^\varepsilon(x)
		:=
		\bigl(
		r_{A,1}^\varepsilon(x),
		r_{A,2}^\varepsilon(x)
		\bigr).
	\end{equation}
	Identifying $r_{A,j}^\varepsilon$ with the corresponding multiplication
	operator, define
	\begin{equation}
		\label{eq:rho-A-definition}
		\begin{aligned}
			\rho_{A,j}^\varepsilon
			&:=
			r_{A,j}^\varepsilon\tau_j,
			&
			\rho_{A,j}^{\varepsilon *}
			&:=
			\tau_j^*\overline{r_{A,j}^\varepsilon},\\
			U_{B,j}^\varepsilon
			&:=
			\exp\bigl(i\varepsilon B_j(\cdot)\bigr),
			&
			U_{B,j}^{\varepsilon *}
			&:=
			U_{-B,j}^\varepsilon.
		\end{aligned}
	\end{equation}
	Moreover,
	\begin{equation}
		\label{eq:covariant-derivative-rho-decomposition}
		\begin{aligned}
			D_{A,j}^\varepsilon
			&=D_j^\varepsilon+i\rho_{A,j}^\varepsilon,
			&
			D_{A,j}^{\varepsilon *}
			&=D_j^{\varepsilon *}-i\rho_{A,j}^{\varepsilon *},\\
			D_{A+B,j}^\varepsilon
			&=D_{B,j}^\varepsilon
			+iU_{B,j}^\varepsilon\rho_{A,j}^\varepsilon,
			&
			D_{A+B,j}^{\varepsilon *}
			&=D_{B,j}^{\varepsilon *}
			-i\rho_{A,j}^{\varepsilon *}U_{B,j}^{\varepsilon *}.
		\end{aligned}
	\end{equation}
	We also use
	\begin{equation}
		\label{eq:covariant-derivative-translation-identities}
		U_{B,j}^{\varepsilon *}D_{B,j}^\varepsilon
		=-\tau_jD_{B,j}^{\varepsilon *},
		\qquad
		D_{B,j}^{\varepsilon *}U_{B,j}^\varepsilon
		=-D_{B,j}^\varepsilon\tau_j^*.
	\end{equation}
	For a vector field $F=(F_1,F_2)$, set
	\begin{equation}
		\label{eq:operator-family-gradient-divergence}
		\begin{aligned}
			\rho_A^\varepsilon\phi
			&:=
			\bigl(
			\rho_{A,1}^\varepsilon\phi,
			\rho_{A,2}^\varepsilon\phi
			\bigr),
			&
			\rho_A^{\varepsilon *}F
			&:=
			\sum_{j=1}^2\rho_{A,j}^{\varepsilon *}F_j,\\
			U_B^\varepsilon F
			&:=
			\bigl(
			U_{B,1}^\varepsilon F_1,
			U_{B,2}^\varepsilon F_2
			\bigr),
			&
			U_B^{\varepsilon *}
			&:=U_{-B}^\varepsilon.
		\end{aligned}
	\end{equation}
	\begin{lemma}[Expansion of the lattice covariant Laplacian]
		\label{lem:lattice-covariant-laplacian-expansion}
		For every
		$A,B\colon\Lambda_\varepsilon\to\mathbb{R}^2$, the following
		identities hold.
		\begin{itemize}
			\item[\textup{(i)}]
			\begin{equation}
				\label{eq:relative-covariant-laplacian-rho-expansion}
				\mathcal{L}_{A+B}^\varepsilon
				-\mathcal{L}_B^\varepsilon
				=
				-i\rho_A^{\varepsilon *}
				U_B^{\varepsilon *}D_B^\varepsilon
				+iD_B^{\varepsilon *}
				U_B^\varepsilon\rho_A^\varepsilon
				+\rho_A^{\varepsilon *}\rho_A^\varepsilon.
			\end{equation}
			In particular,
			\begin{equation}
				\label{eq:covariant-laplacian-rho-expansion}
				\mathcal{L}_A^\varepsilon-\mathcal{L}_0^\varepsilon
				=
				-i\rho_A^{\varepsilon *}D^\varepsilon
				+iD^{\varepsilon *}\rho_A^\varepsilon
				+\rho_A^{\varepsilon *}\rho_A^\varepsilon.
			\end{equation}
			
			\item[\textup{(ii)}]
			For every $f\colon\Lambda_\varepsilon\to\mathbb{C}$,
			\begin{equation}
				\label{eq:relative-covariant-laplacian-quadratic-form}
				\left\langle
				\bigl(
				\mathcal{L}_{A+B}^\varepsilon
				-\mathcal{L}_B^\varepsilon
				\bigr)f,f
				\right\rangle_\varepsilon
				=
				2\left\langle
				iD_B^{\varepsilon *}
				U_B^\varepsilon\rho_A^\varepsilon f,f
				\right\rangle_\varepsilon
				+
				\left\lVert
				\rho_A^\varepsilon f
				\right\rVert_{L^{2,\varepsilon}}^2.
			\end{equation}
			
			\item[\textup{(iii)}]
			For every $\phi\colon\Lambda_\varepsilon\to\mathbb{C}$,
			\begin{equation}
				\label{eq:relative-covariant-laplacian-divergence-rewriting}
				\bigl(
				\mathcal{L}_{A+B}^\varepsilon
				-\mathcal{L}_B^\varepsilon
				\bigr)\phi
				=
				iD_B^{\varepsilon *}
				\Bigl(
				\overline{r_A^\varepsilon}\,\phi
				+U_B^\varepsilon\rho_A^\varepsilon\phi
				\Bigr)
				+\rho_A^{\varepsilon *}\rho_A^\varepsilon\phi
				-i\bigl(
				D^{\varepsilon *}\overline{r_A^\varepsilon}
				\bigr)\phi.
			\end{equation}
		\end{itemize}
	\end{lemma}
	
	\begin{proof}
		Expanding
		\(
		\mathcal{L}_{A+B}^\varepsilon-\mathcal{L}_B^\varepsilon
		=D_{A+B}^{\varepsilon *}D_{A+B}^\varepsilon
		-D_B^{\varepsilon *}D_B^\varepsilon
		\)
		using~\eqref{eq:covariant-derivative-rho-decomposition} yields
		\eqref{eq:relative-covariant-laplacian-rho-expansion}. Taking $B=0$
		gives~\eqref{eq:covariant-laplacian-rho-expansion}.
		The quadratic-form identity
		\eqref{eq:relative-covariant-laplacian-quadratic-form} follows
		immediately because the first two terms in
		\eqref{eq:relative-covariant-laplacian-rho-expansion} are adjoints.
		
		By the discrete covariant Leibniz rule~\eqref{eq:covariant-leibniz-rules},
		\begin{equation}
			\begin{aligned}
				-i\rho_A^{\varepsilon *}
				U_B^{\varepsilon *}D_B^\varepsilon\phi(x)
				&=
				-i\sum_{j=1}^2
				\tau_j^*\Bigl(
				\overline{r_{A,j}^\varepsilon}\,
				U_{B,j}^{\varepsilon *}
				D_{B,j}^\varepsilon\phi
				\Bigr)(x)\\
				&\overset{\eqref{eq:covariant-derivative-translation-identities}}{=}
				i\sum_{j=1}^2
				\bigl(
				\tau_j^*\overline{r_{A,j}^\varepsilon}
				\bigr)(x)
				D_{B,j}^{\varepsilon *}\phi(x)\\
				&=
				i\sum_{j=1}^2
				D_{B,j}^{\varepsilon *}
				\bigl(\overline{r_{A,j}^\varepsilon}\,\phi\bigr)(x)
				-i\sum_{j=1}^2
				\bigl(
				D_j^{\varepsilon *}
				\overline{r_{A,j}^\varepsilon}
				\bigr)(x)\phi(x).
			\end{aligned}
		\end{equation}
		Substituting this identity into
		\eqref{eq:relative-covariant-laplacian-rho-expansion} gives
		\eqref{eq:relative-covariant-laplacian-divergence-rewriting}.
	\end{proof}
	
	\begin{remark}
		For smooth gauge fields
		$A,B\colon\mathbb{T}^2\to\mathbb{R}^2$, the continuum analogue of
		\eqref{eq:relative-covariant-laplacian-divergence-rewriting} is
		\begin{equation}
			\label{eq:continuum-relative-covariant-laplacian-components}
			\bigl(\mathcal{L}_{A+B}-\mathcal{L}_B\bigr)\phi
			=2iD_B^*(A\phi)
			+i(\nabla\mathbin{\cdot}A)\phi
			+\lvert A\rvert^2\phi.
		\end{equation}
		See also \cite[Lemma~2.2]{BC24}.
	\end{remark}
	
	\begin{lemma}[Covariant Laplacian commutator]
		\label{lem:covariant-laplacian-weight-commutator}
		Let $w,u\in\mathbb{C}^{\Lambda_\varepsilon}$.
		Then, pointwise on $\Lambda_\varepsilon$,
		\begin{equation}
			\label{eq:covariant-laplacian-weight-commutator}
			[\mathcal{L}_B^\varepsilon,w]u
			=-(\Delta^\varepsilon w)u
			-\sum_{j=1}^2\bigl[
			(D_j^\varepsilon w)D_{B,j}^\varepsilon u
			+(D_j^{\varepsilon *}w)D_{B,j}^{\varepsilon *}u
			\bigr].
		\end{equation}
		In particular, for $|a|\le l_0$,
		\begin{equation}
			\label{eq:covariant-laplacian-weight-commutator-bound}
			\bigl|[\mathcal{L}_B^\varepsilon,\rho^a]u\bigr|
			\lesssim\rho^a\Bigl(|u|
			+\sum_{j=1}^2\bigl(|D_{B,j}^\varepsilon u|
			+|D_{B,j}^{\varepsilon *}u|\bigr)\Bigr).
		\end{equation}
	\end{lemma}
	\begin{proof}
		For $j=1,2$, the Leibniz rule~\eqref{eq:covariant-leibniz-rules} gives
		\begin{equation}
			\label{eq:covariant-laplacian-weight-product-rule}
			D_{B,j}^\varepsilon(wu)
			=(\tau_jw)D_{B,j}^\varepsilon u+(D_j^\varepsilon w)u.
		\end{equation}
		Applying $D_{B,j}^{\varepsilon *}$ and using the same rule yields
		\begin{equation}
			\label{eq:covariant-laplacian-weight-product-expansion}
			\begin{aligned}
				D_{B,j}^{\varepsilon *}D_{B,j}^\varepsilon(wu)
				=wD_{B,j}^{\varepsilon *}D_{B,j}^\varepsilon u
				+(D_j^{\varepsilon *}\tau_jw)D_{B,j}^\varepsilon u+(\tau_j^*D_j^\varepsilon w)D_{B,j}^{\varepsilon *}u
				+(D_j^{\varepsilon *}D_j^\varepsilon w)u.
			\end{aligned}
		\end{equation}
		Since $D_j^{\varepsilon *}\tau_j=-D_j^\varepsilon$ and
		$\tau_j^*D_j^\varepsilon=-D_j^{\varepsilon *}$, summing over $j=1,2$
		yields~\eqref{eq:covariant-laplacian-weight-commutator}.
		The bound follows from~\eqref{eq:weight-basic-estimates}.
	\end{proof}

	\subsection{Discrete Sobolev and Besov spaces}
	\label{sec:discrete-sobolev-besov}
	
	As in \cite[Proposition~2.10]{BDC11} and \cite[Section~2.3]{BC23},
	fix a smooth, even, radially non-increasing cutoff
	$\varphi\in C_c^\infty(\mathbb{R}^2;[0,1])$ that equals $1$ on
	$[-1,1]^2$ and is supported in $[-9/8,9/8]^2$.
	Set $\varphi_{-1}:=\varphi$ and
	$\varphi_j:=
	\varphi(2^{-j-1}\mathord\cdot)-\varphi(2^{-j}\mathord\cdot)$ for $j\ge0$.
	In particular,
	$\operatorname{supp}\varphi_i\cap\operatorname{supp}\varphi_j
	=\varnothing$ whenever $\lvert i-j\rvert>1$.
	Set
$
		\label{eq:lattice-besov-frequency-cutoff}
		J^\varepsilon
		:=
		\inf\left\{
		j\ge-1:
		\operatorname{supp}\varphi_j
		\not\subset
		\left[-\frac{\pi}{\varepsilon},\frac{\pi}{\varepsilon}\right)^2
		\right\},
$
	and define the truncated multipliers
	\begin{equation}
		\label{eq:truncated-dyadic-multipliers}
		\varphi_j^\varepsilon
		:=
		\begin{cases}
			\varphi_j,
			& -1\le j\le J^\varepsilon-1,\\
			1-\displaystyle\sum_{-1\le\ell\le J^\varepsilon-1}\varphi_\ell,
			& j=J^\varepsilon.
		\end{cases}
	\end{equation}
	The associated Littlewood--Paley blocks are
	$\Delta_j^\varepsilon f
	=\mathcal{F}_\varepsilon^{-1}
	\bigl(\varphi_j^\varepsilon\mathcal F_\varepsilon f\bigr)$.
	With the convention
	$\Delta_j^\varepsilon=0$ for $j>J^\varepsilon$, we define the weighted
	discrete Besov norms by
	\begin{equation}
		\label{eq:discrete-besov-norm}
		\lVert f\rVert_{B_{p,q}^{s,\varepsilon}(\rho)}
		:=
		\left(
		\sum_{j\ge-1}
		2^{s jq}
		\lVert\rho\Delta_j^\varepsilon f\rVert_{L^{p,\varepsilon}(\Lambda_\varepsilon)}^q
		\right)^{1/q}
	\end{equation}
	for $s\in\mathbb{R}$, $p,q\in[1,\infty]$, with the usual
	modifications when $q=\infty$.
	The unweighted norm corresponds to $\rho=1$.
	For $x,y\in\mathbb R^2$, define
	\begin{equation}
		\label{eq:periodized-polynomial-weight}
		w_x^M(y):=
		\sum_{n\in\mathbb Z^2}\rho(x-y+2\pi Mn),
		\qquad
		\rho_x(y):=\rho(x-y).
	\end{equation}
	For $a,s\in\mathbb R$ and $p,q\in[1,\infty]$, set
$
		\label{eq:periodized-weighted-Besov-norm}
		\|f\|_{B_{p,q}^{s}((w_x^M)^a;\Lambda_{M,\varepsilon})}
		:=
		\left(\sum_{j\ge-1}2^{s jq}
		\|(w_x^M)^a\Delta_j^\varepsilon f
		\|_{L^{p,\varepsilon}(\Lambda_{M,\varepsilon})}^q\right)^{1/q}.
$
	
	For $-1\leq j\leq J^\varepsilon$, write
	$\eta_j^\varepsilon:=
	\mathcal{F}_\varepsilon^{-1}(\varphi_j^\varepsilon)$.
	By \cite[Lemma~2.20 and Remark~2.21]{MP19}
	(see also Lemma~\ref{lem:weighted-periodic-LP-kernel} below),
	\begin{equation}
		\label{eq:app-LP-kernel-Lq-bound}
		\lVert\eta_j^\varepsilon\rVert_{L^{q,\varepsilon}}
		\lesssim
		2^{2j(1-1/q)},
		\qquad q\in[1,\infty].
	\end{equation}
	For $2\pi M$-periodic fields, write
	$\eta_j^{M,\varepsilon}(x):=
	\sum_{n\in\mathbb Z^2}\eta_j^\varepsilon(x+2\pi Mn)$.

	We also define the discrete Sobolev norms associated with
	$\mathcal{L}_0^\varepsilon=-\Delta^\varepsilon$ by
	\begin{equation}
		\label{eq:discrete-sobolev-norm}
		\begin{aligned}
		\lVert f\rVert_{W_p^{s,\varepsilon}}
		:=
		\bigl\lVert
		(\mathcal{L}_0^\varepsilon+1)^{s/2}f
		\bigr\rVert_{L^{p,\varepsilon}},\quad
		\lVert f\rVert_{W_p^{s,\varepsilon}(\rho)}
		:=
		\lVert f\rho\rVert_{W_p^{s,\varepsilon}},
		\end{aligned}
	\end{equation}
	and write $H^{s,\varepsilon}(\rho):=W_2^{s,\varepsilon}(\rho)$ and
	$\mathcal{C}^{s,\varepsilon}(\rho):=B_{\infty,\infty}^{s,\varepsilon}(\rho)$.
	Uniformly in $\varepsilon$,
	$\lVert f\rVert_{H^{s,\varepsilon}}
	\simeq\lVert f\rVert_{B_{2,2}^{s,\varepsilon}}$.
	Moreover, by \cite[Lemma~A.1]{GH21},
	$\lVert f\rVert_{B_{p,q}^{s,\varepsilon}(\rho)}
	\simeq\lVert\rho f\rVert_{B_{p,q}^{s,\varepsilon}}$
	for every polynomial weight $\rho$, $s\in\mathbb{R}$, and
	$p,q\in[1,\infty]$.
	
	Define the discrete paraproducts
	\begin{equation}
		\label{eq:discrete-paraproduct-definition}
		f\prec g
		:=
		\sum_{i<j-1}
		\Delta_i^\varepsilon f\,\Delta_j^\varepsilon g,
		\qquad
		f\circ g
		:=
		\sum_{\lvert i-j\rvert\le 1}
		\Delta_i^\varepsilon f\,\Delta_j^\varepsilon g.
	\end{equation}
We also write
$
f\preccurlyeq g:=f\prec g+f\circ g,
f\succ g:=g\prec f.
$
	\begin{lemma}[Weighted discrete paraproduct estimates]
		\label{lem:discrete-paraproduct-estimates}
		Let $\rho_1,\rho_2$ be polynomial weights and let
		$p,p_1,p_2\in[1,\infty]$ satisfy
		$\dfrac{1}{p_1}+\dfrac{1}{p_2}=\dfrac{1}{p}$.
		The constants below may depend on the weights and the displayed
		regularity and integrability exponents.
		\smallskip
		
		\noindent\textup{(i)}
		For any $s_2\in\mathbb{R}$, and, in the second estimate, any
		$s_1<0$,
		\begin{equation}
			\label{eq:discrete-paraproduct-prec-estimate}
			\begin{aligned}
				\lVert f_1\prec f_2\rVert_{B_{p,\infty}^{s_2,\varepsilon}(\rho_1\rho_2)}
				&\lesssim
				\lVert f_1\rVert_{L^{p_1,\varepsilon}(\rho_1)}
				\lVert f_2\rVert_{B_{p_2,\infty}^{s_2,\varepsilon}(\rho_2)},\\
				\lVert f_1\prec f_2\rVert_{B_{p,\infty}^{s_1+s_2,\varepsilon}(\rho_1\rho_2)}
				&\lesssim
				\lVert f_1\rVert_{B_{p_1,\infty}^{s_1,\varepsilon}(\rho_1)}
				\lVert f_2\rVert_{B_{p_2,\infty}^{s_2,\varepsilon}(\rho_2)}.
			\end{aligned}
		\end{equation}
		
		\smallskip
		\noindent\textup{(ii)}
		For any $s_1,s_2\in\mathbb{R}$ with $s_1+s_2>0$,
		\begin{equation}
			\label{eq:discrete-resonant-product-estimate}
			\lVert f_1\circ f_2\rVert_{B_{p,\infty}^{s_1+s_2,\varepsilon}(\rho_1\rho_2)}
			\lesssim
			\lVert f_1\rVert_{B_{p_1,\infty}^{s_1,\varepsilon}(\rho_1)}
			\lVert f_2\rVert_{B_{p_2,\infty}^{s_2,\varepsilon}(\rho_2)}.
		\end{equation}
		
		\smallskip
		\noindent\textup{(iii)}
		Estimates \textup{(i)} and \textup{(ii)} remain valid for
		the periodic Besov spaces defined above,
		with $\rho_i=(w_x^M)^{a_i}$, $a_i\in\mathbb R$.
	\end{lemma}
	
	\begin{proof}
		Apply the weighted Littlewood--Paley argument of
		\cite[Lemma~4.2]{MP19}, using weighted Young's inequality for the
		dyadic kernels and H\"older's inequality with output weight
		$\rho_1\rho_2$. The same argument extends to
		$B_{p,q}^{s,\varepsilon}$ as in~\cite[Theorem~3.1]{MW17a}.
		For \textup{(iii)}, the periodized dyadic kernels satisfy the same
		weighted $L^1$ bounds uniformly in $M,\varepsilon$, so the argument
		applies on $\Lambda_{M,\varepsilon}$.
	\end{proof}
	Following \cite[Lemma~2.4]{GH19}, for $\sigma>0$, $L\ge0$,
	define the spatial localizers on $\mathcal S'(\Lambda_\varepsilon)$ by
	\begin{equation}
		\label{eq:frequency-truncation-operators}
		\begin{alignedat}{2}
			&L_k:=L-\sigma\log_2 \rho(2^k)\ge 0,
			\qquad k\ge-1,\\
			&\mathcal{U}_{>L}^\varepsilon f
			:=
			\sum_{k\ge-1}\varphi_k
			\sum_{L_k<j\le J^\varepsilon}\Delta_j^\varepsilon f,
			&&\qquad\mathcal{U}_{\leq L}^\varepsilon
			:=
			\mathrm{id}-\mathcal{U}_{>L}^\varepsilon.
		\end{alignedat}
	\end{equation}
	The dependence on $\sigma$ is suppressed in the notation.
	
	\begin{lemma}[Spatial localizer estimates]
		\label{lem:frequency-truncation-estimates}
		For every $s,a\in\mathbb{R}$, $\sigma,\delta>0$, and $L\ge0$,
		\begin{equation}
			\label{eq:low-frequency-truncation-estimate}
			\bigl\lVert\mathcal{U}_{\leq L}^\varepsilon f
			\bigr\rVert_{\mathcal C^{s+\delta,\varepsilon}(\rho^{a+\sigma\delta})}
			\lesssim_{s,a,\sigma,\delta}
			2^{\delta L}
			\lVert f\rVert_{\mathcal C^{s,\varepsilon}(\rho^a)},
		\end{equation}
		and
		\begin{equation}
			\label{eq:high-frequency-truncation-estimate}
			\bigl\lVert\mathcal{U}_{>L}^\varepsilon f
			\bigr\rVert_{\mathcal C^{s-\delta,\varepsilon}(\rho^{a-\sigma\delta})}
			\lesssim_{s,a,\sigma,\delta}
			2^{-\delta L}
			\lVert f\rVert_{\mathcal C^{s,\varepsilon}(\rho^a)}.
		\end{equation}
		Moreover, for $j=1,2$,
		\begin{equation}
			\label{eq:high-frequency-truncation-commutator-estimate}
			\begin{aligned}
				\bigl\lVert[D_j^\varepsilon,\mathcal{U}_{>L}^\varepsilon]f
				\bigr\rVert_{\mathcal C^{s-\delta,\varepsilon}(\rho^{a-\sigma\delta})}
				+\bigl\lVert[D_j^{\varepsilon *},\mathcal{U}_{>L}^\varepsilon]f
				\bigr\rVert_{\mathcal C^{s-\delta,\varepsilon}(\rho^{a-\sigma\delta})}
				\lesssim_{s,a,\sigma,\delta}
				2^{-\delta L}\lVert f\rVert_{\mathcal C^{s,\varepsilon}(\rho^a)},
			\end{aligned}
		\end{equation}
		and
		\begin{equation}
			\label{eq:low-frequency-truncation-commutator-estimate}
			\begin{aligned}
				\bigl\lVert[D_j^\varepsilon,\mathcal{U}_{\leq L}^\varepsilon]f
				\bigr\rVert_{\mathcal C^{s+\delta,\varepsilon}(\rho^{a+\sigma\delta})}
				+\bigl\lVert[D_j^{\varepsilon *},\mathcal{U}_{\leq L}^\varepsilon]f
				\bigr\rVert_{\mathcal C^{s+\delta,\varepsilon}(\rho^{a+\sigma\delta})}
				\lesssim_{s,a,\sigma,\delta}
				2^{\delta L}\lVert f\rVert_{\mathcal C^{s,\varepsilon}(\rho^a)}.
			\end{aligned}
		\end{equation}
	\end{lemma}

	\begin{proof}
		The proofs of~\eqref{eq:low-frequency-truncation-estimate}
		and~\eqref{eq:high-frequency-truncation-estimate} are the same as
		that of~\cite[Lemma~2.4]{GH19}. We prove
		\eqref{eq:high-frequency-truncation-commutator-estimate}.
		Since $D_j^\varepsilon$ commutes with the Littlewood--Paley blocks,
		\begin{equation}
			\label{eq:high-frequency-truncation-commutator-identity}
			\begin{aligned}
				[D_j^\varepsilon,\mathcal U_{>L}^\varepsilon]f
			=\sum_{k\ge-1}[D_j^\varepsilon,\varphi_k]
				\sum_{L_k<i\le J^\varepsilon}\Delta_i^\varepsilon f
				=\sum_{k\ge-1}(D_j^\varepsilon\varphi_k)\tau_j
				\sum_{L_k<i\le J^\varepsilon}\Delta_i^\varepsilon f.
			\end{aligned}
		\end{equation}
		Following~\cite[Lemma~2.4]{GH19}, choose a noninteger
		$\beta>|s-\delta|$. By the lattice analogue of~\cite[(2.3)]{GH19},
		finite overlap, the weighted product estimate, and the uniform
		boundedness of $\tau_j$,
		\begin{equation}
			\label{eq:high-frequency-truncation-commutator-localization}
			\begin{aligned}
				\bigl\lVert[D_j^\varepsilon,\mathcal U_{>L}^\varepsilon]f
				\bigr\rVert_{\mathcal C^{s-\delta,\varepsilon}(\rho^{a-\sigma\delta})}
				&\lesssim\sup_{k\ge-1}\biggl\lVert
				(D_j^\varepsilon\varphi_k)\tau_j
				\sum_{L_k<i\le J^\varepsilon}\Delta_i^\varepsilon f
				\biggr\rVert_{\mathcal C^{s-\delta,\varepsilon}(\rho^{a-\sigma\delta})}\\
				&\lesssim\sup_{k\ge-1}\biggl(
				\lVert D_j^\varepsilon\varphi_k
				\rVert_{\mathcal C^{\beta,\varepsilon}(\rho^{-\sigma\delta})}
				\biggl\lVert\sum_{L_k<i\le J^\varepsilon}\Delta_i^\varepsilon f
				\biggr\rVert_{\mathcal C^{s-\delta,\varepsilon}(\rho^a)}
				\biggr)\\
				&\lesssim\sup_{k\ge-1}\bigl(
				\lVert D_j^\varepsilon\varphi_k
				\rVert_{\mathcal C^{\beta,\varepsilon}(\rho^{-\sigma\delta})}
				2^{-L_k\delta}\bigr)
				\lVert f\rVert_{\mathcal C^{s,\varepsilon}(\rho^a)}.
			\end{aligned}
		\end{equation}
		On the support of $D_j^\varepsilon\varphi_k$,
		$\rho(x)\asymp\rho(2^k)$. The weighted H\"older characterization
		and the smooth cutoff derivative bounds yield
		\begin{equation}
			\label{eq:high-frequency-truncation-cutoff-weight-bound}
			\begin{aligned}
				\lVert D_j^\varepsilon\varphi_k
				\rVert_{\mathcal C^{\beta,\varepsilon}(\rho^{-\sigma\delta})}
				\lesssim\rho(2^k)^{-\sigma\delta}
				\lVert D_j^\varepsilon\varphi_k\rVert_{\mathcal C^{\beta,\varepsilon}}
				\lesssim\rho(2^k)^{-\sigma\delta}
				\sum_{\substack{\alpha\in\mathbb N_0^2\\
					|\alpha|\le\lfloor\beta\rfloor+1}}
				\lVert\partial^\alpha D_j^\varepsilon\varphi_k
				\rVert_{L^\infty(\mathbb R^2)}
				\lesssim\rho(2^k)^{-\sigma\delta}.
			\end{aligned}
		\end{equation}
		Since $2^{-L_k\delta}\rho(2^k)^{-\sigma\delta}=2^{-\delta L}$,
		this yields the forward-difference estimate.
		The same argument applies to $D_j^{\varepsilon *}$.
		Estimate~\eqref{eq:low-frequency-truncation-commutator-estimate}
		follows similarly.
	\end{proof}

	For $z_1,z_2\in\mathbb R^2$, define
$
		\label{eq:pinned-periodic-distance}
		d_M(z_1,z_2):=
		\inf_{n\in\mathbb Z^2}|z_1-z_2+2\pi Mn|.
$

	\begin{lemma}[Periodic multiplier kernels]
		\label{lem:periodic-multiplier-kernel-bound}
		Let $N\ge1$ be an integer, $-1\le j\le J^\varepsilon$, and
		$m\in C^N(\widehat\Lambda_\varepsilon;\mathbb R)$.
		Suppose that
		$\operatorname{supp}m\cap Q_{M,\varepsilon}
		\subset\{\xi:|\xi|\le c2^j\}$,
		where $c>0$, and
		\begin{equation*}
			\|m\|_{L^\infty}\le C_0,
			\qquad
			\max_{|\gamma|=N}\|\partial_\xi^\gamma m\|_{L^\infty}\le C_N.
		\end{equation*}
		Set
		\begin{equation}
			\label{eq:periodic-multiplier-kernel-definition}
			K(z):=
			\mathcal F_{M,\varepsilon}^{-1}m(z)
			=\frac1{|\Lambda_{M,\varepsilon}|}
			\sum_{k\in Q_{M,\varepsilon}}m(k)\,\mathrm e_k(z).
		\end{equation}
		Then
		\begin{equation}
			\label{eq:periodic-multiplier-kernel-min-bound}
			|K(z)|\lesssim_{N,c}
			2^{2j}\min\{C_Nd_M(z,0)^{-N},C_0\},
			\qquad z\in\Lambda_{M,\varepsilon}.
		\end{equation}
	\end{lemma}

	\begin{proof}
		For $z\ne0$, choose its representative in $[-\pi M,\pi M)^2$
		and a coordinate $r$ with
		$|z_r|=\max\{|z_1|,|z_2|\}\ge d_M(z,0)/\sqrt2$.
		Set $\delta_r^Mm(k):=
		m(k+M^{-1}e_r)-m(k)$. Discrete summation by parts gives
		\begin{equation}
			\label{eq:periodic-multiplier-summation-by-parts}
			(\mathrm e^{-iz_r/M}-1)^N K(z)
			=\frac1{|\Lambda_{M,\varepsilon}|}
			\sum_{k\in Q_{M,\varepsilon}}
			(\delta_r^M)^Nm(k)\,\mathrm e_k(z).
		\end{equation}
		Since $|z_r|\le\pi M$,
$
			\label{eq:periodic-multiplier-phase-bound}
			|\mathrm e^{-iz_r/M}-1|
			\gtrsim\frac{d_M(z,0)}{M}.
$
		The fundamental theorem of calculus yields
		\begin{equation}
			\label{eq:periodic-multiplier-difference-bound}
			|(\delta_r^M)^Nm(k)|
			\le\int_{[0,M^{-1}]^N}
			|\partial_{\xi_r}^Nm(k+(t_1+\cdots+t_N)e_r)|
			\,dt_1\cdots dt_N
			\le C_NM^{-N}.
		\end{equation}
		The support assumption gives
		\begin{equation}
			\label{eq:periodic-multiplier-frequency-count}
			\frac{\#\{k\in Q_{M,\varepsilon}:(\delta_r^M)^Nm(k)\ne0\}}
			{|\Lambda_{M,\varepsilon}|}
			\lesssim_{N,c}2^{2j}.
		\end{equation}
		Combining these bounds gives
		\begin{equation}
			\label{eq:periodic-multiplier-kernel-decay-proof}
			|K(z)|\lesssim_{N,c}
			C_N2^{2j}d_M(z,0)^{-N}.
		\end{equation}
		Taking $N=0$ gives
		\begin{equation}
			\label{eq:periodic-multiplier-kernel-direct-proof}
			|K(z)|
			\le\frac{C_0}{|\Lambda_{M,\varepsilon}|}
			\#\{k\in Q_{M,\varepsilon}:|k|\le c2^j\}
			\lesssim_c C_0 2^{2j}.
		\end{equation}
The claim follows by combining the two bounds.
	\end{proof}

	\begin{lemma}[Weighted periodic Littlewood--Paley kernels]
		\label{lem:weighted-periodic-LP-kernel}
		For $x\in\Lambda_\varepsilon$ and $-1\le j\le J^\varepsilon$,
		set $\eta_{j,x}^{M,\varepsilon}:=
		\eta_j^{M,\varepsilon}(x-\cdot)$.
		For $a,s\in\mathbb R$ and $p,q\in[1,\infty]$,
		\begin{equation}
			\label{eq:weighted-periodic-LP-kernel-Besov-bound}
			\|\eta_{j,x}^{M,\varepsilon}
			\|_{B_{p,q}^{s}((w_x^M)^a;\Lambda_{M,\varepsilon})}
			\lesssim_{a,s,p,q}
			2^{j(s+2(1-1/p))}.
		\end{equation}
		The same bound holds for
		$\|\eta_j^\varepsilon\|_{B_{p,q}^{s,\varepsilon}(\rho^a)}$.
	\end{lemma}

	\begin{proof}
		In the frequency cell, for some $c_1,c_2>0$,
		\begin{equation}
			\label{eq:weighted-periodic-LP-frequency-support}
			\operatorname{supp}\varphi_j^\varepsilon
			\subset
			\begin{cases}
				\{\xi:c_1 2^j\le|\xi|\le c_2 2^j\},&j\ge0,\\
				\{\xi:|\xi|\le c_2\},&j=-1.
			\end{cases}
		\end{equation}
		By periodization,
		\begin{equation}
			\label{eq:weighted-periodic-LP-Fourier-sum}
			\eta_j^{M,\varepsilon}(z)
			=\frac1{|\Lambda_{M,\varepsilon}|}
			\sum_{k\in Q_{M,\varepsilon}}
			\varphi_j^\varepsilon(k)\,\mathrm e_k(z).
		\end{equation}
		For $0\le j<J^\varepsilon$ and $j=J^\varepsilon$, respectively,
		\begin{equation}
			\label{eq:periodic-LP-cutoff-identities}
			\begin{aligned}
				\varphi_j^\varepsilon(\xi)
				=\varphi(2^{-j-1}\xi)-\varphi(2^{-j}\xi),
				\quad\varphi_{J^\varepsilon}^\varepsilon(\xi)
				=1-\varphi(2^{-J^\varepsilon}\xi).
			\end{aligned}
		\end{equation}
		Differentiating \eqref{eq:periodic-LP-cutoff-identities} gives
		\begin{equation}
			\label{eq:periodic-LP-multiplier-derivatives}
			\|\partial_\xi^\gamma\varphi_j^\varepsilon\|_{L^\infty}
			\lesssim_\gamma2^{-j|\gamma|}.
		\end{equation}
		The same bound holds for $j=-1$, since
		$\varphi_{-1}^\varepsilon=\varphi$.
		By Lemma~\ref{lem:periodic-multiplier-kernel-bound}, for $N\ge1$,
		\begin{equation}
			\label{eq:periodic-LP-kernel-pointwise-decay}
			|\eta_j^{M,\varepsilon}(z)|
			\lesssim_N2^{2j}(1+2^jd_M(z,0))^{-N}.
		\end{equation}
		For $j=-1$, apply the lemma with $j=0$.
		The definition of $w_x^M$ gives
$
			\label{eq:periodized-weight-pointwise-bounds}
			(1+\nu^2d_M(x,y)^2)^{-2}
			\le w_x^M(y)\lesssim1.
$
		Hence
		\begin{equation}
			\label{eq:periodized-weight-kernel-growth}
			(w_x^M(y))^a
			\lesssim_a(1+2^jd_M(x,y))^{4|a|}.
		\end{equation}
		Let $p<\infty$ and choose $N>4|a|+2$.
		Translating $x$ to $0$, we obtain
		\begin{equation}
			\label{eq:weighted-periodic-LP-kernel-lattice-integral}
			\begin{aligned}
				\|(w_x^M)^a\eta_{j,x}^{M,\varepsilon}
				\|_{L^{p,\varepsilon}(\Lambda_{M,\varepsilon})}
				&\lesssim_{a,p}2^{2j}
				\left(\int_{\Lambda_{M,\varepsilon}}
		 (w_0^M(z))^{ap}
				(1+2^jd_M(z,0))^{-Np}\,dz\right)^{1/p}\\
				&\lesssim_{a,p}2^{2j}
				\left(\int_{\Lambda_\varepsilon}
				(1+2^j|z|)^{-(N-4|a|)p}\,dz\right)^{1/p}.
			\end{aligned}
		\end{equation}
		Since $\varepsilon2^j\lesssim1$, for $z\in\Lambda_\varepsilon$
		and $N_0>0$,
		$
			\label{eq:LP-weight-cell-comparison}
			\sup_{\xi\in z+[-\varepsilon/2,\varepsilon/2]^2}
			(1+2^j|\xi|)^{N_0}
			\lesssim_{N_0}(1+2^j|z|)^{N_0}.
	$
		Applying \cite[Lemma~A.2]{MP19} with $r=(N-4|a|)p>2$ gives
		\begin{equation}
			\label{eq:LP-weight-lattice-integral-comparison}
			\int_{\Lambda_\varepsilon}(1+2^j|z|)^{-r}\,dz
			\lesssim_r\int_{\mathbb R^2}(1+2^j|z|)^{-r}\,dz
			\lesssim_r2^{-2j}.
		\end{equation}
		Substituting this into
		\eqref{eq:weighted-periodic-LP-kernel-lattice-integral}, we get
		\begin{equation}
			\label{eq:weighted-periodic-LP-kernel-Lp-bound}
			\|(w_x^M)^a\eta_{j,x}^{M,\varepsilon}
			\|_{L^{p,\varepsilon}(\Lambda_{M,\varepsilon})}
			\lesssim_{a,p}2^{2j(1-1/p)}.
		\end{equation}
		The case $p=\infty$ follows from the pointwise bounds.
		Since $\Delta_\ell^\varepsilon\eta_{j,x}^{M,\varepsilon}=0$
		when $|\ell-j|>1$,
		\begin{equation}
			\label{eq:weighted-periodic-LP-kernel-Besov-proof}
			\begin{aligned}
				\|\eta_{j,x}^{M,\varepsilon}
				\|_{B_{p,q}^{s}((w_x^M)^a;\Lambda_{M,\varepsilon})}
				\lesssim_{s,q}2^{js}
				\max_{|\ell-j|\le1}
				\|(w_x^M)^a\Delta_\ell^\varepsilon\eta_{j,x}^{M,\varepsilon}
				\|_{L^{p,\varepsilon}(\Lambda_{M,\varepsilon})}
				\lesssim_{a,s,p,q}2^{j(s+2(1-1/p))}.
			\end{aligned}
		\end{equation}
	\end{proof}
	
	\begin{lemma}[Periodized weight bounds]
		\label{lem:periodized-weight-power-bounds}
		Let $a\in\mathbb R$ and $\alpha\in\mathbb N_0^2$. The following estimates hold.

		\smallskip
		\noindent\textup{(i)}
		For $x\in\mathbb R^2$ and $y\in\Lambda_\varepsilon$,
		\begin{equation}
			\label{eq:periodized-weight-power-derivatives}
			\bigl|(D^\varepsilon)^\alpha(w_x^M)^a(y)\bigr|
			\lesssim_{a,\alpha}(w_x^M(y))^a.
		\end{equation}

		\smallskip
		\noindent\textup{(ii)}
		For $r>1/2$ and $z\in\mathbb R^2$,
		\begin{equation}
			\label{eq:periodized-weight-power-sum-comparison}
			\sum_{n\in\mathbb Z^2}\rho(z+2\pi Mn)^r
			\asymp_r
			\Bigl(\sup_{n\in\mathbb Z^2}\rho(z+2\pi Mn)\Bigr)^r
			\asymp_r
			\Bigl(\sum_{n\in\mathbb Z^2}\rho(z+2\pi Mn)\Bigr)^r.
		\end{equation}
	\end{lemma}

	\begin{proof}
		\textup{(i)} Write $w_x:=w_x^M$.
		By~\eqref{eq:weight-basic-estimates},
$
			\label{eq:periodized-weight-derivative-bound}
			\bigl|(D^\varepsilon)^\alpha w_x(y)\bigr|
			\le
			\sum_{n\in\mathbb Z^2}
			\bigl|(D^\varepsilon)^\alpha\rho_x(y+2\pi Mn)\bigr|
			\lesssim_\alpha w_x(y).
$
		The chain rule gives~\eqref{eq:periodized-weight-power-derivatives}.

		\smallskip
		\noindent\textup{(ii)} Clearly,
$
			\label{eq:periodized-weight-power-sup-lower-bound}
			\Bigl(\sup_{n\in\mathbb Z^2}\rho(z+2\pi Mn)\Bigr)^r
			\le\sum_{n\in\mathbb Z^2}\rho(z+2\pi Mn)^r.
$
		We now prove the reverse inequality.
		Choose $n_0\in\mathbb Z^2$ maximizing $\rho(z+2\pi Mn)$ and set
		$z_0:=z+2\pi Mn_0$.
		Since $\rho$ is radially decreasing, $z_0\in[-\pi M,\pi M]^2$.
		Reindexing gives
		\begin{equation}
			\label{eq:periodized-weight-power-reindexing}
			\sum_{n\in\mathbb Z^2}\rho(z+2\pi Mn)^r
			=\sum_{n\in\mathbb Z^2}\rho(z_0+2\pi Mn)^r.
		\end{equation}
		For $n\in\mathbb Z^2$ with $|n|_\infty>2$, choose $j\in\{1,2\}$
		with $|n_j|=|n|_\infty$. Then
		\begin{equation}
			\label{eq:periodized-weight-tail-distance}
			|z_0+2\pi Mn|
			\ge|(z_0)_j+2\pi Mn_j|
			\ge2\pi M|n|_\infty-\pi M
			\gtrsim M|n|\gtrsim|n|.
		\end{equation}
		Consequently, since $r>1/2$ and $M\ge1$,
		\begin{equation}
			\label{eq:periodized-weight-power-tail-bound}
			\begin{aligned}
				\sum_{\substack{n\in\mathbb Z^2\\|n|_\infty>2}}
				\rho(z_0+2\pi Mn)^r
				\lesssim_r
				\sum_{\substack{n\in\mathbb Z^2\\|n|_\infty>2}}
				(1+2M|n|)^{-4r}
				\lesssim_r M^{-4r}
				\sum_{\substack{n\in\mathbb Z^2\\|n|_\infty>2}}|n|^{-4r}
				\lesssim_r M^{-4r}\lesssim_r\rho(z_0)^r.
			\end{aligned}
		\end{equation}
		By maximality,
$
			\label{eq:periodized-weight-power-near-bound}
			\sum_{\substack{n\in\mathbb Z^2\\|n|_\infty\le2}}
			\rho(z_0+2\pi Mn)^r\le25\rho(z_0)^r.
$
		This gives the first comparison
		in~\eqref{eq:periodized-weight-power-sum-comparison}.
		Applying it with $r=1$ and raising to the power $r$ gives the second.
	\end{proof}

	\begin{lemma}[Periodized weights]
		\label{lem:periodized-weighted-Besov-equivalences}
		Let $f\in(\mathbb R^2)^{\Lambda_{M,\varepsilon}}$,
		$a,s\in\mathbb R$, and $p,q\in[1,\infty]$. The following estimates hold.

		\smallskip
		\noindent\textup{(i)}
		\begin{equation}
			\label{eq:periodized-weight-Besov-multiplication}
			\|f(w_x^M)^a\|_{B_{p,q}^{s}(\Lambda_{M,\varepsilon})}
			\simeq
			\|f\|_{B_{p,q}^{s}((w_x^M)^a;\Lambda_{M,\varepsilon})}.
		\end{equation}

		\smallskip
		\noindent\textup{(ii)}
		\begin{equation}
			\label{eq:translated-weight-holder-comparison}
			\begin{aligned}
				\rho(x)^{|a|}\|f\|_{B_{p,q}^{s,\varepsilon}(\rho^a)}
				\lesssim
				\|f\|_{B_{p,q}^{s,\varepsilon}(\rho_x^a)}
				\lesssim
				\rho(x)^{-|a|}\|f\|_{B_{p,q}^{s,\varepsilon}(\rho^a)}.
			\end{aligned}
		\end{equation}

		\smallskip
		\noindent\textup{(iii)}
		If $p<\infty$ and $ap>1/2$, or if $p=\infty$ and $a\ge0$, then
		\begin{equation}
			\label{eq:periodized-full-lattice-Besov-equivalence}
			\|f\|_{B_{p,q}^{s}((w_x^M)^a;\Lambda_{M,\varepsilon})}
			\simeq
			\|f\|_{B_{p,q}^{s,\varepsilon}(\rho_x^a)}.
		\end{equation}
	\end{lemma}

	\begin{proof}
		Write $w:=w_x^M$.
		\textup{(i)} Choose an integer $N>|s|+1$.
		Arguing as in \cite[Lemma~A.1]{GH21}, using the paraproduct estimates
		of Lemma~\ref{lem:discrete-paraproduct-estimates}, it suffices to check
		\begin{equation}
			\label{eq:periodized-weight-multiplier-holder-bounds}
			\|w^a\|_{\mathcal C^{N+1}(w^{-a};\Lambda_{M,\varepsilon})}
			+\|w^{-a}\|_{\mathcal C^{N+1}(w^a;\Lambda_{M,\varepsilon})}
			\lesssim_{a,s}1.
		\end{equation}
		Indeed, the weighted H\"older characterization of Besov spaces and
		Lemma~\ref{lem:periodized-weight-power-bounds}\textup{(i)} give
		\begin{equation}
			\label{eq:periodized-weight-multiplier-derivative-control}
			\|w^a\|_{\mathcal C^{N+1}(w^{-a};\Lambda_{M,\varepsilon})}
			\lesssim_{a,N}
			\sum_{\substack{\alpha\in\mathbb N_0^2\\|\alpha|\le N+2}}
			\|(D^\varepsilon)^\alpha w^a
			\|_{L^{\infty,\varepsilon}(w^{-a};\Lambda_{M,\varepsilon})}
			\lesssim_{a,N}1.
		\end{equation}
		The same argument with $a$ replaced by $-a$ controls the second term.

		\smallskip
		\noindent\textup{(ii)} The weight comparison
		in~\eqref{eq:weight-basic-estimates} yields
$
			\label{eq:translated-polynomial-weight-pointwise-comparison}
			\rho(x)^{|a|}\rho(y)^a
			\lesssim_a\rho_x(y)^a
			\lesssim_a\rho(x)^{-|a|}\rho(y)^a,
 y\in\Lambda_\varepsilon.
$
		Apply this to each Littlewood--Paley block to
		obtain~\eqref{eq:translated-weight-holder-comparison}.

		\smallskip
		\noindent\textup{(iii)} Let $g$ be a $2\pi M\mathbb Z^2$-periodic
		lattice field. If $p<\infty$ and $ap>1/2$, then
		Lemma~\ref{lem:periodized-weight-power-bounds}\textup{(ii)} gives
		\begin{equation}
			\label{eq:periodized-full-lattice-Lp-equivalence}
			\begin{aligned}
				\|g\|_{L^{p,\varepsilon}(\rho_x^a)}^p
				&=\int_{\Lambda_\varepsilon}\rho_x(y)^{ap}|g(y)|^p\,dy\\
				&=\int_{\Lambda_{M,\varepsilon}}
				\sum_{n\in\mathbb Z^2}\rho_x(y+2\pi Mn)^{ap}
				|g(y)|^p\,dy\\
				&\asymp_{a,p}\int_{\Lambda_{M,\varepsilon}}
				\Bigl(\sum_{n\in\mathbb Z^2}\rho_x(y+2\pi Mn)\Bigr)^{ap}
				|g(y)|^p\,dy
				=\|g\|_{L^{p,\varepsilon}(w^a;\Lambda_{M,\varepsilon})}^p.
			\end{aligned}
		\end{equation}
		For $p=\infty$ and $a\ge0$, the same norm equivalence follows
		from
		Lemma~\ref{lem:periodized-weight-power-bounds}\textup{(ii)}.
		Apply these equivalences to $g=\Delta_j^\varepsilon f$ and take the
		$\ell^q$ norm over $j\ge-1$ to obtain
		\eqref{eq:periodized-full-lattice-Besov-equivalence}.
	\end{proof}

	Following~\cite[Section~2.4]{MP19}, we introduce a smooth extension
	operator from the full lattice to $\mathbb R^2$. Let
	$(\varphi_j^\varepsilon)_{\mathrm{per}}$ denote the
	$(2\pi/\varepsilon)\mathbb Z^2$-periodic extension of
	$\varphi_j^\varepsilon$ to $\mathbb R^2$.
	Choose a symmetric function
	$\psi\in C_c^\infty((-3\pi/2,3\pi/2)^2;[0,1])$ such that, with
	$\psi_\varepsilon=\psi(\varepsilon\mathord\cdot)$,
	\begin{equation}
		\label{eq:lattice-extension-smear-properties}
		\begin{aligned}
			\textup{(i)}\quad
			&\psi_\varepsilon=1
			\quad\text{on }\operatorname{supp}\varphi_j^\varepsilon
			\quad\text{whenever }j<J^\varepsilon-1,\\
			\textup{(ii)}\quad
			&\left(
			\operatorname{supp}\psi_\varepsilon
			\cap\operatorname{supp}(\varphi_j^\varepsilon)_{\mathrm{per}}
			\right)\setminus
			\left[-\frac{\pi}{\varepsilon},\frac{\pi}{\varepsilon}\right)^2
			=\varnothing
			\quad\text{whenever }j<J^\varepsilon-1,\\
			\textup{(iii)}\quad
			&\sum_{n\in\mathbb Z^2}\psi(\mathord\cdot+2\pi n)=1.
		\end{aligned}
	\end{equation}

	For $f\in\mathcal S'(\Lambda_\varepsilon)$, write
	$\widehat f_{\mathrm{per}}$ for the
	$(2\pi/\varepsilon)\mathbb Z^2$-periodic lift of its discrete Fourier
	transform to $\mathbb R^2$. We define
	\begin{equation}
		\label{eq:lattice-sharp-extension-definition}
		\begin{aligned}
			\mathcal E^\varepsilon\colon
			\mathcal S'(\Lambda_\varepsilon)
			&\longrightarrow C^\infty(\mathbb R^2)\cap\mathcal S'(\mathbb R^2),\\
			f&\longmapsto
			\mathcal E^\varepsilon f
			:=\mathcal F_{\mathbb R^2}^{-1}
			\bigl(\psi_\varepsilon\widehat f_{\mathrm{per}}\bigr).
		\end{aligned}
	\end{equation}
	We use the analogous weighted Besov norms $\lVert\cdot\rVert_{B_{p,q}^s(\rho)}$ on $\mathbb R^2$.
	By \cite[Lemma~2.24]{MP19}, we have the following estimate:
	\begin{lemma}[Uniform estimate for the extension operator]
		\label{lem:sharp-extension-Besov-equivalence}
		\label{lem:lattice-smooth-extension-besov-bound}
		For every $s,a\in\mathbb R$ and $p,q\in[1,\infty]$,
		\begin{equation}
			\label{eq:sharp-extension-Besov-equivalence}
			\left\lVert\mathcal E^\varepsilon f
			\right\rVert_{B_{p,q}^s(\rho^a)}
			\lesssim_{s,p,q,a,\psi}
			\lVert f\rVert_{B_{p,q}^{s,\varepsilon}(\rho^a)}.
		\end{equation}
	\end{lemma}

	\subsection{Covariant Sobolev spaces}
	\label{sec:covariant-sobolev-spaces}

	Throughout this subsection, fix a lattice gauge field
	$B\in(\mathbb{R}^2)^{\Lambda_\varepsilon}$.
	For each fixed $\varepsilon$, the operator
	$\mathcal{L}_B^\varepsilon+1$ is bounded and boundedly invertible on
	$L^{p,\varepsilon}$, $1\le p\le\infty$.
	On $L^{2,\varepsilon}$, it is positive definite and self-adjoint.
	For $t\in\mathbb{R}$ we define
$
	    \label{eq:lattice-covariant-operator-exponential}
	    e^{-t\mathcal{L}_B^\varepsilon}
	    :=\sum_{n=0}^{\infty}\frac{(-t)^n}{n!}
	    (\mathcal{L}_B^\varepsilon)^n.
$
	The series converges in operator norm on each $L^{p,\varepsilon}$.
	For $t\ge0$, the diamagnetic inequality~\eqref{eq:diamagnetic-semigroup}
	gives
	$\|e^{-t(\mathcal{L}_B^\varepsilon+1)}\|_{L^{p,\varepsilon}\to L^{p,\varepsilon}}
	\le e^{-t}$.
	Hence the following integrals converge in operator norm:
	\begin{equation}
	    \label{eq:lattice-covariant-inverse-definition}
	    (\mathcal{L}_B^\varepsilon+1)^{-1}
	    :=
	    \int_0^\infty e^{-t(\mathcal{L}_B^\varepsilon+1)}\,dt.
	\end{equation}
	This operator is the inverse of $\mathcal{L}_B^\varepsilon+1$.
	For $s>0$, define
	\begin{equation}
	    \label{eq:lattice-covariant-negative-powers-definition}
	    (\mathcal{L}_B^\varepsilon+1)^{-s}
	    :=\frac{1}{\Gamma(s)}
	    \int_0^\infty t^{s-1}
	    e^{-t(\mathcal{L}_B^\varepsilon+1)}\,dt.
	\end{equation}
	For $s>0$, choose an integer $m>s$ and define
$
	    \label{eq:lattice-covariant-positive-powers-definition}
	    (\mathcal{L}_B^\varepsilon+1)^s
	    :=
	    (\mathcal{L}_B^\varepsilon+1)^m
	    (\mathcal{L}_B^\varepsilon+1)^{-(m-s)}.
$
The definition is independent of the choice of $m$.
	
	For $s\in\mathbb{R}$, $p\in[1,\infty]$, and a polynomial weight $\rho$,
	we define the weighted covariant Sobolev norms by
	\begin{equation}
		\label{eq:covariant-sobolev-norm}
			\lVert f\rVert_{W_{p,B}^{s,\varepsilon}(\rho)}
			:=
			\bigl\lVert
			(\mathcal{L}_B^\varepsilon+1)^{s/2}(\rho f)
			\bigr\rVert_{L^{p,\varepsilon}},\quad
		H_B^{s,\varepsilon}(\rho)
			:=W_{2,B}^{s,\varepsilon}(\rho).
	\end{equation}
	At fixed $\varepsilon$, these spaces consist of lattice functions with
	$\rho f\in L^{p,\varepsilon}$. We omit the weight when $\rho=1$.

	Each inequality below is understood for $f$ such that its right-hand side is finite.

	\begin{lemma}
		\label{lem:covariant-nonnegative-weight}
		For $\psi\in C^2(\mathbb R^2;\mathbb C)$ with finite $C^2$ norm and
		$f\in H_B^{s,\varepsilon}$,
		\begin{equation}
			\label{eq:covariant-nonnegative-weight}
			\lVert\psi f\rVert_{H_B^{s,\varepsilon}}
			\lesssim_s\lVert\psi\rVert_{C^2}
			\lVert f\rVert_{H_B^{s,\varepsilon}},
			\qquad |s|\le2.
		\end{equation}
	\end{lemma}
	\begin{proof}
		The case $s=0$ is immediate. By the fixed-weight case of
		Lemma~\ref{lem:weighted-operator-interpolation} and duality,
 it suffices to consider $s=2$.
		Lemma~\ref{lem:covariant-laplacian-weight-commutator} with $w=\psi$ gives
		\begin{equation}
			\label{eq:covariant-nonnegative-weight-s-two}
			\begin{aligned}
				\lVert\psi f\rVert_{H_B^{2,\varepsilon}}
				\le\bigl\lVert\psi(\mathcal{L}_B^\varepsilon+1)f
				\bigr\rVert_{L^{2,\varepsilon}}
				+\bigl\lVert[\mathcal{L}_B^\varepsilon,\psi]f
				\bigr\rVert_{L^{2,\varepsilon}}
				\lesssim\lVert\psi\rVert_{C^2}
				\lVert f\rVert_{H_B^{2,\varepsilon}}.
			\end{aligned}
		\end{equation}
		Here the discrete derivatives of $\psi$ are bounded by
		$\lVert\psi\rVert_{C^2}$, and
		$\lVert f\rVert_{L^{2,\varepsilon}}
		+\lVert D_B^\varepsilon f\rVert_{L^{2,\varepsilon}}
		+\lVert D_B^{\varepsilon *}f\rVert_{L^{2,\varepsilon}}
		\lesssim\lVert f\rVert_{H_B^{2,\varepsilon}}$.
	\end{proof}

	\begin{lemma}[Covariant interpolation inequality]
		\label{lem:covariant-interpolation-inequality}
		For $s_1,s_2\in[-2,2]$, $a_1,a_2\in\mathbb{R}$, and
		$\theta\in[0,1]$,
		\begin{equation}
			\label{eq:covariant-interpolation-inequality}
			\lVert f\rVert_{
				H_B^{s_1\theta+s_2(1-\theta),\varepsilon}
				(\rho^{a_1\theta+a_2(1-\theta)})}
			\lesssim
			\lVert f\rVert_{H_B^{s_1,\varepsilon}(\rho^{a_1})}^\theta
			\lVert f\rVert_{H_B^{s_2,\varepsilon}(\rho^{a_2})}^{1-\theta}.
		\end{equation}
		If $a_1=a_2$, the regularities $s_1,s_2$ may be arbitrary real
		numbers, and the constant is one. 
	\end{lemma}
	\begin{proof}
		By density, we may assume that $f$ is finitely supported.
		Choose a nonzero finitely supported $h\in L^{2,\varepsilon}$ and define
$
			\label{eq:covariant-interpolation-rank-one-operator}
			T_f g:=\langle g,h\rangle_\varepsilon f, g\in L^{2,\varepsilon}.
$
		For $s,a\in\mathbb R$,
		\begin{equation}
			\label{eq:covariant-interpolation-rank-one-norm}
			\begin{aligned}
				\|T_f\|_{L^{2,\varepsilon}\to H_B^{s,\varepsilon}(\rho^a)}
				=\sup_{\|g\|_{L^{2,\varepsilon}}=1}
				|\langle g,h\rangle_\varepsilon|\,\|f\|_{H_B^{s,\varepsilon}(\rho^a)}
				=\|h\|_{L^{2,\varepsilon}}\|f\|_{H_B^{s,\varepsilon}(\rho^a)}.
			\end{aligned}
		\end{equation}
		Apply Lemma~\ref{lem:weighted-operator-interpolation} with both input
		spaces equal to $L^{2,\varepsilon}$ and output spaces
		$H_B^{s_1,\varepsilon}(\rho^{a_1})$ and
		$H_B^{s_2,\varepsilon}(\rho^{a_2})$.
		Cancelling $\|h\|_{L^{2,\varepsilon}}>0$ gives
		\eqref{eq:covariant-interpolation-inequality}.
		If $a_1=a_2$, the same lemma allows arbitrary real regularities
		and gives constant one. 
	\end{proof}

	The proofs of Lemmas~\ref{lem:WpB-Wp-comparison} and
	\ref{lem:covariant-derivative-estimates} adapt the arguments in
	\cite[Lemmas~2.2--2.3 and~2.6, and the proof of Proposition~2.8]{BCCKM26},
	using weighted discrete product estimates.

	\begin{lemma}
		\label{lem:WpB-Wp-comparison}
		The following estimates hold.
		\begin{enumerate}
			\item[\textup{(1)}]
			\textit{Diamagnetic inequality.}
			For every $p\in[1,\infty]$, $f\in L^{p,\varepsilon}$,
			$x\in\Lambda_\varepsilon$, $t\ge0$, $s\ge0$, and
			$j\in\{1,2\}$,
			\begin{equation}
				\label{eq:diamagnetic-semigroup}
				\bigl\lvert e^{-t\mathcal{L}_B^\varepsilon}f(x)\bigr\rvert
				\le
				e^{-t\mathcal{L}_0^\varepsilon}\lvert f\rvert(x),
			\end{equation}
			\begin{equation}
				\label{eq:diamagnetic-resolvent}
				\bigl\lvert
				(1+\mathcal{L}_B^\varepsilon)^{-s/2}f(x)
				\bigr\rvert
				\le
				(1+\mathcal{L}_0^\varepsilon)^{-s/2}\lvert f\rvert(x),
			\end{equation}
			\begin{equation}
				\label{eq:diamagnetic-gradient}
				\bigl\lvert D_j^\varepsilon\lvert f\rvert(x)\bigr\rvert
				\le
				\bigl\lvert D_{B,j}^\varepsilon f(x)\bigr\rvert.
			\end{equation}

			\item[\textup{(2)}]
			\textit{Sobolev inequality.}
			If $1<p\le q<\infty$, $s,a\in\mathbb{R}$, and
			$\delta\ge\tfrac{2}{p}-\tfrac{2}{q}$, then
			\begin{equation}
				\label{eq:covariant-sobolev-embedding}
				\lVert f\rVert_{W_{q,B}^{s-\delta,\varepsilon}(\rho^a)}
				\lesssim_{\delta,p,q}
				\lVert f\rVert_{W_{p,B}^{s,\varepsilon}(\rho^a)}.
			\end{equation}
			If $\delta>\tfrac{2}{p}-\tfrac{2}{q}$, the same estimate holds
			for $1\le p\le q\le\infty$.

			\item[\textup{(3)}]
			\textit{Relation between $H^{s,\varepsilon}$ and
				$H_B^{s,\varepsilon}$.}
			If $B\in H^{1,\varepsilon}(\rho)$, then, for
			$s\in[-2,2]$ and $a\in\mathbb{R}$,
			\begin{equation}
				\label{eq:weighted-HB-H-comparison}
				\begin{aligned}
					\lVert f\rVert_{H_B^{s,\varepsilon}(\rho^{a+|s|})}
					&\lesssim
					\bigl(1+\lVert B\rVert_{H^{1,\varepsilon}(\rho)}\bigr)^{|s|}
					\lVert f\rVert_{H^{s,\varepsilon}(\rho^a)},\\
					\lVert f\rVert_{H^{s,\varepsilon}(\rho^{a+|s|})}
					&\lesssim
					\bigl(1+\lVert B\rVert_{H^{1,\varepsilon}(\rho)}\bigr)^{|s|}
					\lVert f\rVert_{H_B^{s,\varepsilon}(\rho^a)}.
				\end{aligned}
			\end{equation}
		\end{enumerate}
	\end{lemma}

	\begin{proof}
		\textup{(1).}
		The diamagnetic estimates~\eqref{eq:diamagnetic-semigroup}--\eqref{eq:diamagnetic-resolvent}
		follow from the Feynman--Kac representation~\eqref{eq:full-lattice-heat-kernels}
		and the integral formula~\eqref{eq:lattice-covariant-negative-powers-definition}.
		For~\eqref{eq:diamagnetic-gradient}, using
		$\lvert e^{i\varepsilon B_j(x)}\rvert=1$ we obtain
		\begin{equation}
			\label{eq:proof-diamagnetic-gradient}
			\begin{aligned}
				\bigl\lvert D_j^\varepsilon\lvert f\rvert(x)\bigr\rvert
				=
				\frac{\bigl\lvert
					\lvert e^{i\varepsilon B_j(x)}
					f(x+\varepsilon e_j)\rvert-\lvert f(x)\rvert
					\bigr\rvert}{\varepsilon}
				\le
				\frac{
					\bigl\lvert
						e^{i\varepsilon B_j(x)}f(x+\varepsilon e_j)-f(x)
					\bigr\rvert
					}{\varepsilon}
				=
				\bigl\lvert D_{B,j}^\varepsilon f(x)\bigr\rvert.
			\end{aligned}
		\end{equation}
		
		\textup{(2).}
		By~\eqref{eq:diamagnetic-resolvent} and the uniform discrete Sobolev
		embedding,
		\begin{equation}
			\label{eq:covariant-sobolev-via-diamagnetic}
			\bigl\lVert(\mathcal L_B^\varepsilon+1)^{-\delta/2}f
			\bigr\rVert_{L^{q,\varepsilon}}
			\le\bigl\lVert(\mathcal L_0^\varepsilon+1)^{-\delta/2}|f|
			\bigr\rVert_{L^{q,\varepsilon}}
			\lesssim_{\delta,p,q}\lVert f\rVert_{L^{p,\varepsilon}}.
		\end{equation}
		Applying~\eqref{eq:covariant-sobolev-via-diamagnetic} to
		$(\mathcal L_B^\varepsilon+1)^{s/2}(\rho^a f)$ yields
		\eqref{eq:covariant-sobolev-embedding}.

		\textup{(3).}
		By Lemma~\ref{lem:weighted-operator-interpolation} and duality,
		it suffices to prove
		\begin{equation}
			\label{eq:weighted-HB-H-s-two-forward}
			\lVert f\rVert_{H_B^{2,\varepsilon}(\rho^2)}
			\lesssim
			\bigl(1+\lVert B\rVert_{H^{1,\varepsilon}(\rho)}\bigr)^2
			\lVert f\rVert_{H^{2,\varepsilon}},
		\end{equation}
		and
		\begin{equation}
			\label{eq:weighted-HB-H-s-two-reverse}
			\lVert f\rVert_{H^{2,\varepsilon}(\rho^2)}
			\lesssim
			\bigl(1+\lVert B\rVert_{H^{1,\varepsilon}(\rho)}\bigr)^2
			\lVert f\rVert_{H_B^{2,\varepsilon}}.
		\end{equation}
		Substitution of $\rho^a f$ then gives arbitrary input weights.
		We first work with finitely supported $f$ and conclude by density.

		For~\eqref{eq:weighted-HB-H-s-two-forward}, the weight bounds give
		\begin{equation}
			\label{eq:weighted-HB-H-first-endpoint-reduction}
			\bigl\lVert(\mathcal L_B^\varepsilon+1)(\rho^2f)
			\bigr\rVert_{L^{2,\varepsilon}}
			\lesssim
			\bigl\lVert(\mathcal L_B^\varepsilon-\mathcal L_0^\varepsilon)
			(\rho^2f)\bigr\rVert_{L^{2,\varepsilon}}
			+\lVert f\rVert_{H^{2,\varepsilon}}.
		\end{equation}
		By~\eqref{eq:covariant-laplacian-rho-expansion}, the discrete product
		rule, H\"older's inequality, and the uniform embedding
		$H^{1,\varepsilon}\hookrightarrow L^{6,\varepsilon}$,
		\begin{equation}
			\label{eq:weighted-HB-H-linear-terms}
			\begin{aligned}
				\bigl\lVert\rho_B^{\varepsilon *}D^\varepsilon(\rho^2f)
				\bigr\rVert_{L^{2,\varepsilon}}
				+\bigl\lVert D^{\varepsilon *}(\rho_B^\varepsilon\rho^2f)
				\bigr\rVert_{L^{2,\varepsilon}}\lesssim
				\lVert B\rVert_{H^{1,\varepsilon}(\rho)}
				\Bigl(\bigl\lVert\rho^{-1}D^\varepsilon(\rho^2f)
				\bigr\rVert_{L^{3,\varepsilon}}
				+\lVert\rho f\rVert_{L^\infty}\Bigr),
			\end{aligned}
		\end{equation}
		and
		\begin{equation}
			\label{eq:weighted-HB-H-quadratic-term}
			\bigl\lVert\rho_B^{\varepsilon *}\rho_B^\varepsilon(\rho^2f)
			\bigr\rVert_{L^{2,\varepsilon}}
			\lesssim
			\lVert B\rVert_{H^{1,\varepsilon}(\rho)}^2
			\lVert f\rVert_{L^\infty}.
		\end{equation}
		The uniform free Sobolev embeddings give
		$\lVert f\rVert_{L^\infty}\lesssim\lVert f\rVert_{H^{2,\varepsilon}}$
		and
		\begin{equation}
			\label{eq:weighted-HB-H-free-embeddings}
			\bigl\lVert\rho^{-1}D^\varepsilon(\rho^2f)
			\bigr\rVert_{L^{3,\varepsilon}}
			\lesssim\lVert f\rVert_{W_3^{1,\varepsilon}}
			\lesssim\lVert f\rVert_{H^{2,\varepsilon}}.
		\end{equation}
		Thus~\eqref{eq:weighted-HB-H-s-two-forward} follows.

		For~\eqref{eq:weighted-HB-H-s-two-reverse},
		Lemma~\ref{lem:covariant-nonnegative-weight} yields
		\begin{equation}
			\label{eq:weighted-HB-H-second-endpoint-reduction}
			\begin{aligned}
				\lVert f\rVert_{H^{2,\varepsilon}(\rho^2)}
				\lesssim\lVert f\rVert_{H_B^{2,\varepsilon}}
				+\bigl\lVert(\mathcal L_B^\varepsilon-\mathcal L_0^\varepsilon)
				(\rho^2f)\bigr\rVert_{L^{2,\varepsilon}}.
			\end{aligned}
		\end{equation}
		By~\eqref{eq:weighted-HB-H-linear-terms} and
		\eqref{eq:weighted-HB-H-quadratic-term},
		\begin{equation}
			\label{eq:weighted-HB-H-perturbation-bound}
			\begin{aligned}
				\bigl\lVert(\mathcal L_B^\varepsilon-\mathcal L_0^\varepsilon)
				(\rho^2f)\bigr\rVert_{L^{2,\varepsilon}}\lesssim
				\bigl(1+\lVert B\rVert_{H^{1,\varepsilon}(\rho)}\bigr)^2
				\lVert f\rVert_{L^\infty}
				+\lVert B\rVert_{H^{1,\varepsilon}(\rho)}
				\bigl\lVert\rho^{-1}D^\varepsilon(\rho^2f)
				\bigr\rVert_{L^{3,\varepsilon}}.
			\end{aligned}
		\end{equation}
		Taking $p=2$, $q\in\{6,\infty\}$, $s=\delta=2$, and $a=0$
		in~\eqref{eq:covariant-sobolev-embedding} gives
$
			\label{eq:weighted-HB-H-covariant-linfty}
			\lVert f\rVert_{L^\infty}
			+\lVert f\rVert_{L^{6,\varepsilon}}
			\lesssim\lVert f\rVert_{H_B^{2,\varepsilon}}.
$
		By the weight bounds and  the weighted Gagliardo--Nirenberg inequality,
		\begin{equation}
			\label{eq:weighted-HB-H-weighted-GN}
			\begin{aligned}
				\bigl\lVert\rho^{-1}D^\varepsilon(\rho^2f)
				\bigr\rVert_{L^{3,\varepsilon}}
			\lesssim\lVert\rho f\rVert_{W_3^{1,\varepsilon}}
				\lesssim\lVert f\rVert_{H^{2,\varepsilon}(\rho^2)}^{1/2}
				\lVert f\rVert_{L^{6,\varepsilon}}^{1/2}
				\lesssim\lVert f\rVert_{H^{2,\varepsilon}(\rho^2)}^{1/2}
				\lVert f\rVert_{H_B^{2,\varepsilon}}^{1/2}.
			\end{aligned}
		\end{equation}
		Combining these estimates and applying Young's inequality,
		\begin{equation}
			\label{eq:weighted-HB-H-absorption}
			\begin{aligned}
				\lVert f\rVert_{H^{2,\varepsilon}(\rho^2)}
				&\le C\bigl(1+\lVert B\rVert_{H^{1,\varepsilon}(\rho)}\bigr)^2
				\lVert f\rVert_{H_B^{2,\varepsilon}}
				+C\lVert B\rVert_{H^{1,\varepsilon}(\rho)}
				\lVert f\rVert_{H^{2,\varepsilon}(\rho^2)}^{1/2}
				\lVert f\rVert_{H_B^{2,\varepsilon}}^{1/2}\\
				&\le\tfrac12\lVert f\rVert_{H^{2,\varepsilon}(\rho^2)}
				+C'\bigl(1+\lVert B\rVert_{H^{1,\varepsilon}(\rho)}\bigr)^2
				\lVert f\rVert_{H_B^{2,\varepsilon}}.
			\end{aligned}
		\end{equation}
		Absorption yields~\eqref{eq:weighted-HB-H-s-two-reverse}.
	\end{proof}

	\begin{lemma}[Covariant derivative estimates]
		\label{lem:covariant-derivative-estimates}
		\label{lem:extended-covariant-derivative-estimates}
		For $\beta\in[-1,2]$, $a\in\mathbb{R}$, and $j\in\{1,2\}$,
		put $d_\beta:=\operatorname{dist}(\beta,[0,1])$.
		If $B\in H^{1,\varepsilon}(\rho)$, then
		\begin{equation}
			\label{eq:extended-covariant-derivative-estimates}
			\csname ltx@label\endcsname{eq:covariant-derivative-estimates}
			\begin{aligned}
				\lVert D_{B,j}^\varepsilon f\rVert_{H_B^{\beta-1,\varepsilon}(\rho^{a+2d_\beta})}
				+\lVert D_{B,j}^{\varepsilon *}f\rVert_{H_B^{\beta-1,\varepsilon}(\rho^{a+2d_\beta})}\lesssim_{a,\beta}
				\bigl(1+\lVert B\rVert_{H^{1,\varepsilon}(\rho)}\bigr)^{2d_\beta}
				\lVert f\rVert_{H_B^{\beta,\varepsilon}(\rho^a)}.
			\end{aligned}
		\end{equation}
	\end{lemma}
	\begin{proof}
		We first work with finitely supported $f$ and give the proof for
		$D_B^\varepsilon$; the adjoint is treated in the same way.
		For $\beta\in[0,1]$, interpolation and duality reduce the proof to
		$\beta=1$. The identity
		$D_{B,j}^{\varepsilon *}D_{B,j}^\varepsilon
		=D_{B,j}^\varepsilon D_{B,j}^{\varepsilon *}
		\le\mathcal L_B^\varepsilon$
		and~\eqref{eq:covariant-leibniz-rules} give
		\begin{equation}
			\label{eq:covariant-derivative-weighted-energy-bound}
			\begin{aligned}
				\lVert D_B^\varepsilon f\rVert_{L^{2,\varepsilon}(\rho^a)}^2
				&\lesssim
				\lVert D_B^\varepsilon(\rho^af)\rVert_{L^{2,\varepsilon}}^2
				+\lVert[D_B^\varepsilon,\rho^a]f\rVert_{L^{2,\varepsilon}}^2\\
				&\lesssim_a
				\lVert\rho^af\rVert_{H_B^{1,\varepsilon}}^2
				+\lVert\rho^af\rVert_{L^{2,\varepsilon}}^2
				\lesssim_a\lVert f\rVert_{H_B^{1,\varepsilon}(\rho^a)}^2.
			\end{aligned}
		\end{equation}
		The same bound holds componentwise for $D_{B,j}^{\varepsilon *}$.

		Interpolation and duality reduce the remaining range to $\beta=2$.
		Put $N:=1+\lVert B\rVert_{H^{1,\varepsilon}(\rho)}$.
		We first take $a=0$.
		The definition of $H_B^{1,\varepsilon}$ and the preceding operator
		identity yield
		\begin{equation}
			\label{eq:covariant-mixed-second-derivative-reduction}
			\begin{aligned}
				\lVert D_{B,j}^\varepsilon f\rVert_{H_B^{1,\varepsilon}(\rho^2)}
				\lesssim
				\lVert D_{B,j}^\varepsilon f\rVert_{L^{2,\varepsilon}(\rho^2)}
				+\sum_{k=1}^2
				\lVert D_{B,k}^{\varepsilon *}D_{B,j}^\varepsilon f
				\rVert_{L^{2,\varepsilon}(\rho^2)}+\sum_{k=1}^2
				\lVert[D_{B,k}^{\varepsilon *},\rho^2]D_{B,j}^\varepsilon f
				\rVert_{L^{2,\varepsilon}}.
			\end{aligned}
		\end{equation}
		Only the second-derivative terms and the weight commutators need estimates.
		For $k,j\in\{1,2\}$, we expand
		\begin{equation}
			\label{eq:covariant-mixed-second-derivative-expansion}
			\begin{aligned}
				D_{B,k}^{\varepsilon *}D_{B,j}^\varepsilon
				&=\bigl(D_k^{\varepsilon *}-i\rho_{B,k}^{\varepsilon *}\bigr)
				\bigl(D_j^\varepsilon+i\rho_{B,j}^\varepsilon\bigr)\\
				&=D_k^{\varepsilon *}D_j^\varepsilon
				-i\rho_{B,k}^{\varepsilon *}D_j^\varepsilon
				+iD_k^{\varepsilon *}\rho_{B,j}^\varepsilon
				+\rho_{B,k}^{\varepsilon *}\rho_{B,j}^\varepsilon.
			\end{aligned}
		\end{equation}
	The discrete product
		rule, H\"older's inequality, and Lemma~\ref{lem:WpB-Wp-comparison} give
		\begin{align}
			\lVert D_k^{\varepsilon *}D_j^\varepsilon f
			\rVert_{L^{2,\varepsilon}(\rho^2)}
			&\lesssim\lVert f\rVert_{H^{2,\varepsilon}(\rho^2)}
			\lesssim N^2\lVert f\rVert_{H_B^{2,\varepsilon}},
			\label{eq:covariant-mixed-second-derivative-free-term}\\
			\lVert\rho_{B,k}^{\varepsilon *}D_j^\varepsilon f
			\rVert_{L^{2,\varepsilon}(\rho^2)}
			&\lesssim\lVert B\rVert_{L^{6,\varepsilon}(\rho)}
			\lVert D_j^\varepsilon f\rVert_{L^{3,\varepsilon}(\rho)}
			\lesssim\lVert B\rVert_{H^{1,\varepsilon}(\rho)}N
			\lVert f\rVert_{H_B^{2,\varepsilon}},
			\label{eq:covariant-mixed-second-derivative-linear-term}\\
			\lVert D_k^{\varepsilon *}(\rho_{B,j}^\varepsilon f)
			\rVert_{L^{2,\varepsilon}(\rho^2)}
			&\lesssim\lVert B\rVert_{H^{1,\varepsilon}(\rho)}\lVert f\rVert_\infty
			+\lVert B\rVert_{L^{6,\varepsilon}(\rho)}
			\lVert D_k^\varepsilon f\rVert_{L^{3,\varepsilon}(\rho)}\notag\\
			&\lesssim\lVert B\rVert_{H^{1,\varepsilon}(\rho)}N
			\lVert f\rVert_{H_B^{2,\varepsilon}},
			\label{eq:covariant-mixed-second-derivative-product-term}\\
			\lVert\rho_{B,k}^{\varepsilon *}\rho_{B,j}^\varepsilon f
			\rVert_{L^{2,\varepsilon}(\rho^2)}
			&\lesssim\lVert B\rVert_{H^{1,\varepsilon}(\rho)}^2
			\lVert f\rVert_{L^{6,\varepsilon}}
			\lesssim\lVert B\rVert_{H^{1,\varepsilon}(\rho)}^2
			\lVert f\rVert_{H_B^{2,\varepsilon}}.
			\label{eq:covariant-mixed-second-derivative-quadratic-term}
		\end{align}
		Here we used the weighted Gagliardo--Nirenberg inequality:
		\begin{equation}
			\label{eq:covariant-derivative-weighted-gagliardo-nirenberg}
			\lVert D_j^\varepsilon f\rVert_{L^{3,\varepsilon}(\rho)}
			\lesssim\lVert f\rVert_{H^{2,\varepsilon}(\rho^2)}^{1/2}
			\lVert f\rVert_{L^{6,\varepsilon}}^{1/2}
			\lesssim N\lVert f\rVert_{H_B^{2,\varepsilon}}.
		\end{equation}
		For the weight commutators, \eqref{eq:covariant-leibniz-rules} gives
		\begin{equation}
			\label{eq:covariant-mixed-second-derivative-weight-commutator}
			\lVert[D_{B,k}^{\varepsilon *},\rho^2]D_{B,j}^\varepsilon f
			\rVert_{L^{2,\varepsilon}}
			\lesssim
			\lVert D_{B,j}^\varepsilon f\rVert_{L^{2,\varepsilon}(\rho^2)}
			\lesssim\lVert f\rVert_{H_B^{1,\varepsilon}}.
		\end{equation}
		Combining~\eqref{eq:covariant-mixed-second-derivative-reduction}
		with the four estimates gives
		\begin{equation}
			\label{eq:covariant-derivative-covariant-upper-endpoint}
			\lVert\rho^2D_{B,j}^\varepsilon f\rVert_{H_B^{1,\varepsilon}}
			+\lVert\rho^2D_{B,j}^{\varepsilon *}f\rVert_{H_B^{1,\varepsilon}}
			\lesssim N^2\lVert f\rVert_{H_B^{2,\varepsilon}}.
		\end{equation}

		For general $a$, the covariant Leibniz rule gives
		\begin{equation}
			\label{eq:covariant-derivative-weight-transfer-reduction}
			\begin{aligned}
				\lVert D_{B,j}^\varepsilon f\rVert_{H_B^{1,\varepsilon}(\rho^{a+2})}
				\le\bigl\lVert D_{B,j}^\varepsilon\bigl((\tau_j^*\rho^a)f\bigr)
				\bigr\rVert_{H_B^{1,\varepsilon}(\rho^2)}
				+\lVert(D_j^{\varepsilon *}\rho^a)f
				\rVert_{H_B^{1,\varepsilon}(\rho^2)}.
			\end{aligned}
		\end{equation}
		For the first term, apply the case $a=0$, since
		$\lVert(\tau_j^*\rho^a)/\rho^a\rVert_{C^2}\lesssim_a1$.
		For the second term, the identity
		\begin{equation}
			\label{eq:covariant-derivative-weight-transfer}
			\rho^aD_{B,j}^\varepsilon f
			=D_{B,j}^\varepsilon\bigl((\tau_j^*\rho^a)f\bigr)
			+(D_j^{\varepsilon *}\rho^a)f
		\end{equation}
		and Lemma~\ref{lem:covariant-nonnegative-weight} give
		\begin{equation}
			\label{eq:covariant-derivative-weight-transfer-remainder}
			\lVert(D_j^{\varepsilon *}\rho^a)f
			\rVert_{H_B^{1,\varepsilon}(\rho^2)}
			\lesssim_a\lVert\rho^af\rVert_{H_B^{1,\varepsilon}(\rho^2)}
			\lesssim_a\lVert f\rVert_{H_B^{2,\varepsilon}(\rho^a)}.
		\end{equation}
	\end{proof}

	\begin{corollary}[Free Sobolev and divergence estimates]
		\label{cor:weighted-covariant-derivative-and-divergence}
		Assume $B\in H^{1,\varepsilon}(\rho)$.
		For $s\in[-2,1]$, put
		$\vartheta_s:=\max\{-3s-2,-s,3s\}$.
		Then, for $a\in\mathbb{R}$ and $j\in\{1,2\}$,
		\begin{equation}
			\label{eq:weighted-covariant-derivative-free-output}
			\begin{aligned}
				\lVert D_{B,j}^\varepsilon f\rVert_{H^{s,\varepsilon}(\rho^{a+\vartheta_s})}
				+\lVert D_{B,j}^{\varepsilon *}f\rVert_{H^{s,\varepsilon}(\rho^{a+\vartheta_s})}\lesssim_{a,s}
				\bigl(1+\lVert B\rVert_{H^{1,\varepsilon}(\rho)}\bigr)^{\vartheta_s}
				\lVert f\rVert_{H_B^{1+s,\varepsilon}(\rho^a)}.
			\end{aligned}
		\end{equation}
		For $s\in[-1,1]$ and vector fields $F,G$,
		\begin{equation}
			\label{eq:weighted-covariant-negative-divergence}
			\begin{aligned}
				\lVert D_B^{\varepsilon *}(F+U_B^\varepsilon G)
				\rVert_{H_B^{-1-s,\varepsilon}(\rho^a)}\lesssim_{a,s}
				\bigl(1+\lVert B\rVert_{H^{1,\varepsilon}(\rho)}\bigr)^{\vartheta_s}
				\bigl(\lVert F\rVert_{H^{-s,\varepsilon}(\rho^{a-\vartheta_s})}
				+\lVert G\rVert_{H^{-s,\varepsilon}(\rho^{a-\vartheta_s})}\bigr).
			\end{aligned}
		\end{equation}
	\end{corollary}
	\begin{proof}
		Put $N:=1+\lVert B\rVert_{H^{1,\varepsilon}(\rho)}$ and $\beta:=1+s$.
		By~\eqref{eq:weighted-HB-H-comparison} and
		\eqref{eq:extended-covariant-derivative-estimates},
		\begin{equation}
			\label{eq:weighted-covariant-derivative-free-output-proof}
			\begin{aligned}
				\lVert D_{B,j}^\varepsilon f\rVert_{H^{s,\varepsilon}(\rho^{a+\vartheta_s})}
				\lesssim N^{|s|}
				\lVert D_{B,j}^\varepsilon f\rVert_{H_B^{s,\varepsilon}(\rho^{a+\vartheta_s-|s|})}
				\lesssim N^{2d_\beta+|s|}
				\lVert f\rVert_{H_B^{1+s,\varepsilon}(\rho^a)}.
			\end{aligned}
		\end{equation}
		Here $\vartheta_s=|s|+2d_\beta$, since
		$d_\beta=\operatorname{dist}(s,[-1,0])$.
		The same argument applies to $D_{B,j}^{\varepsilon *}$.

		For the divergence estimate, put $\tilde\beta:=-s$.
		By~\eqref{eq:covariant-derivative-translation-identities},
		\eqref{eq:extended-covariant-derivative-estimates} and 
		\eqref{eq:weighted-HB-H-comparison} give
		\begin{equation}
			\label{eq:weighted-covariant-negative-divergence-proof}
			\begin{aligned}
				\lVert D_B^{\varepsilon *}(F+U_B^\varepsilon G)
				\rVert_{H_B^{-1-s,\varepsilon}(\rho^a)}
				&\lesssim N^{2d_{\tilde\beta}}\sum_{j=1}^2
				\Bigl(\lVert F_j\rVert_{H_B^{-s,\varepsilon}(\rho^{a-2d_{\tilde\beta}})}
				+\lVert\tau_j^*G_j\rVert_{H_B^{-s,\varepsilon}(\rho^{a-2d_{\tilde\beta}})}\Bigr)\\
				&\lesssim N^{2d_{\tilde\beta}+|s|}
				\Bigl(\lVert F\rVert_{H^{-s,\varepsilon}(\rho^{a-2d_{\tilde\beta}-|s|})}
				+\lVert G\rVert_{H^{-s,\varepsilon}(\rho^{a-2d_{\tilde\beta}-|s|})}\Bigr).
			\end{aligned}
		\end{equation}
		Since $2d_{\tilde\beta}+|s|=\vartheta_s$ for $s\in[-1,1]$,
		this yields~\eqref{eq:weighted-covariant-negative-divergence}.
	\end{proof}
	\begin{lemma}[Weighted energy identity]
		\label{lem:discrete-ims-identity}
		\label{cor:weighted-covariant-coercivity}
		Let $B\in(\mathbb{R}^2)^{\Lambda_\varepsilon}$, let $w$ be a polynomial
		weight, and let $u\in L^{2,\varepsilon}(w)$.
 Then
		\begin{equation}
			\label{eq:discrete-ims-identity}
			\langle\mathcal{L}_B^\varepsilon u,w^2u\rangle_\varepsilon
			=\lVert D_B^\varepsilon(wu)\rVert_{L^{2,\varepsilon}}^2
			-\sum_{j=1}^2
			\bigl\langle |D_j^\varepsilon w|^2,
			\overline{u}\,U_{B,j}^\varepsilon\tau_ju
			\bigr\rangle_\varepsilon.
		\end{equation}
		In particular, if $w=\rho^a$ and $|a|\le l_0$, then by the choice of
		$\nu$ in Section~\ref{subsec:parameters},
		\begin{equation}
			\label{eq:weighted-covariant-coercivity}
			\frac12\lVert\rho^a u\rVert_{H_B^{1,\varepsilon}}^2
			\le\bigl\langle(\mathcal{L}_B^\varepsilon+1)u,
			\rho^{2a}u\bigr\rangle_\varepsilon
			\lesssim\lVert\rho^a u\rVert_{H_B^{1,\varepsilon}}^2.
		\end{equation}
	\end{lemma}
	\begin{proof}
		By Cauchy--Schwarz, all sums below converge absolutely.
		For $j=1,2$, the Leibniz rule~\eqref{eq:covariant-leibniz-rules} gives
		\begin{equation}
			\label{eq:discrete-ims-product-rule}
			\langle D_{B,j}^\varepsilon u,
			D_{B,j}^\varepsilon(w^2u)\rangle_\varepsilon
			=\langle D_{B,j}^\varepsilon u,
			wD_{B,j}^\varepsilon(wu)\rangle_\varepsilon
			+\langle D_{B,j}^\varepsilon u,
			(D_j^\varepsilon w)U_{B,j}^\varepsilon\tau_j(wu)
			\rangle_\varepsilon.
		\end{equation}
		Since $w$ is real-valued, another application
		of~\eqref{eq:covariant-leibniz-rules} yields
		\begin{equation}
			\label{eq:discrete-ims-first-term}
			\langle D_{B,j}^\varepsilon u,
			wD_{B,j}^\varepsilon(wu)\rangle_\varepsilon
			=\lVert D_{B,j}^\varepsilon(wu)\rVert_{L^{2,\varepsilon}}^2
			-\bigl\langle
			(D_j^\varepsilon w)U_{B,j}^\varepsilon\tau_ju,
			D_{B,j}^\varepsilon(wu)
			\bigr\rangle_\varepsilon.
		\end{equation}
		Thus
		\begin{equation}
			\label{eq:discrete-ims-cancellation}
			\begin{aligned}
				\langle D_{B,j}^\varepsilon u,
				D_{B,j}^\varepsilon(w^2u)\rangle_\varepsilon
				-\lVert D_{B,j}^\varepsilon(wu)\rVert_{L^{2,\varepsilon}}^2
				&=\bigl\langle
				(\tau_jw)D_{B,j}^\varepsilon u-D_{B,j}^\varepsilon(wu),
				(D_j^\varepsilon w)U_{B,j}^\varepsilon\tau_ju
				\bigr\rangle_\varepsilon\\
				&\overset{\eqref{eq:covariant-leibniz-rules}}{=}
				-\bigl\langle
				(D_j^\varepsilon w)u,
				(D_j^\varepsilon w)U_{B,j}^\varepsilon\tau_ju
				\bigr\rangle_\varepsilon.
			\end{aligned}
		\end{equation}
		Summing over $j=1,2$ and using~\eqref{eq:lattice-covariant-laplacian}
		yields~\eqref{eq:discrete-ims-identity}.
		For $w=\rho^a$, $|a|\le l_0$, and $j=1,2$,
		\eqref{eq:weight-basic-estimates} and Cauchy--Schwarz give
		\begin{equation}
			\label{eq:ims-polynomial-weight-error}
			\left|
			\Bigl\langle |D_j^\varepsilon\rho^a|^2,
			\overline{u}\,U_{B,j}^\varepsilon\tau_ju
			\Bigr\rangle_\varepsilon
			\right|
			\lesssim\nu\lVert\rho^a u\rVert_{L^{2,\varepsilon}}^2.
		\end{equation}
		Applying~\eqref{eq:discrete-ims-identity} with $w=\rho^a$ gives
		the upper bound in~\eqref{eq:weighted-covariant-coercivity}.
		Absorbing the error into the mass term gives the lower bound
		for $\nu>0$ sufficiently small.
	\end{proof}

	\begin{lemma}[Weighted covariant resolvent estimate]
		\label{lem:weighted-covariant-resolvent}
		For $s\in[0,2]$, $|a|\le l_0$, and
		$f\in H_B^{s-2,\varepsilon}(\rho^a)$,
		\begin{equation}
			\label{eq:weighted-covariant-resolvent}
			\bigl\lVert(\mathcal{L}_B^\varepsilon+1)^{-1}f
			\bigr\rVert_{H_B^{s,\varepsilon}(\rho^a)}
			\lesssim\lVert f\rVert_{H_B^{s-2,\varepsilon}(\rho^a)}.
		\end{equation}
	\end{lemma}
	\begin{proof}
	 Set
		$u:=(\mathcal{L}_B^\varepsilon+1)^{-1}f$.
		For $1\le s\le2$, operator interpolation reduces the proof
		to $s=1,2$.

		For $s=1$, \eqref{eq:weighted-covariant-coercivity}, duality, and
		Young's inequality give
		\begin{equation}
			\label{eq:weighted-covariant-resolvent-s-zero}
			\begin{aligned}
				&\frac12\lVert u\rVert_{H_B^{1,\varepsilon}(\rho^a)}^2
				\le\bigl\langle(\mathcal{L}_B^\varepsilon+1)u,
				\rho^{2a}u\bigr\rangle_\varepsilon
				=\langle f,\rho^{2a}u\rangle_\varepsilon\\
				&\qquad\le\lVert f\rVert_{H_B^{-1,\varepsilon}(\rho^a)}
				\lVert u\rVert_{H_B^{1,\varepsilon}(\rho^a)}
				\le\frac14\lVert u\rVert_{H_B^{1,\varepsilon}(\rho^a)}^2
				+16\lVert f\rVert_{H_B^{-1,\varepsilon}(\rho^a)}^2.
			\end{aligned}
		\end{equation}
		This gives the bound for $s=1$.

		For $s=2$,
		\begin{equation}
			\label{eq:weighted-covariant-resolvent-s-one-reduction}
			\begin{aligned}
				\lVert u\rVert_{H_B^{2,\varepsilon}(\rho^a)}^2
				=\bigl\lVert(\mathcal{L}_B^\varepsilon+1)(\rho^a u)
				\bigr\rVert_{L^{2,\varepsilon}}^2
				\lesssim\lVert f\rVert_{L^{2,\varepsilon}(\rho^a)}^2
				+\bigl\lVert[\mathcal{L}_B^\varepsilon,\rho^a]u
				\bigr\rVert_{L^{2,\varepsilon}}^2.
			\end{aligned}
		\end{equation}
		Since $D_{B,j}^{\varepsilon *}
		=-\tau_j^*U_{B,j}^{\varepsilon *}D_{B,j}^\varepsilon$,
		\eqref{eq:weight-basic-estimates}, Lemma~\ref{lem:covariant-laplacian-weight-commutator},
		and~\eqref{eq:covariant-leibniz-rules} yield
		\begin{equation}
			\label{eq:weighted-covariant-resolvent-commutator}
			\begin{aligned}
				&\bigl\lVert[\mathcal{L}_B^\varepsilon,\rho^a]u
				\bigr\rVert_{L^{2,\varepsilon}}^2
				\lesssim\lVert\rho^a u\rVert_{L^{2,\varepsilon}}^2
				+\lVert\rho^a D_B^\varepsilon u\rVert_{L^{2,\varepsilon}}^2\\
				&\qquad\lesssim\lVert u\rVert_{H_B^{1,\varepsilon}(\rho^a)}^2
				+\sum_{j=1}^2\bigl\lVert(D_j^\varepsilon\rho^a)
				U_{B,j}^\varepsilon\tau_ju\bigr\rVert_{L^{2,\varepsilon}}^2
				\lesssim\lVert u\rVert_{H_B^{1,\varepsilon}(\rho^a)}^2,
			\end{aligned}
		\end{equation}
		where~\eqref{eq:weight-basic-estimates} controls the last sum.
		Using the $s=1$ estimate and
		$\lVert f\rVert_{H_B^{-1,\varepsilon}(\rho^a)}
		\le\lVert f\rVert_{L^{2,\varepsilon}(\rho^a)}$, we conclude
		\begin{equation}
			\label{eq:weighted-covariant-resolvent-s-one}
			\lVert u\rVert_{H_B^{2,\varepsilon}(\rho^a)}^2
			\lesssim\lVert f\rVert_{L^{2,\varepsilon}(\rho^a)}^2
			+\lVert u\rVert_{H_B^{1,\varepsilon}(\rho^a)}^2
			\lesssim\lVert f\rVert_{L^{2,\varepsilon}(\rho^a)}^2.
		\end{equation}
		For $0\le s\le1$, the dual of $H_B^{r,\varepsilon}(\rho^a)$ under
		the unweighted lattice pairing is $H_B^{-r,\varepsilon}(\rho^{-a})$.
		By self-adjointness and the estimate at $2-s$ with weight $\rho^{-a}$,
		\begin{equation}
			\label{eq:weighted-covariant-resolvent-negative-s}
			\begin{aligned}
			\lVert u\rVert_{H_B^{s,\varepsilon}(\rho^a)}
			=\sup_{\lVert g\rVert_{H_B^{-s,\varepsilon}(\rho^{-a})}\le1}
			\bigl|\langle f,(\mathcal L_B^\varepsilon+1)^{-1}g
			\rangle_\varepsilon\bigr|
			\lesssim\lVert f\rVert_{H_B^{s-2,\varepsilon}(\rho^a)}.
			\end{aligned}
		\end{equation}
	\end{proof}

	\section{The lattice Abelian--Higgs model and the associated iterated FBSDEs}
	\label{sec:lattice-abelian-higgs-and-variational-problems}
	
	\subsection{Lattice Abelian--Higgs model}

	In this subsection, we make precise the lattice model and Coulomb gauge
	fixing outlined in the Introduction; the Gaussian reference fields and
	renormalization constants will be specified in
	Subsection~\ref{subsec:gff-pure-abelian-yang-mills}. Using the notation of
	Section~\ref{sec:preliminaries}, we first consider the direct lattice analogue
	of~\eqref{eq:intro-formal-abelian-higgs-measure}. Although the configuration
	space is finite-dimensional, gauge redundancy makes the partition function
	infinite, so the following formal expression does not define a probability
	measure:
	\begin{equation}
		\label{eq:naive-lattice-measure}
		\mu_{M,\varepsilon}^{\mathrm{naive}}(dA,d\phi)
		\overset{\mathrm{formally}}{\propto} \exp\bigl(-S_{M,\varepsilon}(A,\phi)\bigr)
		\prod_{x\in\Lambda_{M,\varepsilon}}dA(x)\,d\phi(x).
	\end{equation}
	Here $dA(x)\,d\phi(x)$ is Lebesgue measure on
	$\mathbb{R}^2\times\mathbb{C}$, and the lattice action is
	\begin{equation}
		\label{eq:lattice-abelian-higgs-action}
		S_{M,\varepsilon}(A,\phi)
		:=
		\int_{\Lambda_{M,\varepsilon}}
		\left(\frac12\lvert F_A^\varepsilon(x)\rvert^2
		+\frac12\lvert D_A^\varepsilon\phi(x)\rvert^2
		+\frac14\lvert\phi(x)\rvert^4\right)\,dx.
	\end{equation}
	Here $F_{A,jk}^\varepsilon
	:=\partial_j^\varepsilon A_k-
	\partial_k^\varepsilon A_j$,
	$\lvert F_A^\varepsilon\rvert^2
	:=\frac12\sum_{j,k=1}^2
	(F_{A,jk}^\varepsilon)^2$, and
	$\lvert D_A^\varepsilon\phi\rvert^2
	:=\sum_{j=1}^2
	\lvert D_{A,j}^\varepsilon\phi\rvert^2$.
	
	The lattice gauge transformation group is
	\begin{equation}
		\label{eq:lattice-gauge-group}
		\mathcal{G}_{M,\varepsilon}
		:=
		\bigl\{g\colon\varepsilon\mathbb{Z}^2\to\mathbb{R}:
		D^\varepsilon g \text{ and } \exp(-ig)
		\text{ are }2\pi M\mathbb{Z}^2\text{-periodic}\bigr\}.
	\end{equation}
	Allowing $g$ to be defined on $\varepsilon\mathbb{Z}^2$ rather than
	$\Lambda_{M,\varepsilon}$ includes the gauge transformations
	$g(x)=n\mathbin{\cdot}x$, $n\in M^{-1}\mathbb{Z}^2$, which are used to remove the
	zero-frequency gauge redundancy. The lattice analogues of
	\eqref{eq:intro-continuum-gauge-transformation} and
	\eqref{eq:intro-continuum-gauge-invariance} are as follows. For
	$g\in\mathcal{G}_{M,\varepsilon}$, define
	$(A^g,\phi^g)
	:=
	(A+\nabla^\varepsilon g,\exp(-ig)\phi)$. Then
	$D_{A^g,j}^\varepsilon\phi^g(x)
	=\exp(-ig(x))D_{A,j}^\varepsilon\phi(x)$,
	$F_{A^g}^\varepsilon=F_A^\varepsilon$, and
	$S_{M,\varepsilon}(A^g,\phi^g)=S_{M,\varepsilon}(A,\phi)$.
	The use of the exponential link $e^{i\varepsilon A_j(x)}$ in
	\eqref{eq:prelim-covariant-derivative} is essential: the linear approximation
	$(\phi(x+\varepsilon e_j)-\phi(x))/\varepsilon+iA_j(x)\phi(x)$ is not gauge
	covariant; see~\cite[Remark~2.2]{S21}.
	
	We next choose a Coulomb-gauge representative of each gauge-field
	orbit. Recall from~\eqref{eq:lattice-pinning-projection} that
	$P_{\chi_\varepsilon}^{M}A=A-\overline A^{\chi_\varepsilon}$.
	For fixed $M$, one may work with the zero-mean Gaussian gauge field and
	treat the zero Fourier modes separately, as in \cite[Section~1.1]{BCCKM26},
	without applying
	$P_{\chi_\varepsilon}^{M}$. We use pinning in preparation for the
	infinite-volume limit, since the unpinned field exhibits an infrared
	divergence near zero frequency; see
	Remark~\ref{rem:pure-yang-mills-infrared-divergence}.
	
	\begin{lemma}[Coulomb-gauge representative]
		\label{lem:canonical-coulomb-gauge-representative}
		Let $A\colon\Lambda_{M,\varepsilon}\to\mathbb{R}^2$ have the pinned
		Helmholtz parametrization
		\begin{equation}
		A=\nabla^\varepsilon f+P_{\chi_\varepsilon}^{M}\nabla_\perp^\varepsilon\eta+h,
		\qquad f,\eta\in\mathbb{R}_0^{\Lambda_{M,\varepsilon}},
		\quad h\in\mathbb{R}^2,
		\end{equation}
		obtained from~\eqref{eq:lattice-helmholtz-decomposition} by translating
		the constant coordinate. There exists a unique
		$g_A\in\mathcal{G}_{M,\varepsilon}$ satisfying $g_A(0)=0$ such that
		\begin{equation}
			\label{eq:canonical-coulomb-representative}
			A^{g_A}=P_{\chi_\varepsilon}^{M}\nabla_\perp^\varepsilon\eta+h_0,
			\qquad h_0\in[0,M^{-1})^2.
		\end{equation}
		Without the normalization $g_A(0)=0$, the gauge function is unique only
		up to an additive constant.
	\end{lemma}

	\begin{proof}
		Write uniquely $h=h_0+n$, where $h_0\in[0,M^{-1})^2$ and
		$n\in M^{-1}\mathbb{Z}^2$. Using the periodic extension of $f$, define
$
			\label{eq:canonical-coulomb-gauge-function-explicit}
			g_A(x):=-f(x)-n\mathbin{\cdot}x+f(0).
$
		The periodicity of $f$ and the fact that $2\pi M n\in2\pi\mathbb{Z}^2$ show that
		$g_A\in\mathcal{G}_{M,\varepsilon}$. Moreover, $g_A(0)=0$ and
		$A^{g_A}=P_{\chi_\varepsilon}^{M}\nabla_\perp^\varepsilon\eta+h_0$, proving existence.
		
		For uniqueness, let $q=g_1-g_2$. Then
		$D^\varepsilon q=h_{0,1}-h_{0,2}\in(-M^{-1},M^{-1})^2$. Since
		$D^\varepsilon q$ is constant and $e^{-iq}$ is periodic, one also has
		$D^\varepsilon q\in M^{-1}\mathbb{Z}^2$, and hence $D^\varepsilon q=0$.
		Thus $q$ is constant, and $q(0)=0$ gives $q=0$. 
	\end{proof}

 Let $F\geq0$ be a bounded measurable
	gauge-invariant observable on
	$(\mathbb{R}^2\times\mathbb{C})^{\Lambda_{M,\varepsilon}}$, that is,
	\begin{equation}
		\label{eq:gauge-invariant-lattice-observable}
		F(A^g,\phi^g)=F(A,\phi),
		\qquad g\in\mathcal{G}_{M,\varepsilon}.
	\end{equation}
	By the uniqueness of the Helmholtz decomposition, every
	$A\in(\mathbb{R}^2)^{\Lambda_{M,\varepsilon}}$ is represented uniquely by
	$(f,\eta,h)\in
	\mathbb{R}^{\Lambda_{M,\varepsilon}}_0
	\times\mathbb{R}^{\Lambda_{M,\varepsilon}}_0
	\times\mathbb{R}^2$ through
$
	A=D^\varepsilon f+P_{\chi_\varepsilon}^{M}\nabla_\perp^\varepsilon\eta+h.
$	Lebesgue measure on each copy of
	$\mathbb{R}^{\Lambda_{M,\varepsilon}}_0$ is understood with
	respect to the inner product in
	~\eqref{eq:lattice-inner-product-and-adjoints}. Writing $h=h_0+n$ uniquely with
	$h_0\in[0,M^{-1})^2$ and $n\in M^{-1}\mathbb{Z}^2$ gives the corresponding coordinates
	$(f,\eta,h_0,n)$.
	Consequently, for a constant $C_{M,\varepsilon}>0$,
	\begin{equation}
		\label{eq:lattice-helmholtz-lebesgue-decomposition}
		dA=C_{M,\varepsilon}\,df\,d\eta\,dh
		=C_{M,\varepsilon}\,df\,d\eta\,dh_0\,dn,
	\end{equation}
	where $dh_0$ is Lebesgue measure on $[0,M^{-1})^2$ and $dn$ is counting
	measure on $M^{-1}\mathbb{Z}^2$.
	
	Let $g_A$ be the gauge function from
	Lemma~\ref{lem:canonical-coulomb-gauge-representative}; thus
	$A^{g_A}=P_{\chi_\varepsilon}^{M}\nabla_\perp^\varepsilon\eta+h_0$.
	Gauge invariance of $F$ and $S_{M,\varepsilon}$ gives
	\begin{equation}
		\label{eq:gauge-invariant-integrand-reduction}
		\bigl(Fe^{-S_{M,\varepsilon}}\bigr)(A,\phi)
		=\bigl(Fe^{-S_{M,\varepsilon}}\bigr)(A^{g_A},\phi^{g_A})
		=\bigl(Fe^{-S_{M,\varepsilon}}\bigr)
		\bigl(P_{\chi_\varepsilon}^{M}\nabla_\perp^\varepsilon\eta+h_0,\phi^{g_A}\bigr).
	\end{equation}
	For fixed $A$, multiplication by $e^{-ig_A(x)}$ is unitary at every
	vertex, so $\phi\mapsto\phi^{g_A}$ preserves the product Lebesgue measure.
	Using~\eqref{eq:lattice-helmholtz-lebesgue-decomposition}, followed by
	this change of variables in $\phi$, and with the integration ranges
	$f,\eta\in\mathbb{R}^{\Lambda_{M,\varepsilon}}_0$,
	$h_0\in[0,M^{-1})^2$, $n\in M^{-1}\mathbb{Z}^2$, and
	$\phi\in\mathbb{C}^{\Lambda_{M,\varepsilon}}$ understood, we obtain
	\begin{equation}
		\label{eq:gauge-volume-factorization}
		\begin{aligned}
			\int\bigl(Fe^{-S_{M,\varepsilon}}\bigr)(A,\phi)\,dA\,d\phi
			&=C_{M,\varepsilon}\int
			\bigl(Fe^{-S_{M,\varepsilon}}\bigr)
			\bigl(P_{\chi_\varepsilon}^{M}\nabla_\perp^\varepsilon\eta+h_0,\phi\bigr)
			\,df\,d\eta\,dh_0\,dn\,d\phi\\
			&=C_{M,\varepsilon}\left(\int df\,dn\right)\int
			\bigl(Fe^{-S_{M,\varepsilon}}\bigr)
			\bigl(P_{\chi_\varepsilon}^{M}\nabla_\perp^\varepsilon\eta+h_0,\phi\bigr)
			\,d\eta\,dh_0\,d\phi.
		\end{aligned}
	\end{equation}
	
	    Recall that
	    $\mathcal A_{\mathrm{gf}}^{M,\varepsilon}
	    :=
	    \bigl(\mathbb{R}^2\bigr)_{\chi_\varepsilon,\mathrm{div}}^{\Lambda_{M,\varepsilon}}+[0,M^{-1})^2$.
	    The Higgs configuration space is 
	    $\mathbb{C}^{\Lambda_{M,\varepsilon}}$.
	   After removing the infinite factor $\int df\,dn$,
	    we define the
	    unrenormalized lattice Abelian--Higgs probability measure on
	    $\mathcal A_{\mathrm{gf}}^{M,\varepsilon}
	    \times\mathbb{C}^{\Lambda_{M,\varepsilon}}$ by
		\begin{equation}
			\label{eq:rigorous-lattice-abelian-higgs-measure}
			\begin{aligned}
				\mu_{M,\varepsilon}^{\mathrm{YMH,unren}}(dA,d\phi)
				&\overset{\mathrm{def}}{\propto}
				\exp\bigl(-S_{M,\varepsilon}(A,\phi)\bigr)
				\,\lambda_{M,\varepsilon}^{\mathrm{GF}}(dA)\,d\phi\\
				&=
				\exp\left(-\int_{\Lambda_{M,\varepsilon}}
				\left(
				\frac12\lvert\nabla^\varepsilon A(x)\rvert^2
				+\frac12\lvert D_A^\varepsilon\phi(x)\rvert^2
				+\frac14\lvert\phi(x)\rvert^4
				\right)\,dx\right)\,\lambda_{M,\varepsilon}^{\mathrm{GF}}(dA)\,d\phi,
			\end{aligned}
		\end{equation}
		where $\lambda_{M,\varepsilon}^{\mathrm{GF}}$ is the pushforward of the
		product Lebesgue measure $d\eta\,dh_0$ on
		$\mathbb{R}^{\Lambda_{M,\varepsilon}}_0\times[0,M^{-1})^2$ under the parametrization
		$A=P_{\chi_\varepsilon}^{M}\nabla_\perp^\varepsilon\eta+h_0$, and $d\phi$ is product Lebesgue
		measure on $\mathbb{C}^{\Lambda_{M,\varepsilon}}$.		Here
		$\lvert\nabla^\varepsilon A(x)\rvert^2
		:=
		\sum_{i,j=1}^2\lvert\partial_i^\varepsilon A_j(x)\rvert^2$,
		and the second equality uses $\operatorname{div}^\varepsilon A=0$.
		By  Lemma~\ref{lem:app-cor-satisfies-hypothesis},
		the corresponding partition function is finite, so the measure can be
		normalized to a probability measure.

	Although lattice regularization makes all nonlinear products in the action
	well defined for each fixed $\varepsilon>0$, the ultraviolet divergence
	reappears as the lattice spacing tends to zero. To obtain a meaningful
	continuum limit, we must renormalize the model by subtracting an
	$M,\varepsilon$-dependent divergent counterterm.

	\begin{remark}[Faddeev--Popov gauge fixing]
		\label{rem:helmholtz-fp-gaussian-averaging}
		The above gauge fixing is equivalent to the Faddeev--Popov gauge fixing
		used in physics; cf.~\cite[Section~9.4]{Pes95}. Formally, up to a
		normalization constant,
		\begin{equation}
			\label{eq:coulomb-slice-delta-representation}
			\lambda_{M,\varepsilon}^{\mathrm{GF}}(dA)
			\overset{\mathrm{formally}}{\propto}
			\delta_0\bigl(\operatorname{div}^\varepsilon A\bigr)
			\mathbf{1}_{\left\{
			\overline A^{\chi_\varepsilon}\in[0,M^{-1})^2
			\right\}}\,dA,
		\end{equation}
		where $\delta_0$ denotes the Dirac delta at zero on
		$\mathbb{R}^{\Lambda_{M,\varepsilon}}_0$.
		Thus~\eqref{eq:rigorous-lattice-abelian-higgs-measure} corresponds to the
		third displayed formula on page~296 of~\cite{Pes95}, with Gaussian
		averaging subsequently performed in their formula~(9.56).
		The correspondence between the Coulomb-gauge and Gaussian-averaged
		measures on the torus is rigorously established
		in~\cite[Theorem~1.5]{BCCKM26}.
	\end{remark}
	
	\subsection{The Gaussian free field and pure Abelian Yang--Mills field}
	\label{subsec:gff-pure-abelian-yang-mills}
	
	Using the quadratic terms in
	~\eqref{eq:rigorous-lattice-abelian-higgs-measure}, let
	$\mu_{M,\varepsilon}^{\mathrm{GFF}}$ be the centered complex Gaussian
	probability measure defined by
	\begin{equation}
		\label{eq:lattice-gff-measure}
		\begin{aligned}
			\mu_{M,\varepsilon}^{\mathrm{GFF}}(d\phi)
			&\overset{\mathrm{def}}{\propto}
			\exp\left(-\frac12
			\left\langle
			(\mathcal{L}_0^\varepsilon+1)\phi,\phi
			\right\rangle_{M,\varepsilon}\right)\,d\phi,
		\end{aligned}
	\end{equation}
	Its covariance on the underlying real Hilbert space is
	$(\mathcal{L}_0^\varepsilon+1)^{-1}$.
	Let $\mathbb{B}^{M,\varepsilon}$ be a standard Brownian motion in
	$(\mathbb{C}^{\Lambda_{M,\varepsilon}},
	\langle\cdot,\cdot\rangle_{M,\varepsilon})$
	on a filtered Wiener space
	$(\Omega^\phi,\mathcal{F}^\phi,\mathbb{P}^\phi,
	(\mathcal{F}_t^\phi)_{0\leq t\leq1})$. Set
	\begin{equation}
		\label{eq:brownian-realizations-of-free-fields}
		\phi^{M,\varepsilon}
		:=
		\bigl(\mathcal{L}_0^\varepsilon+1\bigr)^{-1/2}
		\mathbb{B}_1^{M,\varepsilon},
		\qquad
		\mu_{M,\varepsilon}^{\mathrm{GFF}}
		=\operatorname{Law}_{\mathbb{P}^\phi}(\phi^{M,\varepsilon}).
	\end{equation}
	
	Set
$
	K^{M,\varepsilon}:=\nabla_\perp^\varepsilon(-\Delta^\varepsilon)^{-1}.
$
	For the gauge field $A$, let $\mu_{M,\varepsilon}^{\mathrm{YM}}$ be the
	probability measure on
	$\mathcal A_{\mathrm{gf}}^{M,\varepsilon}$ defined by
	\begin{equation}
		\label{eq:pure-abelian-yang-mills-measure}
		\mu_{M,\varepsilon}^{\mathrm{YM}}(dA)
		\overset{\mathrm{def}}{\propto}
		\exp\left(-\frac12\int_{\Lambda_{M,\varepsilon}}
		\lvert\nabla^\varepsilon A(x)\rvert^2\,dx\right)
		\lambda_{M,\varepsilon}^{\mathrm{GF}}(dA)
	\end{equation}
	This measure admits the random-field representation
	\begin{equation}
		\label{eq:pure-yang-mills-law}
		\begin{aligned}
		\mu_{M,\varepsilon}^{\mathrm{YM}}
		=\operatorname{Law}\bigl(A^{M,\varepsilon,h_0}\bigr),\quad
		A^{M,\varepsilon,h_0}
		:=P_{\chi_\varepsilon}^{M}K^{M,\varepsilon}W_1^{M,\varepsilon}+h_0.
		\end{aligned}
	\end{equation}
	Here $W^{M,\varepsilon}$ is a standard Brownian motion in
	$\mathbb{R}^{\Lambda_{M,\varepsilon}}_0$ with inner product
	$\langle\cdot,\cdot\rangle_{M,\varepsilon}$, on a filtered Wiener space
	$(\Omega^A,\mathcal{F}^A,\mathbb{P}^A,
	(\mathcal{F}_t^A)_{0\leq t\leq1})$. $h_0\sim\operatorname{Unif}([0,M^{-1})^2)$ is independent of
	$W^{M,\varepsilon}$.
	Indeed, under the Coulomb parametrization described above,
	$\lambda_{M,\varepsilon}^{\mathrm{GF}}$ is the pushforward of
	$d\eta\,dh_0$. For
	$A=P_{\chi_\varepsilon}^{M}\nabla_\perp^\varepsilon\eta+h_0$, Parseval's identity and
	~\eqref{eq:perpendicular-gradient-symbol} give
	$\lVert\nabla^\varepsilon A\rVert_{L^{2,\varepsilon}}^2
	=\lVert-\Delta^{M,\varepsilon}\eta\rVert_{L^{2,\varepsilon}}^2$.
	Thus the linear change of variables $z=-\Delta^{M,\varepsilon}\eta$ transforms
	the measure in~\eqref{eq:pure-abelian-yang-mills-measure} into one proportional to
	$\exp(-\frac12\lVert z\rVert_{L^{2,\varepsilon}}^2)\,dz\,dh_0$,
	which yields the representation~\eqref{eq:pure-yang-mills-law}.
	
	Rewriting~\eqref{eq:rigorous-lattice-abelian-higgs-measure} relative to
	$\mu_{M,\varepsilon}^{\mathrm{GFF}}$ and
	$\mu_{M,\varepsilon}^{\mathrm{YM}}$ and introducing the renormalization
	counterterms, we define the renormalized lattice Abelian--Higgs probability
	measure by
	\begingroup
	\renewcommand{\theHequation}{renormalized-lattice-ymh}
	\begin{equation}
		\begin{aligned}
		\mu_{M,\varepsilon}^{\mathrm{YMH}}(dA,d\phi)
		\overset{\mathrm{def}}{\propto}{}
		\exp\Bigg(
		-\frac12
		\left\langle
		\left(
		:\!\mathcal{L}_A^\varepsilon\!:_{M,\varepsilon}
		-\mathcal{L}_0^\varepsilon-1
		\right)\phi,\phi
		\right\rangle_{M,\varepsilon}
		-\frac14\int_{\Lambda_{M,\varepsilon}}
		:\!\lvert\phi(x)\rvert^4\!:_{M,\varepsilon}\,dx
		\Bigg)
		\mu_{M,\varepsilon}^{\mathrm{GFF}}(d\phi)\,
		\mu_{M,\varepsilon}^{\mathrm{YM}}(dA).
		\end{aligned}
		\label{eq:renormalized-lattice-abelian-higgs-measure}
	\end{equation}
	\endgroup
	Here
	$:\!\mathcal{L}_A^\varepsilon\!:_{M,\varepsilon}
	:=\mathcal{L}_A^\varepsilon-a_{M,\varepsilon}$
	and
	$:\!\lvert\phi(x)\rvert^4\!:_{M,\varepsilon}
	:=
	\lvert\phi(x)\rvert^4-4b_{M,\varepsilon}\lvert\phi(x)\rvert^2
	+2b_{M,\varepsilon}^2$.
	The renormalization constants are
	\begin{equation}
		\label{eq:lattice-renormalization-constants}
		a_{M,\varepsilon}
		:=
		\mathbb{E}^A\left[
		\lvert A^{M,\varepsilon,0}(0)\rvert^2
		\right],
		\qquad
		b_{M,\varepsilon}
		:=
		\mathbb{E}^\phi\left[
		\lvert\phi^{M,\varepsilon}(0)\rvert^2
		\right].
	\end{equation}
	
	\subsection{FBSDE representations}
	\label{subsec:variational-representation-and-fbsdes}
	
	Define
	\begin{equation}
		\label{eq:effective-interaction-tilde-V}
		\widetilde{V}^{M,\varepsilon}(A,\phi)
		:=
		\frac12
		\left\langle
		\bigl(:\!\mathcal{L}_A^\varepsilon\!:_{M,\varepsilon}
		-\mathcal{L}_0^\varepsilon-1\bigr)\phi,\phi
		\right\rangle_{M,\varepsilon}
		+\frac14\int_{\Lambda_{M,\varepsilon}}
		:\!\lvert\phi(x)\rvert^4\!:_{M,\varepsilon}\,dx.
	\end{equation}
	For
	$A\in
	\mathcal A_{\mathrm{gf}}^{M,\varepsilon}$,
	we define the effective action by
	\begin{equation}
		\label{eq:effective-gauge-potential}
		\mathcal{W}^{M,\varepsilon}(A)
		:=
		-\log\int
		\exp\bigl(-\widetilde{V}^{M,\varepsilon}(A,\phi)\bigr)\,
		\mu_{M,\varepsilon}^{\mathrm{GFF}}(d\phi).
	\end{equation}
	For $A\in(\mathbb R^2)_{\chi_\varepsilon,\mathrm{div}}^{\Lambda_{M,\varepsilon}}$,
	we integrate out the zero-frequency variable and define
	\begin{equation}
		\label{eq:harmonic-averaged-effective-potential}
		\mathcal W_{\mathrm{av}}^{M,\varepsilon}(A)
		:=
		-\log\left(\int_{[0,M^{-1})^2}
		\exp\bigl(-\mathcal W^{M,\varepsilon}(A+h_0)\bigr)\,dh_0\right).
	\end{equation}
	For
	$B\in
	\bigl(\mathbb{R}^2\bigr)_{\chi_\varepsilon,\mathrm{div}}^{\Lambda_{M,\varepsilon}}$,
	define the covariant Gaussian free field and its law by
	\begin{equation}
		\label{eq:B-covariant-gff-measure}
		\mu_{B,M,\varepsilon}^{\mathrm{GFF}}
		:=\operatorname{Law}_{\mathbb{P}^\phi}
		\bigl(\phi_B^{M,\varepsilon}\bigr),\qquad
		\phi_B^{M,\varepsilon}
		:=
		\bigl(\mathcal{L}_B^\varepsilon+1\bigr)^{-1/2}
		\mathbb{B}_1^{M,\varepsilon}.
	\end{equation}
	
	For $A\in\mathcal A_{\mathrm{gf}}^{M,\varepsilon}$ and
	$B\in(\mathbb R^2)_{\chi_\varepsilon,\mathrm{div}}^{\Lambda_{M,\varepsilon}}$,
	set $C:=A+B$. For $\phi\in\mathbb C^{\Lambda_{M,\varepsilon}}$, define
	\begin{equation}
		\label{eq:relative-covariant-interaction}
		V_{A,B}(\phi)
		:=
		\frac12
		\left\langle
		\bigl(:\!\mathcal{L}_C^\varepsilon\!:_{M,\varepsilon}
		-\mathcal{L}_B^\varepsilon-1\bigr)
		\phi,\phi
		\right\rangle_{M,\varepsilon}
		+\frac14\int_{\Lambda_{M,\varepsilon}}
		:\!\lvert\phi(x)\rvert^4\!:_{M,\varepsilon}\,dx.
	\end{equation}
	Define the conditional Higgs probability measure by
	\begin{equation}
		\label{eq:conditional-renormalized-higgs-measure}
		\begin{aligned}
			\mu_{C,M,\varepsilon}^{\mathrm{Higgs}}(d\phi)
			&\overset{\mathrm{def}}{\propto}_{A,B}
			\exp\bigl(-\widetilde{V}^{M,\varepsilon}(C,\phi)\bigr)
			\,\mu_{M,\varepsilon}^{\mathrm{GFF}}(d\phi)\\
			&\propto_{A,B}
			\exp\bigl(-V_{A,B}(\phi)\bigr)
			\,\mu_{B,M,\varepsilon}^{\mathrm{GFF}}(d\phi).
		\end{aligned}
	\end{equation}
	The gauge-field marginal on
	$(\mathbb R^2)_{\chi_\varepsilon,\mathrm{div}}^{\Lambda_{M,\varepsilon}}$ is
	\begin{equation}
		\label{eq:conditional-ymh-gauge-law}
		\mu_{M,\varepsilon}^{\mathrm{gauge}}(dA)
		\overset{\mathrm{def}}{\propto}
		\exp\bigl(-\mathcal{W}_{\mathrm{av}}^{M,\varepsilon}(A)\bigr)\,
		\operatorname{Law}_{\mathbb{P}^A}
		\bigl(A^{M,\varepsilon,0}\bigr)(dA).
	\end{equation}
	These Gibbs measures are normalizable by
	Lemma~\ref{lem:app-cor-satisfies-hypothesis}.
	
	\begin{remark}[Choice of Gaussian reference measure]
		\label{rem:gaussian-reference-measure-choice}
		In defining the effective gauge potential \eqref{eq:effective-gauge-potential},
		we use the free GFF with covariance $(\mathcal{L}_0^\varepsilon+1)^{-1}$.
		By contrast, \cite{BCCKM26} uses an $A$-covariant GFF, requiring an
		$A$-dependent spectral shift and introducing a determinant factor and
		a compensating interaction term. Our choice avoids these additional terms in the outer gauge-field problem,
		but directly estimating $\mathcal{L}_{A+B}^\varepsilon-\mathcal{L}_0^\varepsilon$
		leads to excessive dependence on the gauge shift $B$.
		Inspired by \cite[(6.19)]{BC24}, we therefore reorganize the conditional
		Higgs problem around the $B$-covariant GFF, as in
		\eqref{eq:conditional-renormalized-higgs-measure}.
		This replaces $\mathcal{L}_{A+B}^\varepsilon-\mathcal{L}_0^\varepsilon$ by
		$\mathcal{L}_{A+B}^\varepsilon-\mathcal{L}_B^\varepsilon$ in the interaction
		and yields sufficiently mild dependence on $B$ to close the gauge-field estimates.
	\end{remark}
	
	We use the canonical Wiener spaces
	$
	\Omega^\phi=C_0\bigl([0,1];\mathbb{C}^{\Lambda_{M,\varepsilon}}\bigr),
	\Omega^A=C_0\bigl([0,1];\mathbb{R}^{\Lambda_{M,\varepsilon}}_0\bigr).
	$
	In what follows, $\mathbb{E}_t^\phi$ and $\mathbb{E}_t^A$ denote
	conditional expectation with respect to $\mathcal{F}_t^\phi$ and
	$\mathcal{F}_t^A$, respectively. Here $\nabla_\phi V_{A,B}$ denotes the gradient in
	$\mathbb C^{\Lambda_{M,\varepsilon}}$, and
	$\nabla_A\mathcal W^{M,\varepsilon}$ and
	$\nabla_A\mathcal W_{\mathrm{av}}^{M,\varepsilon}$ denote gradients in
	$(\mathbb R^2)_{\chi_\varepsilon,\mathrm{div}}^{\Lambda_{M,\varepsilon}}$, all with respect to
	$\langle\cdot,\cdot\rangle_{M,\varepsilon}$.
	
	\begin{proposition}[FBSDE representation]
		\label{prop:maximum-principle-and-fbsdes}
		For every
		$A\in\mathcal A_{\mathrm{gf}}^{M,\varepsilon}$ and
		$B\in\bigl(\mathbb{R}^2\bigr)_{\chi_\varepsilon,\mathrm{div}}^{\Lambda_{M,\varepsilon}}$,
		the shifts $Z^{M,\varepsilon,A,B}$ and $B^{M,\varepsilon}$ associated with
		the Gibbs measures~\eqref{eq:conditional-renormalized-higgs-measure}
		and~\eqref{eq:conditional-ymh-gauge-law}, respectively, in the sense of
		Proposition~\ref{lem:app-abstract-maximum-principle}, satisfy the following assertions.
		
		\smallskip
		\noindent\textup{(i)}
		\emph{(FBSDEs).}
		Define the terminal states by
		\begin{equation}
			\label{eq:optimal-state-processes}
			\Phi^{M,\varepsilon,A,B}:=
			\phi_B^{M,\varepsilon}+Z_1^{M,\varepsilon,A,B},\qquad
			C^{\varepsilon,0}:=A^{M,\varepsilon,0}+B_1^{M,\varepsilon}.
		\end{equation}
		For $t\in[0,1]$, the shifts satisfy
		\begin{align}
			\label{eq:scalar-integral-fbsde}
			&Z_t^{M,\varepsilon,A,B}
			=-\bigl(\mathcal{L}_B^\varepsilon+1\bigr)^{-1}
			\int_0^t
			\mathbb{E}_s^\phi\left[
			\nabla_\phi V_{A,B}
			\bigl(\Phi^{M,\varepsilon,A,B}\bigr)
			\right]ds,\\
			\label{eq:gauge-integral-fbsde}
			&B_t^{M,\varepsilon}
			=-P_{\chi_\varepsilon}^{M}K^{M,\varepsilon}\bigl(P_{\chi_\varepsilon}^{M}K^{M,\varepsilon}\bigr)^*
			\int_0^t\mathbb E_s^A\left[
			\nabla_A\mathcal W_{\mathrm{av}}^{M,\varepsilon}
			\bigl(C^{\varepsilon,0}\bigr)
			\right]ds.
		\end{align}
		
		\smallskip
		\noindent\textup{(ii)}
		\emph{(Terminal law).}
		The terminal states satisfy
		\begin{equation}
			\label{eq:optimal-terminal-laws}
			\operatorname{Law}_{\mathbb{P}^\phi}
			\bigl(\Phi^{M,\varepsilon,A,B}\bigr)
			=\mu_{A+B,M,\varepsilon}^{\mathrm{Higgs}},
			\qquad
			\operatorname{Law}_{\mathbb{P}^A}
			\bigl(C^{\varepsilon,0}\bigr)
			=\mu_{M,\varepsilon}^{\mathrm{gauge}}.
		\end{equation}
		
		\smallskip
		\noindent\textup{(iii)}
		\emph{(Measurability).}
		The map $(A,B,\omega_\phi)\mapsto
		\Phi^{M,\varepsilon,A,B}(\omega_\phi)$ is jointly measurable.
	\end{proposition}
	
	\begin{proof}
		Parts~\textup{(i)} and~\textup{(ii)} follow from
		\eqref{eq:app-integral-fbsde} and~\eqref{eq:app-optimal-terminal-law},
		applied with diffusion coefficients
		$\bigl(\mathcal{L}_B^\varepsilon+1\bigr)^{-1/2}$ and
		$P_{\chi_\varepsilon}^{M}K^{M,\varepsilon}$, and terminal potentials
		$V_{A,B}$ and $\mathcal{W}_{\mathrm{av}}^{M,\varepsilon}$, respectively.
		The hypotheses hold by Lemma~\ref{lem:app-cor-satisfies-hypothesis}.
		
		To prove~\textup{(iii)}, note that the diffusion coefficient
		$\bigl(\mathcal{L}_B^\varepsilon+1\bigr)^{-1/2}$ is independent of time
		and, together with $V_{A,B}$, depends on $A,B$ only through \(e^{i\varepsilon A_j(x)}\) and \(e^{i\varepsilon B_j(x)}\).
		By Lemma~\ref{lem:app-cor-satisfies-hypothesis},
		the feedback drift in~\eqref{eq:app-feedback-sde} is jointly continuous in $(A,B)$
		and locally Lipschitz in the state variable, uniformly in $(t,A,B)$.
 Picard iteration for the truncated equations~\eqref{eq:app-cutoff-sde}
		yields jointly measurable solutions. The non-explosion argument in
		Proposition~\ref{lem:app-abstract-maximum-principle} then provides a jointly
		measurable version of $\Phi^{M,\varepsilon,A,B}$.
	\end{proof}
	
	\begin{remark}[Adjoint identities]
		\label{rem:adjoint-gauge-control-operator}
		$P_{\chi_\varepsilon}^{M}K^{M,\varepsilon}$ is the composition
		\begin{equation}
			\label{eq:gauge-control-operator-mapping}
			\mathbb R_0^{\Lambda_{M,\varepsilon}}
			\xrightarrow{\ K^{M,\varepsilon}\ }
			(\mathbb R^2)_{0,\mathrm{div}}^{\Lambda_{M,\varepsilon}}
			\xrightarrow{\ P_{\chi_\varepsilon}^{M}\ }
			(\mathbb R^2)_{\chi_\varepsilon,\mathrm{div}}^{\Lambda_{M,\varepsilon}}.
		\end{equation}
		Let $P_\perp^{M,\varepsilon}$ be the orthogonal projection onto
		$(\mathbb R^2)_{0,\mathrm{div}}^{\Lambda_{M,\varepsilon}}$. Then
$
		(P_{\chi_\varepsilon}^{M})^*A
		=P_\perp^{M,\varepsilon}
		\Bigl(A-\widehat{A}(0)\chi_\varepsilon\Bigr).
$
		Moreover,
$
		(K^{M,\varepsilon})^*
		=(-\Delta^\varepsilon)^{-1}(\nabla_\perp^\varepsilon)^*,
$
		where, for $k\in\widehat\Lambda_{M,\varepsilon}$,
$
		\bigl((\nabla_\perp^\varepsilon)^*A\bigr)^\wedge(k)
		=-i\bigl\langle(k^\varepsilon)^\perp,\widehat A(k)\bigr\rangle_{\mathbb C^2}.
$
	\end{remark}
	
	\begin{remark}
		\label{rem:maximum-principle-literature}
		The
		FBSDEs~\eqref{eq:scalar-integral-fbsde}
		and~\eqref{eq:gauge-integral-fbsde} are classical;
		see~\cite[Chapter~3, Theorem~3.2]{YZ99}.
In fact, in the present setting, these
		systems take a particularly simple form: their BSDEs have zero drift,
		so the adjoint processes can be eliminated via conditional expectations,
		leaving equations for the state processes alone.
		For the unbounded potential~\eqref{eq:relative-covariant-interaction},
		we establish the FBSDE representation in
		Appendix~\ref{app:maximum-principle-fbsde}.
		For specific models with bounded potentials, analogues of
		Proposition~\ref{prop:maximum-principle-and-fbsdes} have been
		established in~\cite[Proposition~3.8]{BGH23}
		and~\cite[Theorem~2.2]{GM24}.
	\end{remark}

	In what follows, we suppress $M$ in the notation for random fields,
	but retain it for configuration spaces and inner products.

	We now use the terminal states of the two FBSDEs to represent $\mu_{M,\varepsilon}^{\mathrm{YMH}}$. Define the product probability space
	\begin{equation}
		\label{eq:product-probability-space}
		\begin{aligned}
			\Omega:=
			[0,M^{-1})^2\times\Omega^A\times\Omega^\phi,\quad
			\mathcal{F}:=
			\mathcal{B}([0,M^{-1})^2)\otimes\mathcal{F}^A\otimes\mathcal{F}^\phi,\quad
			\mathbb{P}:=M^2\,dh_0
			\otimes\mathbb{P}^A\otimes\mathbb{P}^\phi.
		\end{aligned}
	\end{equation}
	We identify random variables on the two Wiener factors with their canonical
	lifts to $\Omega$. For
	$C\in(\mathbb R^2)_{\chi_\varepsilon,\mathrm{div}}^{\Lambda_{M,\varepsilon}}$,
	define
	\begin{equation}
		\label{eq:harmonic-conditional-law}
		\mu_C^{\mathrm{av}}(dh)
		:=\frac{\exp\bigl(-\mathcal W^{M,\varepsilon}(C+h)\bigr)}
		{\exp\bigl(-\mathcal W_{\mathrm{av}}^{M,\varepsilon}(C)\bigr)}\,dh,
		\qquad h\in[0,M^{-1})^2.
	\end{equation}
	We choose maps
	$H^{\varepsilon,C}\colon[0,M^{-1})^2\to[0,M^{-1})^2$
	such that $(C,h_0)\mapsto H^{\varepsilon,C}(h_0)$ is jointly Borel
	measurable and $(H^{\varepsilon,C})_\#(M^2\,dh_0)=\mu_C^{\mathrm{av}}$,
	where $\#$ denotes pushforward. For
	$\omega=(h_0,\omega_A,\omega_\phi)\in\Omega$, set
	\begin{equation}
		\label{eq:random-harmonic-variable}
		H^\varepsilon(\omega)
		:=H^{\varepsilon,C^{\varepsilon,0}(\omega_A)}(h_0).
	\end{equation}
	Then $\mathbb P(H^\varepsilon\in dh\mid\mathcal F_1^A)
	=\mu_{C^{\varepsilon,0}}^{\mathrm{av}}(dh)$. Define
	\begin{equation}
		\label{eq:random-harmonic-terminal-state-lifts}
		A^\varepsilon(\omega):=
		A^{\varepsilon,0}(\omega_A)+H^\varepsilon(\omega),
		\qquad
		C^\varepsilon:=A^\varepsilon+B_1^\varepsilon.
	\end{equation}
	\begin{equation}
		\label{eq:random-harmonic-field-shorthands}
		\Phi^\varepsilon(\omega):=
		\Phi^{\varepsilon,A^\varepsilon(\omega),B_1^\varepsilon(\omega_A)}(\omega_\phi),
		\qquad
		Z_t^\varepsilon(\omega):=
		Z_t^{\varepsilon,A^\varepsilon(\omega),B_1^\varepsilon(\omega_A)}(\omega_\phi).
	\end{equation}
	By construction and Proposition~\ref{prop:maximum-principle-and-fbsdes}\textup{(iii)},
	the terminal states $C^\varepsilon$ and $\Phi^\varepsilon$ are
	$\mathcal{F}$-measurable.
	We also write $B^\varepsilon:=B_1^\varepsilon$ for
	the terminal gauge-field shift. 
	On $\Omega$, we set
	\begin{equation}
		\label{eq:product-space-conditional-expectations}
		\begin{aligned}
			\mathbb E^\phi[\cdot]
			:=\mathbb E\bigl[\cdot\,\big|\,
			\mathcal B([0,M^{-1})^2)\otimes\mathcal F^A\bigr],\quad
			\mathbb E_t^\phi[\cdot]
			:=\mathbb E\bigl[\cdot\,\big|\,
			\mathcal B([0,M^{-1})^2)\otimes\mathcal F^A\otimes\mathcal F_t^\phi\bigr],\quad
			\mathbb E_t^A[\cdot]
			:=\mathbb E\bigl[\cdot\,\big|\,\mathcal F_t^A\bigr].
		\end{aligned}
	\end{equation}
	Here the conditioning $\sigma$-fields are identified with their canonical lifts to $\Omega$.
	
	\begin{proposition}[FBSDE representation of the renormalized measure]
		\label{thm:stochastic-representation-renormalized-ymh-measure} The renormalized lattice Abelian--Higgs measure satisfies
		\begin{equation}
			\label{eq:law-representation-renormalized-ymh-measure}
			\mu_{M,\varepsilon}^{\mathrm{YMH}}
			=\operatorname{Law}_{\mathbb{P}}
			\bigl(C^\varepsilon,\Phi^\varepsilon\bigr).
		\end{equation}
	\end{proposition}
	
	\begin{proof}
		Let $F\colon
		\mathcal A_{\mathrm{gf}}^{M,\varepsilon}
		\times\mathbb{C}^{\Lambda_{M,\varepsilon}}\to\mathbb{R}$
		be bounded and measurable. Define
$
			\label{eq:higgs-averaged-test-function}
			\widetilde F(C)
			:=\int F(C,\phi)\,\mu_{C,M,\varepsilon}^{\mathrm{Higgs}}(d\phi).
$
		By~\eqref{eq:optimal-terminal-laws},
$
			\label{eq:higgs-test-function-terminal-identity}
			\widetilde F(C^\varepsilon)
			=\mathbb E^\phi\bigl[F(C^\varepsilon,\Phi^\varepsilon)\bigr].
$
		For $C\in(\mathbb R^2)_{\chi_\varepsilon,\mathrm{div}}^{\Lambda_{M,\varepsilon}}$,
		define
$
			\label{eq:harmonic-averaged-test-function}
			\widetilde F_{\mathrm{av}}(C)
			:=\int_{[0,M^{-1})^2}\widetilde F(C+h)\,\mu_C^{\mathrm{av}}(dh).
$
		By the construction of $H^\epsilon$,
		\begin{equation}
			\label{eq:harmonic-test-function-conditional-expectation}
			\begin{aligned}
				\widetilde F_{\mathrm{av}}(C^{\varepsilon,0})
				=\mathbb E\bigl[\widetilde F(C^\varepsilon)\mid\mathcal F_1^A\bigr]
				=\mathbb E\bigl[F(C^\varepsilon,\Phi^\varepsilon)\mid\mathcal F_1^A\bigr].
			\end{aligned}
		\end{equation}
	Using ~\eqref{eq:renormalized-lattice-abelian-higgs-measure}
		and~\eqref{eq:optimal-terminal-laws}, we obtain
		\begin{equation}
			\label{eq:ymh-expectation-product-space}
			\begin{aligned}
				\int F(C,\phi)\,\mu_{M,\varepsilon}^{\mathrm{YMH}}(dC,d\phi)
				&\propto\mathbb E^A\left[\int_{[0,M^{-1})^2}\int
				F(A^{\varepsilon,h},\phi)
				\exp\bigl(-\widetilde V^{M,\varepsilon}(A^{\varepsilon,h},\phi)\bigr)
				\mu_{M,\varepsilon}^{\mathrm{GFF}}(d\phi)\,dh\right]\\
				&=\mathbb E^A\left[\int_{[0,M^{-1})^2}
				\widetilde F(A^{\varepsilon,h})
				\exp\bigl(-\mathcal W^{M,\varepsilon}(A^{\varepsilon,h})\bigr)\,dh\right]\\
				&=\mathbb E^A\left[
				\widetilde F_{\mathrm{av}}(A^{\varepsilon,0})
				\exp\bigl(-\mathcal W_{\mathrm{av}}^{M,\varepsilon}(A^{\varepsilon,0})\bigr)
				\right]\\
				&\overset{\eqref{eq:optimal-terminal-laws}}{\propto}
				\mathbb E\bigl[\widetilde F_{\mathrm{av}}(C^{\varepsilon,0})\bigr]\\
				&\overset{\eqref{eq:harmonic-test-function-conditional-expectation}}{=}
				\mathbb E\bigl[F(C^\varepsilon,\Phi^\varepsilon)\bigr].
			\end{aligned}
		\end{equation}
		This gives the law representation~\eqref{eq:law-representation-renormalized-ymh-measure}.
	\end{proof}
	
	Let $\Pi_{\chi_\varepsilon}^{M}$ denote the orthogonal projection from
	$(\mathbb{R}^2)^{\Lambda_{M,\varepsilon}}$ onto
	$(\mathbb{R}^2)_{\chi_\varepsilon,\mathrm{div}}^{\Lambda_{M,\varepsilon}}$
	with respect to $\langle\cdot,\cdot\rangle_{M,\varepsilon}$.
	
	\begin{lemma}
		\label{lem:projected-lattice-current-identity}
		For every $C\colon\Lambda_{M,\varepsilon}\to\mathbb{R}^2$ and
		$\phi,\psi\colon\Lambda_{M,\varepsilon}\to\mathbb{C}$,
		\begin{equation}
			\label{eq:projected-lattice-current-identity}
			P_\perp^{M,\varepsilon}
			\left(
			\operatorname{Im}\left[
			\overline{\phi(\mathord\cdot+\varepsilon e_j)}
			U_{-C,j}^\varepsilon D_{C,j}^\varepsilon\psi
			\right]
			\right)_{j=1}^2
			=
			P_\perp^{M,\varepsilon}
			\left(
			\operatorname{Im}\left[
			\overline{\psi}\,D_{C,j}^\varepsilon\phi
			\right]
			\right)_{j=1}^2.
		\end{equation}
		The same identity holds with $P_\perp^{M,\varepsilon}$ replaced by
		$\Pi_{\chi_\varepsilon}^{M}$. Moreover, for $j=1,2$,
		\begin{equation}
			\label{eq:pointwise-lattice-current-identity}
			\operatorname{Im}\left[
			\overline{\phi(\mathord\cdot+\varepsilon e_j)}
			U_{-C,j}^\varepsilon D_{C,j}^\varepsilon\phi
			\right]
			=\operatorname{Im}\bigl(\overline\phi D_{C,j}^\varepsilon\phi\bigr).
		\end{equation}
	\end{lemma}
	
	\begin{proof}
		Both $P_\perp^{M,\varepsilon}$ and $\Pi_{\chi_\varepsilon}^{M}$ annihilate
		$\nabla^\varepsilon$, so it is enough to prove,
		for $j=1,2$,
		\begin{equation}
			\label{eq:covariant-discrete-product-rule}
			D_j^\varepsilon(\overline{\phi}\psi)
			=
			\overline{D_{C,j}^\varepsilon\phi}\,\psi
			{}+
			\overline{\phi(\mathord\cdot+\varepsilon e_j)}
			U_{-C,j}^\varepsilon D_{C,j}^\varepsilon\psi.
		\end{equation}
		Combining the ordinary lattice product rule with the definition of the
		covariant derivative, one verifies
		\eqref{eq:covariant-discrete-product-rule} directly. Taking imaginary
		parts and applying either projection yields~\eqref{eq:projected-lattice-current-identity}.
		The left-hand side of~\eqref{eq:pointwise-lattice-current-identity} equals
		\begin{equation}
			\label{eq:pointwise-lattice-current-calculation}
			\begin{aligned}
				\frac1\varepsilon\operatorname{Im}\left[
				\overline{\phi(\mathord\cdot+\varepsilon e_j)}
				\bigl(\phi(\mathord\cdot+\varepsilon e_j)-U_{-C,j}^\varepsilon\phi\bigr)
				\right]
				&=-\frac1\varepsilon\operatorname{Im}\left[
				\overline{\phi(\mathord\cdot+\varepsilon e_j)}U_{-C,j}^\varepsilon\phi
				\right]\\
				&=\frac1\varepsilon\operatorname{Im}\left[
				\overline\phi U_{C,j}^\varepsilon\phi(\mathord\cdot+\varepsilon e_j)
				\right]
				=\operatorname{Im}\bigl(\overline\phi D_{C,j}^\varepsilon\phi\bigr).
			\end{aligned}
		\end{equation}
	\end{proof}
	
	\begin{lemma}[Calculation of gradients]
		\label{lem:calculation-of-gradients}
		The
		FBSDEs~\eqref{eq:scalar-integral-fbsde} and
		\eqref{eq:gauge-integral-fbsde} can be written explicitly as
		\begin{equation}
			\label{eq:explicit-scalar-integral-fbsde}
			Z_t^\varepsilon
			=-\bigl(\mathcal{L}_{B_1^\varepsilon}^\varepsilon+1\bigr)^{-1}
			\int_0^t\mathbb{E}_s^\phi\Bigg[
			\Bigl(
			:\!\mathcal{L}_{C^\varepsilon}^\varepsilon\!:_{M,\varepsilon}
			-\mathcal{L}_{B_1^\varepsilon}^\varepsilon-1
			\Bigr)
			\Phi^\varepsilon
			+:\!\lvert\Phi^\varepsilon\rvert^2
			\Phi^\varepsilon\!:_{M,\varepsilon}
			\Bigg]ds,
		\end{equation}
		and
		\begin{equation}
		\label{eq:explicit-gauge-integral-fbsde}
		B_t^\varepsilon
		=-P_{\chi_\varepsilon}^{M}K^{M,\varepsilon}\bigl(P_{\chi_\varepsilon}^{M}K^{M,\varepsilon}\bigr)^*
		\int_0^t\mathbb E_s^A\Bigl[
		\mathbb E^\phi\bigl[
		\Pi_{\chi_\varepsilon}^{M}
		\operatorname{Im}\bigl(
		\overline{\Phi^\varepsilon}D_{C^\varepsilon}^\varepsilon\Phi^\varepsilon
		\bigr)\bigr]
		\Bigr]ds.
	\end{equation}
	\end{lemma}
	
	\begin{proof}
		Since Wick powers preserve the differentiation rules of polynomials,
		\begin{equation}
			\label{eq:scalar-gradient-relative-interaction}
			\nabla_\phi V_{A,B}(\phi)
			=
			\Bigl(
			:\!\mathcal{L}_{A+B}^\varepsilon\!:_{M,\varepsilon}
			-\mathcal{L}_B^\varepsilon-1
			\Bigr)\phi
			+:\!\lvert\phi\rvert^2\phi\!:_{M,\varepsilon},
			\qquad
			:\!\lvert\phi\rvert^2\phi\!:_{M,\varepsilon}
			=\bigl(\lvert\phi\rvert^2-2b_{M,\varepsilon}\bigr)\phi.
		\end{equation}
		Substitution into~\eqref{eq:scalar-integral-fbsde} with
		$A=A^\varepsilon$ and $B=B_1^\varepsilon$ gives
		\eqref{eq:explicit-scalar-integral-fbsde}.
		
		For the gauge-field gradient, differentiation under the
		integral \eqref{eq:effective-gauge-potential} and the first identity in~\eqref{eq:optimal-terminal-laws} give
		\begin{equation}
			\label{eq:gauge-gradient-effective-potential-expectation}
			\begin{aligned}
				\nabla_A\mathcal{W}^{M,\varepsilon}(C^\varepsilon)
				&=
				\mathbb{E}_{\mu_{C^\varepsilon,M,\varepsilon}^{\mathrm{Higgs}}}
				\left[
				\nabla_A\widetilde{V}^{M,\varepsilon}(C^\varepsilon,\phi)
				\right]
				=
				\mathbb{E}^\phi\left[
				\nabla_A\widetilde{V}^{M,\varepsilon}
				\bigl(C^\varepsilon,\Phi^\varepsilon\bigr)
				\right].
			\end{aligned}
		\end{equation}
		For fixed $C$, only the covariant quadratic form in
		$\widetilde{V}^{M,\varepsilon}(C+\delta h,\phi)$ depends on $\delta$. For
		$h\in(\mathbb{R}^2)^{\Lambda_{M,\varepsilon}}$,
		Lemma~\ref{lem:lattice-covariant-laplacian-expansion}\textup{(ii)},
		applied with $A=\delta h$ and $B=C$, and the expansion
		$(\rho_{\delta h,j}^\varepsilon\phi)(x)
		=\delta h_j(x)\phi(x+\varepsilon e_j)+O(\delta^2)$ as $\delta\to0$
		give, by~\eqref{eq:pointwise-lattice-current-identity}, the identity used below:
		\begin{equation}
			\label{eq:gauge-gradient-effective-interaction-current}
			\begin{aligned}
			\left.\frac{d}{d\delta}\right|_{\delta=0}
			\widetilde{V}^{M,\varepsilon}(C+\delta h,\phi)
			&=\left\langle
			\left(
			\operatorname{Im}\left[
			\overline{\phi(\mathord\cdot+\varepsilon e_j)}
			U_{-C,j}^\varepsilon D_{C,j}^\varepsilon\phi
			\right]
			\right)_{j=1}^2,h
			\right\rangle_{M,\varepsilon}\\
			&=\left\langle
			\operatorname{Im}\bigl(\overline\phi D_C^\varepsilon\phi\bigr),h
			\right\rangle_{M,\varepsilon}.
			\end{aligned}
		\end{equation}
		Now take $h\in(\mathbb{R}^2)_{\chi_\varepsilon,\mathrm{div}}^{\Lambda_{M,\varepsilon}}$.
		Using self-adjointness of $\Pi_{\chi_\varepsilon}^{M}$ and
		$\Pi_{\chi_\varepsilon}^{M}h=h$, we differentiate
		\eqref{eq:harmonic-averaged-effective-potential} and apply
		\eqref{eq:harmonic-conditional-law} to obtain
		\begin{equation}
			\label{eq:harmonic-averaged-gauge-current}
			\nabla_A\mathcal W_{\mathrm{av}}^{M,\varepsilon}(C^{\varepsilon,0})
			=\mathbb E\Bigl[\mathbb E^\phi\bigl[
			\Pi_{\chi_\varepsilon}^{M}\operatorname{Im}
			(\overline{\Phi^\varepsilon}D_{C^\varepsilon}^\varepsilon\Phi^\varepsilon)
			\bigr]\,\big|\,\mathcal F_1^A\Bigr].
		\end{equation}
	 The tower property in
		\eqref{eq:gauge-integral-fbsde} now yields
		\eqref{eq:explicit-gauge-integral-fbsde}.
	\end{proof}

	\section{Estimates of the Higgs field}
	\label{sec:estimate-higgs-field}
	
	Fix $A\in\mathcal A_{\mathrm{gf}}^{M,\varepsilon}$ and
	$B\in(\mathbb R^2)_{\chi_\varepsilon,\mathrm{div}}^{\Lambda_{M,\varepsilon}}$.
	We estimate the Higgs-field equation
	~\eqref{eq:explicit-scalar-integral-fbsde} following the  low–high frequency decomposition strategy of~\cite{GH19}.
	We collect the stochastic objects in $\Xi^{\varepsilon,A,B}$ and set
	\begin{equation}
		\label{eq:stochastic-object-Xi-norm}
		\begin{aligned}
		\lVert\Xi^{\varepsilon,A,B}\rVert_\kappa	&:=1
			+\lVert r_A^\varepsilon
			\rVert_{\mathcal{C}^{-\kappa,\varepsilon}(\rho^\kappa)}
			+\sum_{\substack{
					0\leq k_1\leq2,\ 0\leq k_2\leq1\\
					k_1+k_2\neq0}}
			\left\lVert
			:\!(\phi_B^\varepsilon)^{k_1}
			(\overline{\phi_B^\varepsilon})^{k_2}\!:_\varepsilon
			\right\rVert_{\mathcal{C}^{-\kappa,\varepsilon}(\rho^\kappa)}
			+\lVert\phi^\varepsilon
			\rVert_{\mathcal{C}^{-\kappa,\varepsilon}(\rho^\kappa)}\\
			&\quad+\lVert
			\overline{r_A^\varepsilon}\phi_B^\varepsilon
			\rVert_{\mathcal{C}^{-\kappa,\varepsilon}(\rho^\kappa)}
			+\lVert
			\overline{r_A^\varepsilon}\phi^\varepsilon
			\rVert_{\mathcal{C}^{-\kappa,\varepsilon}(\rho^\kappa)}
			+\lVert\rho_A^\varepsilon\phi_B^\varepsilon
			\rVert_{\mathcal{C}^{-\kappa,\varepsilon}(\rho^\kappa)}
			+\lVert\rho_A^\varepsilon\phi^\varepsilon
			\rVert_{\mathcal{C}^{-\kappa,\varepsilon}(\rho^\kappa)}\\
			&\quad+
			\left\lVert
			:\!\rho_A^{\varepsilon *}\rho_A^\varepsilon\!:_\varepsilon
			\right\rVert_{\mathcal{C}^{-\kappa,\varepsilon}(\rho^\kappa)}
			+
			\left\lVert
			:\!\rho_A^{\varepsilon *}\rho_A^\varepsilon\!:_\varepsilon
			\phi_B^\varepsilon
			\right\rVert_{\mathcal{C}^{-\kappa,\varepsilon}(\rho^\kappa)}
			+
			\left\lVert
			D^{\varepsilon *}\overline{r_A^\varepsilon}
			\right\rVert_{\mathcal{C}^{-\kappa,\varepsilon}(\rho^\kappa)}
			+
			\left\lVert
			\bigl(
			D^{\varepsilon *}\overline{r_A^\varepsilon}
			\bigr)\phi_B^\varepsilon
			\right\rVert_{\mathcal{C}^{-\kappa,\varepsilon}(\rho^\kappa)}\\
			&\quad+
			\frac{
				\left\lVert
				\phi_B^\varepsilon-\phi^\varepsilon
				\right\rVert_{\mathcal{C}^{2\kappa,\varepsilon}(\rho^{16\kappa})}
			}{
				\bigl(1+\lVert B\rVert_{H^{1,\varepsilon}(\rho)}\bigr)^{16\kappa}
			}
			+
			\frac{
				\left\lVert
				\phi_B^\varepsilon-\phi^\varepsilon
				\right\rVert_{\mathcal{C}^{1-\kappa,\varepsilon}(\rho^3)}
			}{
				\bigl(1+\lVert B\rVert_{H^{1,\varepsilon}(\rho)}\bigr)^2
			},
		\end{aligned}
	\end{equation}
	where, for $\phi\in\mathbb{C}^{\Lambda_\varepsilon}$ and
	$k_1,k_2\in\mathbb{N}_0$, the Wick powers are defined pointwise by
	\begin{equation}
		\label{eq:complex-wick-power-definition}
		:\!\phi^{k_1}\overline{\phi}^{k_2}\!:_\varepsilon
		:=
		H_{k_1,k_2}\bigl(\phi,\overline{\phi};b_\varepsilon\bigr).
	\end{equation}
	Here $H_{k_1,k_2}$ are the complex Hermite polynomials
	of~\cite[Section~2.5]{BC24}, with $b_\varepsilon$ as in
	\eqref{eq:lattice-renormalization-constants}.
	For brevity, we temporarily suppress the superscript
	$(\varepsilon,A,B)$ from $Z^{\varepsilon,A,B}$,
	$\Phi^{\varepsilon,A,B}$, and
	$\Xi^{\varepsilon,A,B}$.
	We use \textcolor{singularpink}{pink} for irregular terms,
	\textcolor{energyblue}{blue} for regular terms, and
	\textcolor{remainderorange}{orange} for higher-order remainders
	in~$\varepsilon$, with $(\mathord\cdot)^\prec$ denoting a regular remainder.
	
	\subsection{Low--high frequency decomposition}
	\label{subsec:higgs-low-high-decomposition}
	
	We decompose
	$\nabla_\phi V_{A,B}
	\bigl(\phi_B^\varepsilon+Z_1^\varepsilon\bigr)$
	into low- and high-fre\-quency components. By
	\eqref{eq:scalar-gradient-relative-interaction},
	\begin{equation}
		\label{eq:higgs-gradient-low-high-start}
		\nabla_\phi V_{A,B}
		\bigl(\phi_B^\varepsilon+Z_1^\varepsilon\bigr)
		=
		\bigl(
		\mathcal{L}_{A+B}^\varepsilon
		-\mathcal{L}_B^\varepsilon-a_\varepsilon
		\bigr)\Phi^\varepsilon
		+:\!\lvert\Phi^\varepsilon\rvert^2\Phi^\varepsilon\!:_\varepsilon
		-\Phi^\varepsilon.
	\end{equation}
	
	\noindent\textup{(1)}
	We first treat
	$\bigl(
	\mathcal{L}_{A+B}^\varepsilon
	-\mathcal{L}_B^\varepsilon-a_\varepsilon
	\bigr)\Phi^\varepsilon$.
	We use the notation
	\begin{equation}
		\label{eq:wick-rho-square}
		:\!\rho_A^{\varepsilon *}\rho_A^\varepsilon\!:_\varepsilon
		:=
		\rho_A^{\varepsilon *}\rho_A^\varepsilon-a_\varepsilon,
		\qquad
		:\!\lvert r_A^\varepsilon\rvert^2\!:_\varepsilon
		:=
		\sum_{j=1}^2
		\lvert r_{A,j}^\varepsilon\rvert^2-a_\varepsilon.
	\end{equation}
	By Lemma~\ref{lem:lattice-covariant-laplacian-expansion}\textup{(iii)}
	and~\eqref{eq:wick-rho-square},
	\begin{equation}
		\label{eq:higgs-relative-laplacian-expansion}
		\begin{aligned}
			\Bigl(
			\mathcal{L}_{A+B}^\varepsilon
			-\mathcal{L}_B^\varepsilon-a_\varepsilon
			\Bigr)\Phi^\varepsilon
			&=
			iD_B^{\varepsilon *}
			\bigl(
			\overline{r_A^\varepsilon}\Phi^\varepsilon
			+U_B^\varepsilon\rho_A^\varepsilon\Phi^\varepsilon
			\bigr)
			+:\!\rho_A^{\varepsilon *}\rho_A^\varepsilon\!:_\varepsilon
			\Phi^\varepsilon
			-i\bigl(D^{\varepsilon *}\overline{r_A^\varepsilon}\bigr)
			\Phi^\varepsilon\\
			&=
			iD_B^{\varepsilon *}
			\bigl(
			\overline{r_A^\varepsilon}\phi_B^\varepsilon
			+U_B^\varepsilon\rho_A^\varepsilon
			\phi_B^\varepsilon
			\bigr)
			+iD_B^{\varepsilon *}
			\bigl(
			\overline{r_A^\varepsilon}Z_1^\varepsilon
			+U_B^\varepsilon\rho_A^\varepsilon Z_1^\varepsilon
			\bigr)\\
			&\qquad
			+:\!\rho_A^{\varepsilon *}\rho_A^\varepsilon\!:_\varepsilon
			\phi_B^\varepsilon
			+:\!\rho_A^{\varepsilon *}\rho_A^\varepsilon\!:_\varepsilon
			Z_1^\varepsilon
			\textcolor{remainderorange}{
				-i\bigl(D^{\varepsilon *}\overline{r_A^\varepsilon}\bigr)
				\phi_B^\varepsilon
				-i\bigl(D^{\varepsilon *}\overline{r_A^\varepsilon}\bigr)
				Z_1^\varepsilon
			}.
		\end{aligned}
	\end{equation}
	
	We decompose these terms as follows:
	\begin{equation}
		\label{eq:higgs-first-low-high-splitting}
		iD_B^{\varepsilon *}
		\bigl(
		\overline{r_A^\varepsilon}\phi_B^\varepsilon
		+U_B^\varepsilon\rho_A^\varepsilon\phi_B^\varepsilon
		\bigr)
		=
		iD_B^{\varepsilon *}\Bigl(
		\textcolor{singularpink}{
			\mathcal{U}_{>L_1}^\varepsilon
			\bigl(\overline{r_A^\varepsilon}\phi_B^\varepsilon\bigr)
			+U_B^\varepsilon\mathcal{U}_{>L_1}^\varepsilon
			\bigl(\rho_A^\varepsilon\phi_B^\varepsilon\bigr)}
		\Bigr)
		+iD_B^{\varepsilon *}\Bigl(
		\textcolor{energyblue}{
			\bigl(\overline{r_A^\varepsilon}\phi_B^\varepsilon\bigr)^\prec
			+U_B^\varepsilon\bigl(\rho_A^\varepsilon\phi_B^\varepsilon\bigr)^\prec}
		\Bigr).
	\end{equation}
	Moreover,
	\begin{equation}
		\label{eq:higgs-second-low-high-splitting}
		iD_B^{\varepsilon *}
		\bigl(
		\overline{r_A^\varepsilon}Z_1^\varepsilon
		+U_B^\varepsilon\rho_A^\varepsilon Z_1^\varepsilon
		\bigr)
		=
		iD_B^{\varepsilon *}\Bigl(
		\textcolor{singularpink}{
			\bigl(
			\mathcal{U}_{>L_2}^\varepsilon\overline{r_A^\varepsilon}
			\bigr)\succ Z_1^\varepsilon
			+U_B^\varepsilon
			\Bigl(\bigl(\mathcal{U}_{>L_2}^\varepsilon\rho_A^\varepsilon\bigr)
			\succ Z_1^\varepsilon\Bigr)}
		\Bigr)
		+iD_B^{\varepsilon *}\Bigl(
		\textcolor{energyblue}{
			\bigl(\overline{r_A^\varepsilon}Z_1^\varepsilon\bigr)^\prec
			+U_B^\varepsilon\bigl(\rho_A^\varepsilon Z_1^\varepsilon\bigr)^\prec}
		\Bigr).
	\end{equation}
	In paraproduct expressions involving the operator
	$\rho_A^\varepsilon$, we use the componentwise convention
	\begin{equation}
		\label{eq:rho-A-operator-paraproduct-convention}
		\Bigl(
		(\mathcal{U}_{>L_2}^\varepsilon\rho_A^\varepsilon)
		\succ Z_1^\varepsilon
		\Bigr)_j
		:=
		\mathcal{U}_{>L_2}^\varepsilon
		r_{A,j}^\varepsilon
		\succ
		Z_1^\varepsilon(\,\cdot+\varepsilon e_j).
	\end{equation}
	
	\noindent\textup{(2)}
	We next decompose
	$:\!\lvert\Phi^\varepsilon\rvert^2\Phi^\varepsilon\!:_\varepsilon
	-\Phi^\varepsilon$.
	Set $q=3$. The binomial identity
	\cite[Lemma~2.22]{BC24} gives
	\begin{equation}
		\label{eq:higgs-wick-cubic-expansion}
		:\!\lvert\Phi^\varepsilon\rvert^2\Phi^\varepsilon\!:_\varepsilon
		-\Phi^\varepsilon
		=
		\sum_{\substack{
				0\leq k_1\leq\frac{q+1}{2}\\
				0\leq k_2\leq\frac{q-1}{2}\\
				k_1+k_2\neq0}}
		\binom{\frac{q+1}{2}}{k_1}
		\binom{\frac{q-1}{2}}{k_2}
		:\!(\phi_B^\varepsilon)^{k_1}
		(\overline{\phi_B^\varepsilon})^{k_2}\!:_\varepsilon
		(Z_1^\varepsilon)^{\frac{q+1}{2}-k_1}
		(\overline{Z_1^\varepsilon})^{\frac{q-1}{2}-k_2}
		-Z_1^\varepsilon-\phi_B^\varepsilon
		+\lvert Z_1^\varepsilon\rvert^2Z_1^\varepsilon.
	\end{equation}
	
	For almost every $t\in[0,1]$, the shift satisfies
	\begin{equation}
		\label{eq:higgs-feedback-equation}
		\bigl(\mathcal{L}_B^\varepsilon+1\bigr)\dot Z_t^\varepsilon
		=-
		\mathbb{E}_t^\phi\left[
		\nabla_\phi V_{A,B}(\Phi^\varepsilon)
		\right]
	\end{equation}
	We define $Z^{\mathrm{hi}}$ by the high-frequency equation below with
	$Z_0^{\mathrm{hi}}=0$ and set
	$Z_t^{\mathrm{lo}}:=Z_t^\varepsilon-Z_t^{\mathrm{hi}}$.
	The \textcolor{singularpink}{pink terms} enter the high-frequency equation:
	\begin{equation}
		\label{eq:higgs-high-feedback-equation}
		\begin{aligned}
			\textcolor{singularpink}{
			\bigl(\mathcal{L}_B^\varepsilon+1\bigr)\dot Z_t^{\mathrm{hi}}}
			={}&
			-
			\mathbb{E}_t^\phi\Bigg[
			\textcolor{singularpink}{
				iD_B^{\varepsilon *}\Bigl(
				\mathcal{U}_{>L_1}^\varepsilon
				\bigl(\overline{r_A^\varepsilon}\phi_B^\varepsilon\bigr)
				+U_B^\varepsilon\mathcal{U}_{>L_1}^\varepsilon
				(\rho_A^\varepsilon\phi_B^\varepsilon)
				\Bigr)}\\
			&\quad+
			\textcolor{singularpink}{
				iD_B^{\varepsilon *}\Bigl(
				\bigl(
				\mathcal{U}_{>L_2}^\varepsilon\overline{r_A^\varepsilon}
				\bigr)\succ Z_1^\varepsilon
				+U_B^\varepsilon\Bigl(
				(\mathcal{U}_{>L_2}^\varepsilon\rho_A^\varepsilon)
				\succ Z_1^\varepsilon
				\Bigr)
				\Bigr)}
			\Bigg].
		\end{aligned}
	\end{equation}
	
	The low-frequency equation is
	\begin{equation}
		\label{eq:higgs-low-feedback-equation}
		\begin{aligned}
			\textcolor{energyblue}{
			\bigl(\mathcal{L}_B^\varepsilon+1\bigr)\dot Z_t^{\mathrm{lo}}}
			={}&
			-
			\mathbb{E}_t^\phi\Bigg[
			\lvert Z_1^{\mathrm{lo}}\rvert^2Z_1^{\mathrm{lo}}
			+
			\Bigl(
			\lvert Z_1^\varepsilon\rvert^2Z_1^\varepsilon
			-\lvert Z_1^{\mathrm{lo}}\rvert^2Z_1^{\mathrm{lo}}
			\Bigr)
			+
			:\!\rho_A^{\varepsilon *}\rho_A^\varepsilon\!:_\varepsilon
			\phi_B^\varepsilon
			+
			:\!\rho_A^{\varepsilon *}\rho_A^\varepsilon\!:_\varepsilon
			Z_1^\varepsilon\\
			&\quad+
			\sum_{\substack{
					0\leq k_1\leq\frac{q+1}{2}\\
					0\leq k_2\leq\frac{q-1}{2}\\
					k_1+k_2\neq0}}
			\binom{\frac{q+1}{2}}{k_1}
			\binom{\frac{q-1}{2}}{k_2}
			:\!(\phi_B^\varepsilon)^{k_1}
			(\overline{\phi_B^\varepsilon})^{k_2}\!:_\varepsilon
			(Z_1^\varepsilon)^{\frac{q+1}{2}-k_1}
			(\overline{Z_1^\varepsilon})^{\frac{q-1}{2}-k_2}
			-Z_1^\varepsilon-\phi_B^\varepsilon\\
			&\quad+
			\textcolor{energyblue}{
				iD_B^{\varepsilon *}\Bigl(
				\bigl(
				\overline{r_A^\varepsilon}\phi_B^\varepsilon
				\bigr)^\prec
				+
				\bigl(
				\overline{r_A^\varepsilon}Z_1^\varepsilon
				\bigr)^\prec
				+U_B^\varepsilon
				\Bigl(
				(\rho_A^\varepsilon\phi_B^\varepsilon)^\prec
				+(\rho_A^\varepsilon Z_1^\varepsilon)^\prec
				\Bigr)
				\Bigr)}\\
			&\quad
			\textcolor{remainderorange}{
				-i\bigl(D^{\varepsilon *}\overline{r_A^\varepsilon}\bigr)
				\phi_B^\varepsilon
				-i\bigl(D^{\varepsilon *}\overline{r_A^\varepsilon}\bigr)
				Z_1^\varepsilon}
			\Bigg].
		\end{aligned}
	\end{equation}

	For the remainder of this section, the cutoffs $L_i$ and $\sigma_i$ are
	chosen as in Table~\ref{tab:higgs-high-frequency-cutoff-choices}.
	Recall that $\mathcal U_{>L_i}^\varepsilon$ in
	\eqref{eq:frequency-truncation-operators} also depends on the choice of
	$\sigma_i$.
	Their dependence on $A$, $B$, $\kappa$, $\alpha$, and $\varepsilon$
	is suppressed in the notation.
	\begin{table}[H]
		\centering
		\small
		\renewcommand{\arraystretch}{1.7}
		\setlength{\tabcolsep}{4pt}
		\begin{tabular}{p{0.12\textwidth}p{0.62\textwidth}p{0.20\textwidth}}
			\hline
			$\mathcal{U}_i$ & \textbf{Choice of $2^{L_i}$}
			& \textbf{Choice of $\sigma_i$}\\
			\hline
			$i=1$
			& $\displaystyle
			\left[
			\bigl(1+\lVert B\rVert_{H^{1,\varepsilon}(\rho)}\bigr)^{l_0}
			\lVert\Xi\rVert_\kappa
			\right]^{\frac{1}{\alpha/2-\kappa}}$
			& $\displaystyle\frac{4+\kappa}{\alpha/2-\kappa}$\\
			\hline
			$i=2$
			& $\displaystyle
			\left[
			\bigl(1+\lVert B\rVert_{H^{1,\varepsilon}(\rho)}\bigr)^{l_0}
			\lVert\Xi\rVert_\kappa
			\bigl(1+\lVert Z_1^\varepsilon\rVert_{L^{2,\varepsilon}(\rho^2)}\bigr)
			\right]^{\frac{1}{\alpha/2-\kappa}}$
			& $\displaystyle\frac{4+\kappa}{\alpha/2-\kappa}$\\
			\hline
		\end{tabular}
		\caption{Summary of the high-frequency cutoff choices.}
		\label{tab:higgs-high-frequency-cutoff-choices}
	\end{table}
	
	\begin{proposition}[Conditional estimates for the terminal Higgs-field shift]
		\label{prop:higgs-terminal-shift-summary}
		The following estimates hold:
		\begin{equation}
			\label{eq:higgs-summary-high-frequency-bound}
			\lVert Z_1^{\mathrm{hi}}\rVert_{
				L^\infty(\Omega^\phi;H_B^{1-\alpha,\varepsilon})}
			\lesssim_\alpha1.
		\end{equation}
		\begin{equation}
			\label{eq:higgs-summary-low-energy-bound}
			\mathbb{E}^\phi\lVert Z_1^{\mathrm{lo}}\rVert_{H_B^{1,\varepsilon}(\rho^2)}^2
			+\mathbb{E}^\phi\lVert Z_1^{\mathrm{lo}}\rVert_{L^{4,\varepsilon}(\rho)}^4
			\lesssim
			\bigl(1+\lVert B\rVert_{H^{1,\varepsilon}(\rho)}\bigr)^{\kappa_1}
			+\mathbb{E}^\phi\lVert\Xi\rVert_\kappa^{K_0}.
		\end{equation}
		\begin{equation}
			\label{eq:higgs-summary-low-schauder-bound}
			\mathbb{E}^\phi
			\lVert Z_1^{\mathrm{lo}}\rVert_{H_B^{1+\kappa_1,\varepsilon}(\rho^3)}^{\frac43(1-\nu)}
			\lesssim
			\bigl(1+\lVert B\rVert_{H^{1,\varepsilon}(\rho)}\bigr)^{\kappa_2}
			+
			\mathbb{E}^\phi\lVert\Xi\rVert_\kappa^{K_0}.
		\end{equation}
		Consequently,
		\begin{equation}
			\label{eq:higgs-summary-full-shift-bound}
			\begin{aligned}
				\mathbb{E}^\phi\left[
				\lVert Z_1\rVert_{H_B^{1-\alpha,\varepsilon}(\rho^2)}^2
				+\lVert Z_1\rVert_{L^{4,\varepsilon}(\rho)}^4
			\right]
				\lesssim
				\bigl(1+\lVert B\rVert_{H^{1,\varepsilon}(\rho)}\bigr)^{\kappa_1}
				+\mathbb{E}^\phi\lVert\Xi\rVert_\kappa^{K_0}.
			\end{aligned}
		\end{equation}
	\end{proposition}
	
	\begin{proof}
		The estimates
		\eqref{eq:higgs-summary-high-frequency-bound},
		\eqref{eq:higgs-summary-low-energy-bound}, and
		\eqref{eq:higgs-summary-low-schauder-bound} follow respectively from
		Propositions~\ref{prop:higgs-high-frequency-uniform-bounds},
		\ref{prop:higgs-low-frequency-uniform-bound}, and
		\ref{prop:higgs-low-schauder-bound}.
		Moreover,
		$H_B^{1,\varepsilon}\hookrightarrow
		H_B^{1-\alpha,\varepsilon}$, and
		Proposition~\ref{prop:higgs-high-frequency-uniform-bounds} with
		$p=4$ gives a uniform $L^{4,\varepsilon}$ bound for
		$Z_1^{\mathrm{hi}}$. Hence the decomposition
		$Z_1^\varepsilon=Z_1^{\mathrm{lo}}+Z_1^{\mathrm{hi}}$, the triangle
		inequality, and
		Proposition~\ref{prop:higgs-low-frequency-uniform-bound} yield
		\eqref{eq:higgs-summary-full-shift-bound}.
	\end{proof}
	
	\subsection{Estimate of
		\texorpdfstring{$Z_1^{\mathrm{hi}}$}{Z1 hi}}
	\label{subsec:estimate-higgs-high-component}
	
	In this subsection, we estimate the
	$H_B^{1-}$-norm of $Z_1^{\mathrm{hi}}$.
	Lemma~\ref{lem:weighted-covariant-resolvent} and Jensen's inequality
	give the following estimate.
	
	\begin{lemma}[Covariant FBSDE Schauder estimate]
		\label{lem:covariant-fbsde-schauder-estimate}
		Let $1\le p\le\infty$, $s\in[-1,1]$, $|a|\le l_0$, fix
		$B\in(\mathbb{R}^2)^{\Lambda_\varepsilon}$, and let
		$f\in L^p\bigl(\Omega^\phi;H_B^{-1+s,\varepsilon}(\rho^a)\bigr)$
		be $\mathcal{F}_1^\phi$-measurable. Set
		\begin{equation}
			\label{eq:covariant-fbsde-schauder-system}
			Z_t
			:=-\bigl(\mathcal{L}_B^\varepsilon+1\bigr)^{-1}
			\int_0^t\mathbb{E}_r^\phi[f]\,dr,
			\qquad t\in[0,1].
		\end{equation}
		Then
		\begin{equation}
			\label{eq:covariant-fbsde-schauder-estimate}
			\left\lVert Z_1\right\rVert_{
				L^p\bigl(
				\Omega^\phi;H_B^{1+s,\varepsilon}(\rho^a)
				\bigr)}
			\lesssim
			\left\lVert f\right\rVert_{
				L^p\bigl(
				\Omega^\phi;H_B^{-1+s,\varepsilon}(\rho^a)
				\bigr)}.
		\end{equation}
	\end{lemma}
	
	We now use Lemma~\ref{lem:covariant-fbsde-schauder-estimate} to
	estimate $Z_1^{\mathrm{hi}}$, taking $p=\infty$, $a=0$, and $s=-\alpha$.
	Thus the terminal forcing inside the conditional expectation
	in~\eqref{eq:higgs-high-feedback-equation} will be estimated in
	$H_B^{-1-\alpha,\varepsilon}$.
	
	\smallskip
	\noindent\textup{(1)}
	We first estimate
	\begin{equation}
		\label{eq:higgs-high-first-divergence-term}
		iD_B^{\varepsilon *}\Bigl(
		\mathcal{U}_{>L_1}^\varepsilon
		\bigl(\overline{r_A^\varepsilon}\phi_B^\varepsilon\bigr)
		{}+U_B^\varepsilon\mathcal{U}_{>L_1}^\varepsilon
		(\rho_A^\varepsilon\phi_B^\varepsilon)
		\Bigr).
	\end{equation}
	By~\eqref{eq:weighted-covariant-negative-divergence}, its
	$H_B^{-1-\alpha,\varepsilon}$-norm is bounded by a constant multiple of
	\begin{equation}
		\label{eq:higgs-high-first-divergence-reduction}
		\bigl(1+\lVert B\rVert_{H^{1,\varepsilon}(\rho)}\bigr)^{3\alpha}
		\Bigl(
		\lVert\mathcal U_{>L_1}^\varepsilon
		(\overline{r_A^\varepsilon}\phi_B^\varepsilon)
		\rVert_{H^{-\alpha,\varepsilon}(\rho^{-3\alpha})}
		+\lVert\mathcal U_{>L_1}^\varepsilon
		(\rho_A^\varepsilon\phi_B^\varepsilon)
		\rVert_{H^{-\alpha,\varepsilon}(\rho^{-3\alpha})}\Bigr).
	\end{equation}
	Since $\rho^{4-3\alpha}\in L^{2,\varepsilon}$, we have the uniform
	embedding $\mathcal C^{-\alpha/2,\varepsilon}(\rho^{-4})
	\hookrightarrow H^{-\alpha,\varepsilon}(\rho^{-3\alpha})$.
	Using Lemma~\ref{lem:frequency-truncation-estimates} with
	$\sigma_1=(4+\kappa)/(\alpha/2-\kappa)$ gives
	\begin{equation}
		\label{eq:higgs-high-rho-phi-frequency-estimate}
		\begin{aligned}
			\lVert\mathcal U_{>L_1}^\varepsilon
			(\overline{r_A^\varepsilon}\phi_B^\varepsilon)
			\rVert_{H^{-\alpha,\varepsilon}(\rho^{-3\alpha})}
			\lesssim
			\lVert\mathcal U_{>L_1}^\varepsilon
			(\overline{r_A^\varepsilon}\phi_B^\varepsilon)
			\rVert_{\mathcal C^{-\alpha/2,\varepsilon}(\rho^{-4})}
			\lesssim
			2^{-(\alpha/2-\kappa)L_1}
			\lVert\overline{r_A^\varepsilon}\phi_B^\varepsilon
			\rVert_{\mathcal C^{-\kappa,\varepsilon}(\rho^\kappa)}.
		\end{aligned}
	\end{equation}
	The second term in~\eqref{eq:higgs-high-first-divergence-reduction}
	is estimated in the same way. Hence the
	$H_B^{-1-\alpha,\varepsilon}$-norm of
	\eqref{eq:higgs-high-first-divergence-term} is bounded by a constant multiple of
	\begin{equation}
		\label{eq:higgs-high-first-divergence-estimate}
		\bigl(1+\lVert B\rVert_{H^{1,\varepsilon}(\rho)}\bigr)^{3\alpha}
		2^{-(\alpha/2-\kappa)L_1}\lVert\Xi\rVert_\kappa.
	\end{equation}
	Choose $L_1$ as in Table~\ref{tab:higgs-high-frequency-cutoff-choices}.
	This cutoff is $\mathcal F_1^\phi$-measurable. Since
	$\lVert\Xi\rVert_\kappa\ge1$ and $l_0>3\alpha$, we obtain
	\begin{equation}
		\label{eq:higgs-high-first-divergence-uniform-bound}
		\left\lVert iD_B^{\varepsilon *}\Bigl(
		\mathcal U_{>L_1}^\varepsilon
		(\overline{r_A^\varepsilon}\phi_B^\varepsilon)
		+U_B^\varepsilon\mathcal U_{>L_1}^\varepsilon
		(\rho_A^\varepsilon\phi_B^\varepsilon)
		\Bigr)\right\rVert_{H_B^{-1-\alpha,\varepsilon}}
		\lesssim1.
	\end{equation}
	
	\smallskip
	\noindent\textup{(2)}
	We next estimate
	\begin{equation}
		\label{eq:higgs-high-second-divergence-term}
		iD_B^{\varepsilon *}\Bigl(
		\bigl(
		\mathcal{U}_{>L_2}^\varepsilon\overline{r_A^\varepsilon}
		\bigr)\succ Z_1^\varepsilon
		+U_B^\varepsilon
		\bigl(
		(\mathcal{U}_{>L_2}^\varepsilon\rho_A^\varepsilon)
		\succ Z_1^\varepsilon
		\bigr)
		\Bigr).
	\end{equation}
	By~\eqref{eq:weighted-covariant-negative-divergence}, its
	$H_B^{-1-\alpha,\varepsilon}$-norm is bounded by a constant multiple of
	\begin{equation}
		\label{eq:higgs-high-second-divergence-reduction}
		\bigl(1+\lVert B\rVert_{H^{1,\varepsilon}(\rho)}\bigr)^{3\alpha}
		\Bigl(
		\lVert(\mathcal U_{>L_2}^\varepsilon\overline{r_A^\varepsilon})
		\succ Z_1^\varepsilon
		\rVert_{H^{-\alpha,\varepsilon}(\rho^{-3\alpha})}
		+\lVert(\mathcal U_{>L_2}^\varepsilon\rho_A^\varepsilon)
		\succ Z_1^\varepsilon
		\rVert_{H^{-\alpha,\varepsilon}(\rho^{-3\alpha})}\Bigr).
	\end{equation}
	Using	$B_{2,\infty}^{-\alpha/2,\varepsilon}(\rho^{-2})
	\hookrightarrow H^{-\alpha,\varepsilon}(\rho^{-3\alpha})$, ~\eqref{eq:discrete-paraproduct-prec-estimate} and
	Lemma~\ref{lem:frequency-truncation-estimates} with
	$\sigma_2=(4+\kappa)/(\alpha/2-\kappa)$, we obtain
	\begin{equation}
		\label{eq:higgs-high-A-Z-frequency-estimate}
		\begin{aligned}
			\lVert(\mathcal U_{>L_2}^\varepsilon\overline{r_A^\varepsilon})
			\succ Z_1^\varepsilon
			\rVert_{H^{-\alpha,\varepsilon}(\rho^{-3\alpha})}
			&\lesssim
			\lVert\mathcal U_{>L_2}^\varepsilon\overline{r_A^\varepsilon}
			\rVert_{\mathcal C^{-\alpha/2,\varepsilon}(\rho^{-4})}
			\lVert Z_1^\varepsilon\rVert_{L^{2,\varepsilon}(\rho^2)}\\
			&\lesssim
			2^{-(\alpha/2-\kappa)L_2}
			\lVert\overline{r_A^\varepsilon}
			\rVert_{\mathcal C^{-\kappa,\varepsilon}(\rho^\kappa)}
			\lVert Z_1^\varepsilon\rVert_{L^{2,\varepsilon}(\rho^2)}.
		\end{aligned}
	\end{equation}
	The second term in~\eqref{eq:higgs-high-second-divergence-term}
	is estimated in the same way.
	Hence the $H_B^{-1-\alpha,\varepsilon}$-norm of
	\eqref{eq:higgs-high-second-divergence-term} is bounded by a constant multiple of
	\begin{equation}
		\label{eq:higgs-high-second-divergence-estimate}
		\bigl(1+\lVert B\rVert_{H^{1,\varepsilon}(\rho)}\bigr)^{3\alpha}
		2^{-(\alpha/2-\kappa)L_2}\lVert\Xi\rVert_\kappa
		\lVert Z_1^\varepsilon\rVert_{L^{2,\varepsilon}(\rho^2)}.
	\end{equation}
	Choose $L_2$ as in Table~\ref{tab:higgs-high-frequency-cutoff-choices}.
	We conclude that
	\begin{equation}
		\label{eq:higgs-high-second-divergence-uniform-bound}
		\left\lVert iD_B^{\varepsilon *}\Bigl(
		(\mathcal U_{>L_2}^\varepsilon\overline{r_A^\varepsilon})
		\succ Z_1^\varepsilon
		+U_B^\varepsilon
		\bigl((\mathcal U_{>L_2}^\varepsilon\rho_A^\varepsilon)
		\succ Z_1^\varepsilon\bigr)
		\Bigr)\right\rVert_{H_B^{-1-\alpha,\varepsilon}}
		\lesssim1.
	\end{equation}
	
	\begin{proposition}
		\label{prop:higgs-high-frequency-uniform-bounds}
		With the cutoffs in Table~\ref{tab:higgs-high-frequency-cutoff-choices},
		the following estimates hold:
		\begin{equation}
			\label{eq:higgs-high-uniform-bound}
			\lVert Z_1^{\mathrm{hi}}\rVert_{
				L^\infty(\Omega^\phi;H_B^{1-\alpha,\varepsilon})}
			\lesssim_\alpha1.
		\end{equation}
		In particular, for every $2\le p\le 2/\alpha$,
		\begin{equation}
			\label{eq:higgs-high-Lp-uniform-bound}
			\lVert Z_1^{\mathrm{hi}}\rVert_{
				L^\infty(\Omega^\phi;L^{p,\varepsilon})}
			\lesssim_{\alpha,p}1.
		\end{equation}
	\end{proposition}
	
	\begin{proof}
		Estimates~\eqref{eq:higgs-high-first-divergence-uniform-bound} and
		\eqref{eq:higgs-high-second-divergence-uniform-bound}, together with
		Lemma~\ref{lem:covariant-fbsde-schauder-estimate}, give
		\eqref{eq:higgs-high-uniform-bound}.
		The covariant Sobolev embedding
		\eqref{eq:covariant-sobolev-embedding} gives
		\eqref{eq:higgs-high-Lp-uniform-bound}.
	\end{proof}
	
	\subsection{Estimate of
		\texorpdfstring{$Z_1^{\mathrm{lo}}$}{Z1 lo}}
	\label{subsec:estimate-higgs-low-component}
	
	We rewrite the equation for $Z^{\mathrm{lo}}$ as
	\begin{equation}
		\label{eq:higgs-low-feedback-abbreviated}
		\begin{aligned}
			\textcolor{energyblue}{
			\bigl(\mathcal{L}_B^\varepsilon+1\bigr)\dot Z_t^{\mathrm{lo}}}
			={}&
			-
			\mathbb{E}_t^\phi\left[
			\lvert Z_1^{\mathrm{lo}}\rvert^2Z_1^{\mathrm{lo}}
			+\sum_{j=1}^2 I_j^{\mathrm{blue}}
			+\sum_{k=1}^3 J_k
			+R^{\mathrm{orange}}
			\right].
		\end{aligned}
	\end{equation}
	Here the blue terms are
	\begin{equation}
		\label{eq:higgs-low-blue-terms}
		I_1^{\mathrm{blue}}
		=
		\textcolor{energyblue}{
			iD_B^{\varepsilon *}\Bigl(
			\bigl(
			\overline{r_A^\varepsilon}\phi_B^\varepsilon
			\bigr)^\prec
			+U_B^\varepsilon
			(\rho_A^\varepsilon\phi_B^\varepsilon)^\prec
			\Bigr)},
		\qquad
		I_2^{\mathrm{blue}}
		=
		\textcolor{energyblue}{
			iD_B^{\varepsilon *}\Bigl(
			\bigl(
			\overline{r_A^\varepsilon}Z_1^\varepsilon
			\bigr)^\prec
			+U_B^\varepsilon
			(\rho_A^\varepsilon Z_1^\varepsilon)^\prec
			\Bigr)}.
	\end{equation}
	We also set
	\begin{equation}
		\label{eq:higgs-low-J1}
		J_1
		=
		\sum_{\substack{
				0\leq k_1\leq\frac{q+1}{2}\\
				0\leq k_2\leq\frac{q-1}{2}\\
				k_1+k_2\neq0}}
		\binom{\frac{q+1}{2}}{k_1}
		\binom{\frac{q-1}{2}}{k_2}
		:\!(\phi_B^\varepsilon)^{k_1}
		(\overline{\phi_B^\varepsilon})^{k_2}\!:_\varepsilon
		(Z_1^\varepsilon)^{\frac{q+1}{2}-k_1}
		(\overline{Z_1^\varepsilon})^{\frac{q-1}{2}-k_2}
		-Z_1^\varepsilon-\phi_B^\varepsilon,
	\end{equation}
	\begin{equation}
		\label{eq:higgs-low-J2-J3}
		J_2
		=
		:\!\rho_A^{\varepsilon *}\rho_A^\varepsilon\!:_\varepsilon
		\phi_B^\varepsilon
		+
		:\!\rho_A^{\varepsilon *}\rho_A^\varepsilon\!:_\varepsilon
		Z_1^\varepsilon,
		\qquad
		J_3
		=
		\lvert Z_1^\varepsilon\rvert^2Z_1^\varepsilon
		-\lvert Z_1^{\mathrm{lo}}\rvert^2Z_1^{\mathrm{lo}}.
	\end{equation}
	Finally, the orange remainder is
	\begin{equation}
		\label{eq:higgs-low-orange-remainder}
		R^{\mathrm{orange}}
		=
		\textcolor{remainderorange}{
			-i\bigl(D^{\varepsilon *}\overline{r_A^\varepsilon}\bigr)
			\phi_B^\varepsilon
			-i\bigl(D^{\varepsilon *}\overline{r_A^\varepsilon}\bigr)
			Z_1^\varepsilon}.
	\end{equation}
	
	\begin{lemma}[Energy estimate]
		\label{lem:higgs-low-energy-estimate}
		Fix $B\in(\mathbb R^2)^{\Lambda_\varepsilon}$. Let $Z$ be adapted
		with almost surely absolutely continuous paths in
		$H_B^{1,\varepsilon}(\rho^2)$ and $Z_0=0$, and suppose that
		\begin{equation}
			\label{eq:higgs-low-feedback-control}
			(\mathcal L_B^\varepsilon+1)\dot Z_t
			=-\mathbb E_t^\phi\bigl[|Z_1|^2Z_1+f\bigr]
		\end{equation}
		for almost every $t$, almost surely.
		Assume that $f$ is $\mathcal F_1^\phi$-measurable and
		\begin{equation}
			\label{eq:higgs-low-moment-assumptions}
			\mathbb E^\phi
			\bigl\lVert |Z_1|^3+|f|\bigr\rVert_{L^{2,\varepsilon}(\rho^2)}^2
			<\infty.
		\end{equation}
		Then
		\begin{equation}
			\label{eq:higgs-low-energy-estimate}
			\frac12\mathbb E^\phi\left[
			\lVert Z_1\rVert_{H_B^{1,\varepsilon}(\rho^2)}^2
			+\lVert Z_1\rVert_{L^{4,\varepsilon}(\rho)}^4
			\right]
			\le-\mathbb E^\phi\langle f,\rho^4Z_1\rangle_\varepsilon.
		\end{equation}
	\end{lemma}
	
	\begin{proof}
		By~\eqref{eq:weighted-covariant-resolvent} and conditional Jensen,
		assumption~\eqref{eq:higgs-low-moment-assumptions} justifies the testing below.
		Testing~\eqref{eq:higgs-low-feedback-control} with
		$\rho^4\dot Z_t$, using adaptedness and
		$\int_0^1\dot Z_t\,dt=Z_1$, gives
		\begin{equation}
			\label{eq:higgs-low-testing-identity}
			\begin{aligned}
				\mathbb E^\phi\int_0^1
				\langle(\mathcal L_B^\varepsilon+1)\dot Z_t,
				\rho^4\dot Z_t\rangle_\varepsilon\,dt
=-\mathbb E^\phi\int_0^1
				\left\langle\mathbb E_t^\phi[|Z_1|^2Z_1+f],
				\rho^4\dot Z_t\right\rangle_\varepsilon\,dt
=-\mathbb E^\phi\lVert Z_1\rVert_{L^{4,\varepsilon}(\rho)}^4
				-\mathbb E^\phi\langle f,\rho^4Z_1\rangle_\varepsilon.
			\end{aligned}
		\end{equation}
		By Jensen's inequality and~\eqref{eq:weighted-covariant-coercivity},
		\begin{equation}
			\label{eq:higgs-low-H1-control-bound}
			\frac12\mathbb E^\phi
			\lVert Z_1\rVert_{H_B^{1,\varepsilon}(\rho^2)}^2
			\le\frac12\mathbb E^\phi\int_0^1
			\lVert\dot Z_t\rVert_{H_B^{1,\varepsilon}(\rho^2)}^2\,dt
			\le\mathbb E^\phi\int_0^1
			\langle(\mathcal L_B^\varepsilon+1)\dot Z_t,
			\rho^4\dot Z_t\rangle_\varepsilon\,dt.
		\end{equation}
		Combining these estimates   yields~\eqref{eq:higgs-low-energy-estimate}.
	\end{proof}
	
	Fix $M,\varepsilon,A,B$. The terminal
	law~\eqref{eq:app-optimal-terminal-law}, Lemma~\ref{lem:app-cor-satisfies-hypothesis} and Gaussian
	moments give
	$\mathbb E^\phi\lVert Z_1^\varepsilon\rVert_{L^{2,\varepsilon}(\Lambda_{M,\varepsilon})}^m<\infty$
	for every $1\le m<\infty$.
	Its periodicity, the fixed-lattice embedding
	$L^{2,\varepsilon}(\Lambda_{M,\varepsilon})
	\hookrightarrow L^{6,\varepsilon}(\Lambda_{M,\varepsilon})$,
	and the summability of $\rho^4$ give
	$Z_1^\varepsilon\in L^m(\Omega^\phi;L^{6,\varepsilon}(\rho^{2/3}))$.
	Proposition~\ref{prop:higgs-high-frequency-uniform-bounds} and
	$L^{2,\varepsilon}(\Lambda_\varepsilon)
	\hookrightarrow L^{6,\varepsilon}(\Lambda_\varepsilon)$
	give the same integrability for $Z_1^{\mathrm{hi}}$.
	Hence $Z_1^{\mathrm{lo}}=Z_1^\varepsilon-Z_1^{\mathrm{hi}}$
	belongs to this space for every $1\le m<\infty$.

	In~\eqref{eq:higgs-low-feedback-abbreviated}, set
	$f:=\sum_{j=1}^2 I_j^{\mathrm{blue}}
	+\sum_{k=1}^3 J_k+R^{\mathrm{orange}}$.
	At fixed $\varepsilon$, the equivalence between the \(H_B^{-1-\alpha,\varepsilon}\) and \(L^{2,\varepsilon}\) norms and
	estimates~\eqref{eq:higgs-high-first-divergence-uniform-bound},
	\eqref{eq:higgs-high-second-divergence-uniform-bound} give
	$L^{2,\varepsilon}(\Lambda_\varepsilon)$-bounds for the two high-frequency forcings.
	Since the full gradient
	$\nabla_\phi V_{A,B}(\Phi^\varepsilon)$ is periodic and satisfies
	$\mathbb E^\phi\lVert\nabla_\phi V_{A,B}(\Phi^\varepsilon)
	\rVert_{L^{2,\varepsilon}(\Lambda_{M,\varepsilon})}^m<\infty$
	for every $1\le m<\infty$, subtracting these forcings and
	$|Z_1^{\mathrm{lo}}|^2Z_1^{\mathrm{lo}}$ gives
$	f\in L^2\bigl(\Omega^\phi;L^{2,\varepsilon}(\rho^2)\bigr).$
	Moreover, $f$ is $\mathcal F_1^\phi$-measurable. Applying
	Lemma~\ref{lem:higgs-low-energy-estimate} to
	\eqref{eq:higgs-low-feedback-abbreviated} therefore gives
	\begin{equation}
		\label{eq:higgs-low-final-energy-estimate}
		\begin{aligned}
			\frac12\mathbb E^\phi\left[
			\lVert Z_1^{\mathrm{lo}}\rVert_{H_B^{1,\varepsilon}(\rho^2)}^2
			+\lVert Z_1^{\mathrm{lo}}\rVert_{L^{4,\varepsilon}(\rho)}^4
			\right]
\le\mathbb E^\phi\left[
			\sum_{j=1}^2
			\bigl|\langle I_j^{\mathrm{blue}},\rho^4Z_1^{\mathrm{lo}}\rangle_\varepsilon\bigr|
			+\sum_{k=1}^3
			\bigl|\langle J_k,\rho^4Z_1^{\mathrm{lo}}\rangle_\varepsilon\bigr|
			+\bigl|\langle R^{\mathrm{orange}},\rho^4Z_1^{\mathrm{lo}}\rangle_\varepsilon\bigr|
			\right].
		\end{aligned}
	\end{equation}
	
	\smallskip
	\noindent\textup{(1)} \textit{Estimate of
		$\mathbb E^\phi|\langle I_1^{\mathrm{blue}},
		\rho^4Z_1^{\mathrm{lo}}\rangle_\varepsilon|$.}
	By the definition of $I_1^{\mathrm{blue}}$ and $|U_{B,j}^\varepsilon|=1$,
	\begin{equation}
		\label{eq:higgs-low-I1-pairing-splitting}
			\bigl|\langle I_1^{\mathrm{blue}},\rho^4Z_1^{\mathrm{lo}}
			\rangle_\varepsilon\bigr|
			\lesssim
			\Bigl(
			\lVert(\overline{r_A^\varepsilon}\phi_B^\varepsilon)^\prec
			\rVert_{L^{2,\varepsilon}(\rho^2)}
			+\lVert(\rho_A^\varepsilon\phi_B^\varepsilon)^\prec
			\rVert_{L^{2,\varepsilon}(\rho^2)}
			\Bigr)
			\lVert\rho^{-2}D_B^\varepsilon(\rho^4Z_1^{\mathrm{lo}})
			\rVert_{L^{2,\varepsilon}}.
	\end{equation}
	By~\eqref{eq:covariant-leibniz-rules}
	and~\eqref{eq:covariant-derivative-estimates},
	\begin{equation}
		\label{eq:higgs-low-I1-weighted-derivative-bound}
		\begin{aligned}
			\lVert\rho^{-2}D_B^\varepsilon(\rho^4Z_1^{\mathrm{lo}})
			\rVert_{L^{2,\varepsilon}}
			\le
			\lVert D_B^\varepsilon(\rho^2Z_1^{\mathrm{lo}})
			\rVert_{L^{2,\varepsilon}}
			+\lVert[\rho^{-2},D_B^\varepsilon](\rho^4Z_1^{\mathrm{lo}})
			\rVert_{L^{2,\varepsilon}}
			\lesssim
			\lVert Z_1^{\mathrm{lo}}\rVert_{H_B^{1,\varepsilon}(\rho^2)}.
		\end{aligned}
	\end{equation}
	The parameter choices in Subsection~\ref{subsec:parameters} give
	$K_0\gg1$ and
	$2\kappa\sigma_1\lesssim\kappa/\alpha\ll\kappa_1$,
	so $1/4-2\kappa\sigma_1\ge\kappa$.
	Since $\rho^{7/4}\in L^{2,\varepsilon}$, we have
	$\mathcal C^{\kappa,\varepsilon}(\rho^{1/4})
	\hookrightarrow L^{2,\varepsilon}(\rho^2)$.
	Using~\eqref{eq:low-frequency-truncation-estimate} and
	Table~\ref{tab:higgs-high-frequency-cutoff-choices}, we obtain
	\begin{equation}
		\label{eq:higgs-low-I1-low-frequency-localization}
			\lVert(\overline{r_A^\varepsilon}\phi_B^\varepsilon)^\prec
			\rVert_{L^{2,\varepsilon}(\rho^2)}
			\lesssim 2^{2\kappa L_1}
			\lVert\overline{r_A^\varepsilon}\phi_B^\varepsilon
			\rVert_{\mathcal C^{-\kappa,\varepsilon}
			(\rho^{1/4-2\kappa\sigma_1})}
			\lesssim
			\bigl(1+\lVert B\rVert_{H^{1,\varepsilon}(\rho)}\bigr)^{\kappa_1/4}
			\lVert\Xi\rVert_\kappa^{K_0/4}.
	\end{equation}
	The second term is estimated in the same way.
	Combining these bounds with~\eqref{eq:higgs-low-I1-pairing-splitting}
	and applying Young's inequality gives
	\begin{equation}
		\label{eq:higgs-low-I1-final-bound}
			\mathbb E^\phi
			\bigl|\langle I_1^{\mathrm{blue}},\rho^4Z_1^{\mathrm{lo}}
			\rangle_\varepsilon\bigr|
			\le
			\frac1{100}\mathbb E^\phi
			\lVert Z_1^{\mathrm{lo}}\rVert_{H_B^{1,\varepsilon}(\rho^2)}^2
			+C\Bigl(\lVert B\rVert_{H^{1,\varepsilon}(\rho)}^{\kappa_1}
			+\mathbb E^\phi\lVert\Xi\rVert_\kappa^{K_0}\Bigr).
	\end{equation}

	\smallskip
	\noindent\textup{(2)} \textit{Estimate of
		$\mathbb E^\phi|\langle I_2^{\mathrm{blue}},
		\rho^4Z_1^{\mathrm{lo}}\rangle_\varepsilon|$.}
	By the definition of $I_2^{\mathrm{blue}}$
	and~\eqref{eq:higgs-low-I1-weighted-derivative-bound},
	\begin{equation}
		\label{eq:higgs-low-I2-pairing-bound}
			\mathbb E^\phi
			\bigl|\langle I_2^{\mathrm{blue}},\rho^4Z_1^{\mathrm{lo}}
			\rangle_\varepsilon\bigr|\\
			\lesssim\mathbb E^\phi\left[
			\Bigl(
			\lVert(\overline{r_A^\varepsilon}Z_1^\varepsilon)^\prec
			\rVert_{L^{2,\varepsilon}(\rho^2)}
			+\lVert(\rho_A^\varepsilon Z_1^\varepsilon)^\prec
			\rVert_{L^{2,\varepsilon}(\rho^2)}
			\Bigr)
			\lVert Z_1^{\mathrm{lo}}\rVert_{H_B^{1,\varepsilon}(\rho^2)}
			\right].
	\end{equation}
	Since $(\overline{r_A^\varepsilon}Z_1^\varepsilon)^\prec
	=(\mathcal U_{\le L_2}^\varepsilon\overline{r_A^\varepsilon})
	\succ Z_1^\varepsilon
	+\overline{r_A^\varepsilon}\preccurlyeq Z_1^\varepsilon$,
	Lemma~\ref{lem:discrete-paraproduct-estimates} and
	\eqref{eq:low-frequency-truncation-estimate} give
	\begin{equation}
		\label{eq:higgs-low-AZ-regular-remainder-decomposition}
			\lVert(\overline{r_A^\varepsilon}Z_1^\varepsilon)^\prec
			\rVert_{L^{2,\varepsilon}(\rho^2)}
			\lesssim 2^{2\kappa L_2}
			\lVert r_A^\varepsilon\rVert_{\mathcal C^{-\kappa,\varepsilon}
			(\rho^{1/4-2\kappa\sigma_2})}
			\lVert Z_1^\varepsilon\rVert_{L^{2,\varepsilon}(\rho^{7/4})}\\
			+
			\lVert r_A^\varepsilon\rVert_{\mathcal C^{-\kappa,\varepsilon}(\rho^\kappa)}
			\lVert Z_1^\varepsilon\rVert_{H^{2\kappa,\varepsilon}(\rho^{2-\kappa})}.
	\end{equation}
	The parameter choices give $1/4-2\kappa\sigma_2\ge\kappa$ and
	$\frac{2\kappa}{\alpha/2-\kappa}
	\le\frac{2l_0\kappa}{\alpha/2-\kappa}\le\kappa_1/16$.
	Since $\rho^{3/4}\in L^{4,\varepsilon}$,
	$L^{4,\varepsilon}(\rho)\hookrightarrow L^{2,\varepsilon}(\rho^{7/4})$.
	Thus Table~\ref{tab:higgs-high-frequency-cutoff-choices} and
	Proposition~\ref{prop:higgs-high-frequency-uniform-bounds} yield
	\begin{equation}
		\label{eq:higgs-low-AZ-low-frequency-bound}
			2^{2\kappa L_2}
			\lVert r_A^\varepsilon\rVert_{\mathcal C^{-\kappa,\varepsilon}
			(\rho^{1/4-2\kappa\sigma_2})}
			\lVert Z_1^\varepsilon\rVert_{L^{2,\varepsilon}(\rho^{7/4})}
			\lesssim
			\lVert\Xi\rVert_\kappa^2
			\bigl(1+\lVert Z_1^{\mathrm{lo}}\rVert_{L^{2,\varepsilon}(\rho^{7/4})}\bigr)^{1+\kappa_1/16}
			\bigl(1+\lVert B\rVert_{H^{1,\varepsilon}(\rho)}^{\kappa_1/16}\bigr).
	\end{equation}
	For the second term, \eqref{eq:weighted-HB-H-comparison} and
	Proposition~\ref{prop:higgs-high-frequency-uniform-bounds} give
	\begin{equation}
		\label{eq:higgs-low-AZ-covariant-interpolation-bound}
			\lVert r_A^\varepsilon\rVert_{\mathcal C^{-\kappa,\varepsilon}(\rho^\kappa)}
			\lVert Z_1^\varepsilon\rVert_{H^{2\kappa,\varepsilon}(\rho^{2-\kappa})}
			\lesssim
			\lVert\Xi\rVert_\kappa
			\bigl(1+\lVert B\rVert_{H^{1,\varepsilon}(\rho)}\bigr)^{2\kappa}
			\bigl(1+\lVert Z_1^{\mathrm{lo}}
			\rVert_{H_B^{2\kappa,\varepsilon}(\rho^{2-3\kappa})}\bigr).
	\end{equation}
	Since $(2-7\kappa)/(1-2\kappa)\ge7/4$,
	Lemma~\ref{lem:covariant-interpolation-inequality} and the same embedding yield
	\begin{equation}
		\label{eq:higgs-low-AZ-weighted-interpolation}
			\lVert Z_1^{\mathrm{lo}}
			\rVert_{H_B^{2\kappa,\varepsilon}(\rho^{2-3\kappa})}\lesssim
			\lVert Z_1^{\mathrm{lo}}\rVert_{H_B^{1,\varepsilon}(\rho^2)}^{2\kappa}
			\lVert Z_1^{\mathrm{lo}}
			\rVert_{L^{2,\varepsilon}(\rho^{\frac{2-7\kappa}{1-2\kappa}})}^{1-2\kappa}
			\lesssim
			\lVert Z_1^{\mathrm{lo}}\rVert_{H_B^{1,\varepsilon}(\rho^2)}^{2\kappa}
			\lVert Z_1^{\mathrm{lo}}\rVert_{L^{4,\varepsilon}(\rho)}^{1-2\kappa}.
	\end{equation}
	The term involving $\rho_A^\varepsilon$ is treated in the same way,
	using~\eqref{eq:rho-A-operator-paraproduct-convention} and
	$\rho(x+\varepsilon e_j)\asymp\rho(x)$.
	Combining the preceding estimates and applying Young's inequality
	with the parameter choices in Subsection~\ref{subsec:parameters} gives
	\begin{equation}
		\label{eq:higgs-low-I2-final-bound}
		\begin{aligned}
			&\Bigl(
			\lVert(\overline{r_A^\varepsilon}Z_1^\varepsilon)^\prec
			\rVert_{L^{2,\varepsilon}(\rho^2)}
			+\lVert(\rho_A^\varepsilon Z_1^\varepsilon)^\prec
			\rVert_{L^{2,\varepsilon}(\rho^2)}
			\Bigr)\lVert Z_1^{\mathrm{lo}}\rVert_{H_B^{1,\varepsilon}(\rho^2)}\\
			&\quad\le C\lVert\Xi\rVert_\kappa^2
			\bigl(1+\lVert Z_1^{\mathrm{lo}}\rVert_{L^{4,\varepsilon}(\rho)}\bigr)^{1+\kappa_1/16}
			\bigl(1+\lVert B\rVert_{H^{1,\varepsilon}(\rho)}^{\kappa_1/16}\bigr)
			\lVert Z_1^{\mathrm{lo}}\rVert_{H_B^{1,\varepsilon}(\rho^2)}\\
			&\quad+C\lVert\Xi\rVert_\kappa
			\bigl(1+\lVert B\rVert_{H^{1,\varepsilon}(\rho)}\bigr)^{2\kappa}
			\Bigl(1+
			\lVert Z_1^{\mathrm{lo}}\rVert_{H_B^{1,\varepsilon}(\rho^2)}^{2\kappa}
			\lVert Z_1^{\mathrm{lo}}\rVert_{L^{4,\varepsilon}(\rho)}^{1-2\kappa}
			\Bigr)
			\lVert Z_1^{\mathrm{lo}}\rVert_{H_B^{1,\varepsilon}(\rho^2)}\\
			&\quad\le\frac1{100}\left(
			\lVert Z_1^{\mathrm{lo}}\rVert_{H_B^{1,\varepsilon}(\rho^2)}^2
			+\lVert Z_1^{\mathrm{lo}}\rVert_{L^{4,\varepsilon}(\rho)}^4
			\right)
			+C\Bigl(\lVert B\rVert_{H^{1,\varepsilon}(\rho)}^{\kappa_1}
			+\lVert\Xi\rVert_\kappa^{K_0}\Bigr).
		\end{aligned}
	\end{equation}
	Taking $\mathbb E^\phi$ and using~\eqref{eq:higgs-low-I2-pairing-bound}
	completes the estimate.

	\smallskip
	\noindent\textup{(3)} \textit{Estimate of
		$\mathbb{E}^\phi\bigl[
		\lvert\langle J_1,\rho^4Z_1^{\mathrm{lo}}
		\rangle_\varepsilon\rvert\bigr]$.}
	Recalling that $q=3$, the definition of $J_1$ gives
	\begin{equation}
		\label{eq:higgs-low-J1-pairing-splitting}
		\begin{aligned}
			\mathbb{E}^\phi\left[
			\left\lvert
			\left\langle J_1,\rho^4Z_1^{\mathrm{lo}}
			\right\rangle_\varepsilon
			\right\rvert\right]
			&\lesssim
			\sum_{\substack{
					0\leq k_1\leq\frac{q+1}{2}\\
					0\leq k_2\leq\frac{q-1}{2}\\
					k_1+k_2\neq0}}
			\mathbb{E}^\phi\left[
			\left\lvert
			\left\langle
			:\!(\phi_B^\varepsilon)^{k_1}
			(\overline{\phi_B^\varepsilon})^{k_2}\!:_\varepsilon
			(Z_1^\varepsilon)^{\frac{q+1}{2}-k_1}
			(\overline{Z_1^\varepsilon})^{\frac{q-1}{2}-k_2},
			\rho^4Z_1^{\mathrm{lo}}
			\right\rangle_\varepsilon
			\right\rvert\right]\\
			&\quad+
			\mathbb{E}^\phi\left[
			\left\lvert
			\left\langle Z_1^\varepsilon,\rho^4Z_1^{\mathrm{lo}}
			\right\rangle_\varepsilon
			\right\rvert\right]
			+
			\mathbb{E}^\phi\left[
			\left\lvert
			\left\langle\phi_B^\varepsilon,\rho^4Z_1^{\mathrm{lo}}
			\right\rangle_\varepsilon
			\right\rvert\right].
		\end{aligned}
	\end{equation}
	For a fixed pair $(k_1,k_2)\neq(0,0)$ in the preceding sum, discrete
	Besov duality yields
	\begin{equation}
		\label{eq:higgs-low-J1-Besov-duality}
		\begin{aligned}
			&\mathbb{E}^\phi\left[
			\left\lvert
			\left\langle
			:\!(\phi_B^\varepsilon)^{k_1}
			(\overline{\phi_B^\varepsilon})^{k_2}\!:_\varepsilon
			(Z_1^\varepsilon)^{\frac{q+1}{2}-k_1}
			(\overline{Z_1^\varepsilon})^{\frac{q-1}{2}-k_2},
			\rho^4Z_1^{\mathrm{lo}}
			\right\rangle_\varepsilon
			\right\rvert\right]\\
			&\quad\lesssim
			\mathbb{E}^\phi\left[
			\left\lVert
			:\!(\phi_B^\varepsilon)^{k_1}
			(\overline{\phi_B^\varepsilon})^{k_2}\!:_\varepsilon
			\right\rVert_{\mathcal{C}^{-\kappa,\varepsilon}(\rho^\kappa)}
			\left\lVert
			(\overline{Z_1^\varepsilon})^{\frac{q+1}{2}-k_1}
			(Z_1^\varepsilon)^{\frac{q-1}{2}-k_2}
			Z_1^{\mathrm{lo}}
			\right\rVert_{B_{1,1}^{\kappa,\varepsilon}(\rho^{4-\kappa})}
			\right].
		\end{aligned}
	\end{equation}
	The regular product in~\eqref{eq:higgs-low-J1-Besov-duality} contains
	$q-k_1-k_2+1\leq q$ factors. Since $\rho^{4-\kappa}\in L^{1,\varepsilon}$,
	we may insert $k_1+k_2-1$ constant factors $1$ in
	$B_{q,\infty}^{2\kappa,\varepsilon}(\rho^{(4-\kappa)/q})$.
	Lemma~\ref{lem:discrete-paraproduct-estimates} then gives
	\begin{equation}
		\label{eq:higgs-low-J1-product-Besov-bound}
		\left\lVert
		(\overline{Z_1^\varepsilon})^{\frac{q+1}{2}-k_1}
		(Z_1^\varepsilon)^{\frac{q-1}{2}-k_2}
		Z_1^{\mathrm{lo}}
		\right\rVert_{B_{1,1}^{\kappa,\varepsilon}(\rho^{4-\kappa})}
		\lesssim
		\left\lVert Z_1^\varepsilon
		\right\rVert_{B_{q,\infty}^{2\kappa,\varepsilon}(\rho^{(4-\kappa)/q})}^{
			q-k_1-k_2}
		\left\lVert Z_1^{\mathrm{lo}}
		\right\rVert_{B_{q,\infty}^{2\kappa,\varepsilon}(\rho^{(4-\kappa)/q})}.
	\end{equation}
	Interpolation, weighted H\"older's inequality,
	\eqref{eq:weighted-HB-H-comparison}, and
	Proposition~\ref{prop:higgs-high-frequency-uniform-bounds} imply
	\begin{equation}
		\label{eq:higgs-low-J1-high-Besov-bound}
		\left\lVert Z_1^{\mathrm{hi}}
		\right\rVert_{B_{q,\infty}^{2\kappa,\varepsilon}(\rho^{(4-\kappa)/q})}
		\lesssim
		\left\lVert Z_1^{\mathrm{hi}}
		\right\rVert_{H^{1-\alpha,\varepsilon}(\rho^{1-\alpha})}^{
			\frac{2\kappa}{1-\alpha}}
		\left\lVert Z_1^{\mathrm{hi}}
		\right\rVert_{L^{2/\alpha,\varepsilon}}^{
			1-\frac{2\kappa}{1-\alpha}}
		\lesssim_\alpha
		\bigl(1+\lVert B\rVert_{H^{1,\varepsilon}(\rho)}\bigr)^{2\kappa}.
	\end{equation}
	Here $\alpha/2+\kappa<1/q$, and the remaining weight is integrable.
	Similarly, using $q+1=4$, interpolation and weighted H\"older's
	inequality yield
	\begin{equation}
		\label{eq:higgs-low-J1-low-Besov-bound}
		\begin{aligned}
			\left\lVert Z_1^{\mathrm{lo}}
			\right\rVert_{B_{q,\infty}^{2\kappa,\varepsilon}(\rho^{(4-\kappa)/q})}
			&\lesssim
			\left\lVert Z_1^{\mathrm{lo}}
			\right\rVert_{H^{1,\varepsilon}(\rho^3)}^{2\kappa}
			\left\lVert Z_1^{\mathrm{lo}}
			\right\rVert_{L^{q+1,\varepsilon}(\rho)}^{1-2\kappa}\\
			&\lesssim
			\bigl(1+\lVert B\rVert_{H^{1,\varepsilon}(\rho)}\bigr)^{2\kappa}
			\left\lVert Z_1^{\mathrm{lo}}
			\right\rVert_{H_B^{1,\varepsilon}(\rho^2)}^{2\kappa}
			\left\lVert Z_1^{\mathrm{lo}}
			\right\rVert_{L^{4,\varepsilon}(\rho)}^{1-2\kappa}.
		\end{aligned}
	\end{equation}
	Combining
	\eqref{eq:higgs-low-J1-product-Besov-bound}--%
	\eqref{eq:higgs-low-J1-low-Besov-bound}, we obtain
	\begin{equation}
		\label{eq:higgs-low-J1-regular-product-bound}
		\left\lVert
		(\overline{Z_1^\varepsilon})^{\frac{q+1}{2}-k_1}
		(Z_1^\varepsilon)^{\frac{q-1}{2}-k_2}
		Z_1^{\mathrm{lo}}
		\right\rVert_{B_{1,1}^{\kappa,\varepsilon}(\rho^{4-\kappa})}
		\lesssim_\alpha
		\bigl(1+\lVert Z_1^{\mathrm{lo}}
		\rVert_{H_B^{1,\varepsilon}(\rho^2)}\bigr)^{2\kappa q}
		\bigl(1+\lVert Z_1^{\mathrm{lo}}
		\rVert_{L^{4,\varepsilon}(\rho)}\bigr)^{(1-2\kappa)q}
		\bigl(1+\lVert B\rVert_{H^{1,\varepsilon}(\rho)}\bigr)^{4\kappa q}.
	\end{equation}
	By~\eqref{eq:stochastic-object-Xi-norm}, the first factor on the
	right-hand side of~\eqref{eq:higgs-low-J1-Besov-duality} is bounded by
	$\lVert\Xi\rVert_\kappa$. Hence  Young's inequality gives
	\begin{equation}
		\label{eq:higgs-low-J1-Wick-pairing-bound}
		\begin{aligned}
			&\mathbb{E}^\phi\left[
			\left\lVert
			:\!(\phi_B^\varepsilon)^{k_1}
			(\overline{\phi_B^\varepsilon})^{k_2}\!:_\varepsilon
			\right\rVert_{\mathcal{C}^{-\kappa,\varepsilon}(\rho^\kappa)}
			\left\lVert
			(\overline{Z_1^\varepsilon})^{\frac{q+1}{2}-k_1}
			(Z_1^\varepsilon)^{\frac{q-1}{2}-k_2}
			Z_1^{\mathrm{lo}}
			\right\rVert_{B_{1,1}^{\kappa,\varepsilon}(\rho^{4-\kappa})}
			\right]\\
			&\quad\leq
			\frac{1}{1000}\mathbb{E}^\phi
			\left[
			\left\lVert Z_1^{\mathrm{lo}}
			\right\rVert_{H_B^{1,\varepsilon}(\rho^2)}^2
			+
			\left\lVert Z_1^{\mathrm{lo}}
			\right\rVert_{L^{4,\varepsilon}(\rho)}^4
			\right]
			+
			C\lVert B\rVert_{H^{1,\varepsilon}(\rho)}^{\kappa_1}
			+C\mathbb{E}^\phi\lVert\Xi\rVert_\kappa^{K_0}.
		\end{aligned}
	\end{equation}
	The two linear terms in
	\eqref{eq:higgs-low-J1-pairing-splitting} are simpler and satisfy
	the same estimate.
	Summing over the finitely many pairs $(k_1,k_2)$ and reducing the
	constants in Young's inequality, we conclude that
	\begin{equation}
		\label{eq:higgs-low-J1-final-bound}
		\mathbb{E}^\phi\left[
		\left\lvert
		\left\langle J_1,\rho^4Z_1^{\mathrm{lo}}
		\right\rangle_\varepsilon
		\right\rvert\right]
		\leq
		\frac{1}{100}\mathbb{E}^\phi
		\left\lVert Z_1^{\mathrm{lo}}
		\right\rVert_{H_B^{1,\varepsilon}(\rho^2)}^2
		+
		\frac{1}{100}\mathbb{E}^\phi
		\left\lVert Z_1^{\mathrm{lo}}
		\right\rVert_{L^{4,\varepsilon}(\rho)}^4
		+
		C\lVert B\rVert_{H^{1,\varepsilon}(\rho)}^{\kappa_1}
		+C\mathbb{E}^\phi\lVert\Xi\rVert_\kappa^{K_0}.
	\end{equation}
	
	\smallskip
	\noindent\textup{(4)} \textit{Estimate of
		$\mathbb{E}^\phi\bigl[
		\lvert\langle J_2,\rho^4Z_1^{\mathrm{lo}}
		\rangle_\varepsilon\rvert\bigr]$ and
		$\mathbb{E}^\phi\bigl[
		\lvert\langle R^{\mathrm{orange}},\rho^4Z_1^{\mathrm{lo}}
		\rangle_\varepsilon\rvert\bigr]$.}
By the definitions in~\eqref{eq:higgs-low-J2-J3}
and~\eqref{eq:higgs-low-orange-remainder}, the stochastic factors
appearing in $J_2$ and $R^{\mathrm{orange}}$ are included in $\Xi$.
Estimates analogous to those in (3) therefore yield the same bound
as in~\eqref{eq:higgs-low-J1-final-bound}.

	\smallskip
	\noindent\textup{(5)} \textit{Estimate of
		$\mathbb{E}^\phi\bigl[
		\lvert\langle J_3,\rho^4Z_1^{\mathrm{lo}}
		\rangle_\varepsilon\rvert\bigr]$.}
	Expanding $J_3$ in~\eqref{eq:higgs-low-J2-J3}, we have
	\begin{equation}
		\label{eq:higgs-low-J3-expansion}
		J_3
		=
		\sum_{\substack{
				0\leq k_1\leq\frac{q+1}{2}\\
				0\leq k_2\leq\frac{q-1}{2}\\
				k_1+k_2\neq0}}
		\binom{\frac{q+1}{2}}{k_1}
		\binom{\frac{q-1}{2}}{k_2}
		(Z_1^{\mathrm{hi}})^{k_1}
		(\overline{Z_1^{\mathrm{hi}}})^{k_2}
		(Z_1^{\mathrm{lo}})^{\frac{q+1}{2}-k_1}
		(\overline{Z_1^{\mathrm{lo}}})^{\frac{q-1}{2}-k_2}.
	\end{equation}
	For each term, H\"older's inequality, $\rho\leq1$, and
	\eqref{eq:higgs-high-Lp-uniform-bound} with $p=q+1=4$ give
	\begin{equation}
		\label{eq:higgs-low-J3-term-bound}
		\mathbb{E}^\phi\left[
		\left\lvert
		\left\langle
		(Z_1^{\mathrm{hi}})^{k_1}
		(\overline{Z_1^{\mathrm{hi}}})^{k_2}
		(Z_1^{\mathrm{lo}})^{\frac{q+1}{2}-k_1}
		(\overline{Z_1^{\mathrm{lo}}})^{\frac{q-1}{2}-k_2},
		\rho^4Z_1^{\mathrm{lo}}
		\right\rangle_\varepsilon
		\right\rvert\right]
		\lesssim_\alpha
		\mathbb{E}^\phi
		\left\lVert Z_1^{\mathrm{lo}}
		\right\rVert_{L^{4,\varepsilon}(\rho)}^{4-k_1-k_2}
		\leq
		\frac{1}{100}\mathbb{E}^\phi
		\left\lVert Z_1^{\mathrm{lo}}
		\right\rVert_{L^{4,\varepsilon}(\rho)}^4
		+C.
	\end{equation}
	Summing the finitely many terms in
	\eqref{eq:higgs-low-J3-expansion} and reducing the constants in
	Young's inequality, we conclude that
	\begin{equation}
		\label{eq:higgs-low-J3-final-bound}
		\mathbb{E}^\phi\left[
		\left\lvert
		\left\langle J_3,\rho^4Z_1^{\mathrm{lo}}
		\right\rangle_\varepsilon
		\right\rvert\right]
		\leq
		\frac{1}{100}\mathbb{E}^\phi
		\left\lVert Z_1^{\mathrm{lo}}
		\right\rVert_{L^{4,\varepsilon}(\rho)}^4
		+C.
	\end{equation}
	
	\begin{proposition}
		\label{prop:higgs-low-frequency-uniform-bound}We have 
		\begin{equation}
			\label{eq:higgs-low-frequency-uniform-bound}
				\mathbb{E}^\phi\lVert Z_1^{\mathrm{lo}}\rVert_{H_B^{1,\varepsilon}(\rho^2)}^2
				+\mathbb{E}^\phi\lVert Z_1^{\mathrm{lo}}\rVert_{L^{4,\varepsilon}(\rho)}^4
				\lesssim
				\bigl(1+\lVert B\rVert_{H^{1,\varepsilon}(\rho)}\bigr)^{\kappa_1}
				+\mathbb{E}^\phi\lVert\Xi\rVert_\kappa^{K_0}.
		\end{equation}
	\end{proposition}
	
	\begin{proof}
		Combining \eqref{eq:higgs-low-final-energy-estimate} with estimates
		\textup{(1)}--\textup{(5)} above and absorbing the corresponding
		$\frac{1}{100}$-terms into the left-hand side yields
		\eqref{eq:higgs-low-frequency-uniform-bound}.
	\end{proof}
	
	\subsection{Schauder estimate of
		\texorpdfstring{$Z_1^{\mathrm{lo}}$}{Z1 lo}}
	\label{subsec:higgs-low-schauder-estimate}
	
	We apply Lemma~\ref{lem:covariant-fbsde-schauder-estimate} to
	\eqref{eq:higgs-low-feedback-abbreviated} with
	$s=\alpha$ and $a=16$. Thus, estimates of the source terms in
	$H_B^{-1+\alpha,\varepsilon}(\rho^{16})$ yield an estimate of
	$Z_1^{\mathrm{lo}}$ in $H_B^{1+\alpha,\varepsilon}(\rho^{16})$.
	Interpolation with the $H_B^{1,\varepsilon}(\rho^2)$-estimate in
	Proposition~\ref{prop:higgs-low-frequency-uniform-bound} then gives
	an estimate in $H_B^{1+\kappa_1,\varepsilon}(\rho^3)$. We first estimate the
	two blue terms in~\eqref{eq:higgs-low-blue-terms}.
	
	\smallskip
	\noindent\textup{(1)} \textit{Estimate of $I_1^{\mathrm{blue}}$.}
	By~\eqref{eq:weighted-covariant-negative-divergence} with $s=-\alpha$
	and $\vartheta_{-\alpha}=\alpha$, it suffices to estimate the regular
	remainders in $H^{\alpha,\varepsilon}(\rho^{16-\alpha})$.
	Since $\rho^{7-\alpha}\in L^{2,\varepsilon}$,
	$\mathcal C^{\alpha+\kappa,\varepsilon}(\rho^9)
	\hookrightarrow H^{\alpha,\varepsilon}(\rho^{16-\alpha})$.
	Table~\ref{tab:higgs-high-frequency-cutoff-choices} gives
	$(\alpha+2\kappa)\sigma_i\le9-\kappa$, $i=1,2$.
	Thus~\eqref{eq:low-frequency-truncation-estimate} yields
	\begin{equation}
		\label{eq:higgs-low-schauder-I1-free-bound}
		\begin{aligned}
			\left\lVert
			\bigl(\overline{r_A^\varepsilon}\phi_B^\varepsilon\bigr)^\prec
			\right\rVert_{H^{\alpha,\varepsilon}(\rho^{16-\alpha})}
			&\lesssim 2^{(\alpha+2\kappa)L_1}
			\left\lVert\overline{r_A^\varepsilon}\phi_B^\varepsilon
			\right\rVert_{\mathcal C^{-\kappa,\varepsilon}
				(\rho^{9-(\alpha+2\kappa)\sigma_1})}\\
			&\lesssim
			\left[
			\bigl(1+\lVert B\rVert_{H^{1,\varepsilon}(\rho)}\bigr)^{l_0}
			\lVert\Xi\rVert_\kappa
			\right]^{\frac{\alpha+2\kappa}{\alpha/2-\kappa}}
			\lVert\Xi\rVert_\kappa.
		\end{aligned}
	\end{equation}
	The term involving $\rho_A^\varepsilon$ is estimated in the same way.
	Hence~\eqref{eq:weighted-covariant-negative-divergence} and Young's inequality give
	\begin{equation}
		\label{eq:higgs-low-schauder-I1-bound}
		\begin{aligned}
			\lVert I_1^{\mathrm{blue}}
			\rVert_{H_B^{-1+\alpha,\varepsilon}(\rho^{16})}
			&\lesssim
			\bigl(1+\lVert B\rVert_{H^{1,\varepsilon}(\rho)}\bigr)^\alpha
			\left(
			\left\lVert
			\bigl(\overline{r_A^\varepsilon}\phi_B^\varepsilon\bigr)^\prec
			\right\rVert_{H^{\alpha,\varepsilon}(\rho^{16-\alpha})}
			+
			\left\lVert(\rho_A^\varepsilon\phi_B^\varepsilon)^\prec
			\right\rVert_{H^{\alpha,\varepsilon}(\rho^{16-\alpha})}
			\right)\\
			&\lesssim
			\bigl(1+\lVert B\rVert_{H^{1,\varepsilon}(\rho)}\bigr)^{3l_0}
			+\lVert\Xi\rVert_\kappa^{K_0/2}.
		\end{aligned}
	\end{equation}
	
	\smallskip
	\noindent\textup{(2)} \textit{Estimate of $I_2^{\mathrm{blue}}$.}
	The same derivative estimate gives
	\begin{equation}
		\label{eq:higgs-low-schauder-I2-derivative-bound}
		\begin{aligned}
			\lVert I_2^{\mathrm{blue}}
			\rVert_{H_B^{-1+\alpha,\varepsilon}(\rho^{16})}\lesssim
			\bigl(1+\lVert B\rVert_{H^{1,\varepsilon}(\rho)}\bigr)^\alpha
			\left(
			\left\lVert
			\bigl(\overline{r_A^\varepsilon}Z_1^\varepsilon\bigr)^\prec
			\right\rVert_{H^{\alpha,\varepsilon}(\rho^{16-\alpha})}
			+
			\left\lVert(\rho_A^\varepsilon Z_1^\varepsilon)^\prec
			\right\rVert_{H^{\alpha,\varepsilon}(\rho^{16-\alpha})}
			\right).
		\end{aligned}
	\end{equation}
	The weighted paraproduct estimates,
	\eqref{eq:low-frequency-truncation-estimate},
	\eqref{eq:weighted-HB-H-comparison}, and
	Proposition~\ref{prop:higgs-high-frequency-uniform-bounds} give
	\begin{equation}
		\label{eq:higgs-low-schauder-AZ-remainder-bound}
		\begin{aligned}
			\left\lVert
			\bigl(\overline{r_A^\varepsilon}Z_1^\varepsilon\bigr)^\prec
			\right\rVert_{H^{\alpha,\varepsilon}(\rho^{16-\alpha})}
			&\le
			\left\lVert
			\bigl(\mathcal U_{\leq L_2}^\varepsilon
			\overline{r_A^\varepsilon}\bigr)\succ Z_1^\varepsilon
			\right\rVert_{H^{\alpha,\varepsilon}(\rho^{16-\alpha})}
			+
			\left\lVert\overline{r_A^\varepsilon}\preccurlyeq Z_1^\varepsilon
			\right\rVert_{H^{\alpha,\varepsilon}(\rho^{16-\alpha})}\\
			&\lesssim
			\left\lVert\mathcal U_{\leq L_2}^\varepsilon
			\overline{r_A^\varepsilon}
			\right\rVert_{\mathcal C^{\alpha+\kappa,\varepsilon}(\rho^9)}
			\lVert Z_1^\varepsilon\rVert_{L^{2,\varepsilon}(\rho^{7-\alpha})}
			+
			\lVert r_A^\varepsilon\rVert_{\mathcal C^{-\kappa,\varepsilon}(\rho^\kappa)}
			\lVert Z_1^\varepsilon
			\rVert_{H^{\alpha+2\kappa,\varepsilon}(\rho^{16-\alpha-\kappa})}\\
			&\lesssim
			\left[
			\bigl(1+\lVert B\rVert_{H^{1,\varepsilon}(\rho)}\bigr)^{l_0}
			\lVert\Xi\rVert_\kappa
			\bigl(1+\lVert Z_1^{\mathrm{lo}}\rVert_{L^{4,\varepsilon}(\rho)}\bigr)
			\right]^{\frac{\alpha+2\kappa}{\alpha/2-\kappa}}
			\lVert\Xi\rVert_\kappa
			\bigl(1+\lVert Z_1^{\mathrm{lo}}\rVert_{L^{4,\varepsilon}(\rho)}\bigr)\\
			&+
			\lVert\Xi\rVert_\kappa
			\bigl(1+\lVert B\rVert_{H^{1,\varepsilon}(\rho)}\bigr)^{\alpha+2\kappa}
			\bigl(1+\lVert Z_1^{\mathrm{lo}}\rVert_{H_B^{1,\varepsilon}(\rho^2)}\bigr)\\
			&\lesssim
			\bigl(1+\lVert B\rVert_{H^{1,\varepsilon}(\rho)}\bigr)^{3l_0}
			\lVert\Xi\rVert_\kappa^4
			\left[
			1+\lVert Z_1^{\mathrm{lo}}\rVert_{L^{4,\varepsilon}(\rho)}^{
				1+\frac{\alpha+2\kappa}{\alpha/2-\kappa}}
			+\lVert Z_1^{\mathrm{lo}}\rVert_{H_B^{1,\varepsilon}(\rho^2)}
			\right].
		\end{aligned}
	\end{equation}
	Here we used $L^{4,\varepsilon}(\rho)\hookrightarrow
	L^{2,\varepsilon}(\rho^2)\cap L^{2,\varepsilon}(\rho^{7-\alpha})$ and
	$16-2\alpha-3\kappa\ge2$.
	The term involving $\rho_A^\varepsilon$ satisfies the same bound by
	\eqref{eq:rho-A-operator-paraproduct-convention} and
	$\rho(x+\varepsilon e_j)\asymp\rho(x)$.
	
	By the choice in Subsection~\ref{subsec:parameters},
	\begin{equation}
		\label{eq:higgs-low-schauder-exponent-comparison}
		1+\frac{\alpha+2\kappa}{\alpha/2-\kappa}
		=
		3+\frac{8\kappa}{\alpha-2\kappa}
		\leq 3+\kappa_1.
	\end{equation}
	We use Young's inequality with
	\begin{equation}
		\label{eq:higgs-low-schauder-young-exponents}
		p_2:=
		\frac{3}{(3+\kappa_1)(1-\nu)},
		\qquad
		p_1:=
		\frac{3}{3-(3+\kappa_1)(1-\nu)}
		\leq\frac{3}{2\nu}
		\ll \frac{1}{\alpha}\ll K_0.
	\end{equation}
	In particular, for $a,b\ge0$,
$
		\label{eq:higgs-low-schauder-young-application}
		a b^{(3+\kappa_1)\frac{4}{3}(1-\nu)}
		\lesssim
		a^{\frac{3}{2\nu}}+b^4+1.
$
	Combining
	\eqref{eq:higgs-low-schauder-I2-derivative-bound}--
	\eqref{eq:higgs-low-schauder-young-exponents}, we obtain
	\begin{equation}
		\label{eq:higgs-low-schauder-I2-moment-bound}
		\begin{aligned}
			\mathbb{E}^\phi\left[
			\left\lVert I_2^{\mathrm{blue}}
			\right\rVert_{H_B^{-1+\alpha,\varepsilon}(\rho^{16})}^{
				\frac{4}{3}(1-\nu)}
			\right]
			&\lesssim
			\mathbb{E}^\phi\left[
			\bigl(1+\lVert B\rVert_{H^{1,\varepsilon}(\rho)}\bigr)^{5l_0}
			\lVert\Xi\rVert_\kappa^8
			\left(
			\left\lVert Z_1^{\mathrm{lo}}\right\rVert_{L^{4,\varepsilon}(\rho)}^{
				(3+\kappa_1)\frac{4}{3}(1-\nu)}
			+
			\left\lVert Z_1^{\mathrm{lo}}\right\rVert_{H_B^{1,\varepsilon}(\rho^2)}^{
				\frac{4}{3}(1-\nu)}
			+1
			\right)
			\right]\\
			&\lesssim
			\mathbb{E}^\phi\left[
			\left\lVert Z_1^{\mathrm{lo}}\right\rVert_{L^{4,\varepsilon}(\rho)}^4
			+
			\left\lVert Z_1^{\mathrm{lo}}\right\rVert_{H_B^{1,\varepsilon}(\rho^2)}^2
			\right]
			+1+\lVert B\rVert_{H^{1,\varepsilon}(\rho)}^{1/\alpha}
			+\mathbb{E}^\phi\lVert\Xi\rVert_\kappa^{K_0}\\
			&\lesssim
			1+\lVert B\rVert_{H^{1,\varepsilon}(\rho)}^{1/\alpha}
			+\mathbb{E}^\phi\lVert\Xi\rVert_\kappa^{K_0},
		\end{aligned}
	\end{equation}
	where the last inequality follows from
	Proposition~\ref{prop:higgs-low-frequency-uniform-bound}.
	
	\smallskip
	\noindent\textup{(3)} \textit{Estimate of $J_1$.}
	By~\eqref{eq:higgs-low-J1},
	\begin{equation}
		\label{eq:higgs-low-schauder-J1-splitting}
		\begin{aligned}
			\lVert J_1\rVert_{H_B^{-1+\alpha,\varepsilon}(\rho^{16})}
			&\lesssim
			\sum_{\substack{
					0\leq k_1\leq\frac{q+1}{2}\\
					0\leq k_2\leq\frac{q-1}{2}\\
					k_1+k_2\neq0}}
			\left\lVert
			:\!(\phi_B^\varepsilon)^{k_1}
			(\overline{\phi_B^\varepsilon})^{k_2}\!:_\varepsilon
			(Z_1^\varepsilon)^{\frac{q+1}{2}-k_1}
			(\overline{Z_1^\varepsilon})^{\frac{q-1}{2}-k_2}
			\right\rVert_{H_B^{-1+\alpha,\varepsilon}(\rho^{16})}\\
			&\quad+
			\lVert Z_1^\varepsilon\rVert_{H_B^{-1+\alpha,\varepsilon}(\rho^{16})}
			+
			\lVert\phi_B^\varepsilon\rVert_{H_B^{-1+\alpha,\varepsilon}(\rho^{16})}.
		\end{aligned}
	\end{equation}
	For every pair $(k_1,k_2)$ in the preceding sum, the discrete
	paraproduct estimates give
	\begin{equation}
		\label{eq:higgs-low-schauder-J1-Wick-product-bound}
		\begin{aligned}
			&\left\lVert
			:\!(\phi_B^\varepsilon)^{k_1}
			(\overline{\phi_B^\varepsilon})^{k_2}\!:_\varepsilon
			(Z_1^\varepsilon)^{\frac{q+1}{2}-k_1}
			(\overline{Z_1^\varepsilon})^{\frac{q-1}{2}-k_2}
			\right\rVert_{B_{4/3,\infty}^{-2\kappa,\varepsilon}(\rho^{3+\kappa})}\\
			&\quad\lesssim
			\left\lVert
			:\!(\phi_B^\varepsilon)^{k_1}
			(\overline{\phi_B^\varepsilon})^{k_2}\!:_\varepsilon
			\right\rVert_{\mathcal{C}^{-\kappa,\varepsilon}(\rho^\kappa)}
			\left\lVert
			(Z_1^\varepsilon)^{\frac{q+1}{2}-k_1}
			(\overline{Z_1^\varepsilon})^{\frac{q-1}{2}-k_2}
			\right\rVert_{B_{4/3,\infty}^{2\kappa,\varepsilon}(\rho^3)}.
		\end{aligned}
	\end{equation}
	Since $q=3$ and the second factor contains at most two copies of
	$Z_1^\varepsilon$ and $\overline{Z_1^\varepsilon}$, inserting constant
	factors in $B_{8/3,\infty}^{2\kappa,\varepsilon}(\rho^{3/2})$ as needed gives
	\begin{equation}
		\label{eq:higgs-low-schauder-J1-regular-product-bound}
		\left\lVert
		(Z_1^\varepsilon)^{\frac{q+1}{2}-k_1}
		(\overline{Z_1^\varepsilon})^{\frac{q-1}{2}-k_2}
		\right\rVert_{B_{4/3,\infty}^{2\kappa,\varepsilon}(\rho^3)}
		\lesssim
		1+
		\left\lVert Z_1^\varepsilon
		\right\rVert_{B_{8/3,\infty}^{2\kappa,\varepsilon}(\rho^{3/2})}^2.
	\end{equation}
	Weighted interpolation,
	\eqref{eq:weighted-HB-H-comparison}, and
	Proposition~\ref{prop:higgs-high-frequency-uniform-bounds} yield
	\begin{equation}
		\label{eq:higgs-low-schauder-Z-Besov-bound}
		\begin{aligned}
			\left\lVert Z_1^\varepsilon
			\right\rVert_{B_{8/3,\infty}^{2\kappa,\varepsilon}(\rho^{3/2})}
			&\lesssim
			\left\lVert Z_1^\varepsilon
			\right\rVert_{B_{8/3,4}^{5\kappa/2,\varepsilon}(\rho^{3/2})}\\
			&\lesssim
			\left\lVert Z_1^\varepsilon
			\right\rVert_{H^{5\kappa,\varepsilon}(\rho^2)}^{1/2}
			\left\lVert Z_1^\varepsilon
			\right\rVert_{B_{4,\infty}^{0,\varepsilon}(\rho)}^{1/2}\\
			&\lesssim
			\left\lVert Z_1^\varepsilon
			\right\rVert_{H^{5\kappa,\varepsilon}(\rho^2)}^{1/2}
			\left\lVert Z_1^\varepsilon
			\right\rVert_{L^{4,\varepsilon}(\rho)}^{1/2}\\
			&\lesssim
			\bigl(1+\lVert B\rVert_{H^{1,\varepsilon}(\rho)}\bigr)^{5\kappa/2}
			\bigl(1+\lVert Z_1^{\mathrm{lo}}
			\rVert_{H_B^{1,\varepsilon}(\rho^2)}\bigr)^{5\kappa/2}
			\bigl(1+\lVert Z_1^{\mathrm{lo}}
			\rVert_{L^{4,\varepsilon}(\rho)}\bigr)^{(2-5\kappa)/2}.
		\end{aligned}
	\end{equation}
	Here we used $L^{4,\varepsilon}(\rho)\hookrightarrow
	L^{2,\varepsilon}(\rho^{7/4})$ and
	$2-5\kappa\ge7/4+5\kappa/4$ in
	Lemma~\ref{lem:covariant-interpolation-inequality}.
	Combining
	\eqref{eq:higgs-low-schauder-J1-Wick-product-bound}--
	\eqref{eq:higgs-low-schauder-Z-Besov-bound} with
	\eqref{eq:stochastic-object-Xi-norm}, we obtain
	\begin{equation}
		\label{eq:higgs-low-schauder-J1-Besov-bound}
		\begin{aligned}
			&\left\lVert
			:\!(\phi_B^\varepsilon)^{k_1}
			(\overline{\phi_B^\varepsilon})^{k_2}\!:_\varepsilon
			(Z_1^\varepsilon)^{\frac{q+1}{2}-k_1}
			(\overline{Z_1^\varepsilon})^{\frac{q-1}{2}-k_2}
			\right\rVert_{B_{4/3,\infty}^{-2\kappa,\varepsilon}(\rho^{3+\kappa})}\\
			&\quad\lesssim
			\lVert\Xi\rVert_\kappa
			\bigl(1+\lVert B\rVert_{H^{1,\varepsilon}(\rho)}\bigr)^{5\kappa}
			\bigl(1+\lVert Z_1^{\mathrm{lo}}
			\rVert_{H_B^{1,\varepsilon}(\rho^2)}\bigr)^{5\kappa}
			\bigl(1+\lVert Z_1^{\mathrm{lo}}
			\rVert_{L^{4,\varepsilon}(\rho)}\bigr)^{2-5\kappa}.
		\end{aligned}
	\end{equation}
	The uniform embedding
	$B_{4/3,\infty}^{-2\kappa,\varepsilon}(\rho^{3+\kappa})
	\hookrightarrow H^{-1/2-3\alpha,\varepsilon}(\rho^{3+\kappa})$,
	\eqref{eq:weighted-HB-H-comparison}, and
	$3+\kappa+1/2+3\alpha<16$ imply
	\begin{equation}
		\label{eq:higgs-low-schauder-Besov-to-covariant-bound}
		\lVert f\rVert_{H_B^{-1+\alpha,\varepsilon}(\rho^{16})}
		\lesssim
		\bigl(1+\lVert B\rVert_{H^{1,\varepsilon}(\rho)}\bigr)^{1/2+3\alpha}
		\lVert f\rVert_{B_{4/3,\infty}^{-2\kappa,\varepsilon}(\rho^{3+\kappa})}.
	\end{equation}
	Consequently, Young's inequality gives
	\begin{equation}
		\label{eq:higgs-low-schauder-J1-Wick-moment-bound}
		\begin{aligned}
			&\mathbb{E}^\phi\left[
			\left\lVert
			:\!(\phi_B^\varepsilon)^{k_1}
			(\overline{\phi_B^\varepsilon})^{k_2}\!:_\varepsilon
			(Z_1^\varepsilon)^{\frac{q+1}{2}-k_1}
			(\overline{Z_1^\varepsilon})^{\frac{q-1}{2}-k_2}
			\right\rVert_{H_B^{-1+\alpha,\varepsilon}(\rho^{16})}^{
				\frac{4}{3}(1-\nu)}
			\right]\\
			&\quad\lesssim
			\mathbb{E}^\phi\left[
			\left\lVert Z_1^{\mathrm{lo}}
			\right\rVert_{L^{4,\varepsilon}(\rho)}^4
			+
			\left\lVert Z_1^{\mathrm{lo}}
			\right\rVert_{H_B^{1,\varepsilon}(\rho^2)}^2
			\right]
			+\lVert B\rVert_{H^{1,\varepsilon}(\rho)}^{1/\alpha}
			+\mathbb{E}^\phi\lVert\Xi\rVert_\kappa^{K_0}.
		\end{aligned}
	\end{equation}
	The remaining two linear terms in~\eqref{eq:higgs-low-schauder-J1-splitting}
	are simpler. Proposition~\ref{prop:higgs-high-frequency-uniform-bounds},
	\eqref{eq:stochastic-object-Xi-norm}, and Young's inequality therefore
	yield
	\begin{equation}
		\label{eq:higgs-low-schauder-J1-moment-bound}
		\mathbb{E}^\phi\left[
		\lVert J_1\rVert_{H_B^{-1+\alpha,\varepsilon}(\rho^{16})}^{
			\frac{4}{3}(1-\nu)}
		\right]
		\lesssim
		\lVert B\rVert_{H^{1,\varepsilon}(\rho)}^{1/\alpha}
		+
		\mathbb{E}^\phi\lVert\Xi\rVert_\kappa^{K_0}.
	\end{equation}
	
	\smallskip
	\noindent\textup{(4)} \textit{Estimate of $J_2$ and $R^{\mathrm{orange}}$.}
	The estimates of $J_2$ and $R^{\mathrm{orange}}$ are similar to those
	in~\textup{(3)}. In particular,
	\begin{equation}
		\label{eq:higgs-low-schauder-J2-orange-moment-bound}
		\mathbb{E}^\phi\left[
		\left\lVert J_2+R^{\mathrm{orange}}
		\right\rVert_{H_B^{-1+\alpha,\varepsilon}(\rho^{16})}^{
			\frac{4}{3}(1-\nu)}
		\right]
		\lesssim
		\lVert B\rVert_{H^{1,\varepsilon}(\rho)}^{1/\alpha}
		+
		\mathbb{E}^\phi\lVert\Xi\rVert_\kappa^{K_0}.
	\end{equation}
	
	\smallskip
	\noindent\textup{(5)} \textit{Estimate of
		$J_3+\lvert Z_1^{\mathrm{lo}}\rvert^2Z_1^{\mathrm{lo}}$.}
	By the definition of $J_3$ in~\eqref{eq:higgs-low-J2-J3},
	\begin{equation}
		J_3+\lvert Z_1^{\mathrm{lo}}\rvert^2Z_1^{\mathrm{lo}}
		=\lvert Z_1^\varepsilon\rvert^2Z_1^\varepsilon.
	\end{equation}
	The uniform embeddings
	$L^{4/3,\varepsilon}(\rho^{16})
	\hookrightarrow H_B^{-1/2,\varepsilon}(\rho^{16})
	\hookrightarrow H_B^{-1+\alpha,\varepsilon}(\rho^{16})$ from
	\eqref{eq:covariant-sobolev-embedding}, together with $\rho\le1$, give
	\begin{equation}
		\label{eq:higgs-low-schauder-J3-cubic-moment-bound}
		\begin{aligned}
			\mathbb{E}^\phi\left[
			\left\lVert
			J_3+\lvert Z_1^{\mathrm{lo}}\rvert^2Z_1^{\mathrm{lo}}
			\right\rVert_{H_B^{-1+\alpha,\varepsilon}(\rho^{16})}^{
				\frac{4}{3}(1-\nu)}
			\right]
			&\lesssim
			\mathbb{E}^\phi
			\left\lVert Z_1^\varepsilon\right\rVert_{L^{4,\varepsilon}(\rho)}^{
				4(1-\nu)}
			\lesssim
			1+
			\mathbb{E}^\phi
			\left\lVert Z_1^{\mathrm{lo}}\right\rVert_{L^{4,\varepsilon}(\rho)}^4\\
			&\lesssim
			\lVert B\rVert_{H^{1,\varepsilon}(\rho)}^{1/\alpha}
			+
			\mathbb{E}^\phi\lVert\Xi\rVert_\kappa^{K_0}.
		\end{aligned}
	\end{equation}
	
	\begin{proposition}
		\label{prop:higgs-low-schauder-bound}
		We have
		\begin{equation}
			\label{eq:higgs-low-H1-alpha-moment-bound}
			\mathbb{E}^\phi
			\lVert Z_1^{\mathrm{lo}}\rVert_{H_B^{1+\alpha,\varepsilon}(\rho^{16})}^{\frac43(1-\nu)}
			\lesssim
			\bigl(1+\lVert B\rVert_{H^{1,\varepsilon}(\rho)}\bigr)^{1/\alpha}
			+
			\mathbb{E}^\phi\lVert\Xi\rVert_\kappa^{K_0}.
		\end{equation}
		Moreover,
		\begin{equation}
			\label{eq:higgs-low-H1-kappa1-moment-bound}
			\mathbb{E}^\phi
			\lVert Z_1^{\mathrm{lo}}\rVert_{H_B^{1+\kappa_1,\varepsilon}(\rho^3)}^{\frac43(1-\nu)}
			\lesssim
			\bigl(1+\lVert B\rVert_{H^{1,\varepsilon}(\rho)}\bigr)^{\kappa_2}
			+
			\mathbb{E}^\phi\lVert\Xi\rVert_\kappa^{K_0}.
		\end{equation}
	\end{proposition}
	
	\begin{proof}
		Applying Lemma~\ref{lem:covariant-fbsde-schauder-estimate} to
		\eqref{eq:higgs-low-feedback-abbreviated}, and using
		\eqref{eq:higgs-low-schauder-I1-bound},
		\eqref{eq:higgs-low-schauder-I2-moment-bound}, and
		\eqref{eq:higgs-low-schauder-J1-moment-bound}--
		\eqref{eq:higgs-low-schauder-J3-cubic-moment-bound}, gives
		\eqref{eq:higgs-low-H1-alpha-moment-bound}.
		
		Since $2+14\kappa_1/\alpha\le3$,
		Lemma~\ref{lem:covariant-interpolation-inequality} and the
		polynomial-weight multiplier bound give
		\begin{equation}
			\label{eq:higgs-low-H1-kappa1-interpolation}
			\left\lVert Z_1^{\mathrm{lo}}
			\right\rVert_{H_B^{1+\kappa_1,\varepsilon}(\rho^3)}
			\lesssim
			\left\lVert Z_1^{\mathrm{lo}}
			\right\rVert_{H_B^{1,\varepsilon}(\rho^2)}^{
				1-\frac{\kappa_1}{\alpha}}
			\left\lVert Z_1^{\mathrm{lo}}
			\right\rVert_{H_B^{1+\alpha,\varepsilon}(\rho^{16})}^{
				\frac{\kappa_1}{\alpha}}.
		\end{equation}
		Set $p:=\frac43(1-\nu)$,
		$p_1:=\frac{2}{p(1-\kappa_1/\alpha)}$, and
		$p_2:=(1-1/p_1)^{-1}$.
		Young's and Jensen's inequalities,
		Proposition~\ref{prop:higgs-low-frequency-uniform-bound}, and
		\eqref{eq:higgs-low-H1-alpha-moment-bound} yield
		\begin{equation}
			\label{eq:higgs-low-H1-kappa1-moment-interpolation}
			\begin{aligned}
				\mathbb E^\phi
				\lVert Z_1^{\mathrm{lo}}\rVert_{H_B^{1+\kappa_1,\varepsilon}(\rho^3)}^p
				&\lesssim
				\mathbb E^\phi
				\lVert Z_1^{\mathrm{lo}}\rVert_{H_B^{1,\varepsilon}(\rho^2)}^2
				+
				\left(
				\mathbb E^\phi
				\lVert Z_1^{\mathrm{lo}}\rVert_{H_B^{1+\alpha,\varepsilon}(\rho^{16})}^p
				\right)^{\frac{\kappa_1}{\alpha}p_2}\\
				&\lesssim
				\bigl(1+\lVert B\rVert_{H^{1,\varepsilon}(\rho)}\bigr)^{\kappa_1}
				+
				\bigl(1+\lVert B\rVert_{H^{1,\varepsilon}(\rho)}\bigr)^{
					\frac{\kappa_1}{\alpha^2}p_2}
				+\mathbb E^\phi\lVert\Xi\rVert_\kappa^{K_0}\\
				&\lesssim
				\bigl(1+\lVert B\rVert_{H^{1,\varepsilon}(\rho)}\bigr)^{\kappa_2}
				+\mathbb E^\phi\lVert\Xi\rVert_\kappa^{K_0}.
			\end{aligned}
		\end{equation}
		Here we used $p<2$, $(\kappa_1/\alpha)p_2<1$, and
		$(\kappa_1/\alpha^2)p_2\ll\kappa_2$.
	\end{proof}

	\section{Estimates of the gauge field}
	\label{sec:estimate-gauge-field}
	
	We estimate
	$\lVert B_1^\varepsilon\rVert_{H^{1,\varepsilon}(\rho)}$
	using the preceding results with
	$(A,B)=\bigl(A^\varepsilon,B_1^\varepsilon\bigr)$.
	We use the notation in~\eqref{eq:random-harmonic-field-shorthands} and set
$
		\lVert\Theta^\varepsilon\rVert_\kappa
		:=
		\left\lVert
		\Xi^{\varepsilon,A^\varepsilon,B_1^\varepsilon}
		\right\rVert_\kappa.
$
	Equation~\eqref{eq:gauge-integral-fbsde} gives, for $0\leq t\leq1$,
	\begin{equation}
		\label{eq:gauge-field-control-equations}
		B_t^\varepsilon
		=-P_{\chi_\varepsilon}^{M}K^{M,\varepsilon}
		\bigl(P_{\chi_\varepsilon}^{M}K^{M,\varepsilon}\bigr)^*
		\int_0^t\mathbb E_s^A
		\bigl[\nabla_A\mathcal W_{\mathrm{av}}^{M,\varepsilon}(C^{\varepsilon,0})\bigr]ds.
	\end{equation}
	Here
$
		K^{M,\varepsilon}=\nabla_\perp^\varepsilon(-\Delta^\varepsilon)^{-1},
		\nabla_A\mathcal W_{\mathrm{av}}^{M,\varepsilon}(C^{\varepsilon,0})
		=\mathbb E\Bigl[\mathbb E^\phi\bigl[
		\Pi_{\chi_\varepsilon}^{M}\operatorname{Im}
		(\overline{\Phi^\varepsilon}D_{C^\varepsilon}^\varepsilon\Phi^\varepsilon)
		\bigr]\,\big|\,\mathcal F_1^A\Bigr].
$

		On the periodic lattice,
$
			\label{eq:gauge-covariance-operator}
			K^{M,\varepsilon}(K^{M,\varepsilon})^*
			=(-\Delta^\varepsilon)^{-1}.
$
		The multiplier has operator norm $|k^\varepsilon|^{-2}$ at
		$k\in\widehat\Lambda_{M,\varepsilon}\setminus\{0\}$ and vanishes at zero.
		Since $\min_{k\ne0}|k^\varepsilon|\asymp M^{-1}$,
$
			\label{eq:gauge-infrared-operator-growth}
			\bigl\|(-\Delta^\varepsilon)^{-1}\bigr\|_{
			H^{1,\varepsilon}(\Lambda_{M,\varepsilon})\to
			H^{1,\varepsilon}(\Lambda_{M,\varepsilon})}
			\asymp M^2.
$	This growth in $M$ obstructs a direct uniform estimate of the full
	$B^\varepsilon$ from the FBSDE structure. We therefore use a low–high frequency decomposition of $B^\epsilon$.
	The low-frequency part is  controlled by a Gaussian domination argument adapted from \cite{CC24}, while the high-frequency part is estimated through the FBSDE.
	Fix a real even $\eta\in C_c^\infty(\mathbb R^2)$ such that
	\begin{equation}
		\label{eq:gauge-infrared-cutoff}
		0\leq\eta\leq1,
		\qquad\eta=1\text{ on }[-3/2,3/2]^2,
		\qquad\operatorname{supp}\eta\subset[-2,2]^2,
	\end{equation}
	and set $\eta^c:=1-\eta$.
	We define
	\begin{equation}
		\label{eq:gauge-low-high-decomposition}
		B_t^{\mathrm{lo}}:=\eta(D)B_t^\varepsilon,
		\qquad
		B_t^{\mathrm{hi}}:=B_t^\varepsilon-B_t^{\mathrm{lo}}.
	\end{equation}

	\begin{proposition}
		\label{prop:gauge-low-shift-uniform-moments}
		For every $2\leq p<\infty$,
		\begin{equation}
			\label{eq:gauge-low-shift-uniform-moments}
			\sup_{M,\varepsilon}\mathbb E
			\bigl[\|B_1^{\mathrm{lo}}\|_{H^{1,\varepsilon}(\rho)}^p\bigr]
			\lesssim_p1.
		\end{equation}
	\end{proposition}
	\begin{proof}
		By the definition of $C^\varepsilon$,
		$B_1^{\mathrm{lo}}=\eta(D)C^\varepsilon-\eta(D)A^\varepsilon$.
		We first prove
		\begin{equation}
			\label{eq:gauge-low-field-uniform-moments}
			\sup_{M,\varepsilon}\mathbb E
			\bigl[\|\eta(D)C^\varepsilon\|_{H^{1,\varepsilon}(\rho)}^p\bigr]
			\lesssim_p1.
		\end{equation}
		By~\eqref{eq:optimal-terminal-laws},
		$\operatorname{Law}(C^{\varepsilon,0})=\mu_{M,\varepsilon}^{\mathrm{gauge}}$.
		Apply Corollary~\ref{cor:gauge-gaussian-domination} componentwise to
		$F(A)=(\Delta_j^\varepsilon\eta(D)A)_i(x)$.
		Equations~\eqref{eq:pinned-wick-dyadic-second-moment},
		\eqref{eq:pinned-gauge-low-frequency-variance} and
		\eqref{eq:dyadic-logarithmic-kernel-bound}, with the fixed smooth
		multiplier $\eta(D)$, give
		\begin{equation}
			\label{eq:gauge-low-pinned-block-moments}
			\begin{aligned}
				\mathbb E|\Delta_j^\varepsilon\eta(D)C^{\varepsilon,0}(x)|^p
				\lesssim_p
				\bigl(\mathbb E^A|\Delta_j^\varepsilon\eta(D)
				A^{\varepsilon,0}(x)|^2\bigr)^{p/2}
				\lesssim_p\bigl(\log(2+|x|)\bigr)^{p/2}.
			\end{aligned}
		\end{equation}
		Only finitely many dyadic blocks occur, independently of $M,\varepsilon$.
		The weighted Besov--Sobolev equivalence and Minkowski's inequality
		therefore yield
		\begin{equation}
			\label{eq:gauge-low-pinned-weighted-moments}
			\bigl(\mathbb E\|\eta(D)C^{\varepsilon,0}\|_
			{H^{1,\varepsilon}(\rho)}^p\bigr)^{2/p}
			\lesssim_p\int_{\Lambda_\varepsilon}
			\rho(x)^2\log(2+|x|)\,dx
			\lesssim_p1.
		\end{equation}
		Finally, $\eta(0)=1$ and $|H^\varepsilon|\leq\sqrt2M^{-1}$ imply
$
			\label{eq:gauge-low-harmonic-weighted-bound}
		\|\eta(D)H^\varepsilon\|_{H^{1,\varepsilon}(\rho)}
			=|H^\varepsilon|\|\rho\|_{H^{1,\varepsilon}}
			\lesssim M^{-1}\lesssim1.
$
		The triangle inequality yields~\eqref{eq:gauge-low-field-uniform-moments}.
		
	By Corollary~\ref{cor:standard-gauge-field-stochastic-estimates},
$
			\label{eq:gauge-free-weighted-moments-assumed}
			\sup_{M,\varepsilon}\mathbb E
			\bigl[\|A^\varepsilon\|_{\mathcal C^{-\kappa,\varepsilon}(\rho^\kappa)}^p\bigr]
			\lesssim_{\kappa,p}1.
$
		Since $\rho^{1-\kappa}\in L^{2,\varepsilon}$, weighted Besov embedding
		and the fixed cutoff $\eta(D)$ give
		\begin{equation}
			\label{eq:gauge-free-low-frequency-bound}
			\|\eta(D)A^\varepsilon\|_{H^{1,\varepsilon}(\rho)}
			\lesssim\|A^\varepsilon\|_{H^{-2\kappa,\varepsilon}(\rho)}
			\lesssim\|A^\varepsilon\|_{\mathcal C^{-\kappa,\varepsilon}(\rho^\kappa)}.
		\end{equation}
		Together with~\eqref{eq:gauge-low-field-uniform-moments},
		the triangle inequality concludes the proof.
	\end{proof}

	In the remainder of this section, we use the FBSDE to estimate
	$B_1^{\mathrm{hi}}$.
	\begin{lemma}
		\label{lem:harmonic-averaged-current-constraints}
		The unprojected conditional current satisfies, almost surely,
		\begin{equation}
			\label{eq:harmonic-averaged-current-constraints}
			\mathbb E\Bigl[\mathbb E^\phi\bigl[
			\operatorname{Im}(\overline{\Phi^\varepsilon}
			D_{C^\varepsilon}^\varepsilon\Phi^\varepsilon)
			\bigr]\,\big|\,\mathcal F_1^A\Bigr]
			\in(\mathbb R^2)_{0,\mathrm{div}}^{\Lambda_{M,\varepsilon}}.
		\end{equation}
	\end{lemma}
	\begin{proof}
		In this proof, extend $\mathcal W^{M,\varepsilon}$ and
		$\mathcal W_{\mathrm{av}}^{M,\varepsilon}$ to all periodic real gauge fields
		by the defining formulas~\eqref{eq:effective-gauge-potential} and
		\eqref{eq:harmonic-averaged-effective-potential}.
		The unprojected first variation in
		\eqref{eq:gauge-gradient-effective-interaction-current}, together with
		\eqref{eq:harmonic-conditional-law}, gives, for every
		$h\in(\mathbb R^2)^{\Lambda_{M,\varepsilon}}$,
		\begin{equation}
			\label{eq:harmonic-averaged-current-first-variation}
			\begin{aligned}
				\left.\frac{d}{d\delta}\right|_{\delta=0}
				\mathcal W_{\mathrm{av}}^{M,\varepsilon}
				(C^{\varepsilon,0}+\delta h)
				=\left\langle
				\mathbb E\Bigl[\mathbb E^\phi\bigl[
				\operatorname{Im}(\overline{\Phi^\varepsilon}
				D_{C^\varepsilon}^\varepsilon\Phi^\varepsilon)
				\bigr]\,\big|\,\mathcal F_1^A\Bigr],h
				\right\rangle_{M,\varepsilon}.
			\end{aligned}
		\end{equation}
		Gauge invariance of $\mathcal W^{M,\varepsilon}$ implies the same
		invariance for $\mathcal W_{\mathrm{av}}^{M,\varepsilon}$.
		Taking $h=\nabla^\varepsilon g$ for arbitrary periodic real $g$ in
		\eqref{eq:harmonic-averaged-current-first-variation} and using discrete
		integration by parts shows that the current is divergence-free.

		The large gauge transformations $g(x)=n\mathbin{\cdot}x$,
		$n\in M^{-1}\mathbb Z^2$, also give
		$\mathcal W^{M,\varepsilon}(C+n)=\mathcal W^{M,\varepsilon}(C)$.
		Periodicity of the integrand in
		\eqref{eq:harmonic-averaged-effective-potential} therefore yields
		$\mathcal W_{\mathrm{av}}^{M,\varepsilon}(C+h)
		=\mathcal W_{\mathrm{av}}^{M,\varepsilon}(C)$ for every constant
		$h\in\mathbb R^2$.
		Taking the constant fields $h=e_j$, $j=1,2$, in
		\eqref{eq:harmonic-averaged-current-first-variation} shows that the
		current has zero spatial mean.
	\end{proof}

	\begin{lemma}[High-frequency FBSDE Schauder estimate]
		\label{lem:gauge-fbsde-schauder-estimate}
		For $1\leq p\leq\infty$ and $a\in\mathbb R$, whenever the right-hand side is finite,
		\begin{equation}
			\label{eq:gauge-fbsde-schauder-estimate}
			\|B_1^{\mathrm{hi}}\|_{L^p(\Omega;H^{1,\varepsilon}(\rho^a))}
			\lesssim_a
			\left\|\mathbb E^\phi\left[
			\eta^c(D)P_\perp^{M,\varepsilon}\operatorname{Im}
			(\overline{\Phi^\varepsilon}D_{C^\varepsilon}^\varepsilon\Phi^\varepsilon)
			\right]\right\|_{L^p(\Omega;H^{-1,\varepsilon}(\rho^a))}.
		\end{equation}
	\end{lemma}
	\begin{proof}
		On $(\mathbb R^2)_{0,\mathrm{div}}^{\Lambda_{M,\varepsilon}}$, we have
		\begin{equation}
			\label{eq:gauge-high-control-operator-identities}
			K^{M,\varepsilon}(K^{M,\varepsilon})^*=(-\Delta^\varepsilon)^{-1},
			\qquad \eta^c(D)P_{\chi_\varepsilon}^{M}=\eta^c(D),
			\qquad (P_{\chi_\varepsilon}^{M})^*\Pi_{\chi_\varepsilon}^{M}=\operatorname{Id}.
		\end{equation}
		For the last identity, take $A,B$ in this space. Then
		\begin{equation}
			\label{eq:gauge-pinned-current-adjoint-duality}
			\begin{aligned}
				\left\langle(P_{\chi_\varepsilon}^{M})^*\Pi_{\chi_\varepsilon}^{M}A,B\right\rangle_{M,\varepsilon}
				=\left\langle\Pi_{\chi_\varepsilon}^{M}A,P_{\chi_\varepsilon}^{M}B\right\rangle_{M,\varepsilon}
				=\left\langle A,P_{\chi_\varepsilon}^{M}B\right\rangle_{M,\varepsilon}
				=\langle A,B\rangle_{M,\varepsilon},
			\end{aligned}
		\end{equation}
		where the last equality uses the zero mean of $A$.
		Equations~\eqref{eq:harmonic-averaged-gauge-current}
		and~\eqref{eq:gauge-field-control-equations}, together with
		Lemma~\ref{lem:harmonic-averaged-current-constraints}, give,
		for $0\leq t\leq1$,
		\begin{equation}
			\label{eq:gauge-control-composition-fourier}
			B_t^{\mathrm{hi}}
			=-\eta^c(D)(-\Delta^\varepsilon)^{-1}
			\int_0^t\mathbb E_s^A\Bigl[\mathbb E^\phi\bigl[
			P_\perp^{M,\varepsilon}\operatorname{Im}
			(\overline{\Phi^\varepsilon}D_{C^\varepsilon}^\varepsilon\Phi^\varepsilon)
			\bigr]\Bigr]\,ds.
		\end{equation}
		Here we insert $P_\perp^{M,\varepsilon}$ using~\eqref{eq:harmonic-averaged-current-constraints}
		to apply~\eqref{eq:projected-lattice-current-identity}.
		For mean-zero $g\colon\Lambda_{M,\varepsilon}\to\mathbb R^2$, the remaining estimate is
		\begin{equation}
			\label{eq:gauge-control-composition-smoothing}
			\|\eta^c(D)(-\Delta^\varepsilon)^{-1}g\|_{H^{1,\varepsilon}(\rho^a)}
			\lesssim_a\|\eta^c(D)g\|_{H^{-1,\varepsilon}(\rho^a)}.
		\end{equation}
	Let $K_j^{M,\varepsilon}$ be the
		periodic kernel with multiplier
		$\varphi_j^\varepsilon(k)(1-\eta(2k))/|k^\varepsilon|^2$.
		The factor $1-\eta(2k)$ equals one on $\operatorname{supp}\eta^c$.
		The multiplier and its derivatives of order $N$ are bounded by
		$C2^{-2j}$ and $C_N2^{-(2+N)j}$, respectively.
		By Lemma~\ref{lem:periodic-multiplier-kernel-bound}, with $N>4|a|+2$,
		\begin{equation}
			\label{eq:gauge-high-resolvent-weighted-kernel-bound}
			\int_{\Lambda_{M,\varepsilon}}
			|K_j^{M,\varepsilon}(z)|(1+d_M(z,0))^{4|a|}\,dz
			\lesssim_a 2^{-2j}.
		\end{equation}
		Weighted Young's inequality and  dyadic summation yield
		\eqref{eq:gauge-control-composition-smoothing}.
		Minkowski's inequality and conditional Jensen's inequality
		yield~\eqref{eq:gauge-fbsde-schauder-estimate}.
	\end{proof}

	Let $J$ be the cost functional associated with the Gibbs
	measure~\eqref{eq:conditional-ymh-gauge-law} in the sense
	of~\eqref{eq:app-cost-functional}:
	\begin{equation}
		\label{eq:gauge-control-cost}
		J(v):=
		\mathbb E^A\left[
		\mathcal W_{\mathrm{av}}^{M,\varepsilon}\left(
		P_{\chi_\varepsilon}^{M}K^{M,\varepsilon}
		\left(W_1^{M,\varepsilon}+\int_0^1v_s\,ds\right)\right)
		+\frac12\int_0^1\|v_s\|_{L^{2,\varepsilon}}^2\,ds\right].
	\end{equation}
	Let $v^\varepsilon$ be its optimal control over the admissible class
	$\mathcal U_{\mathrm{ad}}^A$ ,
	as given by Proposition~\ref{lem:app-abstract-maximum-principle} and
	Lemma~\ref{lem:app-cor-satisfies-hypothesis}. Define
	\begin{equation}
		\label{eq:gauge-fbsde-schauder-system}
		\widetilde B_t^\varepsilon:=
		K^{M,\varepsilon}\int_0^t v_s^\varepsilon\,ds,
		\qquad B_t^\varepsilon=P_{\chi_\varepsilon}^{M}\widetilde B_t^\varepsilon.
	\end{equation}
	Then $B_1^{\mathrm{hi}}=\eta^c(D)\widetilde B_1^\varepsilon$.

	\begin{lemma}
		\label{lem:gauge-stationary-weight-transfer}
		The field $B_1^{\mathrm{hi}}$ is spatially stationary:
		\begin{equation}
			\label{eq:gauge-high-shift-stationarity}
			B_1^{\mathrm{hi}}(\mathord\cdot+y)
			\stackrel{\mathrm{law}}{=}B_1^{\mathrm{hi}},
			\qquad y\in\Lambda_{M,\varepsilon}.
		\end{equation}
		Consequently, for $2\leq p<\infty$ and $a>1/4$,
		\begin{equation}
			\label{eq:gauge-stationary-weight-transfer}
			\|B_1^{\mathrm{hi}}\|_{L^p(\Omega;H^{1,\varepsilon}(\rho^a))}
			\lesssim_a
			\|B_1^{\mathrm{hi}}\|_{L^p(\Omega;H^{1,\varepsilon}(\rho^4))}.
		\end{equation}
	\end{lemma}
	\begin{proof}
		We derive the weighted estimate from stationarity.
		Write $f:=B_1^{\mathrm{hi}}$.
		Set $\ell_\varepsilon:=\varepsilon\lceil\varepsilon^{-1}\rceil
		\in[1,2)$, $Q_0:=[0,\ell_\varepsilon)^2\cap\Lambda_\varepsilon$,
		and $Q_k:=Q_0+\ell_\varepsilon k$ for $k\in\mathbb Z^2$.
		The sets $Q_k$ partition the full lattice.
		By the definition of $\rho$,  Minkowski's inequality in $L^{p/2}(\Omega)$ and stationarity,
		\begin{equation}
			\label{eq:gauge-stationary-weight-transfer-proof}
			\begin{aligned}
			\|f\|_{L^p(\Omega;H^{1,\varepsilon}(\rho^a))}^2
			&\lesssim_a
			\left\|\sum_{k\in\mathbb Z^2}\rho(k)^{2a}
			\|f\|_{H^{1,\varepsilon}(Q_k)}^2\right\|_{L^{p/2}(\Omega)}\\
			&\leq\sum_{k\in\mathbb Z^2}\rho(k)^{2a}
			\|f\|_{L^p(\Omega;H^{1,\varepsilon}(Q_k))}^2
			\lesssim_a\|f\|_{L^p(\Omega;H^{1,\varepsilon}(Q_0))}^2.
			\end{aligned}
		\end{equation}
		We used $\sum_k\rho(k)^{2a}<\infty$ for $a>1/4$.
		On $Q_0$, $\rho^4$ is uniformly
		bounded below. The discrete product rule therefore bounds the last
		norm by $C\|f\|_{L^p(\Omega;H^{1,\varepsilon}(\rho^4))}$.
		It remains to prove~\eqref{eq:gauge-high-shift-stationarity}.
		Fix $j\in\{1,2\}$ and write $W = W^{M,\varepsilon}$.
		The proof of Lemma~\ref{lem:harmonic-averaged-current-constraints}
		establishes invariance under constant shifts. Together with spatial
		translation invariance of $\mathcal{W}^{M,\varepsilon}$, this gives
		\begin{equation}
			\label{eq:gauge-averaged-potential-invariances}
			\mathcal W_{\mathrm{av}}^{M,\varepsilon}(P_{\chi_\varepsilon}^{M}C)
			=\mathcal W_{\mathrm{av}}^{M,\varepsilon}(\tau_jC)
			=\mathcal W_{\mathrm{av}}^{M,\varepsilon}(C).
		\end{equation}
		Replacing $W$ by $\tau_jW$ in~\eqref{eq:gauge-control-cost}, define
		\begin{equation}
			\label{eq:gauge-translated-control-cost}
			\widetilde J(v):=
			\mathbb E^A\left[
			\mathcal W_{\mathrm{av}}^{M,\varepsilon}\left(
			P_{\chi_\varepsilon}^{M}K^{M,\varepsilon}
			\left(\tau_j W_1+\int_0^1v_s\,ds\right)\right)
			+\frac12\int_0^1\|v_s\|_{L^{2,\varepsilon}}^2\,ds\right].
		\end{equation}
		The maps $v\mapsto\tau_jv$ and
		$v\mapsto v\circ\tau_j$ are bijections of $\mathcal U_{\mathrm{ad}}^A$,
		where composition acts on the Wiener path.
		Since $K^{M,\varepsilon}$ commutes with $\tau_j$, while $\tau_j$
		preserves the control norm and Wiener measure,
		\eqref{eq:gauge-averaged-potential-invariances} yields
		\begin{equation}
			\label{eq:gauge-translated-cost-identities}
			\widetilde J(\tau_jv)=J(v)
			=\widetilde J(v\circ\tau_j).
		\end{equation}
		By uniqueness of the minimizer of $J$,
		$\widetilde J$ therefore has a unique minimizer, represented both by
		$\tau_jv^\varepsilon$ and $v^\varepsilon\circ\tau_j$. Hence
		\begin{equation}
			\label{eq:gauge-optimal-control-translation-equivariance}
			\tau_jv_t^\varepsilon(\omega_A)
			=v_t^\varepsilon(\tau_j\omega_A),
			\qquad dt\otimes d\mathbb P^A\text{-almost everywhere}.
		\end{equation}
		Integrating in time and using that $\eta^c(D)K^{M,\varepsilon}$
		commutes with translations, translation invariance of Wiener measure gives
		\eqref{eq:gauge-high-shift-stationarity} for $y=\varepsilon e_j$;
		iteration gives every $y\in\Lambda_{M,\varepsilon}$.
	\end{proof}

	Combining Lemmas~\ref{lem:gauge-fbsde-schauder-estimate}
	and~\ref{lem:gauge-stationary-weight-transfer} gives,
	for $2\leq p<\infty$,
	\begin{equation}
		\label{eq:gauge-high-stationary-schauder-bound}
		\|B_1^{\mathrm{hi}}\|_{L^p(\Omega;H^{1,\varepsilon}(\rho))}
		\lesssim
		\left\|\mathbb E^\phi\left[
		\eta^c(D)P_\perp^{M,\varepsilon}\operatorname{Im}
		(\overline{\Phi^\varepsilon}D_{C^\varepsilon}^\varepsilon\Phi^\varepsilon)
		\right]\right\|_{L^p(\Omega;H^{-1,\varepsilon}(\rho^4))}.
	\end{equation}
	We therefore decompose the high-frequency projected current on the right-hand side.
	
	\subsection{Decomposition of the projected current}
	
	For brevity, set $A=A^\varepsilon$, $B=B^\varepsilon$, and $C=A+B$.
	Recall that $P_\perp^{M,\varepsilon}$ is the orthogonal projection onto
	$(\mathbb R^2)_{0,\mathrm{div}}^{\Lambda_{M,\varepsilon}}$.
	For $\varphi,\psi\colon\Lambda_{M,\varepsilon}\to\mathbb{C}$, define the
	real bilinear map
	\begin{equation}
		\label{eq:gauge-current-bilinear-form}
		T_C^\varepsilon(\varphi,\psi)
		:=
		\eta^c(D)P_\perp^{M,\varepsilon}
		\operatorname{Im}\left[
		\overline{\varphi}\,D_C^\varepsilon\psi
		\right].
	\end{equation}
	Lemma~\ref{lem:projected-lattice-current-identity} then gives
	\begin{equation}
		\label{eq:gauge-current-shifted-form}
		T_C^\varepsilon(\varphi,\psi)
		=
		\eta^c(D)P_\perp^{M,\varepsilon}
		\left(
		\operatorname{Im}\left[
		\overline{\psi(\mathord\cdot+\varepsilon e_j)}
		U_{-C,j}^\varepsilon D_{C,j}^\varepsilon\varphi
		\right]
		\right)_{j=1}^2.
	\end{equation}
		Using $\Phi^\varepsilon=\phi_B^\varepsilon+Z_1^\varepsilon$, we obtain
	\begin{equation}
		\label{eq:gauge-gradient-current-decomposition}
		\begin{aligned}
			\mathbb E^\phi\left[
			\eta^c(D)P_\perp^{M,\varepsilon}\operatorname{Im}
			(\overline{\Phi^\varepsilon}D_C^\varepsilon\Phi^\varepsilon)
			\right]
			={}&
			\mathbb{E}^\phi\left[
			T_C^\varepsilon(\Phi^\varepsilon,\Phi^\varepsilon)
			\right]\\
			={}&
			\mathbb{E}^\phi\Big[
			T_C^\varepsilon(\phi_B^\varepsilon,\phi_B^\varepsilon)
			+T_C^\varepsilon(\phi_B^\varepsilon,Z_1^\varepsilon)
			+T_C^\varepsilon(Z_1^\varepsilon,\phi_B^\varepsilon)
			+T_C^\varepsilon(Z_1^\varepsilon,Z_1^\varepsilon)
			\Big].
		\end{aligned}
	\end{equation}
	Using~\eqref{eq:covariant-derivative-rho-decomposition}, we decompose the four terms in
	\eqref{eq:gauge-gradient-current-decomposition} as follows, with the irregular terms again highlighted in \textcolor{singularpink}{pink}.
	
	\textup{(1)} \textit{For the first term,}
	\begin{equation}
		\label{eq:gauge-current-phi-phi-splitting}
		\begin{aligned}
			T_C^\varepsilon
			(\phi_B^\varepsilon,\phi_B^\varepsilon)
			&=
			\textcolor{singularpink}{
				T_B^\varepsilon
				(\phi_B^\varepsilon,\phi_B^\varepsilon)
			}
			+
			\textcolor{singularpink}{
				\eta^c(D)P_\perp^{M,\varepsilon}
				\operatorname{Im}\left[
				i\overline{\phi_B^\varepsilon}\,
				U_B^\varepsilon\rho_A^\varepsilon
				\phi_B^\varepsilon
				\right]
			}.
		\end{aligned}
	\end{equation}
	
	\textup{(2)} \textit{For the first mixed term,}
	\begin{equation}
		\label{eq:gauge-current-phi-Z-splitting}
		\begin{aligned}
			T_C^\varepsilon(\phi_B^\varepsilon,Z_1^\varepsilon)
			&=
			\textcolor{singularpink}{
				T_B^\varepsilon
				(\phi_B^\varepsilon,Z_1^{\mathrm{hi}})
			}
			+
			T_B^\varepsilon
			(\phi_B^\varepsilon,Z_1^{\mathrm{lo}})
			+
			\eta^c(D)P_\perp^{M,\varepsilon}
			\operatorname{Im}\left[
			i\overline{\phi_B^\varepsilon}\,
			U_B^\varepsilon\rho_A^\varepsilon Z_1^\varepsilon
			\right].
		\end{aligned}
	\end{equation}
	Similarly,
	\begin{equation}
		\label{eq:gauge-current-Z-phi-splitting}
		T_C^\varepsilon(Z_1^\varepsilon,\phi_B^\varepsilon)
		=
		\textcolor{singularpink}{
			T_B^\varepsilon
			(Z_1^{\mathrm{hi}},\phi_B^\varepsilon)
		}
		+
		T_B^\varepsilon
		(Z_1^{\mathrm{lo}},\phi_B^\varepsilon)
		+
		\eta^c(D)P_\perp^{M,\varepsilon}
		\operatorname{Im}\left[
		i\overline{Z_1^\varepsilon}\,
		U_B^\varepsilon\rho_A^\varepsilon\phi_B^\varepsilon
		\right].
	\end{equation}
	
	\textup{(3)} \textit{Finally,}
	\begin{equation}
		\label{eq:gauge-current-Z-Z-splitting}
		T_C^\varepsilon(Z_1^\varepsilon,Z_1^\varepsilon)
		=
		T_B^\varepsilon(Z_1^\varepsilon,Z_1^\varepsilon)
		+
		\eta^c(D)P_\perp^{M,\varepsilon}
		\operatorname{Im}\left[
		i\overline{Z_1^\varepsilon}\,
		U_B^\varepsilon\rho_A^\varepsilon Z_1^\varepsilon
		\right].
	\end{equation}
	
	Substituting
	\eqref{eq:gauge-current-phi-phi-splitting},
	\eqref{eq:gauge-current-phi-Z-splitting},
	\eqref{eq:gauge-current-Z-phi-splitting}, and
	\eqref{eq:gauge-current-Z-Z-splitting} into
	\eqref{eq:gauge-gradient-current-decomposition}, we obtain
	\begin{equation}
		\label{eq:gauge-gradient-current-final-splitting}
		\mathbb E^\phi\left[T_C^\varepsilon(\Phi^\varepsilon,\Phi^\varepsilon)\right]
		=
		\mathbb{E}^\phi\left[
		\sum_{j=1}^3 I_j^{\mathrm{Pink}}
		+
		\sum_{k=1}^4 J_k
		\right].
	\end{equation}
	The terms in \eqref{eq:gauge-gradient-current-final-splitting} are as follows:
	\begin{equation}
		\label{eq:gauge-gradient-current-terms}
		\begin{alignedat}{2}
			&I_1^{\mathrm{Pink}}
			=
			\textcolor{singularpink}{T_B^\varepsilon(\phi_B^\varepsilon,\phi_B^\varepsilon)}
			&\qquad
			&I_2^{\mathrm{Pink}}
			=
			\textcolor{singularpink}{\eta^c(D)P_\perp^{M,\varepsilon}\operatorname{Im}\left[i\overline{\phi_B^\varepsilon}\,U_B^\varepsilon\rho_A^\varepsilon\phi_B^\varepsilon\right]}
			\\
			&I_3^{\mathrm{Pink}}
			=
			\textcolor{singularpink}{T_B^\varepsilon(\phi_B^\varepsilon,Z_1^{\mathrm{hi}})+T_B^\varepsilon(Z_1^{\mathrm{hi}},\phi_B^\varepsilon)}
			\\
			&J_1
			=
			T_B^\varepsilon(\phi_B^\varepsilon,Z_1^{\mathrm{lo}})
			+T_B^\varepsilon(Z_1^{\mathrm{lo}},\phi_B^\varepsilon)
			&\qquad
			&J_2
			=
			\eta^c(D)P_\perp^{M,\varepsilon}\operatorname{Im}\left[
			iU_B^\varepsilon\Bigl(
			\overline{\phi_B^\varepsilon}\,\rho_A^\varepsilon Z_1^\varepsilon
			+\overline{Z_1^\varepsilon}\,\rho_A^\varepsilon\phi_B^\varepsilon
			\Bigr)
			\right]
			\\
			&J_3
			=
			\eta^c(D)P_\perp^{M,\varepsilon}\operatorname{Im}\left[
			i\overline{Z_1^\varepsilon}\,
			U_B^\varepsilon\rho_A^\varepsilon Z_1^\varepsilon
			\right]
			&\qquad
			&J_4
			=
			T_B^\varepsilon(Z_1^\varepsilon,Z_1^\varepsilon).
		\end{alignedat}
	\end{equation}
	
	\subsection{Estimates of regular terms}
	\label{subsec:gauge-regular-terms}
	In this subsection, we estimate the terms $J_k$, $1\leq k\leq4$, in
	\eqref{eq:gauge-gradient-current-terms}.
	
	\noindent
	\textup{(1)} \textit{Estimate of
	$\bigl\lVert\mathbb{E}^\phi[J_1]\bigr\rVert_{H^{-1,\varepsilon}(\rho^4)}$.}
	By the definition of $J_1$,
	\begin{equation}
		\label{eq:gauge-regular-J1-expansion}
		\left\lVert\mathbb{E}^\phi[J_1]\right\rVert_{H^{-1,\varepsilon}(\rho^4)}
		=
		\left\lVert\mathbb{E}^\phi\left[
		T_B^\varepsilon(\phi_B^\varepsilon,Z_1^{\mathrm{lo}})
		+T_B^\varepsilon(Z_1^{\mathrm{lo}},\phi_B^\varepsilon)
		\right]\right\rVert_{H^{-1,\varepsilon}(\rho^4)}.
	\end{equation}
	Since $3+7\kappa\le4$, Lemma~\ref{lem:discrete-paraproduct-estimates},
	\eqref{eq:stochastic-object-Xi-norm}, and
	\eqref{eq:weighted-covariant-derivative-free-output} with $s=2\kappa$ and $a=3$ give
	\begin{equation}
		\label{eq:gauge-regular-J1-first-product-bound}
		\begin{aligned}
			\left\lVert\overline{\phi_B^\varepsilon}\,
			D_B^\varepsilon Z_1^{\mathrm{lo}}\right\rVert_{H^{-1,\varepsilon}(\rho^4)}
			&\lesssim
			\lVert\phi_B^\varepsilon\rVert_{\mathcal C^{-\kappa,\varepsilon}(\rho^\kappa)}
			\lVert D_B^\varepsilon Z_1^{\mathrm{lo}}\rVert_{H^{2\kappa,\varepsilon}(\rho^{3+6\kappa})}\\
			&\lesssim
			\lVert\Xi\rVert_\kappa
			\bigl(1+\lVert B\rVert_{H^{1,\varepsilon}(\rho)}\bigr)^{6\kappa}
			\lVert Z_1^{\mathrm{lo}}\rVert_{H_B^{1+2\kappa,\varepsilon}(\rho^3)}.
		\end{aligned}
	\end{equation}
	For fixed $s,a\in\mathbb R$, the multiplier
	$\eta^c(D)P_\perp^{M,\varepsilon}$ is bounded on
	$H^{s,\varepsilon}(\rho^a)$ uniformly in $M,\varepsilon$.
	Since $2\kappa\le\kappa_1$, conditional Jensen's inequality and 
	Young's inequality yield
	\begin{equation}
		\label{eq:gauge-regular-J1-first-expectation-bound}
		\begin{aligned}
			\left\lVert\mathbb E^\phi\left[
			T_B^\varepsilon(\phi_B^\varepsilon,Z_1^{\mathrm{lo}})
			\right]\right\rVert_{H^{-1,\varepsilon}(\rho^4)}
			&\lesssim
			\mathbb E^\phi\left[
			\lVert\Xi\rVert_\kappa^{K_0}
			+\lVert Z_1^{\mathrm{lo}}\rVert_{H_B^{1+\kappa_1,\varepsilon}(\rho^3)}^{\frac43(1-\nu)}
			\right]
			+\lVert B\rVert_{H^{1,\varepsilon}(\rho)}^{\kappa_1}\\
			&\lesssim
			\mathbb E^\phi\left[\lVert\Xi\rVert_\kappa^{K_0}\right]
			+\lVert B\rVert_{H^{1,\varepsilon}(\rho)}^{\kappa_2},
		\end{aligned}
	\end{equation}
	where the last step uses Proposition~\ref{prop:higgs-low-schauder-bound}.

	For the second term, Lemma~\ref{lem:projected-lattice-current-identity} gives
	\begin{equation}
		\label{eq:gauge-regular-J1-second-shifted-form}
		T_B^\varepsilon(Z_1^{\mathrm{lo}},\phi_B^\varepsilon)
		=
		\eta^c(D)P_\perp^{M,\varepsilon}
		\left(\operatorname{Im}\left[
		\overline{\phi_B^\varepsilon(\mathord\cdot+\varepsilon e_j)}\,
		U_{-B,j}^\varepsilon D_{B,j}^\varepsilon Z_1^{\mathrm{lo}}
		\right]\right)_{j=1}^2.
	\end{equation}
	By~\eqref{eq:covariant-derivative-translation-identities},
	\eqref{eq:weighted-covariant-derivative-free-output}, and the boundedness
	of lattice translations on the weighted spaces,
	\begin{equation}
		\label{eq:gauge-regular-covariant-derivative-bound}
		\begin{aligned}
			\left\lVert U_{-B,j}^\varepsilon D_{B,j}^\varepsilon Z_1^{\mathrm{lo}}
			\right\rVert_{H^{2\kappa,\varepsilon}(\rho^{3+6\kappa})}
			&=
			\left\lVert\bigl(D_{B,j}^{\varepsilon *}Z_1^{\mathrm{lo}}\bigr)
			(\mathord\cdot+\varepsilon e_j)
			\right\rVert_{H^{2\kappa,\varepsilon}(\rho^{3+6\kappa})}\\
			&\lesssim
			\bigl(1+\lVert B\rVert_{H^{1,\varepsilon}(\rho)}\bigr)^{6\kappa}
			\lVert Z_1^{\mathrm{lo}}\rVert_{H_B^{1+2\kappa,\varepsilon}(\rho^3)}.
		\end{aligned}
	\end{equation}
	The shifted field $\phi_B^\varepsilon(\mathord\cdot+\varepsilon e_j)$
	satisfies the same $\mathcal C^{-\kappa,\varepsilon}(\rho^\kappa)$ bound.
	Thus the same argument applies to the second term, and
	\eqref{eq:gauge-regular-J1-first-expectation-bound} yields
	\begin{equation}
		\label{eq:gauge-regular-J1-second-expectation-bound}
		\left\lVert\mathbb E^\phi[J_1]\right\rVert_{H^{-1,\varepsilon}(\rho^4)}
		\lesssim
		\mathbb E^\phi\left[\lVert\Xi\rVert_\kappa^{K_0}\right]
		+\lVert B\rVert_{H^{1,\varepsilon}(\rho)}^{\kappa_2}.
	\end{equation}

	\noindent
	\textup{(2)} \textit{Estimate of
	$\bigl\lVert\mathbb{E}^\phi[J_2]\bigr\rVert_{H^{-1,\varepsilon}(\rho^4)}$.}
	
	By~\eqref{eq:gauge-link-multiplier-bound},
	\begin{equation}
		\label{eq:gauge-regular-J2-weighted-link-multiplier}
		\|U_B^\varepsilon f\|_{\mathcal C^{-\kappa,\varepsilon}(\rho^{3\kappa})}
		\lesssim
		\bigl(1+\|B\|_{H^{1,\varepsilon}(\rho)}\bigr)^{2\kappa}
		\|f\|_{\mathcal C^{-\kappa,\varepsilon}(\rho^\kappa)}.
	\end{equation}
	By the definition of $J_2$,
	\begin{equation}
		\label{eq:gauge-regular-J2-bound}
		\begin{aligned}
			\left\lVert\mathbb{E}^\phi[J_2]
			\right\rVert_{H^{-1,\varepsilon}(\rho^4)}
			&\lesssim
			\mathbb{E}^\phi\left[
			\left\lVert
			\left(
			U_{B,j}^\varepsilon\overline{\phi_B^\varepsilon}\,
			r_{A,j}^\varepsilon
			Z_1^\varepsilon(\mathord\cdot+\varepsilon e_j)
			\right)_{j=1}^2
			\right\rVert_{H^{-1,\varepsilon}(\rho^4)}
			\right]
			+
			\mathbb{E}^\phi\left[
			\left\lVert
			\overline{Z_1^\varepsilon}\,
			U_B^\varepsilon\rho_A^\varepsilon\phi_B^\varepsilon
			\right\rVert_{H^{-1,\varepsilon}(\rho^4)}
			\right]\\
			&\lesssim
			\mathbb{E}^\phi\left[
			\left(
			\left\lVert
			U_B^\varepsilon
			\bigl(\overline{\phi_B^\varepsilon}r_A^\varepsilon\bigr)
			\right\rVert_{\mathcal{C}^{-\kappa,\varepsilon}(\rho^{3\kappa})}
			+
			\left\lVert
			U_B^\varepsilon\rho_A^\varepsilon\phi_B^\varepsilon
			\right\rVert_{\mathcal{C}^{-\kappa,\varepsilon}(\rho^{3\kappa})}
			\right)
			\lVert Z_1^\varepsilon\rVert_{H^{2\kappa,\varepsilon}(\rho^{4-3\kappa})}
			\right]\\
			&\overset{\mathclap{\eqref{eq:gauge-regular-J2-weighted-link-multiplier}}}{\lesssim}
			\mathbb{E}^\phi\left[
			\lVert\Xi\rVert_\kappa
			\bigl(1+\lVert B\rVert_{H^{1,\varepsilon}(\rho)}\bigr)^{2\kappa}
			\lVert Z_1^\varepsilon\rVert_{H^{2\kappa,\varepsilon}(\rho^{4-3\kappa})}
			\right]\\
			&\lesssim
			\mathbb{E}^\phi\left[
			\lVert\Xi\rVert_\kappa
			\left(
			\left\lVert Z_1^{\mathrm{lo}}
			\right\rVert_{H_B^{1,\varepsilon}(\rho^2)}+1
			\right)
			\bigl(1+\lVert B\rVert_{H^{1,\varepsilon}(\rho)}\bigr)^{4\kappa}
			\right]\\
			&\overset{\mathclap{
				\text{Young, Proposition~\ref{prop:higgs-terminal-shift-summary}}}}
			{\lesssim}
			\mathbb{E}^\phi\left[
			\lVert\Xi\rVert_\kappa^{K_0}
			\right]
			+
			\lVert B\rVert_{H^{1,\varepsilon}(\rho)}^{\kappa_2}.
		\end{aligned}
	\end{equation}
	
	\noindent
	\textup{(3)} \textit{Estimate of
	$\bigl\lVert\mathbb E^\phi[J_3]\bigr\rVert_{H^{-1,\varepsilon}(\rho^4)}$.}
	Since $3+3\kappa\le4$, the weighted Besov and paraproduct estimates,
	together with~\eqref{eq:gauge-regular-J2-weighted-link-multiplier}, give
	\begin{equation}
		\label{eq:gauge-regular-J3-first-bound}
		\begin{aligned}
			\left\lVert\mathbb E^\phi[J_3]\right\rVert_{H^{-1,\varepsilon}(\rho^4)}
			&\lesssim\mathbb E^\phi\left[
			\lVert U_B^\varepsilon r_A^\varepsilon\rVert_{\mathcal C^{-\kappa,\varepsilon}(\rho^{3\kappa})}
			\left\lVert\overline{Z_1^\varepsilon}
			\bigl(Z_1^\varepsilon(\mathord\cdot+\varepsilon e_j)\bigr)_{j=1}^2
			\right\rVert_{B_{3/2,\infty}^{2\kappa,\varepsilon}(\rho^3)}
			\right]\\
			&\lesssim\mathbb E^\phi\left[
			\lVert\Xi\rVert_\kappa
			\bigl(1+\lVert B\rVert_{H^{1,\varepsilon}(\rho)}\bigr)^{2\kappa}
			\lVert Z_1^\varepsilon\rVert_{B_{3,\infty}^{2\kappa,\varepsilon}(\rho^{3/2})}^2
			\right].
		\end{aligned}
	\end{equation}
	Interpolation and~\eqref{eq:weighted-HB-H-comparison} yield
	\begin{equation}
		\label{eq:gauge-regular-Z-B3-bound}
		\begin{aligned}
			\lVert Z_1^\varepsilon\rVert_{B_{3,\infty}^{2\kappa,\varepsilon}(\rho^{3/2})}
			&\lesssim
			\lVert Z_1^\varepsilon\rVert_{H^{6\kappa,\varepsilon}(\rho^{5/2})}^{1/3}
			\lVert Z_1^\varepsilon\rVert_{L^{4,\varepsilon}(\rho)}^{2/3}\\
			&\lesssim
			\bigl(1+\lVert B\rVert_{H^{1,\varepsilon}(\rho)}\bigr)^{2\kappa}
			\lVert Z_1^\varepsilon\rVert_{H_B^{1-\alpha,\varepsilon}(\rho^2)}^{1/3}
			\lVert Z_1^\varepsilon\rVert_{L^{4,\varepsilon}(\rho)}^{2/3}.
		\end{aligned}
	\end{equation}
	Consequently, Young's inequality and~\eqref{eq:higgs-summary-full-shift-bound} give
	\begin{equation}
		\label{eq:gauge-regular-J3-final-bound}
		\begin{aligned}
			\left\lVert\mathbb E^\phi[J_3]\right\rVert_{H^{-1,\varepsilon}(\rho^4)}
			&\lesssim\mathbb E^\phi\left[
			\lVert\Xi\rVert_\kappa
			\bigl(1+\lVert B\rVert_{H^{1,\varepsilon}(\rho)}\bigr)^{6\kappa}
			\lVert Z_1^\varepsilon\rVert_{H_B^{1-\alpha,\varepsilon}(\rho^2)}^{2/3}
			\lVert Z_1^\varepsilon\rVert_{L^{4,\varepsilon}(\rho)}^{4/3}
			\right]\\
			&\lesssim
			\mathbb E^\phi\left[\lVert\Xi\rVert_\kappa^{K_0}\right]
			+\lVert B\rVert_{H^{1,\varepsilon}(\rho)}^{\kappa_2}.
		\end{aligned}
	\end{equation}

	\noindent
	\textup{(4)} \textit{Estimate of
	$\bigl\lVert\mathbb E^\phi[J_4]\bigr\rVert_{H^{-1,\varepsilon}(\rho^4)}$.}
	By the definition of $J_4$,
	\begin{equation}
		\label{eq:gauge-regular-J4-initial-bound}
		\left\lVert\mathbb E^\phi[J_4]\right\rVert_{H^{-1,\varepsilon}(\rho^4)}
		\le\mathbb E^\phi\left[
		\left\lVert T_B^\varepsilon(Z_1^\varepsilon,Z_1^\varepsilon)
		\right\rVert_{H^{-1,\varepsilon}(\rho^4)}\right].
	\end{equation}

	\begin{lemma}
		\label{lem:gauge-current-bilinear-estimate}
		For every $\sigma\in(\frac12,1)$ and
		$f,g\colon\Lambda_{M,\varepsilon}\to\mathbb C$,
		\begin{equation}
			\label{eq:gauge-current-bilinear-estimate}
			\lVert T_B^\varepsilon(f,g)\rVert_{H^{2\sigma-2,\varepsilon}(\rho^4)}
			\lesssim_\sigma
			\lVert f\rVert_{H_B^{\sigma,\varepsilon}(\rho^2)}
			\lVert g\rVert_{H_B^{\sigma,\varepsilon}(\rho^2)}.
		\end{equation}
		In particular,
		\begin{equation}
			\label{eq:gauge-current-Z-bilinear-estimate}
			\lVert T_B^\varepsilon(Z_1^\varepsilon,Z_1^\varepsilon)
			\rVert_{H^{-2\alpha,\varepsilon}(\rho^4)}
			\lesssim_\alpha
			\lVert Z_1^\varepsilon\rVert_{H_B^{1-\alpha,\varepsilon}(\rho^2)}^2.
		\end{equation}
	\end{lemma}

	\begin{proof}
		By bilinear interpolation (see~\cite[Theorem~4.4.1]{BL79}),
		it suffices to prove
		\begin{equation}
			\label{eq:gauge-current-first-mixed-regularity-bound}
			\lVert T_B^\varepsilon(f,g)\rVert_{H^{2\sigma-2,\varepsilon}(\rho^4)}
			\lesssim_\sigma
			\lVert f\rVert_{H_B^{2\sigma-1,\varepsilon}(\rho^2)}
			\lVert g\rVert_{H_B^{1,\varepsilon}(\rho^2)}
		\end{equation}
		and
		\begin{equation}
			\label{eq:gauge-current-second-mixed-regularity-bound}
			\lVert T_B^\varepsilon(f,g)\rVert_{H^{2\sigma-2,\varepsilon}(\rho^4)}
			\lesssim_\sigma
			\lVert f\rVert_{H_B^{1,\varepsilon}(\rho^2)}
			\lVert g\rVert_{H_B^{2\sigma-1,\varepsilon}(\rho^2)}.
		\end{equation}
		Weighted Sobolev embedding and H\"older's inequality give
		\begin{equation}
			\label{eq:gauge-current-first-mixed-regularity-proof}
			\begin{aligned}
				\lVert T_B^\varepsilon(f,g)\rVert_{H^{2\sigma-2,\varepsilon}(\rho^4)}
				&\lesssim_\sigma
				\lVert\overline f\,D_B^\varepsilon g
				\rVert_{L^{\frac{2}{3-2\sigma},\varepsilon}(\rho^4)}
				\le\lVert f\rVert_{L^{\frac1{1-\sigma},\varepsilon}(\rho^2)}
				\lVert D_B^\varepsilon g\rVert_{L^{2,\varepsilon}(\rho^2)}\\
				&\lesssim_\sigma
				\lVert f\rVert_{H_B^{2\sigma-1,\varepsilon}(\rho^2)}
				\lVert g\rVert_{H_B^{1,\varepsilon}(\rho^2)}.
			\end{aligned}
		\end{equation}
		The shifted form~\eqref{eq:gauge-current-shifted-form} gives
		\eqref{eq:gauge-current-second-mixed-regularity-bound} analogously.
		Taking $\sigma=1-\alpha$ gives~\eqref{eq:gauge-current-Z-bilinear-estimate}.
	\end{proof}

	Since $H^{-2\alpha,\varepsilon}(\rho^4)\hookrightarrow
	H^{-1,\varepsilon}(\rho^4)$, Lemma~\ref{lem:gauge-current-bilinear-estimate}
	and~\eqref{eq:higgs-summary-full-shift-bound} imply
	\begin{equation}
		\label{eq:gauge-regular-J4-final-bound}
		\left\lVert\mathbb E^\phi[J_4]\right\rVert_{H^{-1,\varepsilon}(\rho^4)}
		\lesssim\mathbb E^\phi\left[
		\lVert Z_1^\varepsilon\rVert_{H_B^{1-\alpha,\varepsilon}(\rho^2)}^2\right]
		\lesssim\mathbb E^\phi\left[\lVert\Xi\rVert_\kappa^{K_0}\right]
		+\lVert B\rVert_{H^{1,\varepsilon}(\rho)}^{\kappa_2}.
	\end{equation}

	\smallskip
	\noindent\textup{(5)} \textit{Conclusion.}
	Combining the estimates in \textup{(1)}--\textup{(4)} gives
	\begin{lemma}
		\label{lem:gauge-regular-terms-conclusion}
		We have
		\begin{equation}
			\label{eq:gauge-regular-terms-conclusion}
			\left\lVert\mathbb E^\phi\left[\sum_{j=1}^4J_j\right]
			\right\rVert_{H^{-1,\varepsilon}(\rho^4)}
			\lesssim\lVert B\rVert_{H^{1,\varepsilon}(\rho)}^{\kappa_2}
			+\mathbb E^\phi\left[\lVert\Xi\rVert_\kappa^{K_0}\right].
		\end{equation}
	\end{lemma}

	\subsection{Estimates of irregular terms}
	\label{subsec:gauge-irregular-terms}
	In this subsection, we estimate the terms $I_j^{\mathrm{Pink}}$,
	$1\leq j\leq3$, in
	\eqref{eq:gauge-gradient-current-final-splitting}.

	\subsubsection{\texorpdfstring{Estimate of the
			$I_1^{\mathrm{Pink}}$ and $I_2^{\mathrm{Pink}}$ terms}
		{Reduction of the I1 and I2 terms}}
	\label{subsubsec:gauge-irregular-I1-I2-reduction}

	For $I_1^{\mathrm{Pink}}$, the estimate will be proved in
	Section~\ref{subsec:stochastic-estimate-mean-covariant-current}.
	Together with the boundedness of
	$\eta^c(D)P_\perp^{M,\varepsilon}$ on
	$H^{-1,\varepsilon}(\rho^4)$ and the embedding
	$H^{-1,\varepsilon}(\rho)\hookrightarrow
	H^{-1,\varepsilon}(\rho^4)$, this gives
	\begin{equation}
		\label{eq:gauge-irregular-I1-target-bound}
		\left\lVert\mathbb{E}^\phi\left[I_1^{\mathrm{Pink}}\right]
		\right\rVert_{H^{-1,\varepsilon}(\rho^4)}
		\lesssim\log\left(2+\lVert B\rVert_{H^{1,\varepsilon}(\rho)}\right).
	\end{equation}
	For $I_2^{\mathrm{Pink}}$, the definition of $\rho_A^\varepsilon$ gives
	\begin{equation}
		\label{eq:gauge-irregular-I2-expectation}
		\begin{aligned}
			\mathbb{E}^\phi[I_2^{\mathrm{Pink}}]
			&=
			\mathbb{E}^\phi\left[
			\eta^c(D)P_\perp^{M,\varepsilon}\operatorname{Im}\left[
			i\overline{\phi_B^\varepsilon}\,
			U_B^\varepsilon\rho_A^\varepsilon\phi_B^\varepsilon
			\right]
			\right]\\
			&=
			\eta^c(D)P_\perp^{M,\varepsilon}\operatorname{Im}\left[
			i\left(
			\mathbb{E}^\phi\left[
			\overline{\phi_B^\varepsilon(x)}\,
			U_{B,j}^\varepsilon(x)
			\phi_B^\varepsilon(x+\varepsilon e_j)
			\right]
			r_{A,j}^\varepsilon(x)
			\right)_{j=1}^2
			\right].
		\end{aligned}
	\end{equation}
	\begin{lemma}
		\label{lem:stochastic-estimate-bB-coefficient}
		Let $A\in(\mathbb{R}^2)^{\Lambda_{M,\varepsilon}}$ and
		$B\in(\mathbb{R}^2)_{\chi_\varepsilon,\mathrm{div}}^{\Lambda_{M,\varepsilon}}$.
		We have
		\begin{equation}
			\label{eq:stochastic-estimate-bB-target}
			\begin{aligned}
				\left\lVert \mathbb{E}^\phi\bigl[I_2^{\mathrm{Pink}}\bigr]
				-\eta^c(D)P_\perp^{M,\varepsilon}\operatorname{Im}\left[
				\bigl(i b_\varepsilon r_{A,j}^\varepsilon\bigr)_{j=1}^2\right]
				\right\rVert_{H^{-2\kappa,\varepsilon}(\rho^4)}\lesssim \mathbb{E}^\phi\left[\lVert\Xi\rVert_\kappa^{K_0}\right]
				+\lVert B\rVert_{H^{1,\varepsilon}(\rho)}^{\kappa_2}.
			\end{aligned}
		\end{equation}
	\end{lemma}

	\begin{proof}
		For fixed $j\in\{1,2\}$, we have the exact decomposition
		\begin{equation}
			\label{eq:stochastic-estimate-bB-remainder-decomposition}
			\begin{aligned}
			\mathbb{E}^\phi\left[\overline{\phi_B^\varepsilon}\,
			U_{B,j}^\varepsilon\rho_{A,j}^\varepsilon\phi_B^\varepsilon\right]
			&=\mathbb{E}^\phi\left[
			\bigl(\overline{\phi_B^\varepsilon}-\overline{\phi^\varepsilon}\bigr)
			U_{B,j}^\varepsilon\rho_{A,j}^\varepsilon\phi_B^\varepsilon\right]
			+\mathbb{E}^\phi\left[\overline{\phi^\varepsilon}\,
			U_{B,j}^\varepsilon r_{A,j}^\varepsilon
			\tau_j\bigl(\phi_B^\varepsilon-\phi^\varepsilon\bigr)\right]\\
			&\quad+\mathbb{E}^\phi\left[\overline{\phi^\varepsilon}
			\bigl(U_{B,j}^\varepsilon-1\bigr)
			\rho_{A,j}^\varepsilon\phi^\varepsilon\right]
			+\mathbb{E}^\phi\left[\overline{\phi^\varepsilon}
			\bigl(\tau_j\phi^\varepsilon-\phi^\varepsilon\bigr)
			r_{A,j}^\varepsilon\right]
			+b_\varepsilon r_{A,j}^\varepsilon\\
			&:=\sum_{\nu=1}^4R_{\nu,j}^\varepsilon
			+b_\varepsilon r_{A,j}^\varepsilon.
			\end{aligned}
		\end{equation}
		Here we used
		$\mathbb{E}^\phi[\lvert\phi^\varepsilon\rvert^2]=b_\varepsilon$.

		By Lemma~\ref{lem:discrete-paraproduct-estimates},
		\eqref{eq:stochastic-object-Xi-norm}, and
		\eqref{eq:gauge-link-multiplier-bound},
		\begin{equation}
			\label{eq:stochastic-estimate-bB-R1-bound}
			\begin{aligned}
			\lVert R_{1,j}^\varepsilon\rVert_{\mathcal{C}^{-\kappa,\varepsilon}(\rho^{19\kappa})}
			&\lesssim\mathbb{E}^\phi\left[
			\lVert\phi_B^\varepsilon-\phi^\varepsilon
			\rVert_{\mathcal{C}^{2\kappa,\varepsilon}(\rho^{16\kappa})}
			\left\lVert U_{B,j}^\varepsilon
			\rho_{A,j}^\varepsilon\phi_B^\varepsilon
			\right\rVert_{\mathcal{C}^{-\kappa,\varepsilon}(\rho^{3\kappa})}\right]\\
			&\lesssim\mathbb{E}^\phi\left[
			\lVert\Xi\rVert_\kappa
			\bigl(1+\lVert B\rVert_{H^{1,\varepsilon}(\rho)}\bigr)^{18\kappa}
			\lVert\rho_A^\varepsilon\phi_B^\varepsilon
			\rVert_{\mathcal{C}^{-\kappa,\varepsilon}(\rho^\kappa)}\right]\\
			&\lesssim\mathbb{E}^\phi\left[\lVert\Xi\rVert_\kappa^{K_0}\right]
			+\lVert B\rVert_{H^{1,\varepsilon}(\rho)}^{\kappa_2}.
			\end{aligned}
		\end{equation}
		For the second remainder, conjugating the first factor and using
		the boundedness of $\tau_j$ on the weighted spaces gives
		\begin{equation}
			\label{eq:stochastic-estimate-bB-R2-bound}
			\begin{aligned}
			\lVert R_{2,j}^\varepsilon\rVert_{\mathcal{C}^{-\kappa,\varepsilon}(\rho^{19\kappa})}
			&\lesssim\mathbb{E}^\phi\left[
			\left\lVert U_{-B,j}^\varepsilon
			\overline{r_{A,j}^\varepsilon}\phi^\varepsilon
			\right\rVert_{\mathcal{C}^{-\kappa,\varepsilon}(\rho^{3\kappa})}
			\lVert\phi_B^\varepsilon-\phi^\varepsilon
			\rVert_{\mathcal{C}^{2\kappa,\varepsilon}(\rho^{16\kappa})}\right]\\
			&\lesssim\mathbb{E}^\phi\left[
			\lVert\Xi\rVert_\kappa^2
			\bigl(1+\lVert B\rVert_{H^{1,\varepsilon}(\rho)}\bigr)^{18\kappa}\right]\\
			&\lesssim\mathbb{E}^\phi\left[\lVert\Xi\rVert_\kappa^{K_0}\right]
			+\lVert B\rVert_{H^{1,\varepsilon}(\rho)}^{\kappa_2}.
			\end{aligned}
		\end{equation}
		For the third remainder, the pointwise estimate
		$\lvert e^{i\theta}-1\rvert\lesssim\min\{1,\lvert\theta\rvert\}$
		and the dyadic decomposition give
		\begin{equation}
			\label{eq:stochastic-estimate-bB-R3-bound}
			\begin{aligned}
			\lVert R_{3,j}^\varepsilon\rVert_{L^{2,\varepsilon}(\rho^4)}
			&\lesssim\mathbb{E}^\phi\left[
			\lVert\phi^\varepsilon\rVert_{L^{\infty,\varepsilon}(\rho^\kappa)}^2
			\lVert r_{A,j}^\varepsilon\rVert_{L^{\infty,\varepsilon}(\rho^\kappa)}
			\lVert U_{B,j}^\varepsilon-1\rVert_{L^{2,\varepsilon}(\rho^{4-3\kappa})}\right]\\
			&\lesssim\varepsilon^{-6\kappa}
			\min\left\{1,\varepsilon\lVert B\rVert_{H^{1,\varepsilon}(\rho)}\right\}
			\mathbb{E}^\phi\left[
			\lVert\phi^\varepsilon\rVert_{\mathcal{C}^{-\kappa,\varepsilon}(\rho^\kappa)}^2
			\lVert r_A^\varepsilon\rVert_{\mathcal{C}^{-\kappa,\varepsilon}(\rho^\kappa)}\right]\\
			&\lesssim\lVert B\rVert_{H^{1,\varepsilon}(\rho)}^{6\kappa}
			\mathbb{E}^\phi\left[\lVert\Xi\rVert_\kappa^3\right]\\
			&\lesssim\mathbb{E}^\phi\left[\lVert\Xi\rVert_\kappa^{K_0}\right]
			+\lVert B\rVert_{H^{1,\varepsilon}(\rho)}^{\kappa_2}.
			\end{aligned}
		\end{equation}
		Since
		$\mathbb{E}^\phi[\overline{\phi^\varepsilon}
		(\tau_j\phi^\varepsilon-\phi^\varepsilon)]$ is independent of the
		spatial variable,
		\begin{equation}
			\label{eq:stochastic-estimate-bB-R4-stationarity}
			\lVert R_{4,j}^\varepsilon\rVert_{\mathcal{C}^{-\kappa,\varepsilon}(\rho^\kappa)}
			=\left\lvert\mathbb{E}^\phi\left[
			\overline{\phi^\varepsilon}
			\bigl(\tau_j\phi^\varepsilon-\phi^\varepsilon\bigr)\right]\right\rvert
			\lVert r_{A,j}^\varepsilon\rVert_{\mathcal{C}^{-\kappa,\varepsilon}(\rho^\kappa)}.
		\end{equation}
		By~\eqref{eq:free-gff-fourier-expansion} and
		\eqref{eq:free-gff-fourier-covariance},
		\begin{equation}
			\label{eq:stochastic-estimate-bB-shifted-free-covariance}
			\begin{aligned}
				\left\lvert\mathbb{E}^\phi\left[\overline{\phi^\varepsilon}
				\bigl(\tau_j\phi^\varepsilon-\phi^\varepsilon\bigr)\right]\right\rvert
				=\frac{2}{\lvert\Lambda_{M,\varepsilon}\rvert}
				\left\lvert\sum_{q\in Q_{M,\varepsilon}}
				\frac{e^{i\varepsilon q_j}-1}{1+\lvert q^\varepsilon\rvert^2}\right\rvert
				\lesssim\frac{\varepsilon}{\lvert\Lambda_{M,\varepsilon}\rvert}
				\sum_{q\in Q_{M,\varepsilon}}
				\frac{\lvert q\rvert}{1+\lvert q\rvert^2}\lesssim1.
			\end{aligned}
		\end{equation}
		Consequently,
		\begin{equation}
			\label{eq:stochastic-estimate-bB-R4-bound}
			\lVert R_{4,j}^\varepsilon\rVert_{\mathcal{C}^{-\kappa,\varepsilon}(\rho^\kappa)}
			\lesssim\lVert r_A^\varepsilon\rVert_{\mathcal{C}^{-\kappa,\varepsilon}(\rho^\kappa)}.
		\end{equation}

		Finally,~\eqref{eq:gauge-irregular-I2-expectation} and
		\eqref{eq:stochastic-estimate-bB-remainder-decomposition}, together
		with the uniform boundedness of $\eta^c(D)P_\perp^{M,\varepsilon}$ on
		$H^{-2\kappa,\varepsilon}(\rho^4)$, bound the left-hand side of
		\eqref{eq:stochastic-estimate-bB-target} by the sum of the
		$H^{-2\kappa,\varepsilon}(\rho^4)$ norms of the four remainders.
		Since $\kappa>0$ and $4-19\kappa>1/4$, the weighted embeddings give
		$\mathcal{C}^{-\kappa,\varepsilon}(\rho^{19\kappa})
		\hookrightarrow H^{-2\kappa,\varepsilon}(\rho^4)$ and
		$L^{2,\varepsilon}(\rho^4)\hookrightarrow H^{-2\kappa,\varepsilon}(\rho^4)$.
		Together with~\eqref{eq:stochastic-object-Xi-norm} and the estimates
		above, this yields
		\eqref{eq:stochastic-estimate-bB-target}.
	\end{proof}

	\subsubsection{\texorpdfstring{Decomposition of $I_3^{\mathrm{Pink}}$
			and the $Z_1^{\mathrm{hi},2}$ contribution}
		{Decomposition of I3 and the Zhi2 contribution}}
	\label{subsubsec:gauge-irregular-I3-Zhi2-contribution}

	We next consider $I_3^{\mathrm{Pink}}$. By
	\eqref{eq:gauge-gradient-current-terms},
	\begin{equation}
		\label{eq:gauge-irregular-I3-identity}
		\mathbb{E}^\phi[I_3^{\mathrm{Pink}}]
		=
		\mathbb{E}^\phi\left[
		T_B^\varepsilon(\phi_B^\varepsilon,Z_1^{\mathrm{hi}})
		+T_B^\varepsilon(Z_1^{\mathrm{hi}},\phi_B^\varepsilon)
		\right].
	\end{equation}
	In accordance with the two terms in
	\eqref{eq:higgs-high-feedback-equation}, write
	$Z_1^{\mathrm{hi}}=Z_1^{\mathrm{hi},1}+Z_1^{\mathrm{hi},2}$, where
	\begin{equation}
		\label{eq:gauge-irregular-Zhi1-definition}
		Z_1^{\mathrm{hi},1}
		=-\bigl(\mathcal{L}_B^\varepsilon+1\bigr)^{-1}
		\int_0^1\mathbb{E}_t^\phi\Bigl[
		iD_B^{\varepsilon *}\Bigl(
		\mathcal{U}_{>L_1}^\varepsilon
		(\overline{r_A^\varepsilon}\phi_B^\varepsilon)
		+U_B^\varepsilon\mathcal{U}_{>L_1}^\varepsilon
		(\rho_A^\varepsilon\phi_B^\varepsilon)
		\Bigr)
		\Bigr]\,dt
	\end{equation}
	and
	\begin{equation}
		\label{eq:gauge-irregular-Zhi2-definition}
		Z_1^{\mathrm{hi},2}
		=-\bigl(\mathcal{L}_B^\varepsilon+1\bigr)^{-1}
		\int_0^1\mathbb{E}_t^\phi\Bigl[
		iD_B^{\varepsilon *}\Bigl(
		(\mathcal{U}_{>L_2}^\varepsilon
		\overline{r_A^\varepsilon})\succ Z_1^\varepsilon
		+U_B^\varepsilon\bigl(
		(\mathcal{U}_{>L_2}^\varepsilon\rho_A^\varepsilon)
		\succ Z_1^\varepsilon
		\bigr)
		\Bigr)
		\Bigr]\,dt.
	\end{equation}
	Consequently,
	\begin{equation}
		\label{eq:gauge-irregular-I3-four-term-splitting}
		I_3^{\mathrm{Pink}}
		=T_B^\varepsilon
		(\phi_B^\varepsilon,Z_1^{\mathrm{hi},1})
		+T_B^\varepsilon
		(Z_1^{\mathrm{hi},1},\phi_B^\varepsilon)
		+T_B^\varepsilon
		(\phi_B^\varepsilon,Z_1^{\mathrm{hi},2})
		+T_B^\varepsilon
		(Z_1^{\mathrm{hi},2},\phi_B^\varepsilon).
	\end{equation}

	We first estimate the $Z_1^{\mathrm{hi},2}$ contribution.
	\begin{lemma}
		\label{lem:gauge-irregular-Zhi2-derivative-bound}
		With the choice of $L_2$ in
		Table~\ref{tab:higgs-high-frequency-cutoff-choices},
		\begin{equation}
			\label{eq:gauge-irregular-Zhi2-derivative-bound}
			\mathbb E^\phi\left[
			\sum_{j=1}^2\left(
			\left\lVert D_{B,j}^\varepsilon Z_1^{\mathrm{hi},2}
			\right\rVert_{H^{1-3\alpha,\varepsilon}}^2
			+\left\lVert D_{B,j}^{\varepsilon *}Z_1^{\mathrm{hi},2}
			\right\rVert_{H^{1-3\alpha,\varepsilon}}^2
			\right)\right]
			\lesssim\mathbb E^\phi\left[\lVert\Xi\rVert_\kappa^{K_0}\right].
		\end{equation}
	\end{lemma}

	\begin{proof}
		The discrete product rule and
		$D_B^{\varepsilon *}=D^{\varepsilon *}-i\rho_B^{\varepsilon *}$
		reduce the first term in
		\eqref{eq:gauge-irregular-Zhi2-definition} to the following three types:
		\begin{equation}
			\label{eq:gauge-irregular-Zhi2-three-types}
			\begin{aligned}
				\bigl(D^{\varepsilon *}\mathcal U_{>L_2}^\varepsilon
				\overline{r_A^\varepsilon}\bigr)\succ Z_1^\varepsilon,\qquad
				\sum_{j=1}^2\bigl(\tau_j^*\mathcal U_{>L_2}^\varepsilon
				\overline{r_{A,j}^\varepsilon}\bigr)
				\succ D_j^{\varepsilon *}Z_1^\varepsilon,\quad
				\rho_B^{\varepsilon *}\Bigl(
				(\mathcal U_{>L_2}^\varepsilon\overline{r_A^\varepsilon})
				\succ Z_1^\varepsilon\Bigr).
			\end{aligned}
		\end{equation}
		Since $\kappa-(\alpha-\kappa)\sigma_2<-8$,
		\eqref{eq:high-frequency-truncation-estimate},
		\eqref{eq:high-frequency-truncation-commutator-estimate}
		with $\delta=\alpha-\kappa$, and the
		choice of $L_2$ in Table~\ref{tab:higgs-high-frequency-cutoff-choices} give
		\begin{equation}
			\label{eq:gauge-irregular-L2-gain}
			\begin{aligned}
				&\left\lVert\mathcal U_{>L_2}^\varepsilon r_A^\varepsilon
				\right\rVert_{\mathcal C^{-\alpha,\varepsilon}(\rho^{-6})}
				+\left\lVert\mathcal U_{>L_2}^\varepsilon
				D^{\varepsilon *}\overline{r_A^\varepsilon}
				\right\rVert_{\mathcal C^{-\alpha,\varepsilon}(\rho^{-6})}
				+\left\lVert[D^{\varepsilon *},\mathcal U_{>L_2}^\varepsilon]
				\overline{r_A^\varepsilon}
				\right\rVert_{\mathcal C^{-\alpha,\varepsilon}(\rho^{-6})}\\
				&\qquad\lesssim\lVert\Xi\rVert_\kappa
				2^{-(\alpha-\kappa)L_2}
				\lesssim\frac{1}{
				\bigl(1+\lVert B\rVert_{H^{1,\varepsilon}(\rho)}\bigr)^{l_0}
				\bigl(1+\lVert Z_1^\varepsilon\rVert_{L^{2,\varepsilon}(\rho^2)}\bigr)}.
			\end{aligned}
		\end{equation}

		By the weighted paraproduct estimates,
		\begin{equation}
			\label{eq:gauge-irregular-Zhi2-first-type-bound}
			\begin{aligned}
				\left\lVert
				\bigl(D^{\varepsilon *}\mathcal U_{>L_2}^\varepsilon
				\overline{r_A^\varepsilon}\bigr)\succ Z_1^\varepsilon
				\right\rVert_{H^{-2\alpha,\varepsilon}(\rho^{-4})}
			\lesssim
				\left\lVert D^{\varepsilon *}\mathcal U_{>L_2}^\varepsilon
				\overline{r_A^\varepsilon}
				\right\rVert_{\mathcal C^{-\alpha,\varepsilon}(\rho^{-6})}
				\lVert Z_1^\varepsilon\rVert_{L^{2,\varepsilon}(\rho^2)}
				\lesssim
				\bigl(1+\lVert B\rVert_{H^{1,\varepsilon}(\rho)}\bigr)^{-l_0}.
			\end{aligned}
		\end{equation}
		For the second type, the weighted comparison
		\eqref{eq:weighted-HB-H-comparison} gives
		\begin{equation}
			\label{eq:gauge-irregular-Zhi2-second-type-bound}
			\begin{aligned}
				\left\lVert
				\bigl(\tau_j^*\mathcal U_{>L_2}^\varepsilon
				\overline{r_{A,j}^\varepsilon}\bigr)
				\succ D_j^{\varepsilon *}Z_1^\varepsilon
				\right\rVert_{H^{-3\alpha,\varepsilon}(\rho^{-3})}
				&\lesssim
				\left\lVert\mathcal U_{>L_2}^\varepsilon r_A^\varepsilon
				\right\rVert_{\mathcal C^{-\alpha,\varepsilon}(\rho^{-6})}
				\lVert D_j^{\varepsilon *}Z_1^\varepsilon
				\rVert_{H^{-\alpha,\varepsilon}(\rho^{3-\alpha})}\\
				&\lesssim
				\left\lVert\mathcal U_{>L_2}^\varepsilon r_A^\varepsilon
				\right\rVert_{\mathcal C^{-\alpha,\varepsilon}(\rho^{-6})}
				\bigl(1+\lVert B\rVert_{H^{1,\varepsilon}(\rho)}\bigr)
				\lVert Z_1^\varepsilon\rVert_{H_B^{1-\alpha,\varepsilon}(\rho^2)}\\
				&\lesssim
				\frac{\lVert Z_1^\varepsilon
				\rVert_{H_B^{1-\alpha,\varepsilon}(\rho^2)}}{
				\bigl(1+\lVert B\rVert_{H^{1,\varepsilon}(\rho)}\bigr)^{l_0/2}}.
			\end{aligned}
		\end{equation}

		For the third type, take $p=2/\alpha$. The uniform embeddings
$
			B_{\frac{2p}{p+2},\infty}^{-\alpha,\varepsilon}(\rho^{-3})
			\hookrightarrow B_{2,\infty}^{-2\alpha,\varepsilon}(\rho^{-3})
			\hookrightarrow H^{-3\alpha,\varepsilon}(\rho^{-3})
$
		and $H_B^{1-\alpha,\varepsilon}(\rho^2)
		\hookrightarrow L^{p,\varepsilon}(\rho^2)$ give
		\begin{equation}
			\label{eq:gauge-irregular-Zhi2-third-type-bound}
			\begin{aligned}
				\left\lVert\rho_B^{\varepsilon *}\Bigl(
				(\mathcal U_{>L_2}^\varepsilon\overline{r_A^\varepsilon})
				\succ Z_1^\varepsilon\Bigr)
				\right\rVert_{H^{-3\alpha,\varepsilon}(\rho^{-3})}
				&\lesssim
				\left\lVert\rho_B^{\varepsilon *}\Bigl(
				(\mathcal U_{>L_2}^\varepsilon\overline{r_A^\varepsilon})
				\succ Z_1^\varepsilon\Bigr)
				\right\rVert_{B_{\frac{2p}{p+2},\infty}^{-\alpha,\varepsilon}(\rho^{-3})}\\
				&\lesssim
				\lVert r_B^\varepsilon\rVert_{H^{1,\varepsilon}(\rho)}
				\left\lVert\mathcal U_{>L_2}^\varepsilon r_A^\varepsilon
				\right\rVert_{\mathcal C^{-\alpha,\varepsilon}(\rho^{-6})}
				\lVert Z_1^\varepsilon\rVert_{L^{p,\varepsilon}(\rho^2)}\\
				&\lesssim
				\frac{\lVert Z_1^\varepsilon
				\rVert_{H_B^{1-\alpha,\varepsilon}(\rho^2)}}{
				\bigl(1+\lVert B\rVert_{H^{1,\varepsilon}(\rho)}\bigr)^{l_0/2}}.
			\end{aligned}
		\end{equation}
		For the term containing $U_B^\varepsilon$, apply
		\eqref{eq:covariant-derivative-translation-identities} and
		\eqref{eq:rho-A-operator-paraproduct-convention} componentwise.
		Expanding $D_{B,j}^\varepsilon=D_j^\varepsilon+i\rho_{B,j}^\varepsilon$
		produces the same three types, with translated coefficients and
		$D_j^\varepsilon Z_1^\varepsilon$ in place of
		$D_j^{\varepsilon *}Z_1^\varepsilon$. The preceding estimates apply.
		Consequently,
		\begin{equation}
			\label{eq:gauge-irregular-Zhi2-source-bound}
			\left\lVert iD_B^{\varepsilon *}\Bigl(
			(\mathcal U_{>L_2}^\varepsilon\overline{r_A^\varepsilon})
			\succ Z_1^\varepsilon
			+U_B^\varepsilon\bigl(
			(\mathcal U_{>L_2}^\varepsilon\rho_A^\varepsilon)
			\succ Z_1^\varepsilon\bigr)\Bigr)
			\right\rVert_{H^{-3\alpha,\varepsilon}(\rho^{-3})}
			\lesssim
			\frac{\lVert Z_1^\varepsilon
			\rVert_{H_B^{1-\alpha,\varepsilon}(\rho^2)}}{
			\bigl(1+\lVert B\rVert_{H^{1,\varepsilon}(\rho)}\bigr)^{l_0/2}}.
		\end{equation}
		
		Combining~\eqref{eq:weighted-HB-H-comparison},
		Lemma~\ref{lem:covariant-fbsde-schauder-estimate} with
		$s=1-3\alpha$ and $a=-3+3\alpha$, and
		\eqref{eq:gauge-irregular-Zhi2-source-bound}, we obtain
		\begin{equation}
			\label{eq:gauge-irregular-Zhi2-HB-bound}
			\mathbb E^\phi\left[
			\lVert Z_1^{\mathrm{hi},2}
			\rVert_{H_B^{2-3\alpha,\varepsilon}(\rho^{-3+3\alpha})}^2\right]
			\lesssim
			\frac{\mathbb E^\phi\left[
			\lVert Z_1^\varepsilon
			\rVert_{H_B^{1-\alpha,\varepsilon}(\rho^2)}^2\right]}{
			\bigl(1+\lVert B\rVert_{H^{1,\varepsilon}(\rho)}\bigr)^{l_0-6\alpha}}.
		\end{equation}
		Moreover,~\eqref{eq:weighted-covariant-derivative-free-output}
		with $s=1-3\alpha$ and $a=-3+3\alpha$,  gives
		\begin{equation}
			\label{eq:gauge-irregular-Zhi2-derivative-comparison}
			\lVert D_{B,j}^\varepsilon Z_1^{\mathrm{hi},2}
			\rVert_{H^{1-3\alpha,\varepsilon}}
			\lesssim
			\bigl(1+\lVert B\rVert_{H^{1,\varepsilon}(\rho)}\bigr)^{3-9\alpha}
			\lVert Z_1^{\mathrm{hi},2}
			\rVert_{H_B^{2-3\alpha,\varepsilon}(\rho^{-3+3\alpha})}.
		\end{equation}
		Hence, for $j=1,2$,
		\begin{equation}
			\label{eq:gauge-irregular-Zhi2-pre-moment-bound}
			\begin{aligned}
				\mathbb E^\phi\left[
				\lVert D_{B,j}^\varepsilon Z_1^{\mathrm{hi},2}
				\rVert_{H^{1-3\alpha,\varepsilon}}^2\right]
				\lesssim
				\frac{\mathbb E^\phi\left[
				\lVert Z_1^\varepsilon
				\rVert_{H_B^{1-\alpha,\varepsilon}(\rho^2)}^2\right]}{
				\bigl(1+\lVert B\rVert_{H^{1,\varepsilon}(\rho)}\bigr)^{l_0-6+12\alpha}}
				\lesssim
				\frac{\mathbb E^\phi\left[
				\lVert Z_1^\varepsilon
				\rVert_{H_B^{1-\alpha,\varepsilon}(\rho^2)}^2\right]}{
				\bigl(1+\lVert B\rVert_{H^{1,\varepsilon}(\rho)}\bigr)^{l_0/2}}.
			\end{aligned}
		\end{equation}
		Finally, Proposition~\ref{prop:higgs-terminal-shift-summary} yields
		\begin{equation}
			\label{eq:gauge-irregular-Zhi2-holder}
			\begin{aligned}
				\frac{\mathbb E^\phi\left[
				\lVert Z_1^\varepsilon
				\rVert_{H_B^{1-\alpha,\varepsilon}(\rho^2)}^2\right]}{
				\bigl(1+\lVert B\rVert_{H^{1,\varepsilon}(\rho)}\bigr)^{l_0/2}}
				\lesssim
				\frac{\bigl(1+\lVert B\rVert_{H^{1,\varepsilon}(\rho)}\bigr)^{\kappa_1}
				+\mathbb E^\phi\left[\lVert\Xi\rVert_\kappa^{K_0}\right]}{
				\bigl(1+\lVert B\rVert_{H^{1,\varepsilon}(\rho)}\bigr)^{l_0/2}}
				\lesssim\mathbb E^\phi\left[\lVert\Xi\rVert_\kappa^{K_0}\right].
			\end{aligned}
		\end{equation}
		The same estimates hold for $D_{B,j}^{\varepsilon *}$.
		Together with~\eqref{eq:gauge-irregular-Zhi2-pre-moment-bound},
		summing over $j=1,2$ yields~\eqref{eq:gauge-irregular-Zhi2-derivative-bound}.
	\end{proof}

	\begin{corollary}
		\label{cor:gauge-irregular-I3-Zhi2-bound}
		We have
		\begin{equation}
			\label{eq:gauge-irregular-I3-Zhi2-final-bound}
			\begin{aligned}
				\left\lVert\mathbb E^\phi\left[
				T_B^\varepsilon(\phi_B^\varepsilon,Z_1^{\mathrm{hi},2})
				+T_B^\varepsilon(Z_1^{\mathrm{hi},2},\phi_B^\varepsilon)
				\right]\right\rVert_{H^{-1,\varepsilon}(\rho^4)}
				&\lesssim\mathbb E^\phi\left[
				\lVert\Xi\rVert_\kappa\sum_{j=1}^2\left(
				\lVert D_{B,j}^\varepsilon Z_1^{\mathrm{hi},2}
				\rVert_{H^{\alpha,\varepsilon}}
				+\lVert D_{B,j}^{\varepsilon *}Z_1^{\mathrm{hi},2}
				\rVert_{H^{\alpha,\varepsilon}}
				\right)\right]\\
				&\lesssim\mathbb E^\phi\left[\lVert\Xi\rVert_\kappa^{K_0}\right].
			\end{aligned}
		\end{equation}
	\end{corollary}

	\begin{proof}
	The first inequality follows from
		the weighted product estimate used in
		\eqref{eq:gauge-regular-J1-first-product-bound},
		the shifted form~\eqref{eq:gauge-regular-J1-second-shifted-form},
		and~\eqref{eq:covariant-derivative-translation-identities}.
		The second follows from
		Lemma~\ref{lem:gauge-irregular-Zhi2-derivative-bound} and Cauchy--Schwarz.
	\end{proof}

\subsubsection{\texorpdfstring{The renormalized
			$Z_1^{\mathrm{hi},1}$ contribution}
		{The renormalized Zhi1 contribution}}
	\label{subsubsec:gauge-irregular-Zhi1-renormalized-contribution}

	For $j=1,2$, define $\tau_j^0=\operatorname{Id}$,
	$\tau_j^1=\tau_j$,
	$r_{A,j}^{\varepsilon,(0)}
	:=
	\overline{r_{A,j}^\varepsilon}$,
	and
	$r_{A,j}^{\varepsilon,(1)}
	:=
	r_{A,j}^\varepsilon$.
	Expanding the first two terms in
	\eqref{eq:gauge-irregular-I3-four-term-splitting} by
	\eqref{eq:gauge-irregular-Zhi1-definition} gives
	\begin{equation}
		\label{eq:gauge-irregular-Zhi1-unified-ab-expansion}
		\begin{aligned}
		&T_B^\varepsilon
		(\phi_B^\varepsilon,Z_1^{\mathrm{hi},1})
		+T_B^\varepsilon
		(Z_1^{\mathrm{hi},1},\phi_B^\varepsilon)\\
		={}&-
		\sum_{a,b=0}^1
		\eta^c(D)P_\perp^{M,\varepsilon}\Bigg(
		\operatorname{Im}\Bigg[
		i\overline{\tau_j^a\phi_B^\varepsilon}\,
		U_{-aB,j}^\varepsilon D_{B,j}^\varepsilon
		\bigl(\mathcal{L}_B^\varepsilon+1\bigr)^{-1}
		\int_0^1\mathbb{E}_t^\phi\Bigg[
		\sum_{\ell=1}^2
		D_{B,\ell}^{\varepsilon *}U_{bB,\ell}^\varepsilon
		\mathcal{U}_{>L_1}^\varepsilon
		\bigl(
		r_{A,\ell}^{\varepsilon,(b)}
		\tau_\ell^b\phi_B^\varepsilon
		\bigr)
		\Bigg]dt
		\Bigg]
		\Bigg)_{j=1}^2.
		\end{aligned}
	\end{equation}
	We compare each term in
	\eqref{eq:gauge-irregular-Zhi1-unified-ab-expansion} with its free
	counterpart.
	\begin{lemma}
		\label{lem:gauge-irregular-Zhi1-unified-free-comparison}
		For every $j=1,2$ and $a,b\in\{0,1\}$,
		\begin{equation}
			\label{eq:gauge-irregular-Zhi1-unified-free-comparison}
			\begin{aligned}
			&\Bigg\lVert
			\mathbb{E}^\phi\Bigg[
			i\overline{\tau_j^a\phi_B^\varepsilon}\,
			U_{-aB,j}^\varepsilon D_{B,j}^\varepsilon
			\bigl(\mathcal{L}_B^\varepsilon+1\bigr)^{-1}
			\int_0^1\mathbb{E}_t^\phi\Bigg[
			\sum_{\ell=1}^2
			D_{B,\ell}^{\varepsilon *}U_{bB,\ell}^\varepsilon
			\mathcal{U}_{>L_1}^\varepsilon
			\bigl(
			r_{A,\ell}^{\varepsilon,(b)}
			\tau_\ell^b\phi_B^\varepsilon
			\bigr)
			\Bigg]dt
			\Bigg]\\
			&\qquad
			-\mathbb{E}^\phi\Bigg[
			i\overline{\tau_j^a\phi^\varepsilon}\,
			D_j^\varepsilon
			\bigl(\mathcal{L}_0^\varepsilon+1\bigr)^{-1}
			\int_0^1\mathbb{E}_t^\phi\Bigg[
			\sum_{\ell=1}^2
			D_\ell^{\varepsilon *}
			\mathcal{U}_{>L_1}^\varepsilon
			\bigl(
			r_{A,\ell}^{\varepsilon,(b)}
			\tau_\ell^b\phi^\varepsilon
			\bigr)
			\Bigg]dt
			\Bigg]
			\Bigg\rVert_{H^{-\alpha,\varepsilon}(\rho^4)}
			\lesssim
			\mathbb{E}^\phi\left[
			\lVert\Xi\rVert_\kappa^{K_0}
			\right].
			\end{aligned}
		\end{equation}
	\end{lemma}
	The proof is deferred to
	Appendix~\ref{app:proof-gauge-irregular-Zhi1-unified-free-comparison}.

	\begin{lemma}
		\label{lem:gauge-irregular-Zhi1-unified-low-free-bound}
		For every $j=1,2$ and $a,b\in\{0,1\}$,
		\begin{equation}
			\label{eq:gauge-irregular-Zhi1-unified-low-free-bound}
			\left\lVert
			\mathbb{E}^\phi\Bigg[
			i\overline{\tau_j^a\phi^\varepsilon}\,
			D_j^\varepsilon
			\bigl(\mathcal{L}_0^\varepsilon+1\bigr)^{-1}
			\int_0^1\mathbb{E}_t^\phi\Bigg[
			\sum_{\ell=1}^2
			D_\ell^{\varepsilon *}
			\mathcal{U}_{\leq L_1}^\varepsilon
			\bigl(
			r_{A,\ell}^{\varepsilon,(b)}
			\tau_\ell^b\phi^\varepsilon
			\bigr)
			\Bigg]dt
			\Bigg]
			\right\rVert_{\mathcal C^{-\kappa,\varepsilon}(\rho)}
			\lesssim
			\mathbb{E}^\phi\left[
			\lVert\Xi\rVert_\kappa^{K_0}
			\right]
			+
			\lVert B\rVert_{H^{1,\varepsilon}(\rho)}^{\kappa_1}.
		\end{equation}
	\end{lemma}

	\begin{proof}
		By~\eqref{eq:low-frequency-truncation-estimate} and
		Table~\ref{tab:higgs-high-frequency-cutoff-choices},
		\begin{equation}
			\label{eq:gauge-irregular-Zhi1-unified-low-free-source-bound}
			\begin{aligned}
			\left\lVert
			\mathcal{U}_{\leq L_1}^\varepsilon
			\bigl(
			r_{A,\ell}^{\varepsilon,(b)}
			\tau_\ell^b\phi^\varepsilon
			\bigr)
			\right\rVert_{\mathcal C^{2\kappa,\varepsilon}(\rho^{1-\kappa})}
			&\lesssim
			2^{3\kappa L_1}
			\left\lVert
			r_{A,\ell}^{\varepsilon,(b)}
			\tau_\ell^b\phi^\varepsilon
			\right\rVert_{\mathcal C^{-\kappa,\varepsilon}(\rho^{1-\kappa-3\kappa\sigma_1})}\\
			&\lesssim
			\left[
			\bigl(1+\lVert B\rVert_{H^{1,\varepsilon}(\rho)}\bigr)^{l_0}
			\lVert\Xi\rVert_\kappa
			\right]^{\frac{3\kappa}{\alpha/2-\kappa}}
			\lVert\Xi\rVert_\kappa\\
			&\lesssim
			\lVert\Xi\rVert_\kappa^{K_0/2}
			+
			\lVert B\rVert_{H^{1,\varepsilon}(\rho)}^{\kappa_1/2}.
			\end{aligned}
		\end{equation}
		The free multiplier
		$D^\varepsilon(\mathcal{L}_0^\varepsilon+1)^{-1}D^{\varepsilon *}$
		is bounded on $\mathcal C^{2\kappa,\varepsilon}(\rho^{1-\kappa})$.
		Hence the paraproduct
		estimate and Cauchy--Schwarz yield
		\begin{equation}
			\begin{aligned}
			&\left\lVert
			\mathbb{E}^\phi\Bigg[
			i\overline{\tau_j^a\phi^\varepsilon}\,
			D_j^\varepsilon
			\bigl(\mathcal{L}_0^\varepsilon+1\bigr)^{-1}
			\int_0^1\mathbb{E}_t^\phi\Bigg[
			\sum_{\ell=1}^2
			D_\ell^{\varepsilon *}
			\mathcal{U}_{\leq L_1}^\varepsilon
			\bigl(
			r_{A,\ell}^{\varepsilon,(b)}
			\tau_\ell^b\phi^\varepsilon
			\bigr)
			\Bigg]dt
			\Bigg]
			\right\rVert_{\mathcal C^{-\kappa,\varepsilon}(\rho)}\\
			&\qquad\lesssim
			\left(
			\mathbb{E}^\phi
			\lVert\phi^\varepsilon\rVert_{\mathcal C^{-\kappa,\varepsilon}(\rho^\kappa)}^2
			\right)^{1/2}
			\Bigg(
			\mathbb{E}^\phi
			\Bigg(
			\sum_{\ell=1}^2
			\left\lVert
			\mathcal{U}_{\leq L_1}^\varepsilon
			\bigl(
			r_{A,\ell}^{\varepsilon,(b)}
			\tau_\ell^b\phi^\varepsilon
			\bigr)
			\right\rVert_{\mathcal C^{2\kappa,\varepsilon}(\rho^{1-\kappa})}
			\Bigg)^2
			\Bigg)^{1/2}\\
			&\qquad\lesssim
			\mathbb{E}^\phi\left[
			\lVert\Xi\rVert_\kappa^{K_0}
			\right]
			+
			\lVert B\rVert_{H^{1,\varepsilon}(\rho)}^{\kappa_1},
			\end{aligned}
		\end{equation}
		which gives the desired estimate.
	\end{proof}

	\begin{lemma}
		\label{lem:gauge-irregular-Zhi1-unified-renormalized-reduction}
		We have
		\begin{equation}
			\label{eq:gauge-irregular-Zhi1-unified-renormalized-bound}
			\left\lVert
			\mathbb{E}^\phi\left[
			T_B^\varepsilon
			(\phi_B^\varepsilon,Z_1^{\mathrm{hi},1})
			+T_B^\varepsilon
			(Z_1^{\mathrm{hi},1},\phi_B^\varepsilon)
			\right]
			+\eta^c(D)P_\perp^{M,\varepsilon}\operatorname{Im}\left[
			i b_\varepsilon r_A^\varepsilon
			\right]
			\right\rVert_{H^{-1,\varepsilon}(\rho^4)}
			\lesssim
			\mathbb{E}^\phi\left[
			\lVert\Xi\rVert_\kappa^{K_0}
			\right]
			+
			\lVert B\rVert_{H^{1,\varepsilon}(\rho)}^{\kappa_1}.
		\end{equation}
	\end{lemma}

	\begin{proof}
		Using~\eqref{eq:gauge-irregular-Zhi1-unified-ab-expansion},
		Lemmas~\ref{lem:gauge-irregular-Zhi1-unified-free-comparison}
		and~\ref{lem:gauge-irregular-Zhi1-unified-low-free-bound},
		and the weighted embeddings, we reduce the problem to estimating
		\begin{equation}
			\label{eq:gauge-irregular-Zhi1-unified-free-renormalized-target}
			\Bigg\lVert
			\eta^c(D)P_\perp^{M,\varepsilon}\operatorname{Im}\Bigg\{
			\mathbb{E}^\phi\Bigg[
			i\overline{\tau_j^a\phi^\varepsilon}\,
			D_j^\varepsilon
			\bigl(\mathcal{L}_0^\varepsilon+1\bigr)^{-1}
			\sum_{\ell=1}^2
			D_\ell^{\varepsilon *}
			r_{A,\ell}^{\varepsilon,(b)}
			\tau_\ell^b
			\int_0^1\mathbb{E}_t^\phi\bigl[\phi^\varepsilon\bigr]dt
			\Bigg]
			-\frac{i}{4}b_\varepsilon
			r_{A,j}^{\varepsilon,(b)}
			\Bigg\}_{j=1}^2
			\Bigg\rVert_{H^{-1,\varepsilon}(\rho^4)}.
		\end{equation}
		By the definition
		$\phi^\varepsilon=(\mathcal{L}_0^\varepsilon+1)^{-1/2}
		\mathbb{B}_1^\varepsilon$, for every deterministic linear operator
		$T$ we have
		\begin{equation}
			\label{eq:gauge-irregular-Zhi1-unified-free-conditional-covariance-identity}
			\begin{aligned}
			\mathbb{E}^\phi\left[
			\overline{\phi^\varepsilon}\,
			T\int_0^1\mathbb{E}_t^\phi\bigl[\phi^\varepsilon\bigr]dt
			\right]
			&=
			\int_0^1\mathbb{E}^\phi\left[
			(\mathcal{L}_0^\varepsilon+1)^{-1/2}
			\overline{\mathbb{B}_1^\varepsilon}\,
			T(\mathcal{L}_0^\varepsilon+1)^{-1/2}
			\mathbb{B}_t^\varepsilon
			\right]dt\\
			&=
			\int_0^1t\,dt\;
			\mathbb{E}^\phi\left[
			\overline{\phi^\varepsilon}\,T\phi^\varepsilon
			\right].
			\end{aligned}
		\end{equation}
		Thus it remains to estimate
		\begin{equation}
			\label{eq:gauge-irregular-Zhi1-unified-free-static-renormalized-target}
			\Bigg\lVert
			\operatorname{Im}\Bigg\{
			\mathbb{E}^\phi\Bigg[
			i\overline{\tau_j^a\phi^\varepsilon}\,
			D_j^\varepsilon
			\bigl(\mathcal{L}_0^\varepsilon+1\bigr)^{-1}
			\sum_{\ell=1}^2
			D_\ell^{\varepsilon *}
			r_{A,\ell}^{\varepsilon,(b)}
			\tau_\ell^b\phi^\varepsilon
			\Bigg]
			-\frac{i}{2}b_\varepsilon
			r_{A,j}^{\varepsilon,(b)}
			\Bigg\}_{j=1}^2
			\Bigg\rVert_{H^{-1,\varepsilon}(\rho^4)}.
		\end{equation}
		The estimate of
		\eqref{eq:gauge-irregular-Zhi1-unified-free-static-renormalized-target}
		is deferred to Lemma~\ref{lem:renormalized-free-current-estimate},
		which, together with
		the embedding $\mathcal C^{-2\kappa,\varepsilon}(\rho^\kappa)
		\hookrightarrow H^{-1,\varepsilon}(\rho^4)$ and
		\eqref{eq:stochastic-object-Xi-norm}, gives
		\begin{equation}
			\label{eq:gauge-irregular-Zhi1-unified-free-static-renormalized-bound}
			\eqref{eq:gauge-irregular-Zhi1-unified-free-static-renormalized-target}
			\lesssim_\kappa
			\lVert r_A^\varepsilon
			\rVert_{\mathcal C^{-\kappa,\varepsilon}(\rho^\kappa)}
			\leq
			\mathbb{E}^\phi\left[
			\lVert\Xi\rVert_\kappa^{K_0}
			\right],
		\end{equation}
		and hence
		\eqref{eq:gauge-irregular-Zhi1-unified-renormalized-bound} holds.
	\end{proof}

	\subsubsection{Conclusion}
	\label{subsubsec:gauge-irregular-conclusion}

	Combining
	\eqref{eq:gauge-irregular-I1-target-bound},
	\eqref{eq:stochastic-estimate-bB-target},
	\eqref{eq:gauge-irregular-I3-Zhi2-final-bound}, and
	\eqref{eq:gauge-irregular-Zhi1-unified-renormalized-bound}, and observing
	that the renormalization terms
	$\eta^c(D)P_\perp^{M,\varepsilon}\operatorname{Im}
	[(i b_\varepsilon r_{A,j}^\varepsilon)_{j=1}^2]$
	in~\eqref{eq:stochastic-estimate-bB-target} and
	\eqref{eq:gauge-irregular-Zhi1-unified-renormalized-bound} cancel, we obtain
	the following conclusion.

	\begin{lemma}
		\label{lem:gauge-irregular-terms-conclusion}
		We have
		\begin{equation}
			\label{eq:gauge-irregular-terms-conclusion}
			\left\lVert
			\mathbb{E}^\phi\left[
			\sum_{j=1}^3 I_j^{\mathrm{Pink}}
			\right]
			\right\rVert_{H^{-1,\varepsilon}(\rho^4)}
			\lesssim
			\lVert B\rVert_{H^{1,\varepsilon}(\rho)}^{\kappa_2}
			+\mathbb{E}^\phi\left[
			\lVert\Xi\rVert_\kappa^{K_0}
			\right].
		\end{equation}
	\end{lemma}

	\subsection{\texorpdfstring{Estimate of
			$B^\varepsilon$}
		{Estimate of the terminal gauge-field shift}}
	\label{subsec:gauge-terminal-shift-estimate}

	\begin{proposition}
		\label{lem:gauge-terminal-shift-moment-bound}
		For every $p\geq1$,
		\begin{equation}
			\label{eq:gauge-terminal-shift-moment-bound}
			\mathbb{E}\left[
			\left\lVert B^\varepsilon
			\right\rVert_{H^{1,\varepsilon}(\rho)}^p
			\right]
			\lesssim_p1.
		\end{equation}
	\end{proposition}

	\begin{proof}
		For $p\geq2$, Proposition~\ref{prop:gauge-low-shift-uniform-moments}
		and~\eqref{eq:gauge-high-stationary-schauder-bound} give
		\begin{equation}
			\label{eq:gauge-terminal-shift-schauder-reduction}
			\mathbb{E}\left[
			\left\lVert B^\varepsilon\right\rVert_{H^{1,\varepsilon}(\rho)}^p
			\right]
			\leq C_p\left(1+\mathbb{E}\left[
			\left\lVert
			\mathbb E^\phi\left[
			T_C^\varepsilon(\Phi^\varepsilon,\Phi^\varepsilon)
			\right]
			\right\rVert_{H^{-1,\varepsilon}(\rho^4)}^p
			\right]\right).
		\end{equation}
		Combining this estimate with
		\eqref{eq:gauge-gradient-current-final-splitting},
		Lemmas~\ref{lem:gauge-regular-terms-conclusion}
		and~\ref{lem:gauge-irregular-terms-conclusion}, conditional Jensen's inequality,
		and Proposition~\ref{lem:gauge-stochastic-object-uniform-moments}
		with moment order $pK_0$, we obtain
		\begin{equation}
			\label{eq:gauge-terminal-shift-pre-absorption}
			\begin{aligned}
				\mathbb{E}\left[
				\left\lVert B^\varepsilon\right\rVert_{H^{1,\varepsilon}(\rho)}^p
				\right]
				&\leq C_p\Bigl(
				\mathbb{E}\left[
				\left\lVert\Theta^\varepsilon
				\right\rVert_\kappa^{pK_0}\right]
				+\mathbb{E}\left[
				\left\lVert B^\varepsilon
				\right\rVert_{H^{1,\varepsilon}(\rho)}^{p\kappa_2}
				\right]\Bigr)\\
				&\leq C_p\Bigl(1+
				\mathbb{E}\left[
				\left\lVert B^\varepsilon
				\right\rVert_{H^{1,\varepsilon}(\rho)}^{p\kappa K_0}
				\right]
				+\mathbb{E}\left[
				\left\lVert B^\varepsilon
				\right\rVert_{H^{1,\varepsilon}(\rho)}^{p\kappa_2}
				\right]\Bigr)\\
				&\leq\frac12\mathbb{E}\left[
				\left\lVert B^\varepsilon
				\right\rVert_{H^{1,\varepsilon}(\rho)}^p
				\right]+C_p.
			\end{aligned}
		\end{equation}
		The last inequality follows from Young's inequality, since
		$\kappa K_0<1$ and $\kappa_2<1$. Absorbing the moment term into the
		left-hand side gives the claimed bound.
	\end{proof}

	\section{Tightness of the lattice Abelian--Higgs measures}
	\label{sec:main-result}

	In this section, we prove the tightness of the extended lattice
	Abelian--Higgs measures. Recall the smooth extension operator
	$\mathcal E^\varepsilon$ defined in~\eqref{eq:lattice-sharp-extension-definition}
	and the renormalized measure
	$\mu_{M,\varepsilon}^{\mathrm{YMH}}$ from
	\eqref{eq:renormalized-lattice-abelian-higgs-measure}.
	For every Borel set
	\begin{equation}
		\label{eq:extended-ymh-measure-definition}
		U\in\mathcal B\Bigl(
		\mathcal C^{-\nu}(\rho^{1+\nu};\mathbb R^2)
		\times\mathcal C^{-\nu}(\rho^{2+\nu};\mathbb C)
		\Bigr),
	\end{equation}
	define
	\begin{equation}
		\label{eq:extended-ymh-pushforward-measure}
		\begin{aligned}
			\mathcal E^\varepsilon_\#\mu_{M,\varepsilon}^{\mathrm{YMH}}(U)
			:=
			&\mu_{M,\varepsilon}^{\mathrm{YMH}}\Bigl(
			\bigl\{(A,\phi)\in
			\mathcal A_{\mathrm{gf}}^{M,\varepsilon}
			\times\mathbb C^{\Lambda_{M,\varepsilon}}:
			(\mathcal E^\varepsilon A,\mathcal E^\varepsilon\phi)\in U
			\bigr\}\Bigr).
		\end{aligned}
	\end{equation}

	\begin{proof}[Proof of Theorem~\ref{thm:intro-subsequential-limits-lattice-ymh}]
		By Lemma~\ref{lem:lattice-smooth-extension-besov-bound},
		\begin{equation}
			\label{eq:sharp-extension-negative-holder-loss}
			\left\lVert\mathcal E^\varepsilon f
			\right\rVert_{\mathcal C^{-\alpha_1}(\rho^a)}
			\lesssim
			\lVert f\rVert_{\mathcal C^{-\alpha_1,\varepsilon}(\rho^a)},
			\qquad a\in\{1,2\}.
		\end{equation}
		The embeddings
		$\mathcal C^{-\alpha_1}(\rho^a)\hookrightarrow
		\mathcal C^{-\nu}(\rho^{a+\nu})$, $a\in\{1,2\}$,
		are compact by~\cite[Lemma~2.23]{MP19}.
		It therefore suffices to prove that, for every $\delta>0$, there
		exists $R_\delta>0$ such that, uniformly in $M,\varepsilon$,
		\begin{equation}
			\label{eq:extended-ymh-tightness-tail-target}
			\mathcal E^\varepsilon_\#\mu_{M,\varepsilon}^{\mathrm{YMH}}
			\left(
			\left\{(A,\phi):
			\lVert A\rVert_{\mathcal C^{-\alpha_1}(\rho)}
			+\lVert\phi\rVert_{\mathcal C^{-\alpha_1}(\rho^2)}
			>R_\delta
			\right\}\right)\le\delta.
		\end{equation}
		By Markov's inequality,
		\eqref{eq:sharp-extension-negative-holder-loss}, and the stochastic
		representation~\eqref{eq:law-representation-renormalized-ymh-measure},
		\begin{equation}
			\label{eq:extended-ymh-tightness-markov-bound}
			\begin{aligned}
				\mathcal E^\varepsilon_\#\mu_{M,\varepsilon}^{\mathrm{YMH}}
				\left(\left\{(A,\phi):
				\lVert A\rVert_{\mathcal C^{-\alpha_1}(\rho)}
				+\lVert\phi\rVert_{\mathcal C^{-\alpha_1}(\rho^2)}
				>R_\delta\right\}\right)\lesssim\frac1{R_\delta^2}\,
				\mathbb E\left[
				\left\lVert C^\varepsilon
				\right\rVert_{\mathcal C^{-\alpha_1,\varepsilon}(\rho)}^2
				+\left\lVert\Phi^\varepsilon
				\right\rVert_{\mathcal C^{-\alpha_1,\varepsilon}(\rho^2)}^2
				\right].
			\end{aligned}
		\end{equation}
		By $C^\varepsilon=A^\varepsilon+B^\varepsilon$,
		Propositions~\ref{prop:higgs-terminal-shift-summary},
		\ref{lem:gauge-terminal-shift-moment-bound}, and
		\ref{lem:gauge-stochastic-object-uniform-moments},
		together with~\eqref{eq:standard-gauge-field-negative-holder-moments}
		and the uniform embedding
		$H^{1,\varepsilon}(\rho)\hookrightarrow
		\mathcal C^{-\alpha_1,\varepsilon}(\rho)$, we obtain
		\begin{equation}
			\label{eq:terminal-states-negative-holder-uniform-bound}
			\mathbb E\left[
			\left\lVert C^\varepsilon
			\right\rVert_{\mathcal C^{-\alpha_1,\varepsilon}(\rho)}^2
			+\left\lVert\Phi^\varepsilon
			\right\rVert_{\mathcal C^{-\alpha_1,\varepsilon}(\rho^2)}^2
			\right]\lesssim1.
		\end{equation}
		Taking $R_\delta$ sufficiently large in
		\eqref{eq:extended-ymh-tightness-markov-bound} completes the proof.
	\end{proof}

	\section{Stochastic estimates}
	\label{sec:stochastic-estimate}
	Throughout this section, we use the abbreviations
	in~\eqref{eq:random-harmonic-field-shorthands} and do not adopt the parameter
	conventions of Section~\ref{subsec:parameters}.
	The main objectives of this section are to estimate the quantity in
	\eqref{eq:gauge-irregular-Zhi1-unified-free-static-renormalized-target} and
$\bigl\lVert
\mathbb{E}^\phi[\operatorname{Im}[\overline{\phi_B^\varepsilon}\,
D_B^\varepsilon\phi_B^\varepsilon]]
\bigr\rVert_{H^{-1,\varepsilon}(\rho)}$,
	and to prove the following proposition.

	\begin{proposition}
		\label{lem:gauge-stochastic-object-uniform-moments}
		For $0<\kappa\ll1$ and $p\geq1$, we have
		\begin{equation}
			\label{eq:gauge-stochastic-object-uniform-moments}
			\begin{aligned}
			\mathbb{E}\left[
			\left\lVert
			\Theta^\varepsilon
			\right\rVert_\kappa^p
			\right]
			&\lesssim_{\kappa,p}
			\mathbb{E}\left[
			\left(
			1+\left\lVert
			B^\varepsilon
			\right\rVert_{H^{1,\varepsilon}(\rho)}
			\right)^{\kappa p}
			\right].
			\end{aligned}
		\end{equation}
	\end{proposition}

	\begin{proof}
		By~\eqref{eq:stochastic-object-Xi-norm},
		Lemma~\ref{lem:covariant-free-gff-difference-negative-holder-moment},
		Corollaries~\ref{cor:standard-gauge-field-stochastic-estimates},
		\ref{cor:covariant-gff-wick-power-negative-holder-moment}, and
		\ref{cor:nonlinear-gauge-coefficients-negative-holder-moments},
		and Lemmas~\ref{lem:nonlinear-divergence-negative-holder-moment}--%
		\ref{cor:mixed-gauge-higgs-product-negative-holder-moments}
		termwise,
		taking the outer expectation in the conditional estimates and applying
		Minkowski's inequality yields the claimed estimate.
	\end{proof}

	We establish a regularity criterion, which is adapted
	from~\cite[Proposition~3]{MWX17}.

	\begin{lemma}[Regularity criterion]
		\label{lem:lattice-block-regularity-criterion}
		\label{lem:lattice-stochastic-regularity-criterion}
		Let $\alpha\in\mathbb R$, $\kappa>0$, $m\in\mathbb N$, and $0\leq b<a$.
		Let $\tau\colon\Omega\to\mathbb C^{\Lambda_\varepsilon}$ be a random
		lattice distribution taking values in the sum of Wiener chaoses of
		order at most $m$. Suppose that, for some $K\geq0$ and every $j\geq-1$,
		\begin{equation}
			\label{eq:lattice-block-pointwise-second-moment}
			\sup_{x\in\Lambda_\varepsilon}\rho(x)^b
			\left(\mathbb E\left[\lvert\Delta_j^\varepsilon\tau(x)\rvert^2\right]\right)^{1/2}
			\leq K2^{j\alpha}.
		\end{equation}
		Then, for every $1\leq p<\infty$,
		\begin{equation}
			\label{eq:lattice-block-regularity-estimate}
			\left(\mathbb E\left[
			\lVert\tau\rVert_{\mathcal C^{-\alpha-\kappa,\varepsilon}(\rho^a)}^p
			\right]\right)^{1/p}
			\lesssim_{\kappa,m,a,b}Kp^{m/2}.
		\end{equation}
	\end{lemma}

	\begin{proof}
		For sufficiently large $p$, weighted Besov embedding gives
		$B_{p,p}^{-\alpha-\kappa/2,\varepsilon}(\rho^a)
		\hookrightarrow\mathcal C^{-\alpha-\kappa,\varepsilon}(\rho^a)$,
		and $4(a-b)p>2$. By hypercontractivity,
		\begin{equation}
			\label{eq:lattice-block-weighted-besov-moment}
			\begin{aligned}
				\mathbb E\left[
				\lVert\tau\rVert_{B_{p,p}^{-\alpha-\kappa/2,\varepsilon}(\rho^a)}^p\right]
				&=\sum_{j\geq-1}2^{-j(\alpha+\kappa/2)p}
				\int_{\Lambda_\varepsilon}
				\mathbb E\left[\lvert\Delta_j^\varepsilon\tau(x)\rvert^p\right]
				\rho(x)^{ap}\,dx\\
				&\leq(C_mKp^{m/2})^p
				\sum_{j\geq-1}2^{-j\kappa p/2}
				\int_{\Lambda_\varepsilon}\rho(x)^{(a-b)p}\,dx.
			\end{aligned}
		\end{equation}
		The dyadic sum and the integral have uniformly bounded $p$-th roots
		for sufficiently large $p$, independently of $\varepsilon$.
		This yields~\eqref{eq:lattice-block-regularity-estimate} for large $p$;
		the remaining values follow by moment monotonicity.
	\end{proof}

	\subsection{Estimates of free stochastic objects}
	\label{subsec:stochastic-estimates-free-fields}
	\label{subsec:stochastic-estimates-gauge-field}

	We continue to suppress the superscript $M$ in the notation.
	Recall from~\eqref{eq:brownian-realizations-of-free-fields} that
	\begin{equation}
		\label{eq:free-gff-fourier-expansion}
		\phi^\varepsilon(x)
		=\frac{1}{\lvert\Lambda_{M,\varepsilon}\rvert}
		\sum_{k\in\widehat\Lambda_{M,\varepsilon}}
		\bigl(1+\lvert k^\varepsilon\rvert^2\bigr)^{-1/2}
		\widehat{\mathbb B_1^\varepsilon}(k)\,\mathrm e_k(x).
	\end{equation}
 For deterministic
	$f,g\in\mathbb C^{\Lambda_{M,\varepsilon}}$, we have
	$\mathbb E^\phi[\langle\mathbb B_1^\varepsilon,f\rangle_{M,\varepsilon}
	\langle\mathbb B_1^\varepsilon,g\rangle_{M,\varepsilon}]
	=\langle f,g\rangle_{M,\varepsilon}$, and
	\begin{equation}
		\label{eq:free-gff-complex-noise-covariance}
		\begin{aligned}
			\mathbb E^\phi\left[
			\langle\mathbb B_1^\varepsilon,f\rangle_{\mathbb C,M,\varepsilon}
			\langle\mathbb B_1^\varepsilon,g\rangle_{\mathbb C,M,\varepsilon}\right]
			=0,\quad
			\mathbb E^\phi\left[
			\overline{\langle\mathbb B_1^\varepsilon,f\rangle_{\mathbb C,M,\varepsilon}}
			\langle\mathbb B_1^\varepsilon,g\rangle_{\mathbb C,M,\varepsilon}\right]
			=2\langle f,g\rangle_{\mathbb C,M,\varepsilon}.
		\end{aligned}
	\end{equation}
	Consequently, for $k,k'\in\widehat\Lambda_{M,\varepsilon}$,
	\begin{equation}
		\label{eq:free-gff-fourier-covariance}
		\begin{aligned}
			\mathbb E^\phi\left[
			\overline{\widehat{\mathbb B_1^\varepsilon}(k')}
			\widehat{\mathbb B_1^\varepsilon}(k)\right]
			=2\lvert\Lambda_{M,\varepsilon}\rvert\mathbf1_{\{k=k'\}},\quad
			\mathbb E^\phi\left[
			\overline{\widehat{\phi^\varepsilon}(k')}
			\widehat{\phi^\varepsilon}(k)\right]
			=\frac{2\lvert\Lambda_{M,\varepsilon}\rvert\mathbf1_{\{k=k'\}}}
			{1+\lvert k^\varepsilon\rvert^2}.
		\end{aligned}
	\end{equation}
Equation~\eqref{eq:pure-yang-mills-law} gives, for $h_0\in[0,M^{-1})^2$,
	\begin{equation}
		\label{eq:pure-yang-mills-fourier-expansion}
		A^{\varepsilon,h_0}(x)
		=\frac{1}{\lvert\Lambda_{M,\varepsilon}\rvert}
		\sum_{k\in\widehat\Lambda_{M,\varepsilon}\setminus\{0\}}
		\frac{i(k^\varepsilon)^\perp}{\lvert k^\varepsilon\rvert^2}
		\widehat{W_1^\varepsilon}(k)
		\bigl(\mathrm e_k(x)-\widehat{\chi_\varepsilon}(-k)\bigr)+h_0.
	\end{equation}
	Since $W_1^\varepsilon$ is a real mean-zero Gaussian field,
	$\widehat{W_1^\varepsilon}(0)=0$ and, for nonzero
	$k,k'\in\widehat\Lambda_{M,\varepsilon}$,
	\begin{equation}
		\label{eq:pure-yang-mills-driving-fourier-covariance}
		\mathbb E^A\left[
		\widehat{W_1^\varepsilon}(k)
		\widehat{W_1^\varepsilon}(k')\right]
		=\lvert\Lambda_{M,\varepsilon}\rvert\mathbf1_{\{k+k'=0\}}.
	\end{equation}
	We write $:\!(A_1^{\varepsilon,0})^{m_1}
	(A_2^{\varepsilon,0})^{m_2}\!:_\chi$ for the pointwise Wick product
	with respect to
	$\bigl(\mathbb E^A[A_i^{\varepsilon,0}(x)A_j^{\varepsilon,0}(x)]\bigr)_{i,j=1}^2$.

	\begin{lemma}
		\label{lem:standard-gauge-field-stochastic-estimates}
		For $m_1,m_2,k_1,k_2\in\mathbb N_0$, $\kappa>0$, and $1\leq p<\infty$,
		\begin{equation}
			\label{eq:pinned-gauge-higgs-wick-product-moment}
			\begin{aligned}
				\left(\mathbb E\left[
				\left\lVert
				:\!(A_1^{\varepsilon,0})^{m_1}(A_2^{\varepsilon,0})^{m_2}\!:_\chi
				:\!(\phi^\varepsilon)^{k_1}
				(\overline{\phi^\varepsilon})^{k_2}\!:_\varepsilon
				\right\rVert_{\mathcal C^{-\kappa,\varepsilon}(\rho^\kappa)}^p
				\right]\right)^{1/p}\lesssim_{\chi,m_1,m_2,k_1,k_2,\kappa}
				p^{(m_1+m_2+k_1+k_2)/2}.
			\end{aligned}
		\end{equation}
	\end{lemma}

	\begin{proof}
		Put $d:=m_1+m_2+k_1+k_2$ and
		\begin{equation}
			\label{eq:pinned-gauge-higgs-wick-product-F}
			F:=
			:\!(A_1^{\varepsilon,0})^{m_1}(A_2^{\varepsilon,0})^{m_2}\!:_\chi
			:\!(\phi^\varepsilon)^{k_1}
			(\overline{\phi^\varepsilon})^{k_2}\!:_\varepsilon.
		\end{equation}
		The case $d=0$ is immediate. For $d\geq1$, Lemmas~\ref{lem:pinned-wick-two-point-bound}
		and~\ref{lem:dyadic-logarithmic-kernel-bound} give, for $-1\leq j\leq J^\varepsilon$,
		\begin{equation}
			\label{eq:pinned-wick-dyadic-second-moment}
			\begin{aligned}
				\mathbb E|\Delta_j^\varepsilon F(x)|^2
				&=\iint_{\Lambda_\varepsilon\times\Lambda_\varepsilon}
				\eta_j^\varepsilon(x-y)\overline{\eta_j^\varepsilon(x-z)}
				\mathbb E[F(y)\overline{F(z)}]\,dy\,dz\\
				&\lesssim_{\chi,d}
				\iint_{\Lambda_\varepsilon\times\Lambda_\varepsilon}
				|\eta_j^\varepsilon(x-y)\eta_j^\varepsilon(x-z)|L(y,z)^d\,dy\,dz\\
				&\lesssim_{\chi,d}(2+j+\log(2+|x|))^d
				\lesssim_{\chi,d,\kappa}\rho(x)^{-\kappa/2}2^{j\kappa/2}.
			\end{aligned}
		\end{equation}
		The independent Gaussian fields $A^{\varepsilon,0}$ and $\phi^\varepsilon$
		make $F$ an element of the $d$-th homogeneous Wiener chaos.
		Apply Lemma~\ref{lem:lattice-stochastic-regularity-criterion} with
		$\alpha=b=\kappa/4$, $a=\kappa$, $m=d$, and regularity-loss parameter
		$3\kappa/4$ to obtain~\eqref{eq:pinned-gauge-higgs-wick-product-moment}.
	\end{proof}

	\begin{lemma}[Two-point bound]
		\label{lem:pinned-wick-two-point-bound}
		Write $\log_+r:=\log(r\vee1)$.
		For $x,y\in\Lambda_\varepsilon$, set
		\begin{equation}
			\label{eq:pinned-logarithmic-covariance-envelope}
			L(x,y):=
			\log_+\frac1{\varepsilon\vee d_M(x,y)}
			+\log(2+|x|)+\log(2+|y|).
		\end{equation}
		For $F$ defined in~\eqref{eq:pinned-gauge-higgs-wick-product-F},
		\begin{equation}
			\label{eq:pinned-wick-two-point-bound}
			|\mathbb E[F(x)\overline{F(y)}]|
			\lesssim_{\chi,d}L(x,y)^d,
			\qquad d=m_1+m_2+k_1+k_2.
		\end{equation}
	\end{lemma}

	\begin{proof}
		Set $a_x(k):=\mathrm e_k(x)-\widehat{\chi_\varepsilon}(-k)$.
		By~\eqref{eq:pure-yang-mills-fourier-expansion} and
		\eqref{eq:pure-yang-mills-driving-fourier-covariance},
		\begin{equation}
			\label{eq:pinned-gauge-covariance-Fourier}
			\begin{aligned}
				C_{ij}^\varepsilon(x,y)
				:=
				\mathbb E^A[A_i^{\varepsilon,0}(x)A_j^{\varepsilon,0}(y)]
				=\frac1{|\Lambda_{M,\varepsilon}|}
				\sum_{k\in Q_{M,\varepsilon}\setminus\{0\}}
				\frac{(k^\varepsilon)^\perp_i
				\overline{(k^\varepsilon)^\perp_j}}{|k^\varepsilon|^4}
				a_x(k)\overline{a_y(k)}.
			\end{aligned}
		\end{equation}
		Insert $\eta+\eta^c=1$ to write
		$C_{ij}^\varepsilon=C_{ij}^{\varepsilon,\mathrm{lo}}+
		C_{ij}^{\varepsilon,\mathrm{hi}}$.

		For the low frequencies, normalization and compact support of
		$\chi_\varepsilon$ give
		\begin{equation}
			\label{eq:pinned-Fourier-mode-cancellation}
			|a_x(k)|
			=\left|\int_{\Lambda_{M,\varepsilon}}
			\chi_\varepsilon(z)(\mathrm e_k(x)-\mathrm e_k(z))\,dz\right|
			\lesssim_\chi\min\{1,|k|(1+|x|)\}.
		\end{equation}
		Comparing nonzero frequency cells with the corresponding integral,
		as in~\cite[Lemma~A.2]{MP19}, yields
		\begin{equation}
			\label{eq:pinned-gauge-low-frequency-variance}
			\begin{aligned}
				C_{ii}^{\varepsilon,\mathrm{lo}}(x,x)
				&\lesssim_\chi\frac1{|\Lambda_{M,\varepsilon}|}
				\sum_{\substack{k\in Q_{M,\varepsilon}\\0<|k|\leq3}}
				\frac{\min\{1,|k|^2(1+|x|^2)\}}{|k|^2}\\
				&\lesssim\int_0^4\frac{\min\{1,r^2(1+|x|^2)\}}r\,dr
				\lesssim\log(2+|x|).
			\end{aligned}
		\end{equation}
		Cauchy--Schwarz therefore gives
		\begin{equation}
			\label{eq:pinned-gauge-low-frequency-covariance}
			|C_{ij}^{\varepsilon,\mathrm{lo}}(x,y)|
			\lesssim_\chi\log(2+|x|)+\log(2+|y|).
		\end{equation}

		For the high frequencies, put
		\begin{equation}
			\label{eq:pinned-gauge-high-frequency-kernel}
			K_{ij}^{\mathrm{hi}}(z):=
			\frac1{|\Lambda_{M,\varepsilon}|}
			\sum_{k\in Q_{M,\varepsilon}\setminus\{0\}}
			\eta^c(k)\frac{(k^\varepsilon)^\perp_i
			\overline{(k^\varepsilon)^\perp_j}}{|k^\varepsilon|^4}
			\mathrm e_k(z).
		\end{equation}
		The choice~\eqref{eq:gauge-infrared-cutoff} gives
		$\eta^c\varphi_{-1}^\varepsilon=0$. Thus
		$K_{ij}^{\mathrm{hi}}=\sum_{\ell=0}^{J^\varepsilon}K_{ij,\ell}^{\mathrm{hi}}$,
		where $K_{ij,\ell}^{\mathrm{hi}}=\Delta_\ell^\varepsilon K_{ij}^{\mathrm{hi}}$
		has multiplier
		\begin{equation}
			\label{eq:pinned-gauge-dyadic-multiplier}
			m_{ij,\ell}(\xi):=
			\varphi_\ell^\varepsilon(\xi)\eta^c(\xi)
			\frac{(\xi^\varepsilon)^\perp_i
			\overline{(\xi^\varepsilon)^\perp_j}}{|\xi^\varepsilon|^4},
			\qquad
			|\partial_\xi^\gamma m_{ij,\ell}(\xi)|
			\lesssim_\gamma2^{-(2+|\gamma|)\ell}.
		\end{equation}
		The multiplier is smooth and periodic, also for $\ell=J^\varepsilon$.
		Applying Lemma~\ref{lem:periodic-multiplier-kernel-bound}
	 gives, for every integer $N\geq1$,
		\begin{equation}
			\label{eq:pinned-gauge-dyadic-kernel-decay}
			|K_{ij,\ell}^{\mathrm{hi}}(z)|
			\lesssim_N(1+2^\ell d_M(z,0))^{-N}.
		\end{equation}
		Using $2^{J^\varepsilon}\asymp\varepsilon^{-1}$ and summing over $\ell$,
		\begin{equation}
			\label{eq:pinned-gauge-high-frequency-log-kernel}
			|K_{ij}^{\mathrm{hi}}(z)|
			\lesssim\sum_{\ell=0}^{J^\varepsilon}(1+2^\ell d_M(z,0))^{-N}
			\lesssim1+\log_+\frac1{\varepsilon\vee d_M(z,0)}.
		\end{equation}
		All integrals against $\chi_\varepsilon$ in the next identity are
		taken over one period:
		\begin{equation}
			\label{eq:pinned-gauge-high-frequency-covariance}
			\begin{aligned}
				C_{ij}^{\varepsilon,\mathrm{hi}}(x,y)
				&=K_{ij}^{\mathrm{hi}}(x-y)
				-\int_{\Lambda_{M,\varepsilon}}\chi_\varepsilon(u)K_{ij}^{\mathrm{hi}}(x-u)\,du\\
				&\quad-\int_{\Lambda_{M,\varepsilon}}\chi_\varepsilon(u)K_{ij}^{\mathrm{hi}}(u-y)\,du
				+\iint_{\Lambda_{M,\varepsilon}\times\Lambda_{M,\varepsilon}}\chi_\varepsilon(u)\chi_\varepsilon(v)
				K_{ij}^{\mathrm{hi}}(u-v)\,du\,dv.
			\end{aligned}
		\end{equation}
		A calculation similar to that in the proof of
		Lemma~\ref{lem:dyadic-logarithmic-kernel-bound} shows that the last three
		terms are uniformly bounded. Hence
		\begin{equation}
			\label{eq:pinned-gauge-covariance-log-bound}
			|C_{ij}^{\varepsilon,\mathrm{hi}}(x,y)|
			\lesssim_\chi1+\log_+\frac1{\varepsilon\vee d_M(x,y)},
			\qquad |C_{ij}^\varepsilon(x,y)|\lesssim_\chi L(x,y).
		\end{equation}
		The same dyadic argument, with multiplier
		$2(1+|k^\varepsilon|^2)^{-1}$ and a bounded low-frequency block, gives
		\begin{equation}
			\label{eq:free-gff-periodic-log-covariance-bound}
			|\mathbb E^\phi[\phi^\varepsilon(x)\overline{\phi^\varepsilon(y)}]|
			\lesssim1+\log_+\frac1{\varepsilon\vee d_M(x,y)}.
		\end{equation}
		Finally, Wick's theorem gives
		\begin{equation}
			\label{eq:pinned-wick-cross-contraction-bound}
			|\mathbb E[F(x)\overline{F(y)}]|
			\lesssim_d\left(\max_{i,j}|C_{ij}^\varepsilon(x,y)|
			+|\mathbb E^\phi[\phi^\varepsilon(x)\overline{\phi^\varepsilon(y)}]|\right)^d
			\lesssim_{\chi,d}L(x,y)^d.
		\end{equation}
	\end{proof}

	\begin{lemma}[Dyadic logarithmic kernel bound]
		\label{lem:dyadic-logarithmic-kernel-bound}
		For every integer $d\geq1$, $-1\leq j\leq J^\varepsilon$, and
		$x\in\Lambda_\varepsilon$,
		\begin{equation}
			\label{eq:dyadic-logarithmic-kernel-bound}
			\iint_{\Lambda_\varepsilon\times\Lambda_\varepsilon}
			|\eta_j^\varepsilon(x-y)\eta_j^\varepsilon(x-z)|L(y,z)^d\,dy\,dz
			\lesssim_d(2+j+\log(2+|x|))^d.
		\end{equation}
	\end{lemma}

	\begin{proof}
		Fix an integer $N>2$. By~\eqref{eq:app-LP-kernel-Lq-bound} and
		integration by parts for the smooth dyadic multipliers,
		\begin{equation}
			\label{eq:dyadic-kernel-pointwise-decay}
			\|\eta_j^\varepsilon\|_{L^{q,\varepsilon}}
			\lesssim2^{2j(1-1/q)},\quad 1\leq q\leq\infty,
			\qquad
			|\eta_j^\varepsilon(u)|
			\lesssim_N2^{2j}(1+2^j|u|)^{-N}.
		\end{equation}
		We estimate the three terms in $L(y,z)$ separately.

		\textup{(1)} Since
		$2+|y|\leq(2+|x|)(1+|x-y|)$, integrating first in $z$ gives
		\begin{equation}
			\label{eq:dyadic-spatial-logarithm-bound}
			\begin{aligned}
				\iint_{\Lambda_\varepsilon\times\Lambda_\varepsilon}
				|\eta_j^\varepsilon(x-y)\eta_j^\varepsilon(x-z)|
				\log^d(2+|y|)\,dy\,dz
				&\lesssim\int_{\Lambda_\varepsilon}
				|\eta_j^\varepsilon(x-y)|\log^d(2+|y|)\,dy\\
				&\lesssim_d\int_{\Lambda_\varepsilon}
				|\eta_j^\varepsilon(x-y)|\log^d(1+|x-y|)\,dy
				+\log^d(2+|x|)\\
				&\lesssim_d\log^d(2+|x|).
			\end{aligned}
		\end{equation}
		The last step uses the decay in~\eqref{eq:dyadic-kernel-pointwise-decay}
		and $2^j\geq1/2$. The term $\log^d(2+|z|)$ has the same bound.

		\textup{(2)} Set
		\begin{equation}
			\label{eq:dyadic-absolute-kernel-convolution}
			h_j^\varepsilon(u):=
			(|\eta_j^\varepsilon|*|\eta_j^\varepsilon|)(u)
			=\int_{\Lambda_\varepsilon}
			|\eta_j^\varepsilon(u-v)\eta_j^\varepsilon(v)|\,dv.
		\end{equation}
		Every $v$ satisfies $|v|\geq|u|/2$ or $|u-v|\geq|u|/2$.
	Using \eqref{eq:dyadic-kernel-pointwise-decay}, we have
		\begin{equation}
			\label{eq:dyadic-absolute-convolution-decay}
			\begin{aligned}
				h_j^\varepsilon(u)
				&\leq\int_{\substack{v\in\Lambda_\varepsilon\\|v|\geq|u|/2}}
				|\eta_j^\varepsilon(u-v)\eta_j^\varepsilon(v)|\,dv
				+\int_{\substack{v\in\Lambda_\varepsilon\\|u-v|\geq|u|/2}}
				|\eta_j^\varepsilon(u-v)\eta_j^\varepsilon(v)|\,dv\\
				&\lesssim_N2^{2j}(1+2^j|u|)^{-N}
				\|\eta_j^\varepsilon\|_{L^{1,\varepsilon}}
				\lesssim_N2^{2j}(1+2^j|u|)^{-N}.
			\end{aligned}
		\end{equation}
		For $|u|\leq1$ and $n\ne0$,
		$|u+2\pi Mn|\geq2\pi M|n|-1\gtrsim M|n|$. Since $j\geq-1$ and $M\geq1$,
		we have $1+2^jM|n|\gtrsim(1+2^j)|n|$, and hence
		\begin{equation}
			\label{eq:dyadic-periodic-logarithmic-images}
			\begin{aligned}
				\sum_{n\in\mathbb Z^2\setminus\{0\}}
				h_j^\varepsilon(u+2\pi Mn)
				&\lesssim_N\sum_{n\in\mathbb Z^2\setminus\{0\}}
				2^{2j}(1+2^jM|n|)^{-N}\\
				&\lesssim_N2^{2j}(1+2^j)^{-N}
				\sum_{n\in\mathbb Z^2\setminus\{0\}}|n|^{-N}
				\lesssim_N2^{2j}(1+2^j|u|)^{-N}.
			\end{aligned}
		\end{equation}
		The series converges because $N>2$. Adding the term $n=0$ gives
		\begin{equation}
			\label{eq:dyadic-periodized-absolute-kernel}
			\sum_{n\in\mathbb Z^2}h_j^\varepsilon(u+2\pi Mn)
			\lesssim_N2^{2j}(1+2^j|u|)^{-N},
			\qquad |u|\leq1.
		\end{equation}

		Now let $\widetilde z=z-y$, so that $x-z=x-y-\widetilde z$.
		Since $\eta_j^\varepsilon$ is even and
		$\log_+\frac1{\varepsilon\vee|u|}=0$ for $|u|\geq1$,
		this change of variables yields
		\begin{equation}
			\label{eq:dyadic-logarithmic-convolution-reduction}
			\begin{aligned}
				\iint_{\Lambda_\varepsilon\times\Lambda_\varepsilon}
				|\eta_j^\varepsilon(x-y)\eta_j^\varepsilon(x-z)|
				\log_+^d\frac1{\varepsilon\vee d_M(y,z)}\,dy\,dz\
				&=\int_{\Lambda_\varepsilon}h_j^\varepsilon(\widetilde z)
				\log_+^d\frac1{\varepsilon\vee d_M(\widetilde z,0)}
				\,d\widetilde z\\
				&\overset{\eqref{eq:dyadic-periodized-absolute-kernel}}{=}\int_{\substack{u\in\Lambda_\varepsilon\\|u|\leq1}}
				\left(\sum_{n\in\mathbb Z^2}h_j^\varepsilon(u+2\pi Mn)\right)
				\log_+^d\frac1{\varepsilon\vee|u|}\,du\\
				&\lesssim_N\int_{\substack{u\in\Lambda_\varepsilon\\|u|\leq1}}
				2^{2j}(1+2^j|u|)^{-N}
				\left(1+\log_+\frac1{\varepsilon\vee|u|}\right)^d\,du.
			\end{aligned}
		\end{equation}

		For $u\in\Lambda_\varepsilon$ and
		$\xi\in u+[-\varepsilon/2,\varepsilon/2]^2$, the bound
		$2^j\varepsilon\lesssim1$ gives
$
			\label{eq:dyadic-logarithmic-cell-comparison}
			\varepsilon\vee|\xi|\asymp\varepsilon\vee|u|,
			\qquad
			1+2^j|\xi|\asymp1+2^j|u|.
$
		Hence
$
			\label{eq:dyadic-logarithmic-cell-weight}
			1+\log_+\frac1{\varepsilon\vee|\xi|}
			\asymp1+\log_+\frac1{\varepsilon\vee|u|}.
$
		By the preceding pointwise comparisons and~\cite[Lemma~A.2]{MP19}, we obtain
		\begin{equation}
			\label{eq:dyadic-local-logarithmic-singularity}
			\int_{\substack{u\in\Lambda_\varepsilon\\|u|\leq1}}
			\frac{2^{2j}\left(1+\log_+\frac1{\varepsilon\vee|u|}\right)^d}
			{(1+2^j|u|)^N}\,du
			\lesssim_{N,d}
			\int_{\mathbb R^2}
			\frac{2^{2j}\left(1+\log_+\frac1{|\xi|}\right)^d}
			{(1+2^j|\xi|)^N}\,d\xi.
		\end{equation}
		Since
		$\log_+(1/|\xi|)\leq(2+j)\log2+\log_+(1/(2^j|\xi|))$
		for $j\geq-1$, the change of variables $v=2^j\xi$ bounds the last integral
		by a constant multiple of
		\begin{equation}
			\label{eq:dyadic-continuous-logarithmic-bound}
			\int_{\mathbb R^2}
			\frac{(2+j)^d+\log_+^d(1/|v|)}{(1+|v|)^N}\,dv
			\lesssim_{N,d}(2+j)^d,
			\qquad N>2.
		\end{equation}
		Combining \textup{(1)} and \textup{(2)} gives the bound
		$(2+j)^d+\log^d(2+|x|)$.
	\end{proof}

	\begin{lemma}
		\label{lem:pinned-gauge-diagonal-covariance-difference}
		For $i,j=1,2$ and $x\in\Lambda_\varepsilon$,
		\begin{equation}
			\label{eq:pinned-gauge-diagonal-covariance-difference}
			\left|C_{ii}^\varepsilon(x,x)-C_{ii}^\varepsilon(0,0)\right|
			\lesssim_\chi\log(2+|x|),
		\end{equation}
		and
		\begin{equation}
			\label{eq:pinned-gauge-diagonal-covariance-first-difference}
			\left|D_j^\varepsilon[C_{ii}^\varepsilon(\cdot,\cdot)](x)\right|
			\lesssim_\chi1.
		\end{equation}
		Here the difference acts on the diagonal map
		$x\mapsto C_{ii}^\varepsilon(x,x)$.
		Consequently, for every $\kappa>0$,
		\begin{equation}
			\label{eq:pinned-gauge-diagonal-covariance-weighted-C1}
			\left\lVert C_{ii}^\varepsilon(\cdot,\cdot)-C_{ii}^\varepsilon(0,0)
			\right\rVert_{\mathcal C^{1,\varepsilon}(\rho^{\kappa/2})}
			\lesssim_{\chi,\kappa}1.
		\end{equation}
	\end{lemma}

	\begin{proof}
		The low-frequency contribution to
		\eqref{eq:pinned-gauge-diagonal-covariance-difference} is bounded
		by~\eqref{eq:pinned-gauge-low-frequency-variance}.
		Recall $a_x(k)=\mathrm e_k(x)-\widehat{\chi_\varepsilon}(-k)$.
		Since $\chi_\varepsilon$ is real,
		\begin{equation}
			\label{eq:pinned-gauge-diagonal-mode-cancellation}
			\begin{aligned}
				|a_x(k)|^2-|a_y(k)|^2
				=-(\mathrm e_k(x)-\mathrm e_k(y))\widehat{\chi_\varepsilon}(k)
				-\overline{(\mathrm e_k(x)-\mathrm e_k(y))}
				\widehat{\chi_\varepsilon}(-k).
			\end{aligned}
		\end{equation}
		Taking $y=0$ bounds this difference by
		$4|\widehat{\chi_\varepsilon}(k)|$.
		Thus~\eqref{eq:pinned-gauge-covariance-Fourier}, the rapid decay of
		$\widehat{\chi_\varepsilon}$, and the sum--integral comparison
		in~\cite[Lemma~A.2]{MP19} give, for $N>0$,
		\begin{equation}
			\label{eq:pinned-gauge-high-frequency-diagonal-difference}
			\begin{aligned}
				\left|C_{ii}^{\varepsilon,\mathrm{hi}}(x,x)
				-C_{ii}^{\varepsilon,\mathrm{hi}}(0,0)\right|
				&\lesssim\frac1{|\Lambda_{M,\varepsilon}|}
				\sum_{k\in Q_{M,\varepsilon}\setminus\{0\}}
				\frac{\eta^c(k)}{|k^\varepsilon|^2}
				|\widehat{\chi_\varepsilon}(-k)|\\
				&\lesssim_{\chi,N}\frac1{|\Lambda_{M,\varepsilon}|}
				\sum_{\substack{k\in Q_{M,\varepsilon}\\|k|\geq1}}
				\frac{(1+|k|)^{-N}}{|k|^2}\lesssim_{\chi,N}1.
			\end{aligned}
		\end{equation}
		This yields~\eqref{eq:pinned-gauge-diagonal-covariance-difference}.
		Taking $y=x+\varepsilon e_j$ in
		\eqref{eq:pinned-gauge-diagonal-mode-cancellation} also gives
		\begin{equation}
			\label{eq:pinned-gauge-diagonal-mode-first-difference}
			\left|D_j^\varepsilon[|a_\cdot(k)|^2](x)\right|
			\lesssim |k|\,|\widehat{\chi_\varepsilon}(k)|.
		\end{equation}
		Hence the same sum--integral comparison and polar coordinates in
		two dimensions yield, for $N>1$,
		\begin{equation}
			\label{eq:pinned-gauge-diagonal-covariance-first-difference-proof}
			\begin{aligned}
				\left|D_j^\varepsilon[C_{ii}^\varepsilon(\cdot,\cdot)](x)\right|
				\lesssim\frac1{|\Lambda_{M,\varepsilon}|}
				\sum_{k\in Q_{M,\varepsilon}\setminus\{0\}}
				\frac{|\widehat{\chi_\varepsilon}(k)|}{|k|}
				\lesssim_{\chi,N}\int_0^\infty(1+r)^{-N}\,dr
				\lesssim_{\chi,N}1.
			\end{aligned}
		\end{equation}
	\end{proof}

	\begin{corollary}
		\label{cor:standard-gauge-field-stochastic-estimates}
		\label{lem:standard-stochastic-estimates}
		For every $k_1,k_2\in\mathbb N_0$, $\kappa>0$, and
		$p\in[1,\infty)$, we have
		\begin{equation}
			\label{eq:free-gff-wick-power-negative-holder-moment}
			\left(
			\mathbb{E}^\phi\left[
			\left\lVert
			:\!(\phi^\varepsilon)^{k_1}
			(\overline{\phi^\varepsilon})^{k_2}\!:_\varepsilon
			\right\rVert_{\mathcal{C}^{-\kappa,\varepsilon}(\rho^\kappa)}^p
			\right]
			\right)^{1/p}
			\lesssim_{\kappa,k_1,k_2}
			p^{(k_1+k_2)/2}.
		\end{equation}
		\begin{equation}
			\label{eq:standard-gauge-field-negative-holder-moments}
				\left(
				\mathbb{E}\left[
				\lVert A^\varepsilon\rVert_{\mathcal{C}^{-\kappa,\varepsilon}(\rho^\kappa)}^p
				\right]
				\right)^{1/p}
				\lesssim_\kappa p^{1/2},
				\qquad
				\left(
				\mathbb{E}\left[
				\left\lVert :\!\lvert A^\varepsilon\rvert^2\!:_\varepsilon
				\right\rVert_{\mathcal{C}^{-\kappa,\varepsilon}(\rho^\kappa)}^p
				\right]
				\right)^{1/p}
				\lesssim_\kappa p.
		\end{equation}
		\begin{equation}
			\label{eq:free-mixed-gauge-higgs-product-negative-holder-moments}
			\left(
			\mathbb{E}\left[
			\left\lVert
			A^\varepsilon\phi^\varepsilon
			\right\rVert_{\mathcal{C}^{-\kappa,\varepsilon}(\rho^\kappa)}^p
			\right]
			\right)^{1/p}
			\lesssim_\kappa p,
			\qquad
			\left(
			\mathbb{E}\left[
			\left\lVert
			:\!\left\lvert A^\varepsilon\right\rvert^2\!:_\varepsilon
			\phi^\varepsilon
			\right\rVert_{\mathcal{C}^{-\kappa,\varepsilon}(\rho^\kappa)}^p
			\right]
			\right)^{1/p}
			\lesssim_\kappa p^{3/2}.
		\end{equation}
		Here $:\!\lvert A^\varepsilon\rvert^2\!:_\varepsilon
		=\lvert A^\varepsilon\rvert^2-a_{M, \varepsilon}$, with
		$a_{M, \varepsilon}$ defined in~\eqref{eq:lattice-renormalization-constants}.
		The same estimates hold with $\phi^\varepsilon$ replaced by
		$\tau_j\phi^\varepsilon$, $j=1,2$.
	\end{corollary}

	\begin{proof}
		The pure Higgs estimate is the case $m_1=m_2=0$ of
		Lemma~\ref{lem:standard-gauge-field-stochastic-estimates}.
		By~\eqref{eq:random-harmonic-terminal-state-lifts} and
		\eqref{eq:lattice-renormalization-constants},
		$A^\varepsilon=A^{\varepsilon,0}+H^\varepsilon$ with
		$|H^\varepsilon|\lesssim1$, and
		$a_{M,\varepsilon}=\sum_{i=1}^2C_{ii}^\varepsilon(0,0)$.
		Hence
		\begin{equation}
			\label{eq:local-to-origin-gauge-Wick-ordering}
			\begin{aligned}
				:\!|A^\varepsilon(x)|^2\!:_\varepsilon
				=\sum_{i=1}^2:\!(A_i^{\varepsilon,0}(x))^2\!:_\chi
				+2H^\varepsilon\cdot A^{\varepsilon,0}(x)+|H^\varepsilon|^2
				+\sum_{i=1}^2
				\bigl(C_{ii}^\varepsilon(x,x)-C_{ii}^\varepsilon(0,0)\bigr).
			\end{aligned}
		\end{equation}
		By Lemma~\ref{lem:pinned-gauge-diagonal-covariance-difference}
		the last term has a uniformly bounded
		$\mathcal C^{1,\varepsilon}(\rho^{\kappa/2})$-norm.
		For $0<\kappa<1$, Lemma~\ref{lem:discrete-paraproduct-estimates}
		bounds its product with $\phi^\varepsilon$ in
		$\mathcal C^{-\kappa,\varepsilon}(\rho^\kappa)$ by
		$\|\phi^\varepsilon\|_{\mathcal C^{-\kappa/2,\varepsilon}(\rho^{\kappa/2})}$.
		Applying Lemma~\ref{lem:standard-gauge-field-stochastic-estimates}
		with parameter $\kappa/2$ gives the remaining estimates.
	\end{proof}

	\begin{remark}[Infrared divergence and pinning]
		\label{rem:pure-yang-mills-infrared-divergence}
		Consider the unpinned pure Yang--Mills field
		\begin{equation}
			\label{eq:unpinned-pure-yang-mills-field}
			A^{\mathrm{unpin},M,\varepsilon}
			:=
			K^{M,\varepsilon}W_1^{M,\varepsilon}.
		\end{equation}
		By~\eqref{eq:pure-yang-mills-driving-fourier-covariance}, for
		$k\in Q_{M,\varepsilon}\setminus\{0\}$,
$
			\mathbb E^A\left[
			\left\lvert\widehat{A^{\mathrm{unpin},M,\varepsilon}}(k)\right\rvert^2
			\right]
			=\frac{\lvert\Lambda_{M,\varepsilon}\rvert}{\lvert k^\varepsilon\rvert^2}
			\asymp\frac{\lvert\Lambda_{M,\varepsilon}\rvert}{\lvert k\rvert^2}.
$
Thus for $M$ sufficiently large,
		\begin{equation}
			\label{eq:unpinned-pure-yang-mills-infrared-variance}
			\mathbb E^A\left[
			\left\lvert\Delta_{-1}^\varepsilon
			A^{\mathrm{unpin},M,\varepsilon}(x)\right\rvert^2\right]
			=\frac{1}{\lvert\Lambda_{M,\varepsilon}\rvert}
			\sum_{k\in Q_{M,\varepsilon}\setminus\{0\}}
			\frac{\lvert\varphi_{-1}^\varepsilon(k)\rvert^2}
			{\lvert k^\varepsilon\rvert^2}
			\asymp\log M.
		\end{equation}
		Since $\rho(0)=1$, Gaussian moments imply, for every $\kappa>0$ and
		$p\in[1,\infty)$,
		\begin{equation}
			\label{eq:unpinned-pure-yang-mills-weighted-moment-divergence}
			\mathbb E^A\left[
			\left\lVert A^{\mathrm{unpin},M,\varepsilon}
			\right\rVert_{\mathcal C^{-\kappa,\varepsilon}(\rho^\kappa)}^p
			\right]
			\gtrsim_{\kappa,p}(\log M)^{p/2}.
		\end{equation}
	\end{remark}

	\subsection{\texorpdfstring{Estimate of
			$\phi_B^\varepsilon$}
		{Estimate of the covariant Gaussian free field}}
	\label{subsec:stochastic-estimate-covariant-gff}

	\begin{lemma}
		\label{lem:covariant-gff-negative-holder-moment}
		For every $B\in(\mathbb R^2)^{\Lambda_{M,\varepsilon}}$,
		$\kappa>0$, and $p\in[1,\infty)$,
		\begin{equation}
			\label{eq:covariant-gff-negative-holder-moment}
			\left(
			\mathbb E^\phi\left[
			\|\phi_B^\varepsilon\|_{\mathcal C^{-\kappa,\varepsilon}(\rho^\kappa)}^p
			\right]\right)^{1/p}
			\lesssim_\kappa p^{1/2}.
		\end{equation}
	\end{lemma}

	\begin{proof}
		It suffices to consider $0<\kappa<1$.
		For periodic $f$, the kernel~\eqref{eq:weighted-periodic-LP-Fourier-sum} gives
		\begin{equation}
			\label{eq:lattice-dyadic-block-convolution}
			\begin{aligned}
				\Delta_j^\varepsilon f(x)
				=\int_{\Lambda_\varepsilon}
				\eta_j^\varepsilon(x-y)f(y)\,dy
				=\int_{\Lambda_{M,\varepsilon}}
				\eta_j^{M,\varepsilon}(x-y)f(y)\,dy.
			\end{aligned}
		\end{equation}
		Since $\phi_B^\varepsilon=(\mathcal L_B^\varepsilon+1)^{-1/2}
		\mathbb B_1^\varepsilon$, the Gaussian isometry and
		\eqref{eq:diamagnetic-resolvent} yield
		\begin{equation}
			\label{eq:covariant-gff-dyadic-block-second-moment}
			\begin{aligned}
				\mathbb E^\phi[|\Delta_j^\varepsilon\phi_B^\varepsilon(x)|^2]
				&=2\left\|(\mathcal L_B^\varepsilon+1)^{-1/2}
				\eta_j^{M,\varepsilon}(x-\cdot)
				\right\|_{L^{2,\varepsilon}(\Lambda_{M,\varepsilon})}^2\\
				&\leq2\left\|(\mathcal L_0^\varepsilon+1)^{-1/2}
				|\eta_j^{M,\varepsilon}|
				\right\|_{L^{2,\varepsilon}(\Lambda_{M,\varepsilon})}^2\\
				&=2\left\||\eta_j^{M,\varepsilon}|
				\right\|_{H^{-1,\varepsilon}(\Lambda_{M,\varepsilon})}^2.
			\end{aligned}
		\end{equation}
		Choose $q\in(1,2)$ such that $2(1-q^{-1})=\kappa/2$.
		The embedding
		$L^{q,\varepsilon}(\Lambda_{M,\varepsilon})
		\hookrightarrow H^{-1,\varepsilon}(\Lambda_{M,\varepsilon})$
		is uniform in $M\geq1$ and $\varepsilon$.
		Thus~\eqref{eq:weighted-periodic-LP-kernel-Lp-bound} with $a=0$ gives
		\begin{equation}
			\label{eq:lattice-dyadic-kernel-negative-Sobolev-bound}
			\left\||\eta_j^{M,\varepsilon}|
			\right\|_{H^{-1,\varepsilon}(\Lambda_{M,\varepsilon})}
			\lesssim_\kappa
			\|\eta_j^{M,\varepsilon}\|_{L^{q,\varepsilon}(\Lambda_{M,\varepsilon})}
			\lesssim_\kappa2^{2j(1-q^{-1})}=2^{j\kappa/2}.
		\end{equation}
		Combining these estimates and applying
		Lemma~\ref{lem:lattice-stochastic-regularity-criterion} with
		$\alpha=\kappa/2$, $a=\kappa$, $b=0$, $m=1$, and regularity loss
		$\kappa/2$ yields~\eqref{eq:covariant-gff-negative-holder-moment}
		for every $p\geq1$.
	\end{proof}

	\subsection{\texorpdfstring{Estimate of
			$\phi_B^\varepsilon-\phi^\varepsilon$}
		{Estimate of the difference of Gaussian free fields}}
	\label{subsec:stochastic-estimate-covariant-free-gff-difference}

	We prove the following weighted estimates.
	\begin{lemma}
		\label{lem:covariant-free-gff-difference-negative-holder-moment}
		For every $B\in(\mathbb R^2)^{\Lambda_{M,\varepsilon}}$,
		$\kappa>0$, and $p\in[1,\infty)$,
		\begin{equation}
			\label{eq:covariant-free-gff-difference-negative-holder-moment}
			\left(\mathbb E^\phi\left[
			\|\phi_B^\varepsilon-\phi^\varepsilon\|_{\mathcal C^{1-\kappa,\varepsilon}(\rho^{2+\kappa})}^p
			\right]\right)^{1/p}
			\lesssim_\kappa p^{1/2}\|B\|_{H^{1,\varepsilon}(\rho)}
			\bigl(1+\|B\|_{H^{1,\varepsilon}(\rho)}\bigr).
		\end{equation}
		Moreover,
		\begin{equation}
			\label{eq:covariant-free-gff-difference-interpolated-holder-moment}
			\left(\mathbb E^\phi\left[
			\left(
			\mathbf1_{\{\|B\|_{H^{1,\varepsilon}(\rho)}\ne0\}}
			\sup_{\sigma\in[-\kappa,1-\kappa]}
			\frac{\|\phi_B^\varepsilon-\phi^\varepsilon\|_{\mathcal C^{\sigma,\varepsilon}(\rho^{3\kappa+2\sigma})}}
			{\|B\|_{H^{1,\varepsilon}(\rho)}^{\kappa+\sigma}
			\bigl(1+\|B\|_{H^{1,\varepsilon}(\rho)}\bigr)^{\kappa+\sigma}}
			\right)^p\right]\right)^{1/p}
			\lesssim_\kappa p^{1/2}.
		\end{equation}
	\end{lemma}

We first establish the following auxiliary estimates. Analogously to~\eqref{eq:covariant-sobolev-norm}, for
	$B\in(\mathbb R^2)^{\Lambda_{M,\varepsilon}}$ and $s\in\mathbb R$, define
	the periodic covariant Sobolev norm
	$\|f\|_{H_B^{s,\varepsilon}(\Lambda_{M,\varepsilon})}
	:=\|(1+\mathcal L_B^\varepsilon)^{s/2}f\|_{L^{2,\varepsilon}(\Lambda_{M,\varepsilon})}$.
	The corresponding unweighted Sobolev embedding estimates
	of Section~\ref{sec:covariant-sobolev-spaces} also hold on the periodic lattice.
For $x\in\Lambda_\varepsilon$ and
	$B\in(\mathbb R^2)^{\Lambda_{M,\varepsilon}}$, let $w_x$ be the weight
	defined in~\eqref{eq:periodized-polynomial-weight} and set
	\begin{equation}
		\label{eq:periodic-localization-weight}
		N_x:=
		\|w_xB\|_{H^{1,\varepsilon}(\Lambda_{M,\varepsilon})}.
	\end{equation}
	By Lemma~\ref{lem:periodized-weighted-Besov-equivalences},
	with $p=q=2$ and $a=s=1$,
	\begin{equation}
		\label{eq:localized-gauge-H1-bound}
		N_x\lesssim\rho(x)^{-1}\|B\|_{H^{1,\varepsilon}(\rho)}.
	\end{equation}

	\begin{lemma}
		\label{lem:relative-covariant-laplacian-Hminus1-estimate}
		Let $B\in(\mathbb R^2)^{\Lambda_{M,\varepsilon}}$,
		$x\in\Lambda_\varepsilon$, $0<\beta<1$, and
		$f\in\mathbb C^{\Lambda_{M,\varepsilon}}$.
	We have
		\begin{equation}
			\label{eq:localized-relative-covariant-laplacian-Hminus1-estimate}
			\|(\mathcal L_B^\varepsilon-\mathcal L_0^\varepsilon)f\|_{H_B^{-1,\varepsilon}(\Lambda_{M,\varepsilon})}
			\lesssim_\beta N_x(1+N_x)
			\|w_x^{-2}f\|_{H^{\beta,\varepsilon}(\Lambda_{M,\varepsilon})}.
		\end{equation}
	\end{lemma}

	\begin{proof}
		Taking $A=-B$ in
		\eqref{eq:relative-covariant-laplacian-divergence-rewriting} gives
		\begin{equation}
			\label{eq:free-minus-covariant-laplacian-expansion}
			(\mathcal L_0^\varepsilon-\mathcal L_B^\varepsilon)f
			=iD_B^{\varepsilon *}
			\bigl(\overline{r_{-B}^\varepsilon}f+U_B^\varepsilon\rho_{-B}^\varepsilon f\bigr)
			+\rho_{-B}^{\varepsilon *}\rho_{-B}^\varepsilon f
			-i(D^{\varepsilon *}\overline{r_{-B}^\varepsilon})f.
		\end{equation}
		We use $\|w_x\|_{L^{\infty,\varepsilon}(\Lambda_{M,\varepsilon})}\lesssim1$ and
		$\tau_j^{\pm1}w_x\asymp w_x$, uniformly for $M\geq1$.

		By the Leibniz rule, $|D_j^\varepsilon w_x|\lesssim w_x$,
		and the uniform periodic embedding $H^{1,\varepsilon}(\Lambda_{M,\varepsilon})\hookrightarrow
		L^{q,\varepsilon}(\Lambda_{M,\varepsilon})$,
		\begin{equation}
			\label{eq:localized-gauge-coefficient-Sobolev-bounds}
			\|w_xr_{-B}^\varepsilon\|_{L^{q,\varepsilon}(\Lambda_{M,\varepsilon})}
			+\sum_{j=1}^2
			\|w_xD_j^{\varepsilon *}r_{-B,j}^\varepsilon\|_{L^{2,\varepsilon}(\Lambda_{M,\varepsilon})}
			\lesssim_q N_x,\qquad 2\leq q<\infty.
		\end{equation}

		The periodic Sobolev embeddings in~\eqref{eq:covariant-sobolev-embedding}
		give
		\begin{equation}
			\label{eq:localized-relative-laplacian-Sobolev-bounds}
			\begin{aligned}
				\|w_x^{-2}f\|_{L^{2/(1-\beta),\varepsilon}(\Lambda_{M,\varepsilon})}
				\lesssim_\beta\|w_x^{-2}f\|_{H^{\beta,\varepsilon}(\Lambda_{M,\varepsilon})},\quad
				\|g\|_{H_B^{-1,\varepsilon}(\Lambda_{M,\varepsilon})}
				\lesssim_\beta\|g\|_{L^{2/(2-\beta),\varepsilon}(\Lambda_{M,\varepsilon})}.
			\end{aligned}
		\end{equation}
		The covariant energy identity, H\"older's inequality, and
		\eqref{eq:localized-gauge-coefficient-Sobolev-bounds} yield
		\begin{equation}
			\label{eq:relative-covariant-laplacian-divergence-term-bound}
			\begin{aligned}
				&\left\|iD_B^{\varepsilon *}
				\bigl(\overline{r_{-B}^\varepsilon}f
				+U_B^\varepsilon\rho_{-B}^\varepsilon f\bigr)\right\|_{H_B^{-1,\varepsilon}(\Lambda_{M,\varepsilon})}\\
				&\quad\leq
				\|r_{-B}^\varepsilon f\|_{L^{2,\varepsilon}(\Lambda_{M,\varepsilon})}
				+\|\rho_{-B}^\varepsilon f\|_{L^{2,\varepsilon}(\Lambda_{M,\varepsilon})}\\
				&\quad\lesssim
				\|w_xr_{-B}^\varepsilon\|_{L^{2/\beta,\varepsilon}(\Lambda_{M,\varepsilon})}
				\|w_x^{-2}f\|_{L^{2/(1-\beta),\varepsilon}(\Lambda_{M,\varepsilon})}
				\lesssim_\beta N_x\|w_x^{-2}f\|_{H^{\beta,\varepsilon}(\Lambda_{M,\varepsilon})}.
			\end{aligned}
		\end{equation}
		For the quadratic term, the embedding
		$L^{2,\varepsilon}(\Lambda_{M,\varepsilon})\hookrightarrow H_B^{-1,\varepsilon}(\Lambda_{M,\varepsilon})$
		and H\"older's inequality yield
		\begin{equation}
			\label{eq:relative-covariant-laplacian-quadratic-term-bound}
			\begin{aligned}
				&\|\rho_{-B}^{\varepsilon *}\rho_{-B}^\varepsilon f\|_{H_B^{-1,\varepsilon}(\Lambda_{M,\varepsilon})}
				\leq\|\rho_{-B}^{\varepsilon *}\rho_{-B}^\varepsilon f\|_{L^{2,\varepsilon}(\Lambda_{M,\varepsilon})}\\
				&\quad\lesssim
				\|w_xr_{-B}^\varepsilon\|_{L^{4/\beta,\varepsilon}(\Lambda_{M,\varepsilon})}^2
				\|w_x^{-2}f\|_{L^{2/(1-\beta),\varepsilon}(\Lambda_{M,\varepsilon})}
				\lesssim_\beta N_x^2\|w_x^{-2}f\|_{H^{\beta,\varepsilon}(\Lambda_{M,\varepsilon})}.
			\end{aligned}
		\end{equation}
Using \eqref{eq:localized-gauge-coefficient-Sobolev-bounds} and~\eqref{eq:localized-relative-laplacian-Sobolev-bounds},
		\begin{equation}
			\label{eq:relative-covariant-laplacian-nondivergence-term-bound}
			\begin{aligned}
				&\|(D^{\varepsilon *}\overline{r_{-B}^\varepsilon})f\|_{H_B^{-1,\varepsilon}(\Lambda_{M,\varepsilon})}
				\lesssim_\beta\|(D^{\varepsilon *}\overline{r_{-B}^\varepsilon})f\|_{L^{2/(2-\beta),\varepsilon}(\Lambda_{M,\varepsilon})}\\
				&\quad\lesssim
				\sum_{j=1}^2
				\|w_xD_j^{\varepsilon *}r_{-B,j}^\varepsilon\|_{L^{2,\varepsilon}(\Lambda_{M,\varepsilon})}
				\|w_x^{-1}f\|_{L^{2/(1-\beta),\varepsilon}(\Lambda_{M,\varepsilon})}
				\lesssim_\beta N_x\|w_x^{-2}f\|_{H^{\beta,\varepsilon}(\Lambda_{M,\varepsilon})}.
			\end{aligned}
		\end{equation}
		Summing the last three bounds yields
		\eqref{eq:localized-relative-covariant-laplacian-Hminus1-estimate}.
	\end{proof}

	\begin{lemma}
		\label{lem:covariant-resolvent-smoothing-estimate}
		For $B\in(\mathbb R^2)^{\Lambda_{M,\varepsilon}}$, $t\geq0$,
		$\alpha\in\mathbb R$, and $\delta\in[0,2]$,
		\begin{equation}
			\label{eq:covariant-resolvent-smoothing-estimate}
			\|(\mathcal L_B^\varepsilon+1+t)^{-1}f\|_{H_B^{\alpha,\varepsilon}(\Lambda_{M,\varepsilon})}
			\lesssim_\delta(1+t)^{-1+\delta/2}
			\|f\|_{H_B^{\alpha-\delta,\varepsilon}(\Lambda_{M,\varepsilon})}.
		\end{equation}
	\end{lemma}

	\begin{proof}
		Let $(e_k^B)_k$ be an orthonormal eigenbasis of
		$(1+\mathcal L_B^\varepsilon)^{1/2}$ in
		$L^{2,\varepsilon}(\Lambda_{M,\varepsilon};\mathbb C)$, with eigenvalues
		$\lambda_k^B\geq1$. Spectral calculus gives
		\begin{equation}
			\label{eq:covariant-resolvent-smoothing-spectral-bound}
			\begin{aligned}
				\|(\mathcal L_B^\varepsilon+1+t)^{-1}f\|_{H_B^{\alpha,\varepsilon}(\Lambda_{M,\varepsilon})}^2
				=\sum_k
				\frac{(\lambda_k^B)^{2\alpha}}{((\lambda_k^B)^2+t)^2}
				|\langle f,e_k^B\rangle_{\mathbb C,M,\varepsilon}|^2
				\leq\sup_{\lambda\geq1}\frac{\lambda^{2\delta}}{(\lambda^2+t)^2}
				\|f\|_{H_B^{\alpha-\delta,\varepsilon}(\Lambda_{M,\varepsilon})}^2.
			\end{aligned}
		\end{equation}
		Since $0\leq\delta\leq2$,
		\begin{equation}
			\label{eq:covariant-resolvent-spectral-multiplier-bound}
			\frac{\lambda^{2\delta}}{(\lambda^2+t)^2}
			\leq
			\begin{cases}
				(1+t)^{-2+\delta},&\lambda^2\leq1+t,\\
				\lambda^{2\delta-4}\leq(1+t)^{-2+\delta},&\lambda^2\geq1+t.
			\end{cases}
		\end{equation}
		Taking square roots yields the estimate.
	\end{proof}

	\begin{proof}[Proof of Lemma~\ref{lem:covariant-free-gff-difference-negative-holder-moment}]
		For the first estimate, weighted embeddings reduce the proof to
		$0<\kappa<1$. The Gaussian isometry gives
		\begin{equation}
			\label{eq:covariant-free-gff-difference-block-second-moment}
			\begin{aligned}
				\left(\mathbb E^\phi[
				|\Delta_j^\varepsilon(\phi_B^\varepsilon-\phi^\varepsilon)(x)|^2]
				\right)^{1/2}=\sqrt2\left\|
				\bigl[(\mathcal L_B^\varepsilon+1)^{-1/2}
				-(\mathcal L_0^\varepsilon+1)^{-1/2}\bigr]
				\eta_j^{M,\varepsilon}(x-\cdot)
				\right\|_{L^{2,\varepsilon}(\Lambda_{M,\varepsilon})}.
			\end{aligned}
		\end{equation}
		For a positive self-adjoint operator $T\geq1$, spectral calculus yields
$
			T^{-1/2}=\frac2\pi\int_0^\infty(T+s^2)^{-1}\,ds.
$
		Hence the resolvent identity gives
		\begin{equation}
			\label{eq:covariant-free-inverse-square-root-difference}
			\begin{aligned}
				&(\mathcal L_B^\varepsilon+1)^{-1/2}
				-(\mathcal L_0^\varepsilon+1)^{-1/2}=-\frac2\pi\int_0^\infty
				(\mathcal L_B^\varepsilon+1+s^2)^{-1}
				(\mathcal L_B^\varepsilon-\mathcal L_0^\varepsilon)
				(\mathcal L_0^\varepsilon+1+s^2)^{-1}\,ds.
			\end{aligned}
		\end{equation}

		By Lemma~\ref{lem:relative-covariant-laplacian-Hminus1-estimate}, for $0<\beta<1$,
		\begin{equation}
			\label{eq:relative-laplacian-Hsigma-HBminus1-mapping}
			\|(\mathcal L_B^\varepsilon-\mathcal L_0^\varepsilon)f\|_{H_B^{-1,\varepsilon}(\Lambda_{M,\varepsilon})}
			\lesssim_{\nu,\beta}N_x(1+N_x)
			\|w_x^{-2}f\|_{H^{\beta,\varepsilon}(\Lambda_{M,\varepsilon})}.
		\end{equation}
		Lemma~\ref{lem:covariant-resolvent-smoothing-estimate}, with
		$\alpha=0$, $\delta=1$, and $t=s^2$, gives
		\begin{equation}
			\label{eq:covariant-resolvent-HBminus1-L2-mapping}
			\|(\mathcal L_B^\varepsilon+1+s^2)^{-1}
			\|_{H_B^{-1,\varepsilon}(\Lambda_{M,\varepsilon})\to L^{2,\varepsilon}(\Lambda_{M,\varepsilon})}
			\lesssim(1+s^2)^{-1/2}.
		\end{equation}
		For $\beta\in\mathbb R$, $0<\delta<2$, and $-1\leq j\leq J^\varepsilon$,
		\begin{equation}
			\label{eq:covariant-free-gff-difference-kernel-Sobolev-bound}
			\begin{aligned}
				\left\|w_x^{-2}(\mathcal L_0^\varepsilon+1+s^2)^{-1}
				\eta_{j,x}^{M,\varepsilon}
				\right\|_{H^{\beta,\varepsilon}(\Lambda_{M,\varepsilon})}
				&\lesssim_\beta
				\left\|(\mathcal L_0^\varepsilon+1+s^2)^{-1}
				\eta_{j,x}^{M,\varepsilon}
				\right\|_{B_{2,2}^{\beta,\varepsilon}(w_x^{-2};\Lambda_{M,\varepsilon})}\\
				&\lesssim_\beta
				\frac1{1+s^2+2^{2j}}
				\left\|\eta_{j,x}^{M,\varepsilon}
				\right\|_{B_{2,2}^{\beta,\varepsilon}(w_x^{-2};\Lambda_{M,\varepsilon})}\\
				&\overset{\text{Lemma~\ref{lem:weighted-periodic-LP-kernel}}}{\lesssim_\beta}
				\frac{2^{j(\beta+1)}}{1+s^2+2^{2j}}
				\lesssim_{\beta,\delta}
				2^{j(\beta+\delta-1)}(1+s^2)^{-\delta/2}.
			\end{aligned}
		\end{equation}
		The first inequality follows from
		Lemma~\ref{lem:periodized-weighted-Besov-equivalences}\textup{(i)}
		with $a=-2$ and $p=q=2$.
		The second follows from
		$\operatorname{supp}\widehat{\eta_j^\varepsilon}\cap Q_{M,\varepsilon}
		\subset\{k:|k|\asymp2^j\}$ for $j\geq0$
		(with $|k|\lesssim1$ for $j=-1$).
		Combining these bounds with
		\eqref{eq:covariant-free-inverse-square-root-difference},
		\eqref{eq:localized-gauge-H1-bound}, and
		\eqref{eq:covariant-free-gff-difference-kernel-Sobolev-bound}, we obtain
		\begin{equation}
			\label{eq:covariant-free-gff-difference-operator-bound}
			\begin{aligned}
				\left\|
				\bigl[(\mathcal L_B^\varepsilon+1)^{-1/2}
				-(\mathcal L_0^\varepsilon+1)^{-1/2}\bigr]
				\eta_j^{M,\varepsilon}(x-\cdot)
				\right\|_{L^{2,\varepsilon}(\Lambda_{M,\varepsilon})}
				&\lesssim_{\nu,\beta,\delta}
				N_x(1+N_x)2^{j(\beta+\delta-1)}
				\int_0^\infty(1+s^2)^{-(1+\delta)/2}\,ds\\
				&\lesssim_{\nu,\beta,\delta}
				\rho(x)^{-2}\|B\|_{H^{1,\varepsilon}(\rho)}
				\bigl(1+\|B\|_{H^{1,\varepsilon}(\rho)}\bigr)
				2^{j(\beta+\delta-1)}.
			\end{aligned}
		\end{equation}
		Take $\beta=\delta=\kappa/3$.
		The Gaussian field $\phi_B^\varepsilon-\phi^\varepsilon$ belongs to
		the first Wiener chaos. By the Gaussian isometry above,
		Lemma~\ref{lem:lattice-stochastic-regularity-criterion} applies with
		$\alpha=2\kappa/3-1$, $a=2+\kappa$, $b=2$, $m=1$, and regularity loss
		$\kappa/3$. This yields
		\eqref{eq:covariant-free-gff-difference-negative-holder-moment}.

		For the interpolation estimate, fix any $\kappa>0$.
		Corollary~\ref{cor:standard-gauge-field-stochastic-estimates} and
		Lemma~\ref{lem:covariant-gff-negative-holder-moment} give
		\begin{equation}
			\label{eq:covariant-free-gff-difference-negative-endpoint}
			\left(\mathbb E^\phi\left[
			\|\phi_B^\varepsilon-\phi^\varepsilon\|_{\mathcal C^{-\kappa,\varepsilon}(\rho^\kappa)}^p
			\right]\right)^{1/p}
			\lesssim_\kappa p^{1/2}.
		\end{equation}
		For $B\ne0$ and $\sigma\in[-\kappa,1-\kappa]$, put
		$\theta=\kappa+\sigma\in[0,1]$.
		The dyadic definition of the norms and
		$(1-\theta)\kappa+\theta(2+\kappa)=3\kappa+2\sigma$ yield
		\begin{equation}
			\label{eq:covariant-free-gff-difference-holder-interpolation}
			\begin{aligned}
				\frac{\|\phi_B^\varepsilon-\phi^\varepsilon\|_{\mathcal C^{\sigma,\varepsilon}(\rho^{3\kappa+2\sigma})}}
				{\|B\|_{H^{1,\varepsilon}(\rho)}^\theta
				\bigl(1+\|B\|_{H^{1,\varepsilon}(\rho)}\bigr)^\theta}
				&\leq
				\|\phi_B^\varepsilon-\phi^\varepsilon\|_{\mathcal C^{-\kappa,\varepsilon}(\rho^\kappa)}^{1-\theta}
				\left(
				\frac{\|\phi_B^\varepsilon-\phi^\varepsilon\|_{\mathcal C^{1-\kappa,\varepsilon}(\rho^{2+\kappa})}}
				{\|B\|_{H^{1,\varepsilon}(\rho)}
				\bigl(1+\|B\|_{H^{1,\varepsilon}(\rho)}\bigr)}
				\right)^\theta\\
				&\leq
				\|\phi_B^\varepsilon-\phi^\varepsilon\|_{\mathcal C^{-\kappa,\varepsilon}(\rho^\kappa)}
				+\frac{\|\phi_B^\varepsilon-\phi^\varepsilon\|_{\mathcal C^{1-\kappa,\varepsilon}(\rho^{2+\kappa})}}
				{\|B\|_{H^{1,\varepsilon}(\rho)}
				\bigl(1+\|B\|_{H^{1,\varepsilon}(\rho)}\bigr)}.
			\end{aligned}
		\end{equation}
		Taking the supremum over $\sigma$, followed by the $L^p$ norm, and
		using~\eqref{eq:covariant-free-gff-difference-negative-endpoint} and
		\eqref{eq:covariant-free-gff-difference-negative-holder-moment}
		yields~\eqref{eq:covariant-free-gff-difference-interpolated-holder-moment}.
	\end{proof}

	\subsection{\texorpdfstring{Estimate of
			$:\!(\phi_B^\varepsilon)^{k_1}
			(\overline{\phi_B^\varepsilon})^{k_2}\!:_\varepsilon$}
		{Estimate of Wick powers of the covariant Gaussian free field}}
	\label{subsec:stochastic-estimate-covariant-gff-wick-powers}

	\begin{corollary}
		\label{cor:covariant-gff-wick-power-negative-holder-moment}
		For every $B\in(\mathbb R^2)^{\Lambda_{M,\varepsilon}}$,
		$k_1,k_2\in\mathbb N_0$, $\kappa>0$, and $p\in[1,\infty)$,
		\begin{equation}
			\label{eq:covariant-gff-wick-power-negative-holder-moment}
			\left(\mathbb E^\phi\left[
			\left\|:\!(\phi_B^\varepsilon)^{k_1}
			(\overline{\phi_B^\varepsilon})^{k_2}\!:_\varepsilon
			\right\|_{\mathcal C^{-\kappa,\varepsilon}(\rho^\kappa)}^p
			\right]\right)^{1/p}
			\lesssim_{\kappa,k_1,k_2}p^{(k_1+k_2)/2}
			\bigl(1+\|B\|_{H^{1,\varepsilon}(\rho)}\bigr)^\kappa.
		\end{equation}
	\end{corollary}

	\begin{proof}
		We may assume $\kappa>0$ small.
		The binomial identity for Wick polynomials gives
		\begin{equation}
			\label{eq:covariant-gff-wick-power-shift-expansion}
			:\!(\phi_B^\varepsilon)^{k_1}
			(\overline{\phi_B^\varepsilon})^{k_2}\!:_\varepsilon
			=\sum_{\substack{j_1+j_3=k_1\\j_2+j_4=k_2}}
			\binom{k_1}{j_1}\binom{k_2}{j_2}
			:\!(\phi^\varepsilon)^{j_1}
			(\overline{\phi^\varepsilon})^{j_2}\!:_\varepsilon
			(\phi_B^\varepsilon-\phi^\varepsilon)^{j_3}
			\overline{(\phi_B^\varepsilon-\phi^\varepsilon)}^{j_4}.
		\end{equation}
		Lemma~\ref{lem:discrete-paraproduct-estimates}
		and weighted embeddings imply
		\begin{equation}
			\label{eq:covariant-gff-wick-power-product-bound}
			\begin{aligned}
				&\left\|:\!(\phi^\varepsilon)^{j_1}
				(\overline{\phi^\varepsilon})^{j_2}\!:_\varepsilon
				(\phi_B^\varepsilon-\phi^\varepsilon)^{j_3}
				\overline{(\phi_B^\varepsilon-\phi^\varepsilon)}^{j_4}
				\right\|_{\mathcal C^{-\kappa,\varepsilon}(\rho^{(1+7(k_1+k_2))\kappa})}\\
				&\quad\lesssim_{\kappa,k_1,k_2}
				\left\|:\!(\phi^\varepsilon)^{j_1}
				(\overline{\phi^\varepsilon})^{j_2}\!:_\varepsilon
				\right\|_{\mathcal C^{-\kappa,\varepsilon}(\rho^\kappa)}
				\|\phi_B^\varepsilon-\phi^\varepsilon
				\|_{\mathcal C^{2\kappa,\varepsilon}(\rho^{7\kappa})}^{j_3+j_4}.
			\end{aligned}
		\end{equation}
		H\"older's inequality yields
		\begin{equation}
			\label{eq:covariant-gff-wick-power-term-moment-bound}
			\begin{aligned}
				&\left(\mathbb E^\phi\left[
				\left\|:\!(\phi^\varepsilon)^{j_1}
				(\overline{\phi^\varepsilon})^{j_2}\!:_\varepsilon
				(\phi_B^\varepsilon-\phi^\varepsilon)^{j_3}
				\overline{(\phi_B^\varepsilon-\phi^\varepsilon)}^{j_4}
				\right\|_{\mathcal C^{-\kappa,\varepsilon}(\rho^{(1+7(k_1+k_2))\kappa})}^p
				\right]\right)^{1/p}\\
				&\quad\lesssim_{\kappa,k_1,k_2}
				\left(\mathbb E^\phi\left[
				\left\|:\!(\phi^\varepsilon)^{j_1}
				(\overline{\phi^\varepsilon})^{j_2}\!:_\varepsilon
				\right\|_{\mathcal C^{-\kappa,\varepsilon}(\rho^\kappa)}^{2p}
				\right]\right)^{1/(2p)}
				\left(\mathbb E^\phi\left[
				\|\phi_B^\varepsilon-\phi^\varepsilon
				\|_{\mathcal C^{2\kappa,\varepsilon}(\rho^{7\kappa})}^{2p(j_3+j_4)}
				\right]\right)^{1/(2p)}.
			\end{aligned}
		\end{equation}
		The last factor is one when $j_3+j_4=0$.
		Otherwise, apply
		\eqref{eq:covariant-free-gff-difference-interpolated-holder-moment}
		with $\sigma=2\kappa$ and moment order $2p(j_3+j_4)$.
		Together with Corollary~\ref{cor:standard-gauge-field-stochastic-estimates},
		this gives
		\begin{equation}
			\label{eq:covariant-gff-wick-power-term-final-bound}
			\begin{aligned}
				&\left(\mathbb E^\phi\left[
				\left\|:\!(\phi^\varepsilon)^{j_1}
				(\overline{\phi^\varepsilon})^{j_2}\!:_\varepsilon
				(\phi_B^\varepsilon-\phi^\varepsilon)^{j_3}
				\overline{(\phi_B^\varepsilon-\phi^\varepsilon)}^{j_4}
				\right\|_{\mathcal C^{-\kappa,\varepsilon}(\rho^{(1+7(k_1+k_2))\kappa})}^p
				\right]\right)^{1/p}\\
				&\quad\lesssim_{\kappa,k_1,k_2}
				p^{(k_1+k_2)/2}
				\bigl(1+\|B\|_{H^{1,\varepsilon}(\rho)}\bigr)^{6\kappa(j_3+j_4)}
				\lesssim_{\kappa,k_1,k_2}
				p^{(k_1+k_2)/2}
				\bigl(1+\|B\|_{H^{1,\varepsilon}(\rho)}\bigr)^{(1+7(k_1+k_2))\kappa}.
			\end{aligned}
		\end{equation}
		Summing and replacing $\kappa$ by $\kappa/[1+7(k_1+k_2)]$ gives the claim.
	\end{proof}

	\subsection{\texorpdfstring{Estimate of
			$r_{A^\varepsilon}^\varepsilon$
			and
			$:\!\lvert r_{A^\varepsilon}^\varepsilon\rvert^2\!:_\varepsilon$}
		{Estimate of the nonlinear gauge coefficients}}
	\label{subsec:stochastic-estimate-nonlinear-gauge-coefficients}

	\begin{lemma}
		\label{lem:nonlinear-link-coefficient-approximation-holder-moment}
		For every $\kappa\in(0,1)$ and $p\in[1,\infty)$,
		\begin{equation}
			\label{eq:nonlinear-link-coefficient-approximation-holder-moment}
			\left(\mathbb E\left[
			\|r_{A^\varepsilon}^\varepsilon-A^\varepsilon
			\|_{\mathcal C^{1-\kappa,\varepsilon}(\rho^\kappa)}^p
			\right]\right)^{1/p}
			\lesssim_\kappa p\varepsilon^{\kappa/2}.
		\end{equation}
	\end{lemma}

	\begin{proof}
		Dyadic reconstruction, $2^{J^\varepsilon}\asymp\varepsilon^{-1}$,
		and~\eqref{eq:standard-gauge-field-negative-holder-moments} give
		\begin{equation}
			\label{eq:gauge-field-Linfty-moment}
			\begin{aligned}
				\left(\mathbb E\left[
				\|A^\varepsilon\|_{L^{\infty,\varepsilon}(\rho^{\kappa/2})}^p
				\right]\right)^{1/p}
				\lesssim_\kappa\varepsilon^{-\kappa/4}
				\left(\mathbb E\left[
				\|A^\varepsilon\|_{\mathcal C^{-\kappa/4,\varepsilon}(\rho^{\kappa/4})}^p
				\right]\right)^{1/p}
				\lesssim_\kappa\varepsilon^{-\kappa/4}p^{1/2}.
			\end{aligned}
		\end{equation}
		Taylor's formula gives, for $j=1,2$,
		\begin{equation}
			\label{eq:nonlinear-link-coefficient-approximation-Linfty}
			\|r_{A^\varepsilon,j}^\varepsilon-A_j^\varepsilon
			\|_{L^{\infty,\varepsilon}(\rho^\kappa)}
			\lesssim\varepsilon
			\|A_j^\varepsilon\|_{L^{\infty,\varepsilon}(\rho^{\kappa/2})}^2.
		\end{equation}
	Similarly,
		$\|h\|_{\mathcal C^{1-\kappa,\varepsilon}(\rho^\kappa)}
		\lesssim_\kappa\varepsilon^{\kappa-1}
		\|h\|_{L^{\infty,\varepsilon}(\rho^\kappa)}$.
		Applying this to the remainder and using the preceding estimates,
		\begin{equation}
			\label{eq:nonlinear-link-coefficient-approximation-holder-proof}
			\begin{aligned}
				\left(\mathbb E\left[
				\|r_{A^\varepsilon}^\varepsilon-A^\varepsilon
				\|_{\mathcal C^{1-\kappa,\varepsilon}(\rho^\kappa)}^p
				\right]\right)^{1/p}&\lesssim_\kappa\varepsilon^{\kappa-1}
				\left(\mathbb E\left[
				\|r_{A^\varepsilon}^\varepsilon-A^\varepsilon
				\|_{L^{\infty,\varepsilon}(\rho^\kappa)}^p
				\right]\right)^{1/p}\\
				&\quad\lesssim_\kappa\varepsilon^\kappa
				\sum_{j=1}^2\left(\mathbb E\left[
				\|A_j^\varepsilon\|_{L^{\infty,\varepsilon}(\rho^{\kappa/2})}^{2p}
				\right]\right)^{1/p}
				\lesssim_\kappa p\varepsilon^{\kappa/2}.
			\end{aligned}
		\end{equation}
	\end{proof}

	\begin{corollary}
		\label{cor:nonlinear-gauge-coefficients-negative-holder-moments}
		For every $\kappa>0$ and $p\in[1,\infty)$,
		\begin{equation}
			\label{eq:nonlinear-gauge-coefficients-negative-holder-moments}
			\begin{aligned}
				\left(\mathbb E\left[
				\|r_{A^\varepsilon}^\varepsilon
				\|_{\mathcal C^{-\kappa,\varepsilon}(\rho^\kappa)}^p
				\right]\right)^{1/p}
				\lesssim_\kappa p,\quad
				\left(\mathbb E\left[
				\|:\!|r_{A^\varepsilon}^\varepsilon|^2\!:_\varepsilon
				\|_{\mathcal C^{-\kappa,\varepsilon}(\rho^\kappa)}^p
				\right]\right)^{1/p}
				\lesssim_\kappa p^2.
			\end{aligned}
		\end{equation}
	\end{corollary}

	\begin{proof}
		We may assume $\kappa>0$ small.
		Set $R^\varepsilon:=
		r_{A^\varepsilon}^\varepsilon-A^\varepsilon$.
		The first bound follows from
		\eqref{eq:standard-gauge-field-negative-holder-moments} and
		Lemma~\ref{lem:nonlinear-link-coefficient-approximation-holder-moment}.
		For the second, \eqref{eq:wick-rho-square} gives
		\begin{equation}
			\label{eq:nonlinear-gauge-coefficients-wick-expansion}
			:\!|r_{A^\varepsilon}^\varepsilon|^2\!:_\varepsilon
			=:\!|A^\varepsilon|^2\!:_\varepsilon
			+2\operatorname{Re}\sum_{j=1}^2
			A_j^\varepsilon\overline{R_j^\varepsilon}
			+|R^\varepsilon|^2.
		\end{equation}
		Lemma~\ref{lem:discrete-paraproduct-estimates} bounds
		$A_j^\varepsilon\overline{R_j^\varepsilon}$ and
		$|R^\varepsilon|^2$ in
		$\mathcal C^{-\kappa/2,\varepsilon}(\rho^\kappa)
		\hookrightarrow\mathcal C^{-\kappa,\varepsilon}(\rho^\kappa)$.
		H\"older's inequality and the same moment estimates with parameter
		$\kappa/2$ yield
		\begin{equation}
			\label{eq:nonlinear-gauge-coefficients-weighted-proof}
			\left(\mathbb E\left[
			\|:\!|r_{A^\varepsilon}^\varepsilon|^2\!:_\varepsilon
			\|_{\mathcal C^{-\kappa,\varepsilon}(\rho^\kappa)}^p
			\right]\right)^{1/p}
			\lesssim_\kappa
			p+p^{3/2}\varepsilon^{\kappa/4}+p^2\varepsilon^{\kappa/2}
			\lesssim_\kappa p^2.
		\end{equation}
	\end{proof}

	\subsection{\texorpdfstring{Estimate of
			$D^{\varepsilon *}
			\overline{r_{A^\varepsilon}^\varepsilon}$
			and
			$\bigl(D^{\varepsilon *}
			\overline{r_{A^\varepsilon}^\varepsilon}\bigr)
			\phi_{B^\varepsilon}^\varepsilon$}
		{Estimate of the nonlinear divergence terms}}
	\label{subsec:stochastic-estimate-nonlinear-divergence}

	\begin{lemma}
		\label{lem:nonlinear-divergence-negative-holder-moment}
		For every $\kappa\in(0,1)$ and $p\in[1,\infty)$,
		\begin{equation}
			\label{eq:nonlinear-divergence-negative-holder-moment}
			\left(
			\mathbb{E}\left[
			\left\lVert
			D^{\varepsilon *}
			\overline{r_{A^\varepsilon}^\varepsilon}
			\right\rVert_{\mathcal{C}^{-\kappa,\varepsilon}(\rho^\kappa)}^p
			\right]
			\right)^{1/p}
			\lesssim_\kappa
			p\varepsilon^{\kappa/2}.
		\end{equation}
	\end{lemma}

	\begin{proof}
		Since $A^\varepsilon$ is real-valued and
		divergence-free,
		\begin{equation}
			\label{eq:nonlinear-divergence-coulomb-cancellation}
			D^{\varepsilon *}
			\overline{r_{A^\varepsilon}^\varepsilon}
			=
			D^{\varepsilon *}
			\left(
			\overline{r_{A^\varepsilon}^\varepsilon}
			-A^\varepsilon
			\right).
		\end{equation}
		$\lVert D^{\varepsilon *}f\rVert_{\mathcal{C}^{-\kappa,\varepsilon}(\rho^\kappa)}
		\lesssim_\kappa
		\lVert f\rVert_{\mathcal{C}^{1-\kappa,\varepsilon}(\rho^\kappa)}$ and
		Lemma~\ref{lem:nonlinear-link-coefficient-approximation-holder-moment}
	yield~\eqref{eq:nonlinear-divergence-negative-holder-moment}.
	\end{proof}

	\begin{lemma}
		\label{lem:nonlinear-divergence-covariant-gff-product-negative-holder-moment}
		For every $\kappa>0$ and $p\in[1,\infty)$,
		\begin{equation}
			\label{eq:nonlinear-divergence-covariant-gff-product-negative-holder-moment}
			\left(\mathbb E\left[
			\left\|\bigl(D^{\varepsilon *}\overline{r_{A^\varepsilon}^\varepsilon}\bigr)
			\phi_{B^\varepsilon}^\varepsilon
			\right\|_{\mathcal C^{-\kappa,\varepsilon}(\rho^\kappa)}^p
			\right]\right)^{1/p}
			\lesssim_\kappa p^{3/2}
			\left(\mathbb E\left[
			\bigl(1+\|B^\varepsilon\|_{H^{1,\varepsilon}(\rho)}\bigr)^{\kappa p}
			\right]\right)^{1/p}.
		\end{equation}
	\end{lemma}

	\begin{proof}
		By the weighted embeddings, it suffices to consider $0<\kappa<1$.
		We decompose
		\begin{equation}
			\label{eq:nonlinear-divergence-covariant-gff-free-difference-splitting}
			\bigl(D^{\varepsilon *}\overline{r_{A^\varepsilon}^\varepsilon}\bigr)
			\phi_{B^\varepsilon}^\varepsilon
			=
			\bigl(D^{\varepsilon *}\overline{r_{A^\varepsilon}^\varepsilon}\bigr)
			\phi^\varepsilon
			+\bigl(D^{\varepsilon *}\overline{r_{A^\varepsilon}^\varepsilon}\bigr)
			\bigl(\phi_{B^\varepsilon}^\varepsilon-\phi^\varepsilon\bigr).
		\end{equation}

		We first estimate the free-field term. Fixing the gauge-field variables,
		$\bigl(D^{\varepsilon *}\overline{r_{A^\varepsilon}^\varepsilon}\bigr)
		\phi^\varepsilon$ belongs to the first Wiener chaos over $\Omega^\phi$.
		The Gaussian isometry gives
		\begin{equation}
			\label{eq:nonlinear-divergence-free-gff-block-isometry}
			\begin{aligned}
				\mathbb E^\phi\left[
				\left|\Delta_j^\varepsilon
				\Bigl(\bigl(D^{\varepsilon *}\overline{r_{A^\varepsilon}^\varepsilon}\bigr)
				\phi^\varepsilon\Bigr)(x)\right|^2\right]
				&=2\left\|(\mathcal L_0^\varepsilon+1)^{-1/2}
				\Bigl(\bigl(D^{\varepsilon *}\overline{r_{A^\varepsilon}^\varepsilon}\bigr)
				\eta_{j,x}^{M,\varepsilon}\Bigr)
				\right\|_{L^{2,\varepsilon}(\Lambda_{M,\varepsilon})}^2\\
				&=2\left\|\bigl(D^{\varepsilon *}\overline{r_{A^\varepsilon}^\varepsilon}\bigr)
				\eta_{j,x}^{M,\varepsilon}
				\right\|_{H^{-1,\varepsilon}(\Lambda_{M,\varepsilon})}^2.
			\end{aligned}
		\end{equation}
		Set $q=(1-\kappa/12)^{-1}$. Since $1<q<2$, the uniform Besov embeddings give
		\begin{equation}
			\label{eq:nonlinear-divergence-free-gff-besov-embeddings}
			B_{q,\infty}^{-\kappa/12}(\Lambda_{M,\varepsilon})
			\hookrightarrow
			B_{2,\infty}^{-1+\kappa/12}(\Lambda_{M,\varepsilon})
			\hookrightarrow H^{-1,\varepsilon}(\Lambda_{M,\varepsilon}).
		\end{equation}
		By Lemmas~\ref{lem:discrete-paraproduct-estimates},
		\ref{lem:periodized-weight-power-bounds}, and
		\ref{lem:periodized-weighted-Besov-equivalences}, using
		$-\kappa/12+\kappa/6>0$,
		\begin{equation}
			\label{eq:nonlinear-divergence-free-gff-test-product-bound}
			\begin{aligned}
				\left\|\bigl(D^{\varepsilon *}\overline{r_{A^\varepsilon}^\varepsilon}\bigr)
				\eta_{j,x}^{M,\varepsilon}
				\right\|_{H^{-1,\varepsilon}(\Lambda_{M,\varepsilon})}
				&\lesssim_\kappa
				\left\|\bigl(D^{\varepsilon *}\overline{r_{A^\varepsilon}^\varepsilon}\bigr)
				\eta_{j,x}^{M,\varepsilon}
				\right\|_{B_{q,\infty}^{-\kappa/12}(\Lambda_{M,\varepsilon})}\\
				&\lesssim_\kappa
				\left\|D^{\varepsilon *}\overline{r_{A^\varepsilon}^\varepsilon}
				\right\|_{\mathcal C^{-\kappa/12}((w_x^M)^{\kappa/12};\Lambda_{M,\varepsilon})}
				\left\|\eta_{j,x}^{M,\varepsilon}
				\right\|_{B_{q,\infty}^{\kappa/6}((w_x^M)^{-\kappa/12};\Lambda_{M,\varepsilon})}\\
				&\lesssim_\kappa \rho(x)^{-\kappa/12}
				\left\|D^{\varepsilon *}\overline{r_{A^\varepsilon}^\varepsilon}
				\right\|_{\mathcal C^{-\kappa/12,\varepsilon}(\rho^{\kappa/12})}
				2^{j\kappa/3}.
			\end{aligned}
		\end{equation}
		The last step uses the weight comparison in
		Lemma~\ref{lem:periodized-weighted-Besov-equivalences} and the kernel bound
		from Lemma~\ref{lem:weighted-periodic-LP-kernel}:
		\begin{equation}
			\label{eq:nonlinear-divergence-free-gff-kernel-besov-bound}
			\begin{aligned}
				\left\|\eta_{j,x}^{M,\varepsilon}
				\right\|_{B_{q,\infty}^{\kappa/6}((w_x^M)^{-\kappa/12};\Lambda_{M,\varepsilon})}
				\simeq_\kappa
				\left\|(w_x^M)^{-\kappa/12}\eta_{j,x}^{M,\varepsilon}
				\right\|_{B_{q,\infty}^{\kappa/6}(\Lambda_{M,\varepsilon})}
				\lesssim_\kappa 2^{2j(1-1/q)}2^{j\kappa/6}
				=2^{j\kappa/3}.
			\end{aligned}
		\end{equation}
		Thus
		\begin{equation}
			\label{eq:nonlinear-divergence-free-gff-block-bound}
			\sup_{x\in\Lambda_\varepsilon}\rho(x)^{\kappa/12}
			\left(\mathbb E^\phi\left[
			\left|\Delta_j^\varepsilon
			\Bigl(\bigl(D^{\varepsilon *}\overline{r_{A^\varepsilon}^\varepsilon}\bigr)
			\phi^\varepsilon\Bigr)(x)\right|^2\right]\right)^{1/2}
			\lesssim_\kappa
			\left\|D^{\varepsilon *}\overline{r_{A^\varepsilon}^\varepsilon}
			\right\|_{\mathcal C^{-\kappa/12,\varepsilon}(\rho^{\kappa/12})}
			2^{j\kappa/3}.
		\end{equation}
		Apply Lemma~\ref{lem:lattice-block-regularity-criterion} conditionally
		with $\alpha=\kappa/3$, $b=\kappa/12<a=\kappa/6$, $m=1$,
		and regularity loss $\kappa/12$. We obtain
		\begin{equation}
			\label{eq:nonlinear-divergence-free-gff-conditional-holder-bound}
			\left(\mathbb E^\phi\left[
			\left\|\bigl(D^{\varepsilon *}\overline{r_{A^\varepsilon}^\varepsilon}\bigr)
			\phi^\varepsilon
			\right\|_{\mathcal C^{-5\kappa/12,\varepsilon}(\rho^{\kappa/6})}^p
			\right]\right)^{1/p}
			\lesssim_\kappa p^{1/2}
			\left\|D^{\varepsilon *}\overline{r_{A^\varepsilon}^\varepsilon}
			\right\|_{\mathcal C^{-\kappa/12,\varepsilon}(\rho^{\kappa/12})}.
		\end{equation}
		Taking the outer expectation, the weighted embedding into
		$\mathcal C^{-\kappa,\varepsilon}(\rho^\kappa)$ and
		Lemma~\ref{lem:nonlinear-divergence-negative-holder-moment} give
		\begin{equation}
			\label{eq:nonlinear-divergence-free-gff-full-moment-bound}
			\left(\mathbb E\left[
			\left\|\bigl(D^{\varepsilon *}\overline{r_{A^\varepsilon}^\varepsilon}\bigr)
			\phi^\varepsilon
			\right\|_{\mathcal C^{-\kappa,\varepsilon}(\rho^\kappa)}^p
			\right]\right)^{1/p}
			\lesssim_\kappa p^{3/2}\varepsilon^{\kappa/24}
			\lesssim_\kappa p^{3/2}.
		\end{equation}

		For the second term, apply
		\eqref{eq:covariant-free-gff-difference-interpolated-holder-moment}
		with parameters $\kappa/12$ and $\sigma=\kappa/6$.
		Since $\kappa/6<1-\kappa/12$, this gives
		\begin{equation}
			\label{eq:nonlinear-divergence-gff-difference-holder-bound}
			\left(\mathbb E^\phi\left[
			\left\|\phi_{B^\varepsilon}^\varepsilon-\phi^\varepsilon
			\right\|_{\mathcal C^{\kappa/6,\varepsilon}(\rho^{7\kappa/12})}^p
			\right]\right)^{1/p}
			\lesssim_\kappa p^{1/2}
			\bigl(1+\|B^\varepsilon\|_{H^{1,\varepsilon}(\rho)}\bigr)^{\kappa/2}.
		\end{equation}
		Lemma~\ref{lem:discrete-paraproduct-estimates} and the weighted embeddings give
		\begin{equation}
			\label{eq:nonlinear-divergence-gff-difference-product-bound}
			\left\|\bigl(D^{\varepsilon *}\overline{r_{A^\varepsilon}^\varepsilon}\bigr)
			\bigl(\phi_{B^\varepsilon}^\varepsilon-\phi^\varepsilon\bigr)
			\right\|_{\mathcal C^{-\kappa,\varepsilon}(\rho^\kappa)}
			\lesssim_\kappa
			\left\|D^{\varepsilon *}\overline{r_{A^\varepsilon}^\varepsilon}
			\right\|_{\mathcal C^{-\kappa/12,\varepsilon}(\rho^{\kappa/12})}
			\left\|\phi_{B^\varepsilon}^\varepsilon-\phi^\varepsilon
			\right\|_{\mathcal C^{\kappa/6,\varepsilon}(\rho^{7\kappa/12})}.
		\end{equation}
		Conditioning on the gauge-field variables, then applying H\"older's
		inequality and Lemma~\ref{lem:nonlinear-divergence-negative-holder-moment}
		with moment order $2p$, we obtain
		\begin{equation}
			\label{eq:nonlinear-divergence-gff-difference-full-moment-bound}
			\begin{aligned}
				\left(\mathbb E\left[
				\left\|\bigl(D^{\varepsilon *}\overline{r_{A^\varepsilon}^\varepsilon}\bigr)
				\bigl(\phi_{B^\varepsilon}^\varepsilon-\phi^\varepsilon\bigr)
				\right\|_{\mathcal C^{-\kappa,\varepsilon}(\rho^\kappa)}^p
				\right]\right)^{1/p}
				&\lesssim_\kappa p^{3/2}\varepsilon^{\kappa/24}
				\left(\mathbb E\left[
				\bigl(1+\|B^\varepsilon\|_{H^{1,\varepsilon}(\rho)}\bigr)^{\kappa p}
				\right]\right)^{1/(2p)}\\
				&\lesssim_\kappa p^{3/2}
				\left(\mathbb E\left[
				\bigl(1+\|B^\varepsilon\|_{H^{1,\varepsilon}(\rho)}\bigr)^{\kappa p}
				\right]\right)^{1/p}.
			\end{aligned}
		\end{equation}
		Combining the two estimates yields the claimed estimate.
	\end{proof}
	
	\subsection{Estimates of mixed gauge--Higgs products}
	\label{subsec:stochastic-estimate-mixed-gauge-higgs-products}

	\begin{lemma}
		\label{cor:mixed-gauge-higgs-product-negative-holder-moments}
		For every $\kappa>0$ and $p\in[1,\infty)$,
		\begin{equation}
			\label{eq:mixed-gauge-higgs-product-negative-holder-moments}
			\begin{aligned}
			&\Bigg(
			\mathbb{E}\Bigg[
			\Bigg(
			\left\lVert
			\overline{r_{A^\varepsilon}^\varepsilon}
			\phi_{B^\varepsilon}^\varepsilon
			\right\rVert_{\mathcal{C}^{-\kappa,\varepsilon}(\rho^\kappa)}
			+
			\left\lVert
			\overline{r_{A^\varepsilon}^\varepsilon}
			\phi^\varepsilon
			\right\rVert_{\mathcal{C}^{-\kappa,\varepsilon}(\rho^\kappa)}
			+
			\left\lVert
			\rho_{A^\varepsilon}^\varepsilon
			\phi_{B^\varepsilon}^\varepsilon
			\right\rVert_{\mathcal{C}^{-\kappa,\varepsilon}(\rho^\kappa)}
			+
			\left\lVert
			\rho_{A^\varepsilon}^\varepsilon
			\phi^\varepsilon
			\right\rVert_{\mathcal{C}^{-\kappa,\varepsilon}(\rho^\kappa)}\\
			&\qquad+
			\left\lVert
			:\!\left\lvert
			r_{A^\varepsilon}^\varepsilon
			\right\rvert^2\!:_\varepsilon
			\phi_{B^\varepsilon}^\varepsilon
			\right\rVert_{\mathcal{C}^{-\kappa,\varepsilon}(\rho^\kappa)}
			+
			\left\lVert
			:\!\rho_{A^\varepsilon}^{\varepsilon *}
			\rho_{A^\varepsilon}^\varepsilon\!:_\varepsilon
			\phi_{B^\varepsilon}^\varepsilon
			\right\rVert_{\mathcal{C}^{-\kappa,\varepsilon}(\rho^\kappa)}
			\Bigg)^p
			\Bigg]
			\Bigg)^{1/p}
			\lesssim_\kappa
			p^{5/2}
			\left(
			\mathbb{E}\left[
			1+\lVert B^\varepsilon\rVert_{H^{1,\varepsilon}(\rho)}^{\kappa p}
			\right]
			\right)^{1/p}.
			\end{aligned}
		\end{equation}
	\end{lemma}

	\begin{proof}
		We may assume that $\kappa>0$ is small.
		By Lemma~\ref{lem:covariant-free-gff-difference-negative-holder-moment}
		with $\sigma=2\kappa$,
		Lemma~\ref{lem:nonlinear-link-coefficient-approximation-holder-moment},
		Corollary~\ref{cor:standard-gauge-field-stochastic-estimates},
		Lemma~\ref{lem:discrete-paraproduct-estimates},
		and H\"older's inequality, we obtain
		\begin{equation}
			\label{eq:mixed-gauge-higgs-product-preliminary-moment-bound}
			\text{LHS of~\eqref{eq:mixed-gauge-higgs-product-negative-holder-moments}}
			\lesssim_\kappa p^{5/2}
			\left(\mathbb E\left[
			\bigl(1+\lVert B^\varepsilon\rVert_{H^{1,\varepsilon}(\rho)}\bigr)^{100\kappa p}
			\right]\right)^{1/p}.
		\end{equation}
		Here all norms on the left are taken in
		$\mathcal C^{-\kappa,\varepsilon}(\rho^{100\kappa})$.
		Replacing $\kappa$ by $\kappa/100$ and using the weighted embeddings
		yields~\eqref{eq:mixed-gauge-higgs-product-negative-holder-moments}.
	\end{proof}

	\subsection{Estimate of the renormalized free mixed current}
	\label{subsec:stochastic-estimate-renormalized-free-current}

	We estimate the renormalized free mixed current appearing in
	\eqref{eq:gauge-irregular-Zhi1-unified-free-static-renormalized-target}, which is the elliptic analogue of the stochastic object in the third line of \cite[(4.26)]{BC24}.
	Set
	\begin{equation}
		\label{eq:renormalized-free-current-T-definition}
		T_j^{a,b}:=
		\mathbb E^\phi\left[
		i\overline{\tau_j^a\phi^\varepsilon}\,
		D_j^\varepsilon(\mathcal L_0^\varepsilon+1)^{-1}
		\sum_{\ell=1}^2D_\ell^{\varepsilon *}
		\bigl(r_{A,\ell}^{\varepsilon,(b)}
		\tau_\ell^b\phi^\varepsilon\bigr)\right].
	\end{equation}

	\begin{lemma}
		\label{lem:renormalized-free-current-estimate}
		For $0<\kappa\leq1$, $a,b\in\{0,1\}$, and $j\in\{1,2\}$,
		\begin{equation}
			\label{eq:renormalized-free-current-estimate}
			\left\|T_j^{a,b}-\frac i2b_\varepsilon
			r_{A,j}^{\varepsilon,(b)}\right\|_
			{\mathcal C^{-2\kappa,\varepsilon}(\rho^\kappa)}
			\lesssim_\kappa
			\|r_A^\varepsilon\|_{\mathcal C^{-\kappa,\varepsilon}(\rho^\kappa)}.
		\end{equation}
	\end{lemma}

	\begin{proof}
		The kernel of $(\mathcal L_0^\varepsilon+1)^{-1}$ is
		\begin{equation}
			\label{eq:renormalized-free-current-Green-kernel}
			G^\varepsilon(x):=
			\frac1{|\Lambda_{M,\varepsilon}|}
			\sum_{q\in Q_{M,\varepsilon}}
			\frac{\mathrm e_q(x)}{1+|q^\varepsilon|^2},
			\qquad
			\mathbb E^\phi[
			\overline{\phi^\varepsilon(x)}\phi^\varepsilon(z)]
			=2G^\varepsilon(x-z).
		\end{equation}
		Define
		\begin{equation}
			\label{eq:renormalized-free-current-kernel}
			K_{j,\ell}^{a,b}(y):=
			\bigl(D_\ell^{\varepsilon *}D_j^\varepsilon
			G^\varepsilon\bigr)(y)\,
			G^\varepsilon\bigl(y+\varepsilon(ae_j-be_\ell)\bigr).
		\end{equation}
		Periodic summation by parts and the covariance identity give
		\begin{equation}
			\label{eq:renormalized-free-current-kernel-contraction}
			\begin{aligned}
				T_j^{a,b}(x)
				&=i\sum_{\ell=1}^2\int_{\Lambda_{M,\varepsilon}}
				\bigl(D_\ell^{\varepsilon *}D_j^\varepsilon
				G^\varepsilon\bigr)(x-z)\,
				r_{A,\ell}^{\varepsilon,(b)}(z)\,
				\mathbb E^\phi\left[
				\overline{\phi^\varepsilon(x+a\varepsilon e_j)}
				\phi^\varepsilon(z+b\varepsilon e_\ell)\right]\,dz\\
				&=2i\sum_{\ell=1}^2
				\bigl(K_{j,\ell}^{a,b}*
				r_{A,\ell}^{\varepsilon,(b)}\bigr)(x),
			\end{aligned}
		\end{equation}
		where convolution is taken over $\Lambda_{M,\varepsilon}$.
		Subtracting the kernel mass yields
		\begin{equation}
			\label{eq:renormalized-free-current-kernel-compensation}
			\begin{aligned}
				T_j^{a,b}-\frac i2b_\varepsilon
				r_{A,j}^{\varepsilon,(b)}
				=2i\sum_{\ell=1}^2
				\int_{\Lambda_{M,\varepsilon}}K_{j,\ell}^{a,b}(y)
				\bigl[r_{A,\ell}^{\varepsilon,(b)}(\cdot-y)
				-r_{A,\ell}^{\varepsilon,(b)}\bigr]\,dy
				+2i\sum_{\ell=1}^2
				\left(\int_{\Lambda_{M,\varepsilon}}
				K_{j,\ell}^{a,b}(y)\,dy
				-\frac{b_\varepsilon}{4}\delta_{\ell j}\right)
				r_{A,\ell}^{\varepsilon,(b)}.
			\end{aligned}
		\end{equation}
		Lemma~\ref{lem:weighted-lattice-translation-increment} below,
		with $\alpha=-\kappa$ and $\delta=\kappa$, gives
		\begin{equation}
			\label{eq:renormalized-free-current-increment-bound}
			\left\|r_{A,\ell}^{\varepsilon,(b)}(\cdot-y)
			-r_{A,\ell}^{\varepsilon,(b)}\right\|_
			{\mathcal C^{-2\kappa,\varepsilon}(\rho^\kappa)}
			\lesssim_\kappa |y|^\kappa\rho(y)^{-\kappa}
			\|r_A^\varepsilon\|_{\mathcal C^{-\kappa,\varepsilon}(\rho^\kappa)}.
		\end{equation}
		Apply this estimate to the first term of
		\eqref{eq:renormalized-free-current-kernel-compensation}.
		Lemmas~\ref{lem:renormalized-free-current-kernel-moment}
		and~\ref{lem:renormalized-free-current-kernel-mass} below
		bound the first and second lines, respectively, and together yield
		\eqref{eq:renormalized-free-current-estimate}.
	\end{proof}

	\begin{lemma}[Weighted lattice translations]
		\label{lem:weighted-lattice-translation-increment}
		For $\alpha\in\mathbb R$, $\kappa>0$, $0<\delta\leq1$,
		$f\in\mathcal C^{\alpha,\varepsilon}(\rho^\kappa)$, and
		$y\in\Lambda_\varepsilon$,
		\begin{equation}
			\label{eq:weighted-lattice-translation-increment}
			\|f(\cdot-y)-f\|_
			{\mathcal C^{\alpha-\delta,\varepsilon}(\rho^\kappa)}
			\lesssim_{\alpha,\kappa,\delta}
			|y|^\delta\rho(y)^{-\kappa}
			\|f\|_{\mathcal C^{\alpha,\varepsilon}(\rho^\kappa)}.
		\end{equation}
	\end{lemma}

	\begin{proof}
		Put $f_m:=\Delta_m^\varepsilon f$. Let
		$\widetilde\varphi_m^\varepsilon
		:=
		\sum_{|n-m|\leq1}\varphi_n^\varepsilon$, where
		$-1\leq n\leq J^\varepsilon$, and set
		$\widetilde\eta_m^\varepsilon
		:=
		\mathcal F_\varepsilon^{-1}\widetilde\varphi_m^\varepsilon$.
		Then $\widetilde\varphi_m^\varepsilon=1$ on
		$\operatorname{supp}\varphi_m^\varepsilon$, so
		$f_m=\widetilde\eta_m^\varepsilon*f_m$ on the full lattice.

		As in Lemma~\ref{lem:weighted-periodic-LP-kernel}, applied
		to the enlarged kernels on the full lattice, we obtain
		the following bound; the derivative term is controlled by
		the $B_{1,1}^{1,\varepsilon}(\rho^{-\kappa})$ norm:
		\begin{equation}
			\label{eq:weighted-translation-enlarged-kernel-L1}
			\|\widetilde\eta_m^\varepsilon\|_{L^{1,\varepsilon}(\rho^{-\kappa})}
			+2^{-m}\|D^\varepsilon\widetilde\eta_m^\varepsilon
			\|_{L^{1,\varepsilon}(\rho^{-\kappa})}
			\lesssim_\kappa1.
		\end{equation}
		Let $y_0,\ldots,y_N$ be a shortest lattice path from $0$ to $y$
		in $\Lambda_\varepsilon$, with $y_0=0$, $y_N=y$, and adjacent
		vertices $y_\ell$ and $y_{\ell+1}$.
		Since $|y_\ell|\leq|y|$ and $\varepsilon N\lesssim|y|$, we have
		$\rho(\cdot)^{-\kappa}\lesssim_\kappa
		\rho(y)^{-\kappa}\rho(\cdot-y_\ell)^{-\kappa}$ and
		\begin{equation}
			\label{eq:weighted-translation-kernel-path-bound}
			\begin{aligned}
				\|\widetilde\eta_m^\varepsilon(\cdot-y)
				-\widetilde\eta_m^\varepsilon\|_{L^{1,\varepsilon}(\rho^{-\kappa})}
				&=\left\|\sum_{\ell=0}^{N-1}
				\bigl(\widetilde\eta_m^\varepsilon(\cdot-y_{\ell+1})
				-\widetilde\eta_m^\varepsilon(\cdot-y_\ell)\bigr)
				\right\|_{L^{1,\varepsilon}(\rho^{-\kappa})}\\
				&\lesssim_\kappa\varepsilon\sum_{\ell=0}^{N-1}
				\rho(y)^{-\kappa}\|D^\varepsilon\widetilde\eta_m^\varepsilon
				\|_{L^{1,\varepsilon}(\rho^{-\kappa})}\\
				&\lesssim_\kappa|y|\rho(y)^{-\kappa}
				\|D^\varepsilon\widetilde\eta_m^\varepsilon
				\|_{L^{1,\varepsilon}(\rho^{-\kappa})}.
			\end{aligned}
		\end{equation}
		The triangle inequality gives the bound
		$\rho(y)^{-\kappa}\|\widetilde\eta_m^\varepsilon
		\|_{L^{1,\varepsilon}(\rho^{-\kappa})}$.
		Combining these bounds with
		\eqref{eq:weighted-translation-enlarged-kernel-L1} yields
		\begin{equation}
			\label{eq:weighted-translation-kernel-increment}
			\|\widetilde\eta_m^\varepsilon(\cdot-y)
			-\widetilde\eta_m^\varepsilon\|_
			{L^{1,\varepsilon}(\rho^{-\kappa})}
			\lesssim_\kappa
			\rho(y)^{-\kappa}\min\{1,2^m|y|\}.
		\end{equation}
		Weighted Young's inequality therefore implies
		\begin{equation}
			\label{eq:weighted-translation-dyadic-increment}
			\|f_m(\cdot-y)-f_m\|_{L^{\infty,\varepsilon}(\rho^\kappa)}
			\lesssim_\kappa
			\rho(y)^{-\kappa}\min\{1,2^m|y|\}
			\|f_m\|_{L^{\infty,\varepsilon}(\rho^\kappa)}.
		\end{equation}
		Use $\min\{1,t\}\leq t^\delta$, multiply by
		$2^{m(\alpha-\delta)}$, and take the supremum over $m$.
		This yields~\eqref{eq:weighted-lattice-translation-increment}.
	\end{proof}

\begin{lemma}
	\label{lem:renormalized-free-current-kernel-moment}
	For $\kappa>0$, $a,b\in\{0,1\}$, and $j,\ell\in\{1,2\}$,
	\begin{equation}
		\label{eq:renormalized-free-current-kernel-moment}
		\int_{\Lambda_\varepsilon\cap[-\pi M,\pi M)^2}
		|K_{j,\ell}^{a,b}(y)|\,|y|^\kappa\rho(y)^{-\kappa}\,dy
		\lesssim_\kappa1.
	\end{equation}
\end{lemma}

\begin{proof}
	As in~\eqref{eq:pinned-gauge-dyadic-multiplier}--%
	\eqref{eq:pinned-gauge-dyadic-kernel-decay}, decompose the two
	Green kernels into dyadic blocks:
	\begin{equation}
		\label{eq:renormalized-free-current-dyadic-Green-kernels}
		G^\varepsilon=\sum_{m=-1}^{J^\varepsilon}G_m^\varepsilon,
		\qquad
		G_m^\varepsilon:=
		\mathcal F_{M,\varepsilon}^{-1}
		\left(\frac{\varphi_m^\varepsilon}{1+|k^\varepsilon|^2}\right).
	\end{equation}
	The multipliers of $G_m^\varepsilon$ and
	$D_\ell^{\varepsilon *}D_j^\varepsilon G_m^\varepsilon$
	are smooth and periodic, supported where $|\xi|\lesssim2^m$,
	and satisfy
	\begin{equation}
		\label{eq:renormalized-free-current-dyadic-symbol-bounds}
		\begin{aligned}
			\left|\partial_\xi^\gamma
			\frac{\varphi_m^\varepsilon(\xi)}{1+|\xi^\varepsilon|^2}
			\right|\lesssim_\gamma2^{-(2+|\gamma|)m},\quad
			\left|\partial_\xi^\gamma
			\frac{\varphi_m^\varepsilon(\xi)
				\xi_j^\varepsilon\overline{\xi_\ell^\varepsilon}}
			{1+|\xi^\varepsilon|^2}
			\right|\lesssim_\gamma2^{-|\gamma|m}.
		\end{aligned}
	\end{equation}
	Lemma~\ref{lem:periodic-multiplier-kernel-bound} therefore gives
	\begin{equation}
		\label{eq:renormalized-free-current-dyadic-kernel-bounds}
		|G_m^\varepsilon(y)|
		+2^{-2m}|D_\ell^{\varepsilon *}D_j^\varepsilon
		G_m^\varepsilon(y)|
		\lesssim_N(1+2^m d_M(y,0))^{-N}.
	\end{equation}
	Since $2^n\varepsilon\lesssim1$, the same decay bound holds
	for $G_n^\varepsilon(y+\varepsilon(ae_j-be_\ell))$.
	
	Choose $N>2+5\kappa$. Since $|y|=d_M(y,0)$, the preceding bounds imply
	\begin{equation}
		\label{eq:renormalized-free-current-dyadic-kernel-moment}
		\begin{aligned}
			&\int_{\Lambda_\varepsilon\cap[-\pi M,\pi M)^2}
			|D_\ell^{\varepsilon *}D_j^\varepsilon G_m^\varepsilon(y)|
			|G_n^\varepsilon(y+\varepsilon(ae_j-be_\ell))|\,
			|y|^\kappa\rho(y)^{-\kappa}\,dy\\
			&\qquad\lesssim_\kappa
			2^{2m}\int_{\Lambda_\varepsilon}
			\frac{|y|^\kappa\langle y\rangle^{4\kappa}}
			{(1+2^{m\vee n}|y|)^N}\,dy
			\lesssim_\kappa 2^{2m-(2+\kappa)(m\vee n)}.
		\end{aligned}
	\end{equation}
	Expanding both factors in
	\eqref{eq:renormalized-free-current-kernel} and splitting the
	scale sum into $n\leq m$ and $n>m$ gives
	\begin{equation}
		\label{eq:renormalized-free-current-dyadic-scale-sum}
		\begin{aligned}
			\sum_{m,n=-1}^{J^\varepsilon}
			2^{2m-(2+\kappa)(m\vee n)}
			&\leq
			\sum_{m=-1}^{\infty}(m+2)2^{-\kappa m}
			+\sum_{n=-1}^{\infty}2^{-(2+\kappa)n}
			\sum_{m=-1}^{n-1}2^{2m}\\
			&\lesssim
			\sum_{m=-1}^{\infty}(m+2)2^{-\kappa m}
			\lesssim_\kappa1.
		\end{aligned}
	\end{equation}
	This gives~\eqref{eq:renormalized-free-current-kernel-moment}.
\end{proof}

	\begin{lemma}
		\label{lem:renormalized-free-current-kernel-mass}
		For $a,b\in\{0,1\}$ and $j,\ell\in\{1,2\}$,
		\begin{equation}
			\label{eq:renormalized-free-current-zero-mode-bound}
			\left|\int_{\Lambda_{M,\varepsilon}}K_{j,\ell}^{a,b}(y)\,dy
			-\frac{b_\varepsilon}{4}\delta_{\ell j}\right|\lesssim1.
		\end{equation}
	\end{lemma}

	\begin{proof}
		The Fourier representation of $G^\varepsilon$ gives
		\begin{equation}
			\label{eq:renormalized-free-current-kernel-mass-Fourier}
			\int_{\Lambda_{M,\varepsilon}}K_{j,\ell}^{a,b}(y)\,dy
			=\frac1{|\Lambda_{M,\varepsilon}|}
			\sum_{q\in Q_{M,\varepsilon}}
			\frac{e^{i\varepsilon(bq_\ell-aq_j)}
			q_j^\varepsilon\overline{q_\ell^\varepsilon}}
			{(1+|q^\varepsilon|^2)^2}.
		\end{equation}
		Since $|e^{i\varepsilon(bq_\ell-aq_j)}-1|
		\lesssim\varepsilon|q|$ and $|q^\varepsilon|\asymp|q|$
		on $Q_{M,\varepsilon}$,
		\begin{equation}
			\label{eq:renormalized-free-current-phase-error-bound}
			\left|\int_{\Lambda_{M,\varepsilon}}
			(K_{j,\ell}^{a,b}-K_{j,\ell}^{0,0})(y)\,dy\right|
			\lesssim
			\frac{\varepsilon}{|\Lambda_{M,\varepsilon}|}
			\sum_{q\in Q_{M,\varepsilon}}\frac{|q|}{1+|q|^2}
			\lesssim1.
		\end{equation}
		It remains to consider $a=b=0$.
		For $j\ne\ell$, write
		\begin{equation}
			\label{eq:lattice-frequency-real-imaginary-decomposition}
			q_j^\varepsilon
			=\frac{\sin(\varepsilon q_j)}{\varepsilon}
			+i\frac{1-\cos(\varepsilon q_j)}{\varepsilon},
			\qquad
			\left|\frac{1-\cos(\varepsilon q_j)}{\varepsilon}\right|
			\leq\frac{\varepsilon}{2}|q_j|^2.
		\end{equation}
		The admissible frequency set is invariant under reflection of
		each coordinate. Thus all terms containing a sine factor vanish
		in the sum, and
		\begin{equation}
			\label{eq:renormalized-free-current-off-diagonal-zero-mode-bound}
			\left|\int_{\Lambda_{M,\varepsilon}}K_{j,\ell}^{0,0}(y)\,dy\right|
			\lesssim
			\frac{\varepsilon^2}{|\Lambda_{M,\varepsilon}|}
			\sum_{q\in Q_{M,\varepsilon}}\frac{|q|^4}{(1+|q|^2)^2}
			\lesssim1.
		\end{equation}
		For $j=\ell$, recall that
$
			\label{eq:lattice-b-epsilon-fourier-sum}
			b_\varepsilon
			=\frac2{|\Lambda_{M,\varepsilon}|}
			\sum_{q\in Q_{M,\varepsilon}}
			\frac1{1+|q^\varepsilon|^2}.
$
		Symmetry between the coordinate directions gives
		\begin{equation}
			\label{eq:renormalized-free-current-diagonal-zero-mode-bound}
			\begin{aligned}
				\int_{\Lambda_{M,\varepsilon}}K_{j,j}^{0,0}(y)\,dy
				-\frac{b_\varepsilon}{4}
				&=\frac1{2|\Lambda_{M,\varepsilon}|}
				\sum_{q\in Q_{M,\varepsilon}}
				\left(\frac{|q^\varepsilon|^2}{(1+|q^\varepsilon|^2)^2}
				-\frac1{1+|q^\varepsilon|^2}\right)\\
				&=-\frac1{2|\Lambda_{M,\varepsilon}|}
				\sum_{q\in Q_{M,\varepsilon}}
				\frac1{(1+|q^\varepsilon|^2)^2}.
			\end{aligned}
		\end{equation}
		The last normalized sum is uniformly bounded.
		Together with~\eqref{eq:renormalized-free-current-phase-error-bound},
		this yields ~\eqref{eq:renormalized-free-current-zero-mode-bound}.
	\end{proof}

	\subsection{\texorpdfstring{Estimate of
			$\mathbb{E}^\phi[\operatorname{Im}[\overline{\phi_B^\varepsilon}\,
			D_B^\varepsilon\phi_B^\varepsilon]]$}
		{Estimate of the mean covariant current}}
	\label{subsec:stochastic-estimate-mean-covariant-current}

	\begin{proposition}
		\label{lem:mean-covariant-current-logarithmic-bound}
		For every $B\in(\mathbb{R}^2)^{\Lambda_{M,\varepsilon}}$,
		\begin{equation}
			\label{eq:mean-covariant-current-logarithmic-bound}
			\left\lVert\mathbb{E}^\phi\left[
				\operatorname{Im}\left[\overline{\phi_B^\varepsilon}\,
				D_B^\varepsilon\phi_B^\varepsilon\right]\right]
			\right\rVert_{H^{-1,\varepsilon}(\rho)}
			\lesssim\log\left(2+\lVert B\rVert_{H^{1,\varepsilon}(\rho)}\right).
		\end{equation}
	\end{proposition}

	\begin{proof}
		\noindent\textit{Step~1: Calculation of the expectation.}
		Set $w^M(y):=w_0^M(y)=\sum_{n\in\mathbb{Z}^2}\rho(y+2\pi Mn)$,
		as in~\eqref{eq:periodized-polynomial-weight}.
		\label{eq:mean-current-periodized-weight}
		By duality and Lemma~\ref{lem:periodized-weighted-Besov-equivalences}\textup{(i), (iii)},
		it suffices to prove that, for every
		$h\in(\mathbb{R}^2)^{\Lambda_{M,\varepsilon}}$,
		\begin{equation}
			\label{eq:mean-covariant-current-dual-target}
			\left\lvert\left\langle
				\mathbb{E}^\phi\left[\operatorname{Im}\left[
					\overline{\phi_B^\varepsilon}\,D_B^\varepsilon\phi_B^\varepsilon\right]\right],
				h\right\rangle_{M,\varepsilon}\right\rvert
			\lesssim\log\left(2+\lVert B\rVert_{H^{1,\varepsilon}(\rho)}\right)
				\lVert h\rVert_{H^{1,\varepsilon}((w^M)^{-1};\Lambda_{M,\varepsilon})}.
		\end{equation}

		By~\eqref{eq:gauge-gradient-effective-interaction-current}, we have
		\begin{equation}
			\label{eq:mean-covariant-current-trace-identity}
			\begin{aligned}
				\left\langle\mathbb{E}^\phi\left[\operatorname{Im}\left[
					\overline{\phi_B^\varepsilon}\,D_B^\varepsilon\phi_B^\varepsilon\right]\right],
					h\right\rangle_{M,\varepsilon}
				&=\frac12\left.\frac{d}{d\delta}\right\rvert_{\delta=0}
					\mathbb{E}^\phi\left[
					\lVert D_{B+\delta h}^\varepsilon\phi_B^\varepsilon
						\rVert_{L^{2,\varepsilon}(\Lambda_{M,\varepsilon})}^2\right]\\
				&=\frac12\left.\frac{d}{d\delta}\right\rvert_{\delta=0}
					\mathbb{E}^\phi\left[
					\left\langle\mathcal{L}_{B+\delta h}^\varepsilon\phi_B^\varepsilon,
						\phi_B^\varepsilon\right\rangle_{M,\varepsilon}\right]\\
				&=\frac12\left.\frac{d}{d\delta}\right\rvert_{\delta=0}
					\operatorname{Tr}_{\mathbb{R}}^{M}\left(
					\mathcal{L}_{B+\delta h}^\varepsilon
					(\mathcal{L}_B^\varepsilon+1)^{-1}\right)\\
				&=\left.\frac{d}{d\delta}\right\rvert_{\delta=0}
					\operatorname{Tr}_{\mathbb{C}}^{M}\left(
					\mathcal{L}_{B+\delta h}^\varepsilon
					(\mathcal{L}_B^\varepsilon+1)^{-1}\right).
			\end{aligned}
		\end{equation}
		Here $\operatorname{Tr}_{\mathbb{R}}^{M}$ and
		$\operatorname{Tr}_{\mathbb{C}}^{M}$ are the traces on
		$L^{2,\varepsilon}(\Lambda_{M,\varepsilon})$, viewed as a real or complex
		Hilbert space, respectively.
		The real covariance operator of $\phi_B^\varepsilon$ is
		$(\mathcal{L}_B^\varepsilon+1)^{-1}$.
		For the last equality in~\eqref{eq:mean-covariant-current-trace-identity},
		we used $\operatorname{Tr}_{\mathbb{R}}^{M}T
		=2\operatorname{Re}\operatorname{Tr}_{\mathbb{C}}^{M}T$
		for a complex-linear operator $T$.

		Writing
		$D_h\operatorname{Tr}_{\mathbb{C}}^{M}(e^{-t\mathcal{L}_B^\varepsilon})
		:=\left.\frac{d}{d\delta}\right\rvert_{\delta=0}
		\operatorname{Tr}_{\mathbb{C}}^{M}(e^{-t\mathcal{L}_{B+\delta h}^\varepsilon})$,
		\label{eq:mean-current-heat-trace-directional-derivative}
		the heat-kernel representation of the resolvent and the identity
		$\frac{d}{d\delta}\operatorname{Tr}(e^{A(\delta)})
		=\operatorname{Tr}(A'(\delta)e^{A(\delta)})$ give

		\begin{equation}
			\label{eq:mean-covariant-current-resolvent-heat-trace}
			\begin{aligned}
				\left.\frac{d}{d\delta}\right\rvert_{\delta=0}
					\operatorname{Tr}_{\mathbb{C}}^{M}\left(
						\mathcal{L}_{B+\delta h}^\varepsilon
						\int_0^\infty e^{-t}e^{-t\mathcal{L}_B^\varepsilon}\,dt\right)=-\int_0^\infty\frac{e^{-t}}t\,
					D_h\operatorname{Tr}_{\mathbb{C}}^{M}
					(e^{-t\mathcal{L}_B^\varepsilon})\,dt.
			\end{aligned}
		\end{equation}

		\smallskip
		\noindent\textit{Step~2: Explicit representation of the heat-trace
		directional derivative.}
		With $X$ the random walk on $\Lambda_{M,\varepsilon}$ and $I_B^X(t)$
		as in Lemma~\ref{lem:FK-LB}, set
		$\mathbb{E}_{x,x}[\mathord\cdot]
		=\mathbb{E}_x[\mathord\cdot\mid X_t=x]$.
		By the same lemma and the self-adjointness of $\mathcal{L}_B^\varepsilon$,
		\begin{equation}
			\label{eq:mean-covariant-current-heat-trace-FK}
			\operatorname{Tr}_{\mathbb{C}}^{M}
				(e^{-t\mathcal{L}_B^\varepsilon})
			=p_t^{M,\varepsilon}(0)\int_{\Lambda_{M,\varepsilon}}
				\mathbb{E}_{x,x}\left[\cos\bigl(I_B^X(t)\bigr)\right]\,dx.
		\end{equation}
		Since $I_{B+\delta h}^X(t)=I_B^X(t)+\delta I_h^X(t)$, differentiation in
		\eqref{eq:mean-covariant-current-heat-trace-FK} yields
		\begin{equation}
			\label{eq:mean-covariant-current-heat-trace-derivative-FK}
			D_h\operatorname{Tr}_{\mathbb{C}}^{M}
				(e^{-t\mathcal{L}_B^\varepsilon})
			=-p_t^{M,\varepsilon}(0)\int_{\Lambda_{M,\varepsilon}}
				\mathbb{E}_{x,x}\left[\sin\bigl(I_B^X(t)\bigr)I_h^X(t)\right]\,dx.
		\end{equation}

		\smallskip
		\noindent\textit{Step~3: Estimate of the heat-trace derivative.}
		For a periodic real vector field $f$ and $w=(w^M)^\alpha$, define
		\begin{equation}
			\label{eq:magnetic-loop-phase-second-moment-definition}
			\mathcal{V}_t^w(f):=
			p_t^{M,\varepsilon}(0)\int_{\Lambda_{M,\varepsilon}}
				 w(x)^2\mathbb{E}_{x,x}\left[\lvert I_f^X(t)\rvert^2\right]\,dx.
		\end{equation}
		Using respectively $\lvert\sin\theta\rvert\leq1$ and
		$\lvert\sin\theta\rvert\leq\lvert\theta\rvert$, followed by
		Cauchy--Schwarz with weights $w^M$ and $(w^M)^{-1}$, we obtain
		\begin{equation}
			\label{eq:mean-covariant-current-two-heat-trace-bounds}
			\left\lvert D_h\operatorname{Tr}_{\mathbb{C}}^{M}
				(e^{-t\mathcal{L}_B^\varepsilon})\right\rvert
			\lesssim
			\min\left\{p_t^{M,\varepsilon}(0)^{1/2},
				\mathcal{V}_t^{w^M}(B)^{1/2}\right\}
				(\mathcal{V}_t^{(w^M)^{-1}}(h))^{1/2}.
		\end{equation}
		By Lemma~\ref{lem:magnetic-loop-phase-second-moment} below,
		with $a=1,-1$, and
		\eqref{eq:free-lattice-heat-kernel-diagonal-estimate},
		\begin{equation}
			\label{eq:mean-covariant-current-piecewise-heat-trace-bound}
			\begin{aligned}
				\left\lvert D_h\operatorname{Tr}_{\mathbb{C}}^{M}
					(e^{-t\mathcal{L}_B^\varepsilon})\right\rvert\lesssim
					\lVert h\rVert_{H^{1,\varepsilon}((w^M)^{-1};\Lambda_{M,\varepsilon})}
					\begin{cases}
						\min\{1,t\lVert B\rVert_{H^{1,\varepsilon}(\rho)}\},&0<t\leq1,\\
						t^{1/2}(1+t)^{R/2},&t\geq1.
					\end{cases}
			\end{aligned}
		\end{equation}
		Substituting this estimate into
		\eqref{eq:mean-covariant-current-resolvent-heat-trace}, we obtain
		\begin{equation}
			\label{eq:mean-covariant-current-logarithmic-integral-bound}
			\begin{aligned}
				&\left\lvert\left\langle
					\mathbb{E}^\phi\left[\operatorname{Im}\left[
						\overline{\phi_B^\varepsilon}\,D_B^\varepsilon\phi_B^\varepsilon\right]\right],
					h\right\rangle_{M,\varepsilon}\right\rvert\\
				&\quad\lesssim
					\left(\int_0^1\frac{\min\{1,t\lVert B\rVert_{H^{1,\varepsilon}(\rho)}\}}t\,dt
						+\int_1^\infty e^{-t}t^{-1/2}(1+t)^{R/2}\,dt\right)
					\lVert h\rVert_{H^{1,\varepsilon}((w^M)^{-1};\Lambda_{M,\varepsilon})}\\
				&\quad\lesssim
					\log\left(2+\lVert B\rVert_{H^{1,\varepsilon}(\rho)}\right)
					\lVert h\rVert_{H^{1,\varepsilon}((w^M)^{-1};\Lambda_{M,\varepsilon})}.
			\end{aligned}
		\end{equation}
		Together with~\eqref{eq:mean-covariant-current-trace-identity}, this yields
		\eqref{eq:mean-covariant-current-dual-target}.
	\end{proof}

	\begin{lemma}
		\label{lem:magnetic-loop-phase-second-moment}
		Let $a\in\mathbb{R}$ and $w=(w^M)^a$.
		There exists $R>0$, depending only on $a$, such that, for every
		$h\in(\mathbb{R}^2)^{\Lambda_{M,\varepsilon}}$ and $t>0$,
		the quantity $\mathcal{V}_t^w(h)$ defined in
		\eqref{eq:magnetic-loop-phase-second-moment-definition} satisfies
		\begin{equation}
			\label{eq:magnetic-loop-phase-second-moment}
			\mathcal{V}_t^w(h)\lesssim_a
				t(1+t)^R\lVert h\rVert_{H^{1,\varepsilon}(w;\Lambda_{M,\varepsilon})}^2.
		\end{equation}
	\end{lemma}

	\begin{proof}
		\noindent\textit{Step~1: Midpoint argument.}
		We obtain the small-time gain by constructing, for each
		$x\in\Lambda_{M,\varepsilon}$, a one-form $g_x$ which approximates
		$h-h(x)$ near $x$. Recall that $I_h^X(t)$ can be viewed as the discrete
		line integral of the one-form $h$ along $X$, which is a loop on
		$\{X_t=X_0=x\}$. To ensure that $I_h^X(t)=I_{g_x}^X(t)$ on this event,
		we seek a $2\pi M$-periodic scalar function $\psi_x$ such that
		\begin{equation}
			\label{eq:small-time-corrector-gradient-relation}
			h-g_x=D^\varepsilon\psi_x.
		\end{equation}
		\begin{remark}
			The naive choice $g_x=h-h(x)$ does not work on the torus, since the
			loop integral of the constant one-form $h(x)$ need not vanish.
		\end{remark}
		By Lemma~\ref{lem:periodic-affine-corrector} below, for every
		$c\in\mathbb{R}^2$ there exists a $2\pi M$-periodic scalar function
		$\psi_{x,c}$ satisfying~\eqref{eq:periodic-affine-corrector-properties}.
		Define $g_x:=h-D^\varepsilon\psi_{x,h(x)}$.
		Since $\psi_{x,h(x)}$ is $2\pi M$-periodic,
		\begin{equation}
			\label{eq:gradient-one-form-loop-integral}
			I_{D^\varepsilon\psi_{x,h(x)}}^X(t)
				=\psi_{x,h(x)}(X_t)-\psi_{x,h(x)}(x)=0
				\quad\text{on }\{X_t=x\}.
		\end{equation}
		Thus $I_h^X(t)=I_{g_x}^X(t)$ on this event.
		Using $\mathbb{P}_x(X_t=x)=\varepsilon^2p_t^{M,\varepsilon}(0)$,
		we first write
		\begin{equation}
			\label{eq:magnetic-loop-second-moment-return-loop}
			\mathcal{V}_t^w(h)
				=\sum_{x\in\Lambda_{M,\varepsilon}}
					w(x)^2\mathbb{E}_x\left[\lvert I_{g_x}^X(t)\rvert^2;X_t=x\right].
		\end{equation}
		Set $s=t/2$. By time reversal of the random walk bridge, we obtain
		\begin{equation}
			\label{eq:magnetic-loop-midpoint-time-reversal}
			\begin{aligned}
				\mathcal{V}_t^w(h)
				&\leq2\sum_{x\in\Lambda_{M,\varepsilon}}w(x)^2\Bigl(
					\mathbb{E}_x\left[\lvert I_{g_x}^X(s)\rvert^2;X_{2s}=x\right]
					+\mathbb{E}_x\left[
					\lvert I_{g_x}^X(2s)-I_{g_x}^X(s)\rvert^2;X_{2s}=x\right]\Bigr)\\
				&=4\sum_{x\in\Lambda_{M,\varepsilon}}
						w(x)^2\mathbb{E}_x\left[\lvert I_{g_x}^X(s)\rvert^2;X_{2s}=x\right].
			\end{aligned}
		\end{equation}
		By the Markov property and~\eqref{eq:free-random-walk-kernel-mass},
		\begin{equation}
			\label{eq:magnetic-loop-midpoint-Markov-property}
			\begin{aligned}
				\mathbb{E}_x\left[\lvert I_{g_x}^X(s)\rvert^2;X_{2s}=x\right]
				=\mathbb{E}_x\left[\lvert I_{g_x}^X(s)\rvert^2
					\mathbb{P}_x(X_{2s}=x\mid\mathcal{F}_s)\right]
				=\varepsilon^2\mathbb{E}_x\left[
					\lvert I_{g_x}^X(s)\rvert^2p_s^{M,\varepsilon}(x;X_s)\right].
			\end{aligned}
		\end{equation}
		By~\eqref{eq:free-lattice-heat-kernel},
		$\lvert p_s^{M,\varepsilon}(x;y)\rvert\leq p_s^{M,\varepsilon}(0)$.
		Consequently,
		\begin{equation}
			\label{eq:magnetic-loop-midpoint-bound}
			\begin{aligned}
				\mathcal{V}_t^w(h)
				&\leq4\varepsilon^2p_s^{M,\varepsilon}(0)
					\sum_{x\in\Lambda_{M,\varepsilon}}
					w(x)^2\mathbb{E}_x\left[\lvert I_{g_x}^X(s)\rvert^2\right]\\
				&=4p_s^{M,\varepsilon}(0)
						\int_{\Lambda_{M,\varepsilon}}
						w(x)^2\mathbb{E}_x\left[\lvert I_{g_x}^X(s)\rvert^2\right]\,dx.
			\end{aligned}
		\end{equation}

		\smallskip
		\noindent\textit{Step~2: Calculation of
		$\mathbb{E}_x[\lvert I_g^X(s)\rvert^2]$.}
		By~\eqref{eq:magnetic-phase-functional}, the processes
		$\widetilde N_j^\pm(r):=N_j^\pm(r)-\varepsilon^{-2}r$ are martingales.
		Hence, for every periodic real vector field $g$,
		\begin{equation}
			\label{eq:path-phase-drift-martingale-decomposition}
			\begin{aligned}
				&I_g^X(s)
				=\int_0^s\operatorname{div}^\varepsilon g(X_r)\,dr+M_s^g,\\
				&M_s^g
				=\varepsilon\sum_{j=1}^2\left(
					\int_0^s g_j(X_{r-})\,d\widetilde N_j^+(r)
					-\int_0^s g_j(X_{r-}-\varepsilon e_j)\,
						d\widetilde N_j^-(r)\right).
			\end{aligned}
		\end{equation}
		Its predictable quadratic variation is
		\begin{equation}
			\label{eq:path-phase-martingale-bracket}
			\langle M^g\rangle_s
			=\int_0^s\sum_{j=1}^2\left(
				\lvert g_j(X_r)\rvert^2
				+\lvert g_j(X_r-\varepsilon e_j)\rvert^2\right)dr
			:=\int_0^s\Gamma_g^\varepsilon(X_r)\,dr.
		\end{equation}
		Here we used
		$\langle\widetilde N_j^\pm\rangle_r=\varepsilon^{-2}r$
		and $g_j(X_{r-})=g_j(X_r)$ for Lebesgue-a.e.\ $r$.
		Consequently, Cauchy--Schwarz and the martingale isometry give
		\begin{equation}
			\label{eq:path-phase-unconditioned-second-moment}
			\mathbb{E}_x\lvert I_g^X(s)\rvert^2
			\lesssim s\int_0^s\mathbb{E}_x\left[
				\lvert\operatorname{div}^\varepsilon g(X_r)\rvert^2\right]dr
				+\int_0^s\mathbb{E}_x\left[\Gamma_g^\varepsilon(X_r)\right]dr.
		\end{equation}
		Applying this estimate with $g=g_x$ in
		\eqref{eq:magnetic-loop-midpoint-bound}, we obtain
		\begin{equation}
			\label{eq:magnetic-loop-small-time-preliminary-bound}
			\begin{aligned}
				\mathcal{V}_t^w(h)
				\lesssim p_s^{M,\varepsilon}(0)\Biggl(
					s\int_0^s\int_{\Lambda_{M,\varepsilon}}w(x)^2
					\mathbb{E}_x\left[
					\lvert\operatorname{div}^\varepsilon g_x(X_r)\rvert^2\right]dx\,dr+\int_0^s\int_{\Lambda_{M,\varepsilon}}w(x)^2
					\mathbb{E}_x\left[\Gamma_{g_x}^\varepsilon(X_r)\right]dx\,dr
				\Biggr).
			\end{aligned}
		\end{equation}

		\smallskip
		\noindent\textit{Step~3: Preparatory estimates.}
		Set $L:=4\lvert a\rvert+1$. We prove that, for every
		$z\in\Lambda_\varepsilon$,
		\begin{equation}
			\label{eq:weighted-loop-spatial-integrand-bounds}
			\begin{aligned}
				&\int_{\Lambda_{M,\varepsilon}}w(x)^2\Gamma_{g_x}^\varepsilon(x+z)\,dx
				\lesssim_a(\lvert z\rvert^2+\varepsilon^2)\langle z\rangle^{2L}
					\lVert h\rVert_{H^{1,\varepsilon}(w;\Lambda_{M,\varepsilon})}^2,\\
				&\int_{\Lambda_{M,\varepsilon}}w(x)^2
					\lvert\operatorname{div}^\varepsilon g_x(x+z)\rvert^2\,dx
				\lesssim_a\langle z\rangle^{2L}
					\lVert h\rVert_{H^{1,\varepsilon}(w;\Lambda_{M,\varepsilon})}^2.
			\end{aligned}
		\end{equation}
		By~\eqref{eq:periodic-affine-corrector-properties}, we have
		\begin{equation}
			\label{eq:small-time-gx-pointwise-decomposition}
			\lvert g_x(y)\rvert
			\lesssim \lvert h(y)-h(x)\rvert
				+\lvert h(x)\rvert d_M(x,y).
		\end{equation}
		\begin{equation}
			\label{eq:small-time-gx-reverse-edge-bound}
			\begin{aligned}
				\lvert g_{x,j}(y-\varepsilon e_j)\rvert
				&\lesssim \lvert h(y)-h(x)\rvert
					+\lvert h(x)\rvert d_M(x,y)
				+\varepsilon\left(
					\lvert D_j^{\varepsilon *}h_j(y)\rvert
					+\lvert h(x)\rvert\right).
			\end{aligned}
		\end{equation}
		\begin{equation}
			\label{eq:small-time-gx-divergence-bound}
			\lvert\operatorname{div}^\varepsilon g_x(y)\rvert
			\lesssim\lvert\operatorname{div}^\varepsilon h(y)\rvert
				+\lvert h(x)\rvert d_M(x,y).
		\end{equation}
		By the definition of $\Gamma_{g_x}^\varepsilon$, we obtain
		\begin{equation}
			\label{eq:weighted-loop-spatial-bracket-bound}
			\begin{aligned}
				\int_{\Lambda_{M,\varepsilon}}w(x)^2
					\Gamma_{g_x}^\varepsilon(x+z)\,dx
				&\lesssim\int_{\Lambda_{M,\varepsilon}}w(x)^2
					\lvert h(x+z)-h(x)\rvert^2\,dx\\
				&\quad+(\lvert z\rvert^2+\varepsilon^2)
					\lVert wh\rVert_{L^{2,\varepsilon}}^2
					+\varepsilon^2\sum_{j=1}^2\int_{\Lambda_{M,\varepsilon}}
					w(x)^2\lvert D_j^{\varepsilon *}h_j(x+z)\rvert^2\,dx\\
				&\lesssim_a
					\lvert z\rvert^2\langle z\rangle^{8\lvert a\rvert}
					\lVert wD^\varepsilon h\rVert_{L^{2,\varepsilon}}^2
					+\lvert z\rvert^2\lVert wh\rVert_{L^{2,\varepsilon}}^2
					+\varepsilon^2\langle z\rangle^{8\lvert a\rvert}
					\lVert wD^\varepsilon h\rVert_{L^{2,\varepsilon}}^2
					+\varepsilon^2\lVert wh\rVert_{L^{2,\varepsilon}}^2.
			\end{aligned}
		\end{equation}
		This gives the first bound in~\eqref{eq:weighted-loop-spatial-integrand-bounds}.
		Here we used 
$
			\label{eq:loop-periodic-weight-ratios}
			w(x+z)\lesssim_a\langle z\rangle^{4\lvert a\rvert}w(x).
$
		For the first term in~\eqref{eq:weighted-loop-spatial-bracket-bound},
		we used a similar argument to that in~\eqref{eq:weighted-translation-kernel-path-bound}.
		
		For the divergence term, we obtain
		\begin{equation}
			\label{eq:weighted-loop-spatial-divergence-bound}
			\begin{aligned}
				\int_{\Lambda_{M,\varepsilon}}w(x)^2
					\lvert\operatorname{div}^\varepsilon g_x(x+z)\rvert^2\,dx
				&\lesssim
					\int_{\Lambda_{M,\varepsilon}}w(x)^2
					\lvert\operatorname{div}^\varepsilon h(x+z)\rvert^2\,dx
					+\lvert z\rvert^2\lVert wh\rVert_{L^{2,\varepsilon}}^2\\
				&\lesssim_a
					\left(\langle z\rangle^{8\lvert a\rvert}
					+\lvert z\rvert^2\right)
					\lVert h\rVert_{H^{1,\varepsilon}
					(w;\Lambda_{M,\varepsilon})}^2.
			\end{aligned}
		\end{equation}
		This gives the second bound in~\eqref{eq:weighted-loop-spatial-integrand-bounds}.

		\smallskip
		\noindent\textit{Step~4: Conclusion.}
		Let $\widetilde{X}$ be the full-lattice lift of $X$,
		and set $m:=\lceil L+1\rceil$.
		By periodicity, the two terms in
		\eqref{eq:magnetic-loop-small-time-preliminary-bound} can be written
		with $x+\widetilde{X}_r$ under $\mathbb{P}_0$.
		Since $m\geq L+1$,~\eqref{eq:weighted-loop-spatial-divergence-bound},
		Lemma~\ref{lem:lifted-random-walk-moment}, and
		\eqref{eq:free-lattice-heat-kernel-diagonal-estimate} yield for the first term
		\begin{equation}
			\label{eq:weighted-loop-displacement-moments}
			\begin{aligned}
				p_s^{M,\varepsilon}(0)s\int_{\Lambda_{M,\varepsilon}}w(x)^2
					\int_0^s\mathbb{E}_0\left[
					\lvert\operatorname{div}^\varepsilon g_x(x+\widetilde{X}_r)\rvert^2
					\right]dr\,dx
				&\lesssim_a p_s^{M,\varepsilon}(0)s
					\int_0^s\mathbb{E}_0\langle\widetilde{X}_r\rangle^{2L}\,dr\,
					\lVert h\rVert_{H^{1,\varepsilon}(w;\Lambda_{M,\varepsilon})}^2\\
				&\lesssim_a p_s^{M,\varepsilon}(0)s^2(1+s)^m
					\lVert h\rVert_{H^{1,\varepsilon}(w;\Lambda_{M,\varepsilon})}^2\\
				&\lesssim_a s(1+s)^{m+1}
					\lVert h\rVert_{H^{1,\varepsilon}(w;\Lambda_{M,\varepsilon})}^2.
			\end{aligned}
			\csname ltx@label\endcsname{eq:weighted-loop-drift-time-integral-bound}
		\end{equation}
		For the second term, using~\eqref{eq:weighted-loop-spatial-integrand-bounds},
		Lemma~\ref{lem:lifted-random-walk-moment} at orders $1$ and $m$, and
		\eqref{eq:free-lattice-heat-kernel-diagonal-estimate}, we obtain
		\begin{equation}
			\label{eq:weighted-loop-averaged-integrand-bounds}
			\begin{aligned}
				&p_s^{M,\varepsilon}(0)\int_{\Lambda_{M,\varepsilon}}w(x)^2
					\int_0^s\mathbb{E}_0\left[
					\Gamma_{g_x}^\varepsilon(x+\widetilde{X}_r)\right]dr\,dx\\
				&\quad\lesssim_a
					\left(1+\frac1{s\vee\varepsilon^2}\right)
					\int_0^s(r+\varepsilon^2)(1+r)^m\,dr\,
					\lVert h\rVert_{H^{1,\varepsilon}(w;\Lambda_{M,\varepsilon})}^2.
			\end{aligned}
		\end{equation}
		Finally,
		\begin{equation}
			\label{eq:weighted-loop-midpoint-heat-kernel-factor}
			\left(1+\frac1{s\vee\varepsilon^2}\right)
			\int_0^s(r+\varepsilon^2)(1+r)^m\,dr
			\lesssim s\left(1+\frac1{s\vee\varepsilon^2}\right)
			(s\vee\varepsilon^2)(1+s)^m
			\lesssim s(1+s)^{m+1}.
		\end{equation}
		Combining these estimates in~\eqref{eq:magnetic-loop-small-time-preliminary-bound}
		yields~\eqref{eq:magnetic-loop-phase-second-moment} with $R=m+1$,
		since $s=t/2$ and $m$ depends only on $a$.
	\end{proof}

	\begin{lemma}
		\label{lem:periodic-affine-corrector}
		For every $x\in\Lambda_{M,\varepsilon}$ and $c\in\mathbb{R}^2$,
		there exists a $2\pi M$-periodic scalar function
		$\psi_{x,c}\colon\Lambda_\varepsilon\to\mathbb{R}$ such that
		\begin{equation}
			\label{eq:periodic-affine-corrector-properties}
			\begin{aligned}
				D^\varepsilon\psi_{x,c}(x)=c,\quad
				\lvert D^\varepsilon\psi_{x,c}(y)-c\rvert
					+\lvert\operatorname{div}^\varepsilon D^\varepsilon\psi_{x,c}(y)\rvert
					\lesssim \lvert c\rvert d_M(x,y).
			\end{aligned}
		\end{equation}
	\end{lemma}

	\begin{proof}
		Define
		\begin{equation}
			\label{eq:periodic-affine-corrector-definition}
			\psi_{x,c}(y):=
				\sum_{j=1}^2c_j\frac{\varepsilon}{\sin(\varepsilon/M)}
					\sin\left(\frac{y_j-x_j}{M}\right).
		\end{equation}
		Then $\psi_{x,c}$ is $2\pi M$-periodic. Moreover,
		\begin{equation}
			\label{eq:periodic-affine-corrector-forward-difference}
			D_j^\varepsilon\psi_{x,c}(y)
				=c_j\frac{\sin((y_j-x_j+\varepsilon)/M)
										-\sin((y_j-x_j)/M)}{\sin(\varepsilon/M)}.
		\end{equation}
		In particular, $D^\varepsilon\psi_{x,c}(x)=c$.
		Since the derivative of the quotient in
		\eqref{eq:periodic-affine-corrector-forward-difference}
		is bounded by $C/M$ for $0<\varepsilon/M<1$, we have
		\begin{equation}
			\label{eq:periodic-affine-corrector-gradient-bound}
			\lvert D_j^\varepsilon\psi_{x,c}(y)-c_j\rvert
				\lesssim \lvert c_j\rvert M^{-1}d_M(x,y).
		\end{equation}
		On the other hand,
		\begin{equation}
			\label{eq:periodic-affine-corrector-second-difference-bound}
			D_j^{\varepsilon *}D_j^\varepsilon\psi_{x,c}(y)
				=\frac{2c_j(1-\cos(\varepsilon/M))}
					{\varepsilon\sin(\varepsilon/M)}
					\sin\left(\frac{y_j-x_j}{M}\right).
		\end{equation}
		Its absolute value is bounded by $C\lvert c_j\rvert M^{-2}d_M(x,y)$.
		Therefore, summing over $j$ and using $M\geq1$ yields
		\eqref{eq:periodic-affine-corrector-properties}.
	\end{proof}
\begin{remark}[Comparison with the SPDE approach]
	\label{rem:intro-mean-covariant-current-comparison-BC24}
	Proposition~\ref{lem:mean-covariant-current-logarithmic-bound}
	is an analogue of \cite[Theorem~1.5(iii)]{BC24}, but estimates
	the expectation of the covariant current in a spatial Sobolev norm,
	rather than the stochastic current itself in space--time
	H\"older--Besov norms. Our gauge-covariant lattice regularization
	requires no counterterm for this estimate. Our proof uses direct random-walk bridge estimates, whereas \cite{BC24} relies on covariant parabolic kernel estimates.
\end{remark}
	\appendix
	
	\section{Interpolation theorems}
	\label{app:interpolation-theorems}
	
	For each lattice gauge field $B$, write $A_B:=1+\mathcal L_B^\varepsilon$.
	Choose a nonnegative $\eta\in C_c^\infty((-2,2)^2)$ equal to one on
	$[-1,1]^2$, and set
	\begin{equation}
		\label{eq:app-interpolation-partition}
		\chi_k(x):=
		\frac{\eta(x-k)}
		{\bigl(\sum_{\ell\in\mathbb Z^2}\eta(x-\ell)^2\bigr)^{1/2}},
		\qquad k\in\mathbb Z^2.
	\end{equation}
 Their supports have uniformly finite overlap, and
	\begin{equation}
		\label{eq:app-interpolation-partition-properties}
		\sum_{k\in\mathbb Z^2}\chi_k^2=1,\qquad
		\operatorname{supp}\chi_k\subset k+[-2,2]^2,\qquad
		\sup_{k\in\mathbb Z^2}\|\partial^\beta\chi_k\|_\infty<\infty
	\end{equation}
	for every multiindex $\beta$.

	\begin{lemma}[Weighted operator interpolation]
		\label{lem:weighted-operator-interpolation}
		Let $B,B'\in(\mathbb R^2)^{\Lambda_\varepsilon}$ be arbitrary.
		Let $s_j,t_j\in[-2,2]$ and $a_j,b_j\in\mathbb R$, $j=0,1$.
		Suppose that a linear operator $T$, initially defined on finitely
		supported lattice functions, has consistent bounded endpoint
		extensions, with
		\begin{equation}
			\label{eq:app-weighted-interpolation-endpoints}
			M_j:=
			\|T\|_{H_B^{s_j,\varepsilon}(\rho^{a_j})
			\to H_{B'}^{t_j,\varepsilon}(\rho^{b_j})}<\infty,
			\qquad j=0,1.
		\end{equation}
		For $\theta\in[0,1]$, set
		\begin{equation}
			\label{eq:app-weighted-interpolation-parameters}
			\begin{aligned}
				s_\theta&:=(1-\theta)s_0+\theta s_1,
				& t_\theta&:=(1-\theta)t_0+\theta t_1,\\
				a_\theta&:=(1-\theta)a_0+\theta a_1,
				& b_\theta&:=(1-\theta)b_0+\theta b_1.
			\end{aligned}
		\end{equation}
		Then $T$ has a bounded extension satisfying
		\begin{equation}
			\label{eq:app-weighted-operator-interpolation}
			\|T\|_{H_B^{s_\theta,\varepsilon}(\rho^{a_\theta})
			\to H_{B'}^{t_\theta,\varepsilon}(\rho^{b_\theta})}
			\le C M_0^{1-\theta}M_1^\theta.
		\end{equation}
		The constant depends only on the endpoint regularities and weight
		exponents, and is independent of $B,B',\varepsilon,\theta,T$.
		If $a_0=a_1$ and $b_0=b_1$, all four endpoint regularities may be
		arbitrary real numbers, and one can take $C=1$.
	\end{lemma}

	\begin{proof}
		Denote
$
			\label{eq:app-interpolation-increments}
			\delta s:=s_1-s_0, \delta t:=t_1-t_0,
			\delta a:=a_1-a_0, \delta b:=b_1-b_0.
$

		\noindent\textit{Fixed weights.}
		Suppose first that $\delta a=\delta b=0$. Conjugating $T$ by the
		fixed weights preserves both endpoint norms:
$			\label{eq:app-fixed-weight-conjugation}
			\widetilde T f:=\rho^{b_0}T(\rho^{-a_0}f).
$
		Relabeling $\widetilde T$ as $T$, we may assume $a_0=b_0=0$.
		Set $\widehat V:=H_B^{s_0,\varepsilon}\times
		H_{B'}^{t_0,\varepsilon}$, with its product Hilbert norm, and define
		\begin{equation}
			\label{eq:app-fixed-weight-block}
			\widehat\Gamma(f,g):=(A_B^{\delta s}f,A_{B'}^{\delta t}g),
			\qquad
			\widehat{\mathcal T}(f,g):=(0,Tf).
		\end{equation}
		Commutation of powers and the spectral bounds for $A_B,A_{B'}$
		show that $\widehat\Gamma$ is positive self-adjoint on $\widehat V$
		and has a strictly positive lower bound.
		For $r\in\mathbb R$, let $\widehat V^r$ denote $\widehat V$ equipped with the norm
		$\|\widehat\Gamma^{r/2}\cdot\|_{\widehat V}$.
		Then, isometrically,
		\begin{equation}
			\label{eq:app-fixed-weight-scale}
			\widehat V^r
			=H_B^{s_0+r\delta s,\varepsilon}
			\times H_{B'}^{t_0+r\delta t,\varepsilon},
			\qquad r\in\mathbb R.
		\end{equation}
		The endpoint assumptions are exactly
		\begin{equation}
			\label{eq:app-fixed-weight-block-endpoints}
			\|\widehat{\mathcal T}\|_{\widehat V^j\to\widehat V^j}
			=\|T\|_{H_B^{s_j,\varepsilon}\to H_{B'}^{t_j,\varepsilon}}
			=M_j,\qquad j=0,1.
		\end{equation}
		Thus \cite[Lemma~A.2]{BCCKM26} gives
		\begin{equation}
			\label{eq:app-fixed-weight-block-interpolation}
			\|T\|_{H_B^{s_\theta,\varepsilon}\to H_{B'}^{t_\theta,\varepsilon}}
			=\|\widehat{\mathcal T}\|_{\widehat V^\theta\to\widehat V^\theta}
			\le M_0^{1-\theta}M_1^\theta.
		\end{equation}
		This gives \eqref{eq:app-weighted-operator-interpolation} with
		$C=1$ for arbitrary real $s_j,t_j$.

		\noindent\textit{Variable weights.}
		For $r,c\in\mathbb R$, introduce the sequence spaces
		\begin{equation}
			\label{eq:app-localized-sequence-space}
			E_B^{r,c}:=\left\{\mathbf u=(u_k)_{k\in\mathbb Z^2}:
			\|\mathbf u\|_{E_B^{r,c}}^2
			:=\sum_{k\in\mathbb Z^2}\rho(k)^{2c}
			\|u_k\|_{H_B^{r,\varepsilon}}^2<\infty\right\}.
		\end{equation}
		Define $E_{B'}^{r,c}$ in the same way, using $H_{B'}^{r,\varepsilon}$.
		Define  maps
		(see diagram~\eqref{eq:app-weighted-interpolation-diagram} below) by
		\begin{equation}
			\label{eq:app-localization-maps}
			Rf:=(\chi_kf)_k,\qquad
			S\mathbf u:=\sum_{k\in\mathbb Z^2}\chi_ku_k,\qquad Q:=RTS.
		\end{equation}
		Since $\sum_{k\in\mathbb Z^2}\chi_k^2=1$, we have
		\begin{equation}
			\label{eq:app-localization-retraction}
			SRf=\sum_{k\in\mathbb Z^2}\chi_k^2f=f,\qquad T=SQR.
		\end{equation}
 Lemma~\ref{lem:weighted-covariant-localization}
		and the endpoint assumptions imply
		\begin{equation}
			\label{eq:app-localized-operator-endpoints}
			\|Q\|_{E_B^{s_j,a_j}\to E_{B'}^{t_j,b_j}}
			\le C_jM_j,\qquad j=0,1,
		\end{equation}
		where $C_j$ depends only on the corresponding regularities and
		weight exponents.

		Set $X_0:=E_B^{s_0,a_0}$ and $Y_0:=E_{B'}^{t_0,b_0}$.
		Choose the two signs independently so that
		\begin{equation}
			\label{eq:app-interpolation-signs}
			\sigma_a,\sigma_b\in\{-1,1\},\qquad
			\sigma_a\delta a\le0,\qquad \sigma_b\delta b\le0.
		\end{equation}
		On the base Hilbert spaces $X_0,Y_0$, respectively, define
		\begin{equation}
			\label{eq:app-interpolation-generators}
			\begin{aligned}
				(\Gamma_X\mathbf u)_k
				:=\rho(k)^{2\sigma_a\delta a}A_B^{\sigma_a\delta s}u_k,\quad
				(\Gamma_Y\mathbf v)_k
				:=\rho(k)^{2\sigma_b\delta b}A_{B'}^{\sigma_b\delta t}v_k.
			\end{aligned}
		\end{equation}
		Since  $\rho(k)\le1$ and
		\eqref{eq:app-interpolation-signs} holds, the  operators
		$\Gamma_X,\Gamma_Y$ are positive self-adjoint on $X_0,Y_0$,
		respectively, and bounded below by strictly positive constants.

		For $r\in\mathbb R$, let $X_0^r,Y_0^r$ be the completions of
		\begin{equation}
			\label{eq:app-interpolation-power-domains}
			\begin{aligned}
				D(\Gamma_X^{r/2}):=\{\mathbf u\in X_0:\Gamma_X^{r/2}\mathbf u\in X_0\},\quad
				D(\Gamma_Y^{r/2}):=\{\mathbf v\in Y_0:\Gamma_Y^{r/2}\mathbf v\in Y_0\},
			\end{aligned}
		\end{equation}
		in the norms $\|\mathbf u\|_{X_0^r}:=\|\Gamma_X^{r/2}\mathbf u\|_{X_0}$
		and $\|\mathbf v\|_{Y_0^r}:=\|\Gamma_Y^{r/2}\mathbf v\|_{Y_0}$,
		respectively. These norms satisfy
		\begin{equation}
			\label{eq:app-interpolation-power-norms}
			\begin{aligned}
				\|\Gamma_X^{\theta\sigma_a/2}\mathbf u\|_{X_0}^2
				=\sum_{k\in\mathbb Z^2}\rho(k)^{2a_\theta}
				\|A_B^{s_\theta/2}u_k\|_{L^{2,\varepsilon}}^2,\quad
				\|\Gamma_Y^{\theta\sigma_b/2}\mathbf v\|_{Y_0}^2
				=\sum_{k\in\mathbb Z^2}\rho(k)^{2b_\theta}
				\|A_{B'}^{t_\theta/2}v_k\|_{L^{2,\varepsilon}}^2.
			\end{aligned}
		\end{equation}
		Thus the identifications are isometric:
		\begin{equation}
			\label{eq:app-interpolation-scale-identifications}
			X_0^{\theta\sigma_a}=E_B^{s_\theta,a_\theta},\qquad
			Y_0^{\theta\sigma_b}=E_{B'}^{t_\theta,b_\theta}.
		\end{equation}
		Now set
		\begin{equation}
			\label{eq:app-interpolation-direct-sum}
			V:=X_0\oplus Y_0,\qquad
			\Gamma:=\Gamma_X\oplus\Gamma_Y,\qquad
			\mathcal T(\mathbf u,\mathbf v):=(0,Q\mathbf u).
		\end{equation}
		Then $V^r=X_0^r\oplus Y_0^r$.
		The following diagram illustrates the above construction.
		\begin{equation}
			\label{eq:app-weighted-interpolation-diagram}
			\begin{tikzcd}[ampersand replacement=\&,
				column sep=3.5em, row sep=1.8em]
				H_B^{s_\theta,\varepsilon}(\rho^{a_\theta})
				\arrow[r, shift left=1.1ex, "R"]
				\arrow[dd, "T"']
				\arrow[ddr, phantom, "\scalebox{1.8}{$\displaystyle\circlearrowright$}" description]
				\& E_B^{s_\theta,a_\theta}
				\arrow[l, shift left=1.1ex, "S"]
				\arrow[r, hook, "\theta=0"]
				\arrow[dd, dashed, "Q"]
				\arrow[ddr, phantom, "\scalebox{1.8}{$\displaystyle\circlearrowright$}" description]
				\& V^0=X_0\oplus Y_0
				\arrow[r, "\Gamma^{-\sigma_a/2}"]
				\arrow[dd, dashed, "\mathcal T"']
				\& V^{\sigma_a}
				\arrow[dd, dashed, "\mathcal T"]\\
				\& \& {}
				\arrow[r, leftrightarrow, "\text{interpolation}"]
				\& {}\\
				H_{B'}^{t_\theta,\varepsilon}(\rho^{b_\theta})
				\arrow[r, shift left=1.1ex, "R"]
				\& E_{B'}^{t_\theta,b_\theta}
				\arrow[l, shift left=1.1ex, "S"]
				\arrow[r, hook, "\theta=0"]
				\& V^0
				\arrow[r, "\Gamma^{-\sigma_b/2}"]
				\& V^{\sigma_b}
			\end{tikzcd}
		\end{equation}
 By \eqref{eq:app-localized-operator-endpoints}, $\mathcal T$ satisfies
		\begin{equation}
			\label{eq:app-interpolation-block-endpoints}
			\|\mathcal T\|_{V^0\to V^0}\le C_0M_0,\qquad
			\|\mathcal T\|_{V^{\sigma_a}\to V^{\sigma_b}}\le C_1M_1.
		\end{equation}
		Applying \cite[Lemma~A.2]{BCCKM26} with input indices
		$0,\sigma_a$ and output indices $0,\sigma_b$, we obtain
		\begin{equation}
			\label{eq:app-interpolation-block-bound}
			\|\mathcal T\|_{V^{\theta\sigma_a}\to V^{\theta\sigma_b}}
			\le (C_0M_0)^{1-\theta}(C_1M_1)^\theta
			\lesssim M_0^{1-\theta}M_1^\theta.
		\end{equation}
		Restricting the input to $(\mathbf u,0)$ and using
		\eqref{eq:app-interpolation-scale-identifications} gives the same
		bound for $Q:E_B^{s_\theta,a_\theta}\to E_{B'}^{t_\theta,b_\theta}$.
		Finally, $T=SQR$ and Lemma~\ref{lem:weighted-covariant-localization} yield
		\begin{equation}
			\label{eq:app-interpolation-reconstruction}
			\begin{aligned}
				\|Tf\|_{H_{B'}^{t_\theta,\varepsilon}(\rho^{b_\theta})}
				\lesssim \|QRf\|_{E_{B'}^{t_\theta,b_\theta}}
				\lesssim M_0^{1-\theta}M_1^\theta
				\|f\|_{H_B^{s_\theta,\varepsilon}(\rho^{a_\theta})}.
			\end{aligned}
		\end{equation}
	\end{proof}

	\begin{lemma}[Weighted covariant localization]
		\label{lem:weighted-covariant-localization}
		For any lattice gauge field $B\in(\mathbb R^2)^{\Lambda_\varepsilon}$, $r\in[-2,2]$, and $c\in\mathbb R$, the maps
		in~\eqref{eq:app-localization-maps} satisfy
		\begin{equation}
			\label{eq:app-localization-and-reconstruction}
			\|Rf\|_{E_B^{r,c}}\lesssim_{r,c}
			\|f\|_{H_B^{r,\varepsilon}(\rho^c)},\qquad
			\|S\mathbf u\|_{H_B^{r,\varepsilon}(\rho^c)}
			\lesssim_{r,c}\|\mathbf u\|_{E_B^{r,c}}.
		\end{equation}
		Consequently,
		\begin{equation}
			\label{eq:app-weighted-localization-equivalence}
			\|f\|_{H_B^{r,\varepsilon}(\rho^c)}
			\asymp_{r,c}
			\left(\sum_{k\in\mathbb Z^2}\rho(k)^{2c}
			\|\chi_kf\|_{H_B^{r,\varepsilon}}^2\right)^{1/2}.
		\end{equation}
	\end{lemma}

	\begin{proof}
		It suffices to prove the two bounds
		in~\eqref{eq:app-localization-and-reconstruction}, since $SR=I$.
		First take $c=0$. At $r=0$, $R$ is an isometry and $S=R^*$ with
		respect to the unweighted $L^2$ pairings, so
		\begin{equation}
			\label{eq:app-unweighted-localization-zero}
			\sum_{k\in\mathbb Z^2}\|\chi_kf\|_2^2=\|f\|_2^2,\qquad
			\left\|\sum_{k\in\mathbb Z^2}\chi_ku_k\right\|_2
			\le\left(\sum_{k\in\mathbb Z^2}\|u_k\|_2^2\right)^{1/2}.
		\end{equation}
		At $r=2$, the commutator identity
		in Lemma~\ref{lem:covariant-laplacian-weight-commutator},
		the bounded derivatives of $\chi_k$, and finite overlap give
		\begin{equation}
			\label{eq:app-unweighted-localization-two}
			\begin{aligned}
				\sum_{k\in\mathbb Z^2}\|A_B(\chi_kf)\|_2^2
				\lesssim\sum_{k\in\mathbb Z^2}\|\chi_kA_Bf\|_2^2
				+\sum_{k\in\mathbb Z^2}\|[\mathcal L_B^\varepsilon,\chi_k]f\|_2^2
				\lesssim\|A_Bf\|_2^2.
			\end{aligned}
		\end{equation}
Using Lemma~\ref{lem:covariant-nonnegative-weight} at \(r=2\) and finite overlap, we get
		\begin{equation}
			\label{eq:app-unweighted-reconstruction-two}
			\left\|A_B\sum_{k\in\mathbb Z^2}\chi_ku_k\right\|_2^2
			\lesssim\sum_{k\in\mathbb Z^2}\|A_B(\chi_ku_k)\|_2^2
			\lesssim\sum_{k\in\mathbb Z^2}\|A_Bu_k\|_2^2.
		\end{equation}
		As in the proof of Lemma~\ref{lem:weighted-operator-interpolation},
		apply \cite[Appendix~A, Lemma~A.2]{BCCKM26} on a direct-sum Hilbert
		space to obtain both bounds for $r\in[0,2]$.
		Duality gives $r\in[-2,0]$ as well.

		For $c\ne 0$, choose $\zeta\in C_c^\infty((-3,3)^2)$ equal
		to one on $[-2,2]^2$, and write $\zeta_k(x):=\zeta(x-k)$.
		Lemma~\ref{lem:covariant-nonnegative-weight} and the unweighted
		bounds already proved imply
		\begin{equation}
			\label{eq:app-weighted-localization-bound}
			\begin{aligned}
				\|Rf\|_{E_B^{r,c}}^2
				&=\sum_{k\in\mathbb Z^2}\rho(k)^{2c}
				\|\chi_kf\|_{H_B^{r,\varepsilon}}^2\\
				&\lesssim_r\sum_{k\in\mathbb Z^2}\rho(k)^{2c}
				\|\rho^{-c}\zeta_k\|_{C^2}^2
				\|\chi_k(\rho^cf)\|_{H_B^{r,\varepsilon}}^2\\
				&\lesssim_{r,c}\sum_{k\in\mathbb Z^2}
				\|\chi_k(\rho^cf)\|_{H_B^{r,\varepsilon}}^2
				\lesssim_{r,c}\|\rho^cf\|_{H_B^{r,\varepsilon}}^2.
			\end{aligned}
		\end{equation}
		Similarly,
		\begin{equation}
			\label{eq:app-weighted-reconstruction-bound}
			\begin{aligned}
				\|S\mathbf u\|_{H_B^{r,\varepsilon}(\rho^c)}
				=\left\|\sum_{k\in\mathbb Z^2}
				\chi_k(\zeta_k\rho^cu_k)\right\|_{H_B^{r,\varepsilon}}
			\lesssim_r
				\left(\sum_{k\in\mathbb Z^2}
				\|\zeta_k\rho^cu_k\|_{H_B^{r,\varepsilon}}^2\right)^{1/2}
				\lesssim_{r,c}
				\left(\sum_{k\in\mathbb Z^2}\rho(k)^{2c}
				\|u_k\|_{H_B^{r,\varepsilon}}^2\right)^{1/2}.
			\end{aligned}
		\end{equation}
	\end{proof}

	\section{Representation of heat kernels}
	\label{app:heat-kernel-representations}
	
		\begin{lemma}[Heat kernel of $\mathcal{L}_0^\varepsilon$]
			\label{lem:free-heat-kernel}
			The semigroup $e^{-t\mathcal{L}_0^\varepsilon}$ admits an integral
			kernel $p_t^{M,\varepsilon}(y;x)$ such that, for every
			$f\in\mathbb C^{\Lambda_{M,\varepsilon}}$ and $t\geq0$,
			\begin{equation}
				\label{eq:free-heat-kernel-representation}
				e^{-t\mathcal{L}_0^\varepsilon}f(x)
				=\int_{\Lambda_{M,\varepsilon}}
				f(y)\,p_t^{M,\varepsilon}(y;x)\,dy.
			\end{equation}
			Moreover, $p_t^{M,\varepsilon}(y;x)=p_t^{M,\varepsilon}(x-y)$ with
			\begin{equation}
				\label{eq:free-lattice-heat-kernel}
				p_t^{M,\varepsilon}(z)
				=\frac{1}{\lvert\Lambda_{M,\varepsilon}\rvert}
				\sum_{k\in Q_{M,\varepsilon}}
				\exp\bigl(-t\lvert k^\varepsilon\rvert^2\bigr)\mathrm e_k(z).
			\end{equation}
			For every $t>0$,
			\begin{equation}
				\label{eq:free-lattice-heat-kernel-diagonal-estimate}
				p_t^{M,\varepsilon}(0)
				\lesssim\min\bigl\{\varepsilon^{-2},M^{-2}+t^{-1}\bigr\}.
			\end{equation}
		\end{lemma}
		
		\begin{proof}
			Finite-volume Fourier inversion gives
			\begin{equation}
				\label{eq:proof-free-heat-kernel}
				\begin{aligned}
					e^{-t\mathcal L_0^\varepsilon}f(x)
					=\frac1{\lvert\Lambda_{M,\varepsilon}\rvert}
					\sum_{k\in Q_{M,\varepsilon}}
					e^{-t|k^\varepsilon|^2}\widehat f(k)\mathrm e_k(x)
					=\int_{\Lambda_{M,\varepsilon}}
					f(y)p_t^{M,\varepsilon}(x-y)\,dy.
				\end{aligned}
			\end{equation}
			Since $|k^\varepsilon|\simeq|k|$ on $Q_{M,\varepsilon}$,
			comparison with the Gaussian sum on $M^{-1}\mathbb Z^2$ gives
			\begin{equation}
				\label{eq:free-lattice-heat-kernel-diagonal-large-time}
				p_t^{M,\varepsilon}(0)
				\lesssim\frac1{\lvert\Lambda_{M,\varepsilon}\rvert}
				\sum_{k\in M^{-1}\mathbb Z^2}e^{-ct|k|^2}
				\lesssim M^{-2}+t^{-1}.
			\end{equation}
			On the other hand,
			\begin{equation}
				\label{eq:free-lattice-heat-kernel-diagonal-small-time-upper}
				p_t^{M,\varepsilon}(0)
				\leq\frac{\#Q_{M,\varepsilon}}{\lvert\Lambda_{M,\varepsilon}\rvert}
				=\varepsilon^{-2}.
			\end{equation}
			These bounds yield~\eqref{eq:free-lattice-heat-kernel-diagonal-estimate}.
		\end{proof}
		
		\begin{lemma}[Feynman--Kac representation]
			\label{lem:FK-LB}
			Let $B\colon\Lambda_{M,\varepsilon}\to\mathbb R^2$.
			The covariant semigroup $e^{-t\mathcal{L}_B^\varepsilon}$ admits
			an integral kernel $p_t^{B,M,\varepsilon}(y;x)$ such that, for every
			$f\in\mathbb C^{\Lambda_{M,\varepsilon}}$ and $t\geq0$,
			\begin{equation}
				\label{eq:magnetic-heat-kernel-representation}
				e^{-t\mathcal{L}_B^\varepsilon}f(x)
				=\int_{\Lambda_{M,\varepsilon}}
				f(y)\,p_t^{B,M,\varepsilon}(y;x)\,dy.
			\end{equation}
			For $t>0$, it admits the Feynman--Kac representation
			\begin{equation}
				\label{eq:FK-representation}
				p_t^{B,M,\varepsilon}(y;x)
				=p_t^{M,\varepsilon}(y;x)\,
				\mathbb{E}_x\Bigl[e^{iI_B^X(t)}\,\Big|\,X_t=y\Bigr],
			\end{equation}
			where $X$ is the continuous-time random walk on $\Lambda_{M,\varepsilon}$,
			\begin{equation}
				\label{eq:CTRW}
				X_t:=
				x+\varepsilon\sum_{j=1}^2
				\bigl(N_j^+(t)-N_j^-(t)\bigr)e_j
				\pmod{2\pi M\mathbb Z^2},
			\end{equation}
			and $\{N_j^+,N_j^-\}_{j=1}^2$ are i.i.d.\@ Poisson processes of
			intensity $\varepsilon^{-2}$.
			Here $I_B^X(t)$ is the discrete line integral of $B$ along $X$
			up to time $t$:
			\begin{equation}
				\label{eq:magnetic-phase-functional}
				I_B^X(t):=\varepsilon\sum_{j=1}^2\left(
				\int_0^t B_j(X_{s-})\,dN_j^+(s)
				-\int_0^t B_j(X_{s-}-\varepsilon e_j)\,dN_j^-(s)
				\right).
			\end{equation}
		\end{lemma}

		\begin{proof}
			It is enough to prove the semigroup identity
			\begin{equation}
				\label{eq:FK-semigroup-expectation}
				e^{-t\mathcal{L}_B^\varepsilon}f(x)
				=\mathbb{E}_x\left[f(X_t)e^{iI_B^X(t)}\right].
			\end{equation}
			Indeed, assuming \eqref{eq:FK-semigroup-expectation}, its case $B=0$
			and Lemma~\ref{lem:free-heat-kernel} give
			\begin{equation}
				\label{eq:free-random-walk-kernel-mass}
				\mathbb{E}_x\left[f(X_t)\right]
				=e^{-t\mathcal{L}_0^\varepsilon}f(x)
				=\varepsilon^2\sum_{y\in\Lambda_{M,\varepsilon}}
				f(y)p_t^{M,\varepsilon}(y;x).
			\end{equation}
			Hence
			$\mathbb{P}_x(X_t=y)=\varepsilon^2p_t^{M,\varepsilon}(y;x)$.
			Conditioning on $X_t$ in \eqref{eq:FK-semigroup-expectation}, we obtain
			\begin{equation}
				\label{eq:FK-kernel-identification}
				\begin{aligned}
					e^{-t\mathcal{L}_B^\varepsilon}f(x)
					&=\sum_{y\in\Lambda_{M,\varepsilon}}f(y)
					\mathbb{E}_x\left[e^{iI_B^X(t)}\mid X_t=y\right]
					\mathbb{P}_x(X_t=y)\\
					&=\int_{\Lambda_{M,\varepsilon}}f(y)p_t^{M,\varepsilon}(y;x)
					\mathbb{E}_x\left[e^{iI_B^X(t)}\mid X_t=y\right]dy.
				\end{aligned}
			\end{equation}
			Thus \eqref{eq:magnetic-heat-kernel-representation} holds with the
			kernel given by \eqref{eq:FK-representation}.

			It remains to prove \eqref{eq:FK-semigroup-expectation}. Set
			\begin{equation}
				\label{eq:magnetic-random-walk-semigroup-definition}
				Z_t:=e^{iI_B^X(t)},
				\qquad
				P_t^{B,M,\varepsilon}f(x)
				:=\mathbb{E}_x\left[f(X_t)Z_t\right].
			\end{equation}
			At a jump time $s$, definition
			\eqref{eq:prelim-covariant-derivative} gives
			\begin{equation}
				\label{eq:magnetic-random-walk-jump-increment}
				f(X_s)Z_s-f(X_{s-})Z_{s-}
				=\varepsilon Z_{s-}
				\begin{cases}
				D_{B,j}^\varepsilon f(X_{s-}),&\text{at a jump of }N_j^+,\\
				D_{B,j}^{\varepsilon *}f(X_{s-}),&\text{at a jump of }N_j^-.
				\end{cases}
			\end{equation}
			Summing the jump increments up to time $t$ yields
			\begin{equation}
				\label{eq:magnetic-random-walk-Poisson-integral}
				\begin{aligned}
					f(X_t)Z_t&=f(x)+\varepsilon\sum_{j=1}^2
					\int_0^t D_{B,j}^\varepsilon f(X_{s-})Z_{s-}\,dN_j^+(s)
					+\varepsilon\sum_{j=1}^2
					\int_0^t D_{B,j}^{\varepsilon *}f(X_{s-})Z_{s-}\,dN_j^-(s).
				\end{aligned}
			\end{equation}
			A direct calculation from
			\eqref{eq:prelim-covariant-derivative}--
			\eqref{eq:lattice-covariant-laplacian} gives
			\begin{equation}
				\label{eq:covariant-laplacian-first-difference-generator}
				\frac1\varepsilon\sum_{j=1}^2
				\left(D_{B,j}^\varepsilon+D_{B,j}^{\varepsilon *}\right)f
				=-\mathcal{L}_B^\varepsilon f.
			\end{equation}
			Since
			$\widetilde N_j^\pm(s)=N_j^\pm(s)-\varepsilon^{-2}s$
			are martingales, taking expectations in
			\eqref{eq:magnetic-random-walk-Poisson-integral} and using
			\eqref{eq:covariant-laplacian-first-difference-generator} gives
			\begin{equation}
				\label{eq:magnetic-random-walk-semigroup-integral-equation}
				\begin{aligned}
					P_t^{B,M,\varepsilon}f(x)
					&=f(x)+\frac1\varepsilon\sum_{j=1}^2\mathbb{E}_x\int_0^t
					\left(D_{B,j}^\varepsilon+D_{B,j}^{\varepsilon *}\right)
					f(X_{s-})Z_{s-}\,ds\\
					&=f(x)-\int_0^t
					P_s^{B,M,\varepsilon}\bigl(\mathcal{L}_B^\varepsilon f\bigr)(x)\,ds.
				\end{aligned}
			\end{equation}
			Consequently,
			\begin{equation}
				\label{eq:magnetic-random-walk-semigroup-backward-equation}
				\partial_tP_t^{B,M,\varepsilon}f
				=-P_t^{B,M,\varepsilon}\mathcal{L}_B^\varepsilon f,
				\qquad P_0^{B,M,\varepsilon}f=f.
			\end{equation}
			Since $\mathcal L_B^\varepsilon$ is bounded on
			$\mathbb C^{\Lambda_{M,\varepsilon}}$, uniqueness for this finite-dimensional ODE gives
			$P_t^{B,M,\varepsilon}=e^{-t\mathcal L_B^\varepsilon}$
			as linear operators on $\mathbb C^{\Lambda_{M,\varepsilon}}$.
			This gives the semigroup identity~\eqref{eq:FK-semigroup-expectation}.
		\end{proof}
		
		\begin{remark}[Full-lattice heat kernels]
			\label{rem:full-lattice-FK}
			For $B\colon\Lambda_\varepsilon\to\mathbb R^2$ and $t>0$, the
			full-lattice heat kernels satisfy
			\begin{equation}
				\label{eq:full-lattice-heat-kernels}
				\begin{aligned}
					p_t^\varepsilon(y;x)=\varepsilon^{-2}\mathbb P_x(X_t=y),\quad
					p_t^{B,\varepsilon}(y;x)=p_t^\varepsilon(y;x)\,
					\mathbb E_x\bigl[e^{iI_B^X(t)}\mid X_t=y\bigr].
				\end{aligned}
			\end{equation}
			Here $X$ and $I_B^X$ are given by~\eqref{eq:CTRW}
			and~\eqref{eq:magnetic-phase-functional} without reduction modulo
			$2\pi M\mathbb Z^2$. The full-lattice versions
			of~\eqref{eq:free-heat-kernel-representation}
			and~\eqref{eq:magnetic-heat-kernel-representation} hold for
			$f\in L^{p,\varepsilon}$, $1\le p\le\infty$.
		\end{remark}
		
		\begin{lemma}[Moments of the lifted random walk]
			\label{lem:lifted-random-walk-moment}
			Let $\widetilde{X}$ be the full-lattice lift of the random walk in
			\eqref{eq:CTRW}. For every integer $m\geq1$ and $r\geq0$,
			\begin{equation}
				\label{eq:random-walk-torus-distance-second-moment}
				\mathbb{E}_x\lvert\widetilde{X}_r-x\rvert^{2m}
				\lesssim_m\sum_{\ell=1}^m
				\varepsilon^{2m-2\ell}r^\ell
				\lesssim_m r(1+r)^{m-1}.
			\end{equation}
		\end{lemma}
		
		\begin{proof}
			Fix $j\in\{1,2\}$ and set
			$Y_r:=N_j^+(r)-N_j^-(r)$ and
			$\lambda:=r\varepsilon^{-2}$.
			The moment-generating function of $Y_r$ is
			\begin{equation}
				\label{eq:lifted-random-walk-coordinate-mgf}
				\begin{aligned}
					\mathbb{E}[e^{\theta Y_r}]
					=\exp\left(\lambda(e^\theta-1)
						+\lambda(e^{-\theta}-1)\right)
						=\exp\left(2\lambda(\cosh\theta-1)\right)
					=\sum_{\ell\geq0}\frac{(2\lambda)^\ell}{\ell!}
						\left(\sum_{k\geq1}\frac{\theta^{2k}}{(2k)!}\right)^\ell.
				\end{aligned}
			\end{equation}
			Comparing the coefficients of $\theta^{2m}$ gives
			\begin{equation}
				\label{eq:lifted-random-walk-coordinate-moment}
				\begin{aligned}
					\mathbb{E}[\lvert Y_r\rvert^{2m}]
					=(2m)!\sum_{\ell=1}^m\frac{(2\lambda)^\ell}{\ell!}
						\sum_{\substack{k_1+\cdots+k_\ell=m\\k_q\geq1}}
						\prod_{q=1}^\ell\frac1{(2k_q)!}
				\lesssim_m\sum_{\ell=1}^m\lambda^\ell.
				\end{aligned}
			\end{equation}
			Since $(a_1^2+a_2^2)^m
			\lesssim_m\lvert a_1\rvert^{2m}+\lvert a_2\rvert^{2m}$,
			\eqref{eq:lifted-random-walk-coordinate-moment} yields
			\begin{equation}
				\begin{aligned}
					\mathbb{E}_x\lvert\widetilde{X}_r-x\rvert^{2m}
					\lesssim_m\varepsilon^{2m}
						\sum_{j=1}^2\mathbb{E}[\lvert N_j^+(r)-N_j^-(r)\rvert^{2m}]
					\lesssim_m\sum_{\ell=1}^m
						\varepsilon^{2m-2\ell}r^\ell.
				\end{aligned}
			\end{equation}
			Finally, $0<\varepsilon<1$ and
			$\sum_{\ell=1}^m r^\ell\lesssim_m r(1+r)^{m-1}$ yield
			\eqref{eq:random-walk-torus-distance-second-moment}.
		\end{proof}

	\section{Bou\'e--Dupuis formula and maximum principle}
	\label{app:maximum-principle-fbsde}
	In this appendix we give a self-contained proof of
	Proposition~\ref{prop:maximum-principle-and-fbsdes}.
	Throughout, we work on the Wiener space
$
		\label{eq:app-wiener-space}
		\bigl(\Omega,\mathcal{F},\{\mathcal{F}_t\}_{t\in[0,1]},\mathbb{P}\bigr),
		\Omega=C_0\bigl([0,1];\mathbb{R}^N\bigr),
 d,N\in\mathbb{N}_+,
$
	where $W\colon\Omega\to C_0([0,1];\mathbb{R}^N)$ denotes the canonical
	process.
	We consider the $\mathbb{R}^d$-valued controlled dynamics
	\begin{equation}
		\label{eq:app-controlled-dynamics}
		X_t
		=\int_0^t A_s u_s\,ds
		+\int_0^t A_s\,dW_s,
		\qquad t\in[0,1].
	\end{equation}
	The cost functional is
	\begin{equation}
		\label{eq:app-cost-functional}
		J(u)
		:=\mathbb{E}\left[
		V(X_1)
		+\frac12\int_0^1\lVert u_s\rVert_{\mathbb{R}^N}^2\,ds
		\right],
	\end{equation}
	where $A\colon[0,1]\to\mathbb{R}^{d\times N}$ is deterministic and continuous,
	$V\colon\mathbb{R}^d\to\mathbb{R}$, and $\mathcal{U}_{\mathrm{ad}}$
	denotes the progressively measurable controls
	$u\colon[0,1]\times\Omega\to\mathbb{R}^N$ satisfying
	$\mathbb{E}\int_0^1\lVert u_s\rVert_{\mathbb{R}^N}^2\,ds<\infty$.
		The following result associates this control problem with the Gibbs measure
	\begin{equation}
		\label{eq:app-gibbs-measure}
		\mu(dx)\overset{\mathrm{def}}{=}
		\frac{1}{\mathcal{Z}}e^{-V(x)}\,\nu(dx),
	\end{equation}
	where $\nu:=\operatorname{Law}(\int_0^1 A_s\,dW_s)$ and
$
		\mathcal{Z}:=\int e^{-V(x)}\,\nu(dx).
$
	\begin{proposition}[Maximum principle and FBSDE]
		\label{lem:app-abstract-maximum-principle}
		Assume that $e^{-V}\in C_b^2(\mathbb{R}^d)$ and that there exist
		constants $C>0$ and $m\in\mathbb{N}$ such that
		\begin{equation}
			\label{eq:app-potential-hypothesis}
			\bigl\lvert\nabla V(x)\bigr\rvert
			\leq C\bigl(1+\lvert x\rvert^m\bigr),
			\qquad x\in\mathbb{R}^d.
		\end{equation}
		Then $\mathcal{Z}\in(0,\infty)$, and the control problem admits a unique
		optimal control $u\in\mathcal{U}_{\mathrm{ad}}$.
		Let $X$ denote the corresponding optimal state.
		The optimal control satisfies the optimality condition
		\begin{equation}
			\label{eq:app-maximum-principle}
			u_t
			=-A_t^*
			\mathbb{E}_t\bigl[\nabla V(X_1)\bigr]
		\end{equation}
		for $dt\otimes d\mathbb{P}$-almost every $(t,\omega)$.
		Consequently, the associated shift $Z_t:=\int_0^t A_s u_s\,ds$
		satisfies the FBSDE
		\begin{equation}
			\label{eq:app-integral-fbsde}
			Z_t
			=-\int_0^t
			A_s A_s^*
			\mathbb{E}_s\bigl[\nabla V(X_1)\bigr]\,ds,
			\qquad t\in[0,1].
		\end{equation}
		The terminal law of $X$ is given by
		\begin{equation}
			\label{eq:app-optimal-terminal-law}
			\operatorname{Law}(X_1)=\mu.
		\end{equation}
		Finally, the Bou\'e--Dupuis formula holds:
		\begin{equation}
			\label{eq:app-boue-dupuis-identity}
			-\log\mathcal{Z}
			=
			\inf_{\widetilde{u}\in\mathcal{U}_{\mathrm{ad}}}
			J(\widetilde{u}).
		\end{equation}
	\end{proposition}
	
	\begin{proof}
		\textit{Step~1: Construction of the feedback control.}
		For $0\leq t\leq1$ and $x\in\mathbb{R}^d$, define
		\begin{equation}
			\label{eq:app-value-function-psi}
			\Psi(t,x)
			:=
			\mathbb{E}\left[
			\exp\Bigl(
			-V\Bigl(
			x+\int_t^1 A_s\,dW_s
			\Bigr)
			\Bigr)
			\right].
		\end{equation}
		In particular, $\Psi(0,0)=\mathcal{Z}$.
		Under these assumptions, one has
		$\Psi(t,x)\in(0,\infty)$ for all
		$(t,x)\in[0,1]\times\mathbb{R}^d$.
		By the Feynman--Kac formula
		(see~\cite[Theorem~5.7.6]{KS91}),
		$\Psi$ solves the linear PDE
		\begin{equation}
			\label{eq:app-feynman-kac-pde}
			\begin{cases}
				\partial_t\Psi(t,x)
				+\dfrac12
				\operatorname{Tr}\bigl(A_t A_t^* D^2\Psi(t,x)\bigr)
				=0,
				&
				(t,x)\in[0,1)\times\mathbb{R}^d,\\[0.4em]
				\Psi(1,x)=\exp\bigl(-V(x)\bigr),
				&
				x\in\mathbb{R}^d.
			\end{cases}
		\end{equation}
		$-\log\Psi$ satisfies the HJB equation associated with the control
		problem. Set $b(t,x):=A_t^*\nabla\log\Psi(t,x)$.
		Since $\Psi$, $\nabla\Psi$, and $D^2\Psi$ are continuous and
		$\Psi$ has a positive lower bound on compact subsets of
		$[0,1]\times\mathbb{R}^d$, $b$ is continuous and locally
		Lipschitz in $x$, uniformly for $t\in[0,1]$.
		(In general $b$ need not be globally bounded or globally
		Lipschitz, since $V$ need not be bounded.)
		
		Consider the SDE
		\begin{equation}
			\label{eq:app-feedback-sde}
			X_t
			=\int_0^t A_s b(s,X_s)\,ds
			+\int_0^t A_s\,dW_s,
			\qquad 0\leq t\leq1.
		\end{equation}
		
		\smallskip
		\noindent\textit{Step~2: Existence of a weak solution.}
		Define
		\begin{equation}
			\label{eq:app-free-field-and-density}
			Y_t
			:=
			\int_0^t A_s\,dW_s,
			\qquad
			M_t
			:=
			\frac{\Psi(t,Y_t)}{\Psi(0,0)}
			=
			\frac{
				\mathbb{E}\bigl[e^{-V(Y_1)}\bigm|\mathcal{F}_t\bigr]
			}{\mathcal{Z}}.
		\end{equation}
		Then $M_t>0$ almost surely,
		$M_1=e^{-V(Y_1)}/\mathcal{Z}$, and It\^o's formula yields
$
			dM_t
			=M_t
			\bigl\langle b(t,Y_t),\,dW_t\bigr\rangle.
$
		Define a probability measure $\mathbb{Q}$ on $\mathcal{F}_1$ by
		\begin{equation}
			\label{eq:app-girsanov-measure}
			\frac{d\mathbb{Q}}{d\mathbb{P}}
			:=M_1
			=
			\mathcal{E}
			\biggl(
			\int_0^\cdot
			\bigl\langle b(s,Y_s),\,dW_s\bigr\rangle
			\biggr)_1.
		\end{equation}
		Set
$
			\widetilde{W}_t
			:=
			W_t-\int_0^t b(s,Y_s)\,ds.
$
		By Girsanov's theorem, $\widetilde{W}$ is a Brownian motion under
		$\mathbb{Q}$. Consequently, under $\mathbb{Q}$,
		\begin{equation}
			\label{eq:app-y-under-q}
			Y_t
			=\int_0^t A_s b(s,Y_s)\,ds
			+\int_0^t A_s\,d\widetilde{W}_s.
		\end{equation}
		Thus
		$\bigl(\Omega,\mathcal{F},\mathbb{Q},\widetilde{W},Y\bigr)$
		is a weak solution of~\eqref{eq:app-feedback-sde}.
		
		\smallskip
		\noindent\textit{Step~3: Construction of a strong solution under
			$\mathbb{P}$ and terminal law.}
		Let $\chi_R$ be a smooth cut-off equal to $1$ on the ball
		$\{\lvert x\rvert\leq R\}$ and vanishing outside
		$\{\lvert x\rvert\leq R+1\}$.
		Since
		$(t,x)\mapsto A_t\chi_R(x)b(t,x)$
		is bounded and globally Lipschitz in $x$ on $[0,1]$,
		the cut-off equation
		\begin{equation}
			\label{eq:app-cutoff-sde}
			X_t^R
			=\int_0^t
			A_s\chi_R(X_s^R)\,b(s,X_s^R)\,ds
			+\int_0^t A_s\,dW_s
		\end{equation}
		admits a unique strong solution.
		Let $\tau_R^X$ and $\tau_R^Y$ denote the first exit times of
		$X^R$ and $Y$, respectively, from the ball of radius $R$.
		Under $\mathbb{Q}$, the stopped process
		$Y_{\cdot\wedge\tau_R^Y}$ solves the same stopped equation
		driven by the $\mathbb{Q}$-Brownian motion $\widetilde{W}$.
		By pathwise uniqueness for the cut-off dynamics,
		\begin{equation}
			\label{eq:app-stopped-law-identification}
			\operatorname{Law}_{\mathbb{P}}
			\bigl(X_{\cdot\wedge\tau_R^X}^R\bigr)
			=
			\operatorname{Law}_{\mathbb{Q}}
			\bigl(Y_{\cdot\wedge\tau_R^Y}\bigr).
		\end{equation}
		In particular,
		$\mathbb{P}(\tau_R^X\leq1)=\mathbb{Q}(\tau_R^Y\leq1)$.
		Since $Y$ is continuous on $[0,1]$ almost surely under
		$\mathbb{Q}$, one has
		$\mathbb{Q}(\tau_R^Y\leq1)\to0$ as $R\to\infty$, and therefore
		\begin{equation}
			\label{eq:app-no-explosion}
			\mathbb{P}(\tau_R^X\leq1)\to0
			\qquad\text{as }R\to\infty.
		\end{equation}
		Equivalently, $\tau_R^X>1$ for all sufficiently large $R$,
		$\mathbb{P}$-almost surely.
		
		Define a process $X$ by setting
		$X_t=X_t^R$ whenever $t\leq\tau_R^X$.
		This is unambiguous: if $R_1<R_2$, then
		$X_t^{R_1}=X_t^{R_2}$ for all
		$t\leq\tau_{R_1}^X$, by uniqueness of the cut-off equation.
		Thus $X$ satisfies~\eqref{eq:app-feedback-sde} on
		$[0,\tau_R^X\wedge1]$ for every $R\geq1$.
		Combining this with~\eqref{eq:app-no-explosion}, we conclude that
		$X$ solves~\eqref{eq:app-feedback-sde} on the full interval
		$[0,1]$.
		Moreover,~\eqref{eq:app-stopped-law-identification} yields
		\begin{equation}
			\label{eq:app-path-law-identification}
			\operatorname{Law}_{\mathbb{P}}
			\bigl((X_s)_{0\leq s\leq1}\bigr)
			=
			\operatorname{Law}_{\mathbb{Q}}
			\bigl((Y_s)_{0\leq s\leq1}\bigr).
		\end{equation}
		Indeed, for every Borel set $A\subset\mathbb{R}^d$ and every
		$s\in[0,1]$,
		\begin{equation}
			\label{eq:app-marginal-law-limit}
			\begin{aligned}
				\mathbb{P}(X_s\in A)
				=\lim_{R\to\infty}
				\mathbb{P}\bigl(X_s^R\in A,\,s\leq\tau_R^X\bigr)
				=\lim_{R\to\infty}
				\mathbb{Q}\bigl(
				Y_{s\wedge\tau_R^Y}\in A,\,
				s\leq\tau_R^Y
				\bigr)
				=\mathbb{Q}(Y_s\in A),
			\end{aligned}
		\end{equation}
		and the same argument identifies all finite-dimensional
		distributions.
		
		For every bounded Borel function
		$f\colon\mathbb{R}^d\to\mathbb{R}$,
		\begin{equation}
			\label{eq:app-terminal-law-computation}
			\mathbb{E}_{\mathbb{P}}\bigl[f(X_1)\bigr]
			=
			\mathbb{E}_{\mathbb{Q}}\bigl[f(Y_1)\bigr]
			=
			\frac{1}{\mathcal{Z}}
			\mathbb{E}_{\mathbb{P}}
			\bigl[f(Y_1)\,e^{-V(Y_1)}\bigr],
		\end{equation}
		which gives~\eqref{eq:app-optimal-terminal-law}.
		
		\smallskip
		\noindent\textit{Step~4: Admissibility.}
		Define $u_t:=b(t,X_t)$.
		We claim that $u\in\mathcal{U}_{\mathrm{ad}}$ and
		$J(u)=-\log\mathcal{Z}$.
		First,
		\begin{equation}
			\label{eq:app-cost-under-q}
			\begin{aligned}
				J(u)
				&=
				\mathbb{E}_{\mathbb{P}}\left[
				V(X_1)
				+\frac12\int_0^1
				\lVert u_s\rVert_{\mathbb{R}^N}^2\,ds
				\right]\\
				&=
				\mathbb{E}_{\mathbb{P}}\left[
				V(X_1)
				+\frac12\int_0^1
				\lVert b(s,X_s)\rVert_{\mathbb{R}^N}^2\,ds
				\right]\\
				&=
				\mathbb{E}_{\mathbb{Q}}\left[
				V(Y_1)
				+\frac12\int_0^1
				\lVert b(s,Y_s)\rVert_{\mathbb{R}^N}^2\,ds
				\right].
			\end{aligned}
		\end{equation}
		By~\eqref{eq:app-girsanov-measure},
		\begin{equation}
			\label{eq:app-relative-entropy-identity}
			\begin{aligned}
				\mathcal{H}(\mathbb{Q}\mid\mathbb{P})
				&:=
				\mathbb{E}_{\mathbb{Q}}
				\left[\log\frac{d\mathbb{Q}}{d\mathbb{P}}\right]
				=
				\frac12
				\mathbb{E}_{\mathbb{Q}}\left[
				\int_0^1
				\lVert b(s,Y_s)\rVert_{\mathbb{R}^N}^2\,ds
				\right];
			\end{aligned}
		\end{equation}
		see~\cite[(3.3)]{Fol85}.
		On the other hand, by~\eqref{eq:app-free-field-and-density},
		\begin{equation}
			\label{eq:app-entropy-via-potential}
			\begin{aligned}
				\mathcal{H}(\mathbb{Q}\mid\mathbb{P})
				=-\mathbb{E}_{\mathbb{Q}}\bigl[V(Y_1)\bigr]
				-\log\mathcal{Z}
				=-\frac{1}{\mathcal{Z}}
				\mathbb{E}_{\mathbb{P}}\bigl[V(Y_1)e^{-V(Y_1)}\bigr]
				-\log\mathcal{Z}<+\infty.
			\end{aligned}
		\end{equation}
		By~\eqref{eq:app-relative-entropy-identity}--\eqref{eq:app-entropy-via-potential}
		and the path-law identity~\eqref{eq:app-path-law-identification},
		the feedback control $u_t=b(t,X_t)$ belongs to
		$\mathcal{U}_{\mathrm{ad}}$.
		Moreover,
		\begin{equation}
			\label{eq:app-optimal-cost-value}
			\begin{aligned}
				J(u)
				&=
				\mathbb{E}_{\mathbb{Q}}\bigl[V(Y_1)\bigr]
				+\mathcal{H}(\mathbb{Q}\mid\mathbb{P})
				=
				-\log\mathcal{Z}.
			\end{aligned}
		\end{equation}
		
		\smallskip
		\noindent\textit{Step~5: Bou\'e--Dupuis formula and optimality.}
		We use a verification argument;
		see~\cite[Chapter~5, Theorems~3.7 and~5.1]{YZ99}.
		Set $v(t,x):=-\log\Psi(t,x)$ and write
		$g:=e^{-V}$.
		Define the approximations
$
			g_n(x)
			:=
			e^{-V(x)}+\frac1n,
			V_n(x)
			:=
			-\log g_n(x),
$
		and set
		\begin{equation}
			\label{eq:app-psi-approximation}
			\begin{aligned}
				\Psi_n(t,x)
				:=
				\mathbb{E}\left[
				\exp\Bigl(
				-V_n\Bigl(
				x+\int_t^1 A_s\,dW_s
				\Bigr)
				\Bigr)
				\right]
			&=
				\mathbb{E}\left[
				g\Bigl(
				x+\int_t^1 A_s\,dW_s
				\Bigr)
				\right]
				+\frac1n
				\geq\frac1n.
			\end{aligned}
		\end{equation}
		Let $v_n(t,x):=-\log\Psi_n(t,x)$.
		Each $v_n$ satisfies the HJB equation associated with the terminal
		cost $V_n$; consequently, for every
		$\widetilde{u}\in\mathcal{U}_{\mathrm{ad}}$ with controlled state
		$\widetilde{X}$,
		\begin{equation}
			\label{eq:app-ito-vn}
			\begin{aligned}
				dv_n\bigl(t,\widetilde{X}_t\bigr)
				=
				\Biggl(
				\frac12
				\bigl\lVert
				A_t^*\nabla v_n\bigl(t,\widetilde{X}_t\bigr)
				+\widetilde{u}_t
				\bigr\rVert_{\mathbb{R}^N}^2
				-
				\frac12
				\lVert\widetilde{u}_t\rVert_{\mathbb{R}^N}^2
				\Biggr)\,dt
				+
				\bigl\langle
				\nabla v_n\bigl(t,\widetilde{X}_t\bigr),
				\,A_t\,dW_t
				\bigr\rangle.
			\end{aligned}
		\end{equation}
		Integrating over $[0,1]$ and taking expectations yields
		\begin{equation}
			\label{eq:app-verification-approx}
			\mathbb{E}\left[
			V_n\bigl(\widetilde{X}_1\bigr)
			+\frac12\int_0^1
			\lVert\widetilde{u}_s\rVert_{\mathbb{R}^N}^2\,ds
			\right]
			-v_n(0,0)
			=
			\mathbb{E}\left[
			\frac12\int_0^1
			\bigl\lVert
			A_s^*\nabla v_n\bigl(s,\widetilde{X}_s\bigr)
			+\widetilde{u}_s
			\bigr\rVert_{\mathbb{R}^N}^2
			ds
			\right]
			\geq0.
		\end{equation}
		Here we have used the representation
		\begin{equation}
			\label{eq:app-grad-vn}
			\nabla v_n(t,x)
			=
			-\frac{1}{\Psi_n(t,x)}
			\mathbb{E}\left[
			\nabla g\Bigl(
			x+\int_t^1 A_s\,dW_s
			\Bigr)
			\right],
		\end{equation}
		together with $g\in C_b^2(\mathbb{R}^d)$ and
		$\Psi_n\geq\frac1n$, which ensure that
		$\nabla v_n\bigl(\cdot,\widetilde{X}_\cdot\bigr)$ is bounded for each
		fixed $n$; thus the stochastic integral in~\eqref{eq:app-ito-vn}
		is a martingale.
		Rearranging~\eqref{eq:app-verification-approx} gives
		\begin{equation}
			\label{eq:app-verification-inequality-n}
			v_n(0,0)
			\leq
			\mathbb{E}\left[
			V_n\bigl(\widetilde{X}_1\bigr)
			+\frac12\int_0^1
			\lVert\widetilde{u}_s\rVert_{\mathbb{R}^N}^2\,ds
			\right].
		\end{equation}
		Since $v_n(0,0)\uparrow v(0,0)$ by~\eqref{eq:app-psi-approximation},
		$V_n\uparrow V$, and
		\begin{equation}
			\label{eq:app-Vn-lower-bound}
			V_n(x)
			=
			-\log\Bigl(g(x)+\frac1n\Bigr)
			\geq
			-\log\bigl(\lVert g\rVert_{L^\infty}+1\bigr),
		\end{equation}
		the monotone convergence theorem yields, upon sending
		$n\to\infty$,
		\begin{equation}
			\label{eq:app-verification-inequality}
			v(0,0)
			=-\log\mathcal{Z}
			\leq
			\mathbb{E}\left[
			V\bigl(\widetilde{X}_1\bigr)
			+\frac12\int_0^1
			\lVert\widetilde{u}_s\rVert_{\mathbb{R}^N}^2\,ds
			\right]
			=
			J(\widetilde{u}).
		\end{equation}
		Since $\widetilde{u}$ was arbitrary,~\eqref{eq:app-optimal-cost-value}
		yields optimality and the Bou\'e--Dupuis formula:
		\begin{equation}
			\label{eq:app-verification-optimality-inf}
			-\log\mathcal{Z}
			=
			J(u)
			=
			\inf_{\widetilde{u}\in\mathcal{U}_{\mathrm{ad}}}
			J(\widetilde{u}).
		\end{equation}
		\smallskip
		\noindent\textit{Step~6: Uniqueness.}
		Equation~\eqref{eq:app-grad-vn} gives
		$A_t^*\nabla v_n(t,x)\to-b(t,x)$ pointwise.
		If $\widetilde{u}$ is optimal, Fatou's lemma applied to
		~\eqref{eq:app-verification-approx} yields
		\begin{equation}
			\label{eq:app-optimal-control-uniqueness}
			\mathbb{E}\left[
			\frac12\int_0^1
			\bigl\lVert\widetilde{u}_s-b(s,\widetilde{X}_s)\bigr\rVert_{\mathbb{R}^N}^2\,ds
			\right]
			\leq J(\widetilde{u})-J(u)=0.
		\end{equation}
		Thus $\widetilde{u}_t=b(t,\widetilde{X}_t)$ for
		$dt\otimes d\mathbb{P}$-almost every $(t,\omega)$, so
		$\widetilde{X}$ solves~\eqref{eq:app-feedback-sde} driven by the same $W$.
		The local Lipschitz property of $b$ yields pathwise uniqueness;
		hence $\widetilde{X}$ and $X$ are indistinguishable, and
		$\widetilde{u}=u$ $dt\otimes d\mathbb{P}$-almost everywhere.
		
		\smallskip
		\noindent\textit{Step~7: The maximum principle.}
		The integral FBSDE~\eqref{eq:app-integral-fbsde} follows at once
		from the maximum principle~\eqref{eq:app-maximum-principle}
		upon writing
		$Z_t=\int_0^t A_s u_s\,ds$.
		It remains to establish~\eqref{eq:app-maximum-principle}.
		We first consider bounded variations
		$h\in\mathcal{U}_{\mathrm{ad}}$.
		Optimality of $u$ yields $J(u+\delta h)\geq J(u)$ for every
		sufficiently small $\delta\in\mathbb{R}$, and therefore
		\begin{equation}
			\label{eq:app-gateaux-vanishing}
			\left.
			\frac{d}{d\delta}J(u+\delta h)
			\right|_{\delta=0}
			=0.
		\end{equation}
		For the running cost,
		\begin{equation}
			\label{eq:app-running-cost-derivative}
			\left.
			\frac{d}{d\delta}
			\right|_{\delta=0}
			\mathbb{E}\left[
			\frac12\int_0^1
			\lVert u_s+\delta h_s\rVert_{\mathbb{R}^N}^2\,ds
			\right]
			=
			\mathbb{E}\left[
			\int_0^1\langle u_s,h_s\rangle\,ds
			\right],
		\end{equation}
		which is justified by
		$\mathbb{E}\bigl[\int_0^1\lVert h_s\rVert_{\mathbb{R}^N}^2\,ds\bigr]
		<\infty$.
		For the terminal cost, the difference quotient reads
		\begin{equation}
			\label{eq:app-terminal-difference-quotient}
			\frac1\delta
			\Biggl(
			V\Bigl(
			X_1+\delta\int_0^1 A_s h_s\,ds
			\Bigr)
			-V(X_1)
			\Biggr)
			=
			\int_0^1
			\left\langle
			\nabla V\Bigl(
			X_1
			+\delta\theta\int_0^1 A_s h_s\,ds
			\Bigr),
			\int_0^1 A_s h_s\,ds
			\right\rangle
			d\theta.
		\end{equation}
		Using~\eqref{eq:app-potential-hypothesis} and the boundedness of $h$,
		\begin{equation}
			\label{eq:app-terminal-dominator}
			\begin{aligned}
				\mathbb{E}
				\int_0^1
				\left\lvert
				\left\langle
				\nabla V\Bigl(
				X_1
				+\delta\theta\int_0^1 A_s h_s\,ds
				\Bigr),
				\int_0^1 A_s h_s\,ds
				\right\rangle
				\right\rvert
				d\theta
				&\lesssim
				\mathbb{E}\left[
				\Bigl(
				1
				+\lVert X_1\rVert^m
				+\Bigl(\int_0^1\lVert h_s\rVert\,ds\Bigr)^m
				\Bigr)
				\int_0^1\lVert h_s\rVert\,ds
				\right]\\
				&\lesssim
				1+\mathbb{E}_{\mathbb{P}}\bigl[\lVert X_1\rVert^m\bigr]
				=
				1
				+\frac{1}{\mathcal{Z}}\mathbb{E}_{\mathbb{P}}\bigl[
				\lVert Y_1\rVert^m e^{-V(Y_1)}
				\bigr]
				<+\infty.
			\end{aligned}
		\end{equation}
		The dominated convergence theorem therefore allows us to
		differentiate under the expectation and obtain
		\begin{equation}
			\label{eq:app-terminal-cost-derivative}
			\left.
			\frac{d}{d\delta}
			\right|_{\delta=0}
			\mathbb{E}\left[
			V\Bigl(
			X_1+\delta\int_0^1 A_s h_s\,ds
			\Bigr)
			\right]
			=
			\mathbb{E}\left[
			\left\langle
			\nabla V(X_1),
			\int_0^1 A_s h_s\,ds
			\right\rangle
			\right].
		\end{equation}
		Consequently,~\eqref{eq:app-gateaux-vanishing} is equivalent to
		\begin{equation}
			\label{eq:app-first-variation}
			\mathbb{E}\left[
			\left\langle
			\nabla V(X_1),
			\int_0^1 A_s h_s\,ds
			\right\rangle
			+\int_0^1
			\langle u_s,h_s\rangle\,ds
			\right]
			=0,
		\end{equation}
		and hence
		\begin{equation}
			\label{eq:app-first-variation-conditional}
			\mathbb{E}\left[
			\int_0^1
			\left\langle
			u_s
			+A_s^*
			\mathbb{E}_s\bigl[\nabla V(X_1)\bigr],
			\,h_s
			\right\rangle
			ds
			\right]
			=0.
		\end{equation}
		Set
$
			\label{eq:app-residual-control}
			R_s
			:=
			u_s
			+A_s^*
			\mathbb{E}_s\bigl[\nabla V(X_1)\bigr],
$
		and define the truncated admissible controls
$
			h_s^{(n)}
			:=
			R_s\mathbf{1}_{\{\lvert R_s\rvert\leq n\}}
			\in\mathcal{U}_{\mathrm{ad}}.
$
		Each $h^{(n)}$ is bounded.
		Substituting into~\eqref{eq:app-first-variation-conditional} yields
		\begin{equation}
			\label{eq:app-truncated-orthogonality}
			\mathbb{E}\left[
			\int_0^1
			\lVert R_s\rVert_{\mathbb{R}^N}^2
			\mathbf{1}_{\{\lvert R_s\rvert\leq n\}}
			\,ds
			\right]
			=0.
		\end{equation}
		Hence $R_s=0$ for $dt\otimes d\mathbb{P}$-almost every
		$(t,\omega)$ with $\lvert R_s(\omega)\rvert\leq n$, for every
		$n\in\mathbb{N}$.
		Letting $n\to\infty$, we conclude that $R\equiv0$
		$dt\otimes d\mathbb{P}$-almost everywhere, which is
		~\eqref{eq:app-maximum-principle}.
	\end{proof}

	\begin{lemma}
		\label{lem:app-cor-satisfies-hypothesis}
		The potentials in the Gibbs measures
		\eqref{eq:conditional-renormalized-higgs-measure} and
		\eqref{eq:conditional-ymh-gauge-law}
		satisfy the hypotheses of
		Proposition~\ref{lem:app-abstract-maximum-principle}.
	\end{lemma}
	
	\begin{proof}
		Fix $M$ and $\varepsilon$. All calculations below take place in finite-dimensional spaces. Constants
		in this proof may depend on $M$, $\varepsilon$.
		It suffices to check that the potentials
		$\mathcal{W}_{\mathrm{av}}^{M,\varepsilon}$ and
		$V_{A,B}(\cdot)$
		satisfy~\eqref{eq:app-potential-hypothesis}
		(together with $e^{-V}\in C_b^2$).
		
		\smallskip
		\noindent\textit{The Higgs potential
			$V_{A,B}(\cdot)$.}
		By~\eqref{eq:relative-covariant-interaction}, 
		\begin{equation}
			\label{eq:app-VA-coercivity}
			\begin{aligned}
				V_{A,B}(\phi)
				&=
				\frac12
				\left\langle
				\bigl(:\!\mathcal{L}_{A+B}^\varepsilon\!:_{M,\varepsilon}
				-\mathcal{L}_B^\varepsilon-1\bigr)\phi,\phi
				\right\rangle_{M,\varepsilon}
				+\frac14\int_{\Lambda_{M,\varepsilon}}
				:\!\lvert\phi(x)\rvert^4\!:_{M,\varepsilon}\,dx\\
				&\geq
				\frac18\lVert\phi\rVert_{L^{4,\varepsilon}}^4
				-c_1^{A,B,M,\varepsilon}
				\geq
				\frac{c_2^{M,\varepsilon}}{8}
				\lVert\phi\rVert_{L^{2,\varepsilon}}^4
				-c_1^{A,B,M,\varepsilon},
			\end{aligned}
		\end{equation}
		for constants $c_1^{A,B,M,\varepsilon}<\infty$ and
		$c_2^{M,\varepsilon}>0$.
		Moreover,
		\begin{equation}
			\label{eq:app-VA-derivatives}
			\begin{aligned}
				\lVert
				\nabla_\phi V_{A,B}(\phi)
				\rVert_{L^{2,\varepsilon}}
				\lesssim_{A,B,M,\varepsilon}
				\bigl(1+\lVert\phi\rVert_{L^{2,\varepsilon}}^3\bigr),\quad
				\bigl\lVert
				D_\phi^2 V_{A,B}(\phi)
				\bigr\rVert_{L^{2,\varepsilon}\otimes L^{2,\varepsilon}}
				\lesssim_{A,B,M,\varepsilon}
				\bigl(1+\lVert\phi\rVert_{L^{2,\varepsilon}}^2\bigr).
			\end{aligned}
		\end{equation}
		In particular,
		\begin{equation}
			\label{eq:app-exp-VA-bound}
			\exp\bigl(-V_{A,B}(\phi)\bigr)
			\leq
			\exp\bigl(c_1^{A,B,M,\varepsilon}\bigr),
		\end{equation}
		and
		\begin{equation}
			\label{eq:app-grad-exp-VA}
			\begin{aligned}
				\bigl\lVert
				\nabla_\phi
				\exp\bigl(-V_{A,B}(\phi)\bigr)
				\bigr\rVert_{L^{2,\varepsilon}}
				&\lesssim_{M,\varepsilon}
				\lVert
				\nabla_\phi V_{A,B}(\phi)
				\rVert_{L^{2,\varepsilon}}
				\,
				\exp\bigl(-V_{A,B}(\phi)\bigr)\\
				&\lesssim_{A,B,M,\varepsilon}
				\bigl(1+\lVert\phi\rVert_{L^{2,\varepsilon}}^3\bigr)
				\exp\Bigl(
				-\frac{c_2^{M,\varepsilon}}{8}
				\lVert\phi\rVert_{L^{2,\varepsilon}}^4
				\Bigr)
				\lesssim_{A,B,M,\varepsilon}1.
			\end{aligned}
		\end{equation}
		Similarly,
		\begin{equation}
			\label{eq:app-hess-exp-VA}
			\bigl\lVert
			D_\phi^2
			\exp\bigl(-V_{A,B}(\phi)\bigr)
			\bigr\rVert_{L^{2,\varepsilon}\otimes L^{2,\varepsilon}}
			\lesssim_{A,B,M,\varepsilon}1.
		\end{equation}
		Combined with~\eqref{eq:app-VA-derivatives}, this yields
		$\exp(-V_{A,B}(\cdot))\in C_b^2$ and the gradient
		growth~\eqref{eq:app-potential-hypothesis}.
		
		\smallskip
		\noindent\textit{The gauge potential $\mathcal{W}_{\mathrm{av}}^{M,\varepsilon}$.}
		Recall from~\eqref{eq:effective-interaction-tilde-V} that
		\begin{equation}
			\label{eq:app-recall-tilde-V}
			\begin{aligned}
				\widetilde{V}^{M,\varepsilon}(A,\phi)
				&=
				\frac12
				\left\langle
				\bigl(:\!\mathcal{L}_A^\varepsilon\!:_{M,\varepsilon}
				-\mathcal{L}_0^\varepsilon-1\bigr)\phi,\phi
				\right\rangle_{M,\varepsilon}
				+\frac14\int_{\Lambda_{M,\varepsilon}}
				:\!\lvert\phi(x)\rvert^4\!:_{M,\varepsilon}\,dx.
			\end{aligned}
		\end{equation}
		By definition ~\eqref{eq:effective-gauge-potential},
		\begin{equation}
			\label{eq:app-exp-W-as-expectation}
			\begin{aligned}
				\exp\bigl(-\mathcal{W}^{M,\varepsilon}(A)\bigr)
				&=
				\mathbb{E}_{\mu_{M,\varepsilon}^{\mathrm{GFF}}}
				\bigl[
				\exp\bigl(-\widetilde{V}^{M,\varepsilon}(A,\phi)\bigr)
				\bigr].
			\end{aligned}
		\end{equation}
		Using
		$\langle\mathcal{L}_A^\varepsilon\phi,\phi\rangle_{M,\varepsilon}\geq0$,
		$\bigl\lvert
		\langle(\mathcal{L}_0^\varepsilon+1)\phi,\phi\rangle_{M,\varepsilon}
		\bigr\rvert
		\lesssim_{M,\varepsilon}
		\lVert\phi\rVert_{L^{2,\varepsilon}}^2$,
		and Young's inequality, one obtains the coercivity bound
		\begin{equation}
			\label{eq:app-tilde-V-coercivity}
			\widetilde{V}^{M,\varepsilon}(A,\phi)
			\geq
			\frac18\lVert\phi\rVert_{L^{4,\varepsilon}}^4
			-c_1^{M,\varepsilon},
		\end{equation}
		and therefore
		\begin{equation}
			\label{eq:app-exp-W-finite}
			\exp\bigl(-\mathcal{W}^{M,\varepsilon}(A)\bigr)
			\leq
			\exp\bigl(c_1^{M,\varepsilon}\bigr)
			<+\infty.
		\end{equation}
		We also have, using that
		$\widetilde{V}^{M,\varepsilon}(A,\phi)$ depends on $A$ only through 
		$e^{i\varepsilon A_j(x)}$,
		\begin{equation}
			\label{eq:app-tilde-V-A-derivatives}
			\begin{aligned}
				\lVert
				\nabla_A\widetilde{V}^{M,\varepsilon}(A,\phi)
				\rVert_{L^{2,\varepsilon}}
				&\lesssim_{M,\varepsilon}
				\bigl(1+\lVert\phi\rVert_{L^{2,\varepsilon}}^3\bigr),\\
				\bigl\lVert
				D_A^2\widetilde{V}^{M,\varepsilon}(A,\phi)
				\bigr\rVert_{L^{2,\varepsilon}\otimes L^{2,\varepsilon}}
				&\lesssim_{M,\varepsilon}
				\bigl(1+\lVert\phi\rVert_{L^{2,\varepsilon}}^2\bigr).
			\end{aligned}
		\end{equation}
		The same dependence implies that
		$\mathcal{W}^{M,\varepsilon}(A)$ depends on $A$ only through
		$e^{i\varepsilon A_j(x)}$.
		Consequently, once the derivatives are shown to exist and be continuous,
		we obtain
		\begin{equation}
			\label{eq:app-exp-W-Cb2}
			\exp(-\mathcal{W}^{M,\varepsilon})
			\in
			C_b^2
			\Bigl(
			\mathcal A_{\mathrm{gf}}^{M,\varepsilon}
			\Bigr)
		\end{equation}
		and
		$\lVert\nabla_A\mathcal{W}^{M,\varepsilon}\rVert_{L^{2,\varepsilon}}
		<+\infty$,
		which is~\eqref{eq:app-potential-hypothesis} for
		$\mathcal{W}^{M,\varepsilon}$.
		Finally,~\eqref{eq:app-tilde-V-A-derivatives}
		and~\eqref{eq:app-tilde-V-coercivity}
		provide integrable dominants for the integrands defining
		$\nabla_A\mathcal{W}^{M,\varepsilon}$ and $D_A^2\mathcal{W}^{M,\varepsilon}$,
		so the dominated convergence theorem yields their continuity.
		Moreover, $e^{-\mathcal W^{M,\varepsilon}}$ is bounded away from zero
		by compactness of the link phases 	$e^{i\varepsilon A_j(x)}$. Differentiating the integral
		in~\eqref{eq:harmonic-averaged-effective-potential} therefore gives
		$\mathcal W_{\mathrm{av}}^{M,\varepsilon}\in C_b^2$,
		which implies the required bounds for the outer potential.
	\end{proof}
		\section{Positive expansion and Gaussian domination}
	\label{app:positive-expansion-gaussian-domination}
	
	The arguments in this appendix are adapted from
	\cite[Sections~3--4]{CC24}.
	Let $E$ be the set of oriented  edges of the periodic lattice:
	\begin{equation}
		\label{eq:periodic-oriented-edge-set}
		E:=\bigl\{(x,x+\varepsilon e_j),(x+\varepsilon e_j,x):
		x\in\Lambda_{M,\varepsilon},\ j=1,2\bigr\},
	\end{equation}For
	$A\in(\mathbb R^2)^{\Lambda_{M,\varepsilon}}$, set
	$
	A_{(x,x+\varepsilon e_j)}:=A_j(x),
	A_{(x+\varepsilon e_j,x)}:=-A_j(x).
	$
	Thus the link variable on $e\in E$ is $e^{i\varepsilon A_e}$. Define
	\begin{equation}
		\label{eq:balanced-edge-multiplicities}
		\begin{aligned}
			&\mathcal N:=\left\{n\in\mathbb N_0^E:
			\sum_{y\sim x}n_{(x,y)}=\sum_{y\sim x}n_{(y,x)}
			\text{ for every }x\in\Lambda_{M,\varepsilon}\right\},\\
			&\ell_n(A):=\varepsilon\sum_{e\in E}n_eA_e,
			\qquad n\in\mathbb N_0^E.
		\end{aligned}
	\end{equation}
	We extend $\mathcal W^{M,\varepsilon}$ and
	$\mathcal W_{\mathrm{av}}^{M,\varepsilon}$ to all real periodic gauge
	fields by~\eqref{eq:effective-gauge-potential} and
	\eqref{eq:harmonic-averaged-effective-potential}.
	
	\begin{lemma}[Positive expansion]
		\label{lem:harmonic-averaged-positive-expansion}
		There exist coefficients $c_n\geq0$, $n\in\mathbb N_0^E$,
		with $\sum_nc_n<\infty$, such that, for every
		$A\in(\mathbb R^2)^{\Lambda_{M,\varepsilon}}$,
		\begin{equation}
			\label{eq:harmonic-averaged-positive-expansion}
			e^{-\mathcal W_{\mathrm{av}}^{M,\varepsilon}(A)}
			=\sum_{n\in\mathbb N_0^E}c_n\cos(\ell_n(A)).
		\end{equation}
	\end{lemma}
	\begin{proof}
		Following~\cite[Theorem~3.17 and Lemma~4.2]{CC24}, apply the positive
		loop expansion to the periodic graph $(\Lambda_{M,\varepsilon},E)$.
		We get
		\begin{equation}
			\label{eq:periodic-positive-loop-expansion}
			e^{-\mathcal W^{M,\varepsilon}(A)}
			=\sum_{n\in\mathcal N}\widetilde c_n e^{i\ell_n(A)},
			\qquad \widetilde c_n\geq0.
		\end{equation}
		At $A=0$, $\sum_{n\in\mathcal N}\widetilde c_n
		=e^{-\mathcal W^{M,\varepsilon}(0)}<\infty$. Thus
		\begin{equation}
			\label{eq:positive-expansion-harmonic-average}
			\begin{aligned}
				e^{-\mathcal W_{\mathrm{av}}^{M,\varepsilon}(A)}
				=\int_{[0,M^{-1})^2}
				e^{-\mathcal W^{M,\varepsilon}(A+h)}\,dh
				=\sum_{n\in\mathcal N}\widetilde c_n e^{i\ell_n(A)}
				\int_{[0,M^{-1})^2}e^{i\ell_n(h)}\,dh.
			\end{aligned}
		\end{equation}
		For $n\in\mathcal N$, write
		\begin{equation}
			\label{eq:positive-expansion-edge-current}
			\begin{aligned}
				J_j(x)&:=n_{(x,x+\varepsilon e_j)}
				-n_{(x+\varepsilon e_j,x)},\\
				\ell_n(h)&=\varepsilon
				\sum_{x\in\Lambda_{M,\varepsilon}}\sum_{j=1}^2J_j(x)h_j.
			\end{aligned}
		\end{equation}
		The balance condition $n\in\mathcal N$ gives, for any $x\in \Lambda_{M,\epsilon}$
		\begin{equation}
			\label{eq:positive-expansion-current-conservation}
			\sum_{y\sim x}(n_{(x,y)}-n_{(y,x)})
			=\sum_{j=1}^2\bigl(J_j(x)-J_j(x-\varepsilon e_j)\bigr)=0.
		\end{equation}
		Set $S(x_1):=\sum_{x_2\in\varepsilon\mathbb Z/(2\pi M\mathbb Z)}
		J_1(x_1,x_2)$. Summing the preceding identity over $x_2$ at fixed $x_1$ yields
		\begin{equation}
			\label{eq:positive-expansion-transverse-flux}
			S(x_1)-S(x_1-\varepsilon)
			=\sum_{x_2\in\varepsilon\mathbb Z/(2\pi M\mathbb Z)}
			\bigl(J_2(x_1,x_2-\varepsilon )-J_2(x_1,x_2)\bigr)=0.
		\end{equation}
		Hence $S$ is constant and
		\begin{equation}
			\label{eq:positive-expansion-integral-flux}
			\varepsilon\sum_{x\in\Lambda_{M,\varepsilon}}J_1(x)
			=\varepsilon\sum_{x_1\in\varepsilon\mathbb Z/(2\pi M\mathbb Z)}
			S(x_1)=2\pi M S(0)\in2\pi M\mathbb Z.
		\end{equation}
		The same argument applies to the second component. Consequently,
		\begin{equation}
			\label{eq:positive-expansion-harmonic-character}
			M^2\int_{[0,M^{-1})^2}e^{i\ell_n(h)}\,dh
			=\prod_{j=1}^2\mathbf1_{\{\sum_xJ_j(x)=0\}}\in\{0,1\}.
		\end{equation}
		Define $c_n:=\widetilde c_n
		\int_{[0,M^{-1})^2}e^{i\ell_n(h)}\,dh$ for $n\in\mathcal N$,
		and $c_n:=0$ otherwise. These coefficients are nonnegative and
		summable. Taking real parts in
		\eqref{eq:positive-expansion-harmonic-average} yields the stated expansion.
	\end{proof}
	
	\begin{lemma}[Gaussian domination; adapted from {\cite[Lemma~4.7]{CC24}}]
		\label{lem:gaussian-domination}
		Let $(X,Y)$ be a centered real Gaussian vector. For every $t\in\mathbb R$,
		\begin{equation}
			\label{eq:gaussian-domination}
			\mathbb E[e^{tX}\cos Y]
			\leq e^{t^2\mathbb E[X^2]/2}\mathbb E[\cos Y].
		\end{equation}
	\end{lemma}
	\begin{proof}
		The Gaussian moment formula gives
		\begin{equation}
			\label{eq:gaussian-domination-proof}
			\begin{aligned}
				\mathbb E[e^{tX}\cos Y]
				&=\operatorname{Re}\mathbb E[e^{tX+iY}]
				=e^{t^2\mathbb E[X^2]/2-\mathbb E[Y^2]/2}
				\cos(t\mathbb E[XY])\\
				&\leq e^{t^2\mathbb E[X^2]/2-\mathbb E[Y^2]/2}
				=e^{t^2\mathbb E[X^2]/2}\mathbb E[\cos Y].
			\end{aligned}
		\end{equation}
	\end{proof}
	
	\begin{corollary}[adapted from {\cite[Theorem~4.6 and Corollary~4.8]{CC24}}]
		\label{cor:gauge-gaussian-domination}
		Let $C^{\varepsilon,0}$ be defined by
		\eqref{eq:optimal-state-processes}. For every real linear functional
		$F\colon(\mathbb R^2)^{\Lambda_{M,\varepsilon}}\to\mathbb R$,
		$t\in\mathbb R$ and $2\leq p<\infty$,
		\begin{equation}
			\label{eq:gauge-linear-exponential-domination}
			\mathbb E[e^{tF(C^{\varepsilon,0})}]
			\leq\exp\left(\frac{t^2}{2}
			\mathbb E^A[F(A^{\varepsilon,0})^2]\right)
		\end{equation}
		and
		\begin{equation}
			\label{eq:gauge-linear-moment-domination}
			\bigl(\mathbb E|F(C^{\varepsilon,0})|^p\bigr)^{1/p}
			\lesssim\sqrt p\,
			\bigl(\mathbb E^A[F(A^{\varepsilon,0})^2]\bigr)^{1/2}.
		\end{equation}
	\end{corollary}
	\begin{proof}
		Put $X:=F(A^{\varepsilon,0})$ and
		$Y_n:=\ell_n(A^{\varepsilon,0})$.
		By~\eqref{eq:optimal-terminal-laws},
		$\operatorname{Law}(C^{\varepsilon,0})=\mu_{M,\varepsilon}^{\mathrm{gauge}}$.
		Lemmas~\ref{lem:harmonic-averaged-positive-expansion}
		and~\ref{lem:gaussian-domination} yield
		\begin{equation}
			\label{eq:gauge-linear-exponential-domination-proof}
			\begin{aligned}
				\mathbb E[e^{tF(C^{\varepsilon,0})}]
				&=\frac{\mathbb E^A[e^{tX}
					e^{-\mathcal W_{\mathrm{av}}^{M,\varepsilon}(A^{\varepsilon,0})}]}
				{\mathbb E^A[e^{-\mathcal W_{\mathrm{av}}^{M,\varepsilon}
						(A^{\varepsilon,0})}]}\\
				&=\frac{\sum_n c_n\mathbb E^A[e^{tX}\cos Y_n]}
				{\sum_n c_n\mathbb E^A[\cos Y_n]}
				\leq e^{t^2\mathbb E^A[X^2]/2}.
			\end{aligned}
		\end{equation}
		Set $Z:=F(C^{\varepsilon,0})$ and
		$\sigma^2:=\mathbb E^A[X^2]$. If $\sigma=0$,
		\eqref{eq:gauge-linear-exponential-domination} implies $Z=0$ almost surely.
		Otherwise, for $r,t>0$,
		\begin{equation}
			\label{eq:gauge-linear-chernoff-bound}
			\mathbb P(Z\geq r)\leq e^{-tr}\mathbb E[e^{tZ}]
			\leq\exp\left(-tr+\frac{t^2\sigma^2}{2}\right).
		\end{equation}
		Choosing $t=r/\sigma^2$ and applying the same argument to $-Z$ gives
		\begin{equation}
			\label{eq:gauge-linear-two-sided-tail}
			\mathbb P(|Z|\geq r)\leq2e^{-r^2/(2\sigma^2)}.
		\end{equation}
		Therefore,
		\begin{equation}
			\label{eq:gauge-linear-tail-moment-integral}
			\begin{aligned}
				\mathbb E|Z|^p
				&=p\int_0^\infty r^{p-1}\mathbb P(|Z|\geq r)\,dr\\
				&\leq2p\int_0^\infty r^{p-1}e^{-r^2/(2\sigma^2)}\,dr
				=p\,2^{p/2}\sigma^p\Gamma(p/2).
			\end{aligned}
		\end{equation}
		Stirling's formula yields~\eqref{eq:gauge-linear-moment-domination}.
	\end{proof}
	
	\section{Proof of
	Lemma~\ref{lem:gauge-irregular-Zhi1-unified-free-comparison}}
	\label{app:proof-gauge-irregular-Zhi1-unified-free-comparison}

	We first record the two multiplication estimates needed below.
	\begin{lemma}
		\label{lem:gauge-link-multiplier}
		Let $B\colon\Lambda_\varepsilon\to\mathbb R^2$ satisfy
		$B\in H^{1,\varepsilon}(\rho)$. For every
		$\sigma\in[0,\frac12]$, with $1/0=\infty$,
		\begin{equation}
			\label{eq:gauge-link-besov-bound}
			\left\lVert U_B^\varepsilon
			\right\rVert_{B_{1/\sigma,\infty}^{2\sigma,\varepsilon}(\rho^{2\sigma})}
			\lesssim_\sigma
			\bigl(1+\lVert B\rVert_{H^{1,\varepsilon}(\rho)}\bigr)^{2\sigma}.
		\end{equation}
		Moreover, for every
		$f\in\mathcal C^{-\kappa,\varepsilon}(\rho^\kappa)$,
		\begin{equation}
			\label{eq:gauge-link-multiplier-bound}
			\left\lVert U_B^\varepsilon f
			\right\rVert_{\mathcal{C}^{-\kappa,\varepsilon}(\rho^{3\kappa})}
			\lesssim
			\bigl(1+\lVert B\rVert_{H^{1,\varepsilon}(\rho)}\bigr)^{2\kappa}
			\lVert f\rVert_{\mathcal{C}^{-\kappa,\varepsilon}(\rho^\kappa)}.
		\end{equation}
	\end{lemma}

	\begin{proof}
		Since $|U_{B,j}^\varepsilon|=1$ and
		$|D_\ell^\varepsilon U_{B,j}^\varepsilon|
		\le\varepsilon|D_\ell^\varepsilon B_j|$, we have
		\begin{equation}
			\label{eq:appE-gauge-link-H1-bound}
			\left\lVert U_B^\varepsilon\right\rVert_{H^{1,\varepsilon}(\rho)}
			\lesssim
			1+\lVert B\rVert_{H^{1,\varepsilon}(\rho)}.
		\end{equation}
		Weighted interpolation gives
		\begin{equation}
			\label{eq:appE-gauge-link-interpolation-bound}
			\begin{aligned}
			\left\lVert U_B^\varepsilon
			\right\rVert_{B_{1/\sigma,\infty}^{2\sigma,\varepsilon}(\rho^{2\sigma})}
			\lesssim_\sigma
			\left\lVert U_B^\varepsilon
			\right\rVert_{L^{\infty,\varepsilon}}^{1-2\sigma}
			\left\lVert U_B^\varepsilon
			\right\rVert_{H^{1,\varepsilon}(\rho)}^{2\sigma}
			\lesssim_\sigma
			\bigl(1+\lVert B\rVert_{H^{1,\varepsilon}(\rho)}\bigr)^{2\sigma}.
			\end{aligned}
		\end{equation}
		Decompose $U_B^\varepsilon f=U_B^\varepsilon\prec f
		+U_B^\varepsilon\circ f+f\prec U_B^\varepsilon$.
		Using $0<\rho\le1$ and the weighted embedding
		$B_{1/\kappa,\infty}^{\kappa,\varepsilon}(\rho^{3\kappa})
		\hookrightarrow\mathcal{C}^{-\kappa,\varepsilon}(\rho^{3\kappa})$,
		the weighted paraproduct estimates yield
		\begin{equation}
			\label{eq:appE-gauge-link-product-bound}
			\begin{aligned}
			\left\lVert U_B^\varepsilon f
			\right\rVert_{\mathcal{C}^{-\kappa,\varepsilon}(\rho^{3\kappa})}
			&\lesssim
			\left(
			\left\lVert U_B^\varepsilon\right\rVert_{L^{\infty,\varepsilon}}
			+\left\lVert U_B^\varepsilon
			\right\rVert_{B_{1/\kappa,\infty}^{2\kappa,\varepsilon}(\rho^{2\kappa})}
			\right)
			\lVert f\rVert_{\mathcal{C}^{-\kappa,\varepsilon}(\rho^\kappa)}\\
			&\lesssim
			\bigl(1+\lVert B\rVert_{H^{1,\varepsilon}(\rho)}\bigr)^{2\kappa}
			\lVert f\rVert_{\mathcal{C}^{-\kappa,\varepsilon}(\rho^\kappa)}.
			\end{aligned}
		\end{equation}
	\end{proof}

	\begin{lemma}
		\label{lem:appE-gauge-link-B-multiplier}
		Let $B\in H^{1,\varepsilon}(\rho)$ and $c\in\mathbb R$.
		For every $f\colon\Lambda_\varepsilon\to\mathbb C$ with finite
		right-hand side, the following estimates hold.
		\begin{equation}
			\label{eq:appE-gauge-link-B-multiplier}
			\left\lVert U_B^\varepsilon(Bf)
			\right\rVert_{H^{2\kappa,\varepsilon}(\rho^{c+1+4\kappa})}
			\lesssim_c
			\bigl(1+\lVert B\rVert_{H^{1,\varepsilon}(\rho)}\bigr)^{1+4\kappa}
			\lVert f\rVert_{H^{1-\alpha,\varepsilon}(\rho^c)}.
		\end{equation}
		\begin{equation}
			\label{eq:appE-gauge-BU-holder-multiplier}
			\left\lVert U_B^\varepsilon(Bf)
			\right\rVert_{H^{-4\alpha,\varepsilon}(\rho^{c+1+2\alpha})}
			\lesssim_c
			\bigl(1+\lVert B\rVert_{H^{1,\varepsilon}(\rho)}\bigr)^{1+2\alpha}
			\lVert f\rVert_{H^{-\alpha,\varepsilon}(\rho^c)}.
		\end{equation}
	\end{lemma}

	\begin{proof}
		The discrete product estimate gives
		\begin{equation}
			\label{eq:appE-gauge-link-B-product-estimate}
			\left\lVert U_B^\varepsilon(Bf)
			\right\rVert_{H^{2\kappa,\varepsilon}(\rho^{c+1+4\kappa})}
			\lesssim
			\left\lVert U_B^\varepsilon B
			\right\rVert_{B_{4,\infty}^{4\kappa,\varepsilon}(\rho^{1+4\kappa})}
			\lVert f\rVert_{B_{4,\infty}^{4\kappa,\varepsilon}(\rho^c)}.
		\end{equation}
		By Lemma~\ref{lem:gauge-link-multiplier}, applied with
		$\sigma=2\kappa$, and the weighted Besov embedding for $B$,
		with $p=4/(1-8\kappa)$,
		\begin{equation}
			\label{eq:appE-gauge-link-B-Besov-estimate}
			\begin{aligned}
			\left\lVert U_B^\varepsilon B
			\right\rVert_{B_{4,\infty}^{4\kappa,\varepsilon}(\rho^{1+4\kappa})}
			&\lesssim
			\left\lVert U_B^\varepsilon
			\right\rVert_{B_{1/(2\kappa),\infty}^{4\kappa,\varepsilon}(\rho^{4\kappa})}
			\lVert B\rVert_{B_{p,\infty}^{4\kappa,\varepsilon}(\rho)}\\
			&\lesssim
			\bigl(1+\lVert B\rVert_{H^{1,\varepsilon}(\rho)}\bigr)^{4\kappa}
			\lVert B\rVert_{H^{1,\varepsilon}(\rho)}.
			\end{aligned}
		\end{equation}
		Here $2\kappa+1/p=1/4$ and $4\kappa<2/p$.
		Moreover, since $4\kappa<1/2-\alpha$, the weighted Besov embedding gives
		\begin{equation}
			\label{eq:appE-gauge-link-B-input-estimate}
			\lVert f\rVert_{B_{4,\infty}^{4\kappa,\varepsilon}(\rho^c)}
			\lesssim
			\lVert f\rVert_{H^{1-\alpha,\varepsilon}(\rho^c)}.
		\end{equation}
		Combining \eqref{eq:appE-gauge-link-B-product-estimate}--%
		\eqref{eq:appE-gauge-link-B-input-estimate}
		yields~\eqref{eq:appE-gauge-link-B-multiplier}.

		For~\eqref{eq:appE-gauge-BU-holder-multiplier}, set
		$\frac1p=\frac12+\alpha$.
		By the discrete paraproduct estimates and the embedding of
		$B_{p,2}^{\alpha,\varepsilon}(\rho^{c+2\alpha})$ into
		$H^{-\alpha,\varepsilon}(\rho^{c+2\alpha})$, we obtain
		\begin{equation}
			\label{eq:appE-gauge-BU-holder-multiplier-proof}
			\begin{aligned}
			\left\lVert U_B^\varepsilon f
			\right\rVert_{H^{-\alpha,\varepsilon}(\rho^{c+2\alpha})}
			&\lesssim
			\left(
			\left\lVert U_B^\varepsilon\right\rVert_{L^{\infty,\varepsilon}}
			+\left\lVert U_B^\varepsilon
			\right\rVert_{B_{1/\alpha,\infty}^{2\alpha,\varepsilon}(\rho^{2\alpha})}
			\right)
			\lVert f\rVert_{H^{-\alpha,\varepsilon}(\rho^c)}\\
			&\lesssim
			\bigl(1+\lVert B\rVert_{H^{1,\varepsilon}(\rho)}\bigr)^{2\alpha}
			\lVert f\rVert_{H^{-\alpha,\varepsilon}(\rho^c)}.
			\end{aligned}
		\end{equation}
		Here we used $\lVert U_B^\varepsilon\rVert_{L^{\infty,\varepsilon}}=1$
		and~\eqref{eq:gauge-link-besov-bound} with $\sigma=\alpha$.
		By the product estimate and the embeddings of
		$B_{p,\infty}^{-\alpha,\varepsilon}(\rho^{c+1})$ into
		$B_{2,\infty}^{-3\alpha,\varepsilon}(\rho^{c+1})$ and then into
		$H^{-4\alpha,\varepsilon}(\rho^{c+1})$, we obtain
		\begin{equation}
			\label{eq:appE-gauge-Bf-H-4alpha-bound}
			\begin{aligned}
			\lVert Bf\rVert_{H^{-4\alpha,\varepsilon}(\rho^{c+1})}
			&\lesssim
			\lVert Bf\rVert_{B_{p,\infty}^{-\alpha,\varepsilon}(\rho^{c+1})}
			\lesssim
			\lVert B\rVert_{B_{1/\alpha,\infty}^{2\alpha,\varepsilon}(\rho)}
			\lVert f\rVert_{H^{-\alpha,\varepsilon}(\rho^c)}\\
			&\lesssim
			\lVert B\rVert_{H^{1,\varepsilon}(\rho)}
			\lVert f\rVert_{H^{-\alpha,\varepsilon}(\rho^c)}.
			\end{aligned}
		\end{equation}
		Applying this estimate with $c+2\alpha$ in place of $c$ to
		$U_B^\varepsilon f$ and using
		\eqref{eq:appE-gauge-BU-holder-multiplier-proof} yields
		\eqref{eq:appE-gauge-BU-holder-multiplier}.
	\end{proof}

	\noindent\textit{Proof of
	Lemma~\ref{lem:gauge-irregular-Zhi1-unified-free-comparison}.}
	Fix $j\in\{1,2\}$ and $a,b\in\{0,1\}$. Define the vector fields
	\begin{equation}
		\label{eq:appE-unified-comparison-sources}
		\mathcal{A}^{A,B}
		:=
		\left(
		\int_0^1\mathbb{E}_t^\phi\left[
		\mathcal{U}_{>L_1}^\varepsilon
		\bigl(
		r_{A,\ell}^{\varepsilon,(b)}
		\tau_\ell^b\phi_B^\varepsilon
		\bigr)
		\right]dt
		\right)_{\ell=1}^2,
		\qquad
		\mathcal{A}^{A,0}
		:=
		\left(
		\int_0^1\mathbb{E}_t^\phi\left[
		\mathcal{U}_{>L_1}^\varepsilon
		\bigl(
		r_{A,\ell}^{\varepsilon,(b)}
		\tau_\ell^b\phi^\varepsilon
		\bigr)
		\right]dt
		\right)_{\ell=1}^2.
	\end{equation}
	 We first decompose the difference between the covariant and free expressions:
	\begin{equation}
		\label{eq:appE-unified-comparison-algebraic-decomposition}
		\begin{aligned}
		&
		\overline{\tau_j^a\phi_B^\varepsilon}\,
		U_{-aB,j}^\varepsilon D_{B,j}^\varepsilon
		\bigl(\mathcal{L}_B^\varepsilon+1\bigr)^{-1}
		D_B^{\varepsilon *}U_{bB}^\varepsilon
		\mathcal{A}^{A,B}
		-
		\overline{\tau_j^a\phi^\varepsilon}\,
		D_j^\varepsilon
		\bigl(\mathcal{L}_0^\varepsilon+1\bigr)^{-1}
		D^{\varepsilon *}\mathcal{A}^{A,0}\\
		={}&
		\int_0^1
		\overline{\tau_j^a\phi_B^\varepsilon}\,
		\frac{d}{d\theta}\Bigl[
		U_{-a\theta B,j}^\varepsilon D_{\theta B,j}^\varepsilon
		\bigl(\mathcal{L}_{\theta B}^\varepsilon+1\bigr)^{-1}
		D_{\theta B}^{\varepsilon *}U_{b\theta B}^\varepsilon
		\Bigr]
		\mathcal{A}^{A,B}\,d\theta
		+S_1+S_2.
		\end{aligned}
	\end{equation}
	Here
	\begin{equation}
		\label{eq:appE-unified-comparison-S-terms}
		\begin{aligned}
		S_1
		&:=
		\overline{\tau_j^a
		(\phi_B^\varepsilon-\phi^\varepsilon)}\,
		D_j^\varepsilon
		\bigl(\mathcal{L}_0^\varepsilon+1\bigr)^{-1}
		D^{\varepsilon *}\mathcal{A}^{A,B},\quad
		S_2
		:=
		\overline{\tau_j^a\phi^\varepsilon}\,
		D_j^\varepsilon
		\bigl(\mathcal{L}_0^\varepsilon+1\bigr)^{-1}
		D^{\varepsilon *}
		\bigl(\mathcal{A}^{A,B}-\mathcal{A}^{A,0}\bigr).
		\end{aligned}
	\end{equation}

	For $j=1,2$, direct differentiation, together with
	\eqref{eq:covariant-derivative-translation-identities}, gives
	\begin{equation}
		\label{eq:appE-shifted-covariant-derivative-derivative}
		\frac{d}{d\theta}
		\bigl(
		U_{-a\theta B,j}^\varepsilon D_{\theta B,j}^\varepsilon
		\bigr)
		=
		\begin{cases}
		iB_jU_{\theta B,j}^\varepsilon\tau_j,
		& a=0,\\
		iB_jU_{-\theta B,j}^\varepsilon,
		& a=1.
		\end{cases}
	\end{equation}
	Moreover,
	\begin{equation}
		\label{eq:appE-interpolated-laplacian-derivative}
		\frac{d}{d\theta}\mathcal{L}_{\theta B}^\varepsilon
		=
		\left(
		\frac{d}{d\theta}D_{\theta B}^\varepsilon
		\right)^*D_{\theta B}^\varepsilon
		+D_{\theta B}^{\varepsilon *}
		\frac{d}{d\theta}D_{\theta B}^\varepsilon,
	\end{equation}
	and hence
	\begin{equation}
		\label{eq:appE-interpolated-resolvent-derivative}
		\frac{d}{d\theta}
		\bigl(\mathcal{L}_{\theta B}^\varepsilon+1\bigr)^{-1}
		=
		-
		\bigl(\mathcal{L}_{\theta B}^\varepsilon+1\bigr)^{-1}
		\left(
		\frac{d}{d\theta}\mathcal{L}_{\theta B}^\varepsilon
		\right)
		\bigl(\mathcal{L}_{\theta B}^\varepsilon+1\bigr)^{-1}.
	\end{equation}

	By~\eqref{eq:high-frequency-truncation-estimate} and the choice of $L_1$ in
	Table~\ref{tab:higgs-high-frequency-cutoff-choices}, for both $b=0$ and $b=1$,
	\begin{equation}
		\label{eq:appE-unified-comparison-source-bound}
		\begin{aligned}
		\left\lVert\mathcal{A}^{A,B}
		\right\rVert_{L^\infty(\Omega^\phi;
			\mathcal{C}^{-\alpha/2,\varepsilon}(\rho^{-4}))}
		&\lesssim
		\left\lVert
		2^{-L_1(\alpha/2-\kappa)}
		\lVert\Xi\rVert_\kappa
		\right\rVert_{L^\infty(\Omega^\phi)}
		\lesssim
		\frac{1}{
		\bigl(1+\lVert B\rVert_{H^{1,\varepsilon}(\rho)}\bigr)^{l_0}}.
		\end{aligned}
	\end{equation}
	Since all polynomial losses in $B$ below are of order strictly less than
	$l_0$, the cutoff gain in~\eqref{eq:appE-unified-comparison-source-bound}
	absorbs them uniformly.

	The derivative in the first term on the right-hand side of
	\eqref{eq:appE-unified-comparison-algebraic-decomposition} can be
	written as
	\begin{equation}
		\label{eq:appE-unified-comparison-R-expansion}
		\overline{\tau_j^a\phi_B^\varepsilon}\,
		\frac{d}{d\theta}\Bigl[
		U_{-a\theta B,j}^\varepsilon D_{\theta B,j}^\varepsilon
		\bigl(\mathcal{L}_{\theta B}^\varepsilon+1\bigr)^{-1}
		D_{\theta B}^{\varepsilon *}U_{b\theta B}^\varepsilon
		\Bigr]
		\mathcal{A}^{A,B}
		=
		\overline{\tau_j^a\phi_B^\varepsilon}\,
		\bigl(
		R_{1,j}^{\theta,a,b}
		+R_{2,j}^{\theta,a,b}
		+R_{3,j}^{\theta,a,b}
		\bigr)\mathcal{A}^{A,B},
	\end{equation}
	where
	\begin{equation}
		\label{eq:appE-unified-comparison-R1}
		R_{1,j}^{\theta,a,b}
		:=
		\frac{d}{d\theta}
		\bigl(
		U_{-a\theta B,j}^\varepsilon D_{\theta B,j}^\varepsilon
		\bigr)
		\bigl(\mathcal{L}_{\theta B}^\varepsilon+1\bigr)^{-1}
		D_{\theta B}^{\varepsilon *}U_{b\theta B}^\varepsilon,
	\end{equation}
	\begin{equation}
		\label{eq:appE-unified-comparison-R2}
		R_{2,j}^{\theta,a,b}
		:=
		U_{-a\theta B,j}^\varepsilon D_{\theta B,j}^\varepsilon
		\bigl(\mathcal{L}_{\theta B}^\varepsilon+1\bigr)^{-1}
		\frac{d}{d\theta}
		\bigl(
		D_{\theta B}^{\varepsilon *}U_{b\theta B}^\varepsilon
		\bigr),
	\end{equation}
	and
	\begin{equation}
		\label{eq:appE-unified-comparison-R3}
		R_{3,j}^{\theta,a,b}
		=-
		R_{2,j}^{\theta,a,0}
		D_{\theta B}^\varepsilon
		\bigl(\mathcal{L}_{\theta B}^\varepsilon+1\bigr)^{-1}
		D_{\theta B}^{\varepsilon *}U_{b\theta B}^\varepsilon
		-
		U_{-a\theta B,j}^\varepsilon D_{\theta B,j}^\varepsilon
		\bigl(\mathcal{L}_{\theta B}^\varepsilon+1\bigr)^{-1}
		D_{\theta B}^{\varepsilon *}
		\bigl(R_{1,\ell}^{\theta,0,b}\bigr)_{\ell=1}^2.
	\end{equation}

	We use the full-lattice covariant operators, with $B$ periodically
	extended. The following estimates hold uniformly in
	$\theta\in[0,1]$ and $a,b\in\{0,1\}$.

	\smallskip
	\noindent\textup{Step 1.} \textit{Estimate of the term containing
	$R_{1,j}^{\theta,a,b}$.}
	We prove that
	\begin{equation}
		\label{eq:appE-unified-comparison-R1-mapping}
		\csname ltx@label\endcsname{eq:appE-unified-comparison-R1-operator-norm}
		\left\lVert R_{1,j}^{\theta,a,b}\right\rVert_{
			\mathcal{C}^{-\alpha/2,\varepsilon}(\rho^{-4})
			\to H^{2\kappa,\varepsilon}(\rho^{-1/2})}
		\lesssim \bigl(1+\lVert B\rVert_{H^{1,\varepsilon}(\rho)}\bigr)^{2+2\alpha+4\kappa}.
	\end{equation}
	By~\eqref{eq:appE-unified-comparison-R1} and
	\eqref{eq:appE-shifted-covariant-derivative-derivative},
	\begin{equation}
		\label{eq:appE-unified-comparison-R1-recalled}
		R_{1,j}^{\theta,a,b}
		=
		\begin{cases}
			iB_jU_{\theta B,j}^\varepsilon\tau_j,& a=0,\\
			iB_jU_{-\theta B,j}^\varepsilon,& a=1
		\end{cases}
		(\mathcal{L}_{\theta B}^\varepsilon+1)^{-1}
		D_{\theta B}^{\varepsilon *}U_{b\theta B}^\varepsilon.
	\end{equation}
	Indeed, since $\rho^{3/2-2\alpha-4\kappa}\in L^{2,\varepsilon}$,
	we have the uniform embedding
	\begin{equation}
		\label{eq:appE-unified-comparison-C-to-HB-embedding}
		\mathcal{C}^{-\alpha/2,\varepsilon}(\rho^{-4})
		\hookrightarrow H^{-\alpha,\varepsilon}(\rho^{-5/2-2\alpha-4\kappa}).
	\end{equation}
	By~\eqref{eq:weighted-covariant-negative-divergence} with $s=\alpha$,
	\begin{equation}
		\label{eq:appE-unified-comparison-divergence-link-mapping}
		\csname ltx@label\endcsname{eq:appE-unified-comparison-divergence-link-norm}
		\left\lVert D_{\theta B}^{\varepsilon *}U_{b\theta B}^\varepsilon
		\right\rVert_{H^{-\alpha,\varepsilon}(\rho^{-5/2-2\alpha-4\kappa})
			\to H_{\theta B}^{-1-\alpha,\varepsilon}(\rho^{-5/2+\alpha-4\kappa})}
		\lesssim \bigl(1+\lVert B\rVert_{H^{1,\varepsilon}(\rho)}\bigr)^{3\alpha}.
	\end{equation}
	Moreover, by Lemma~\ref{lem:weighted-covariant-resolvent} with $s=1-\alpha$,
	\begin{equation}
		\label{eq:appE-unified-comparison-resolvent-mapping}
		\left\lVert(\mathcal{L}_{\theta B}^\varepsilon+1)^{-1}
		\right\rVert_{H_{\theta B}^{-1-\alpha,\varepsilon}(\rho^{-5/2+\alpha-4\kappa})
			\to H_{\theta B}^{1-\alpha,\varepsilon}(\rho^{-5/2+\alpha-4\kappa})}
		\lesssim 1.
	\end{equation}
	Finally, \eqref{eq:weighted-HB-H-comparison} and
	\eqref{eq:appE-gauge-link-B-multiplier}, with $c=-3/2-4\kappa$,
	together with the boundedness of lattice translations on the free
	weighted Sobolev spaces, give
	\begin{equation}
		\label{eq:appE-unified-comparison-outer-derivative-mapping}
		\csname ltx@label\endcsname{eq:appE-unified-comparison-outer-derivative-norm}
		\left\lVert\frac{d}{d\theta}
		\bigl(U_{-a\theta B,j}^\varepsilon D_{\theta B,j}^\varepsilon\bigr)
		\right\rVert_{H_{\theta B}^{1-\alpha,\varepsilon}(\rho^{-5/2+\alpha-4\kappa})
			\to H^{2\kappa,\varepsilon}(\rho^{-1/2})}
		\lesssim \bigl(1+\lVert B\rVert_{H^{1,\varepsilon}(\rho)}\bigr)^{2-\alpha+4\kappa}.
	\end{equation}
	Combining
	\eqref{eq:appE-unified-comparison-C-to-HB-embedding}--%
	\eqref{eq:appE-unified-comparison-outer-derivative-norm}
	yields
	\eqref{eq:appE-unified-comparison-R1-operator-norm}. Therefore, by
	\eqref{eq:appE-unified-comparison-source-bound} and the weighted product estimate,
	\begin{equation}
		\label{eq:appE-unified-comparison-R1-estimate}
		\begin{aligned}
		\left\lVert\overline{\tau_j^a\phi_B^\varepsilon}\,
		R_{1,j}^{\theta,a,b}\mathcal{A}^{A,B}
		\right\rVert_{H^{-2\kappa,\varepsilon}(\rho^4)}
		&\lesssim
		\left\lVert\tau_j^a\phi_B^\varepsilon
		\right\rVert_{\mathcal{C}^{-\kappa,\varepsilon}(\rho^\kappa)}
		\left\lVert R_{1,j}^{\theta,a,b}\mathcal{A}^{A,B}
		\right\rVert_{H^{2\kappa,\varepsilon}(\rho^{-1/2})}\\
		&\lesssim\lVert\Xi\rVert_\kappa \bigl(1+\lVert B\rVert_{H^{1,\varepsilon}(\rho)}\bigr)^{2+2\alpha+4\kappa}
		\left\lVert\mathcal{A}^{A,B}
		\right\rVert_{\mathcal{C}^{-\alpha/2,\varepsilon}(\rho^{-4})}\\
		&\lesssim\lVert\Xi\rVert_\kappa.
		\end{aligned}
	\end{equation}

	\smallskip
	\noindent\textup{Step 2.} \textit{Estimate of the term containing
	$R_{2,j}^{\theta,a,b}$.}
	We prove that
	\begin{equation}
		\label{eq:appE-unified-comparison-R2-mapping}
		\csname ltx@label\endcsname{eq:appE-unified-comparison-R2-operator-norm}
		\left\lVert R_{2,j}^{\theta,a,b}\right\rVert_{
			H_{\theta B}^{-\alpha,\varepsilon}(\rho^{-2})
			\to H^{\alpha,\varepsilon}(\rho^{-1/2})}
		\lesssim \bigl(1+\lVert B\rVert_{H^{1,\varepsilon}(\rho)}\bigr)^{1+10\alpha}.
	\end{equation}
	Recall from~\eqref{eq:appE-unified-comparison-R2} that
	\begin{equation}
		\label{eq:appE-unified-comparison-R2-recalled}
		R_{2,j}^{\theta,a,b}
		=
		U_{-a\theta B,j}^\varepsilon D_{\theta B,j}^\varepsilon
		\bigl(\mathcal{L}_{\theta B}^\varepsilon+1\bigr)^{-1}
		\frac{d}{d\theta}
		\bigl(D_{\theta B}^{\varepsilon *}
		U_{b\theta B}^\varepsilon\bigr).
	\end{equation}
	Componentwise, direct differentiation gives
	\begin{equation}
		\label{eq:appE-unified-comparison-inner-derivative}
		\frac{d}{d\theta}
		\bigl(D_{\theta B,\ell}^{\varepsilon *}
		U_{b\theta B,\ell}^\varepsilon\bigr)
		=
		\begin{cases}
			-i\tau_\ell^*B_\ell U_{-\theta B,\ell}^\varepsilon,
			& b=0,\\
			-iB_\ell U_{\theta B,\ell}^\varepsilon,
			& b=1.
		\end{cases}
	\end{equation}
	By~\eqref{eq:appE-gauge-BU-holder-multiplier} with $c=-2+\alpha$ and
	\eqref{eq:weighted-HB-H-comparison}, applied before and after multiplication,
	the operator in
	\eqref{eq:appE-unified-comparison-inner-derivative} satisfies
	\begin{equation}
		\label{eq:appE-unified-comparison-inner-derivative-mapping}
		\csname ltx@label\endcsname{eq:appE-unified-comparison-inner-derivative-norm}
		\left\lVert\frac{d}{d\theta}
		\bigl(D_{\theta B}^{\varepsilon *}U_{b\theta B}^\varepsilon\bigr)
		\right\rVert_{H_{\theta B}^{-\alpha,\varepsilon}(\rho^{-2})
			\to H_{\theta B}^{-4\alpha,\varepsilon}(\rho^{-1+7\alpha})}
		\lesssim \bigl(1+\lVert B\rVert_{H^{1,\varepsilon}(\rho)}\bigr)^{1+7\alpha}.
	\end{equation}
	Moreover, \eqref{eq:weighted-covariant-resolvent} with $s=2-4\alpha$ gives
	\begin{equation}
		\label{eq:appE-unified-comparison-R2-resolvent-mapping}
		\left\lVert(\mathcal{L}_{\theta B}^\varepsilon+1)^{-1}
		\right\rVert_{H_{\theta B}^{-4\alpha,\varepsilon}(\rho^{-1+7\alpha})
			\to H_{\theta B}^{1+\alpha,\varepsilon}(\rho^{-1+7\alpha})}
		\lesssim 1.
	\end{equation}
	Here we used the embedding from regularity $2-4\alpha$ to $1+\alpha$,
	since $\alpha<1/5$. Finally,
	\eqref{eq:weighted-covariant-derivative-free-output} with $s=\alpha$ and
	\eqref{eq:covariant-derivative-translation-identities} give
	\begin{equation}
		\label{eq:appE-unified-comparison-outer-operator-mapping}
		\csname ltx@label\endcsname{eq:appE-unified-comparison-outer-operator-norm}
		\left\lVert U_{-a\theta B,j}^\varepsilon D_{\theta B,j}^\varepsilon
		\right\rVert_{H_{\theta B}^{1+\alpha,\varepsilon}(\rho^{-1+7\alpha})
			\to H^{\alpha,\varepsilon}(\rho^{-1+10\alpha})}
		\lesssim \bigl(1+\lVert B\rVert_{H^{1,\varepsilon}(\rho)}\bigr)^{3\alpha}.
	\end{equation}
	Since $-1+10\alpha\le-1/2$, this yields
	\eqref{eq:appE-unified-comparison-R2-operator-norm}. By
	\eqref{eq:appE-unified-comparison-C-to-HB-embedding} and
	\eqref{eq:weighted-HB-H-comparison}, we also have
	\begin{equation}
		\label{eq:appE-unified-comparison-source-covariant-embedding}
		\left\lVert\mathrm{Id}\right\rVert_{
			\mathcal C^{-\alpha/2,\varepsilon}(\rho^{-4})
			\to H_{\theta B}^{-\alpha,\varepsilon}(\rho^{-2})}
		\lesssim\bigl(1+\lVert B\rVert_{H^{1,\varepsilon}(\rho)}\bigr)^\alpha.
	\end{equation}
	Consequently, by~\eqref{eq:appE-unified-comparison-source-bound}
	and the weighted product estimate,
	\begin{equation}
		\label{eq:appE-unified-comparison-R2-estimate}
		\begin{aligned}
		\left\lVert\overline{\tau_j^a\phi_B^\varepsilon}\,
		R_{2,j}^{\theta,a,b}\mathcal{A}^{A,B}
		\right\rVert_{H^{-2\kappa,\varepsilon}(\rho^4)}
		&\lesssim
		\left\lVert\tau_j^a\phi_B^\varepsilon
		\right\rVert_{\mathcal{C}^{-\kappa,\varepsilon}(\rho^\kappa)}
		\left\lVert R_{2,j}^{\theta,a,b}\mathcal{A}^{A,B}
		\right\rVert_{H^{\alpha,\varepsilon}(\rho^{-1/2})}\\
		&\lesssim\lVert\Xi\rVert_\kappa \bigl(1+\lVert B\rVert_{H^{1,\varepsilon}(\rho)}\bigr)^{1+11\alpha}
		\left\lVert\mathcal{A}^{A,B}
		\right\rVert_{\mathcal{C}^{-\alpha/2,\varepsilon}(\rho^{-4})}\\
		&\lesssim\lVert\Xi\rVert_\kappa.
		\end{aligned}
	\end{equation}

	\smallskip
	\noindent\textup{Step 3.} \textit{Estimate of the term containing
	$R_{3,j}^{\theta,a,b}$.}
	Recall from~\eqref{eq:appE-unified-comparison-R3} that
	\begin{equation}
		\label{eq:appE-unified-comparison-R3-recalled}
		\begin{aligned}
		R_{3,j}^{\theta,a,b}
		={}-
		R_{2,j}^{\theta,a,0}
		D_{\theta B}^\varepsilon
		\bigl(\mathcal{L}_{\theta B}^\varepsilon+1\bigr)^{-1}
		D_{\theta B}^{\varepsilon *}U_{b\theta B}^\varepsilon
		-
		U_{-a\theta B,j}^\varepsilon D_{\theta B,j}^\varepsilon
		\bigl(\mathcal{L}_{\theta B}^\varepsilon+1\bigr)^{-1}
		D_{\theta B}^{\varepsilon *}
		\bigl(R_{1,\ell}^{\theta,0,b}\bigr)_{\ell=1}^2.
		\end{aligned}
	\end{equation}
	For the first term, the estimates
	\eqref{eq:appE-unified-comparison-C-to-HB-embedding},
	\eqref{eq:appE-unified-comparison-divergence-link-mapping}, and
	\eqref{eq:appE-unified-comparison-resolvent-mapping} imply
	\begin{equation}
		\label{eq:appE-unified-comparison-projected-resolvent-neg-alpha}
		\csname ltx@label\endcsname{eq:appE-unified-comparison-projected-resolvent-neg-alpha-norm}
		\begin{aligned}
		\left\lVert D_{\theta B}^\varepsilon
		(\mathcal{L}_{\theta B}^\varepsilon+1)^{-1}
		D_{\theta B}^{\varepsilon *}U_{b\theta B}^\varepsilon
		\right\rVert_{\mathcal{C}^{-\alpha/2,\varepsilon}(\rho^{-4})
			\to H_{\theta B}^{-\alpha,\varepsilon}(\rho^{-2})}
		\lesssim \bigl(1+\lVert B\rVert_{H^{1,\varepsilon}(\rho)}\bigr)^{3\alpha}.
		\end{aligned}
	\end{equation}
	Indeed, after the embedding
	\eqref{eq:appE-unified-comparison-C-to-HB-embedding}, we use
	\begin{equation}
		\label{eq:appE-unified-comparison-projected-resolvent-chain}
		\begin{aligned}
		D_{\theta B}^{\varepsilon *}U_{b\theta B}^\varepsilon
		&\colon H^{-\alpha,\varepsilon}(\rho^{-5/2-2\alpha-4\kappa})
		\longrightarrow H_{\theta B}^{-1-\alpha,\varepsilon}(\rho^{-5/2+\alpha-4\kappa}),\\
		(\mathcal{L}_{\theta B}^\varepsilon+1)^{-1}
		&\colon H_{\theta B}^{-1-\alpha,\varepsilon}(\rho^{-5/2+\alpha-4\kappa})
		\longrightarrow H_{\theta B}^{1-\alpha,\varepsilon}(\rho^{-5/2+\alpha-4\kappa}),\\
		D_{\theta B}^\varepsilon
		&\colon H_{\theta B}^{1-\alpha,\varepsilon}(\rho^{-5/2+\alpha-4\kappa})
		\longrightarrow H_{\theta B}^{-\alpha,\varepsilon}(\rho^{-5/2+\alpha-4\kappa}).
		\end{aligned}
	\end{equation}
	The first operator has norm bounded by the right-hand side of
	\eqref{eq:appE-unified-comparison-divergence-link-norm}, whereas the
	last two bounds are uniform by
	\eqref{eq:weighted-covariant-resolvent} and
	\eqref{eq:covariant-derivative-estimates}.
	Since $-5/2+\alpha-4\kappa\le-2$, this yields
	\eqref{eq:appE-unified-comparison-projected-resolvent-neg-alpha}.

	For the second term in
	\eqref{eq:appE-unified-comparison-R3-recalled}, we have
	\begin{equation}
		\label{eq:appE-unified-comparison-projected-resolvent-2kappa}
		\csname ltx@label\endcsname{eq:appE-unified-comparison-projected-resolvent-2kappa-norm}
		\begin{aligned}
		&\left\lVert U_{-a\theta B,j}^\varepsilon D_{\theta B,j}^\varepsilon
		(\mathcal{L}_{\theta B}^\varepsilon+1)^{-1}D_{\theta B}^{\varepsilon *}
		\right\rVert_{H^{2\kappa,\varepsilon}(\rho^{-1/2})
			\to H^{2\kappa,\varepsilon}(\rho^{-1/2+8\kappa})}
		\lesssim \bigl(1+\lVert B\rVert_{H^{1,\varepsilon}(\rho)}\bigr)^{8\kappa}.
		\end{aligned}
	\end{equation}
	Indeed,
	\begin{equation}
		\label{eq:appE-unified-comparison-projected-resolvent-2kappa-chain}
		\begin{aligned}
		D_{\theta B}^{\varepsilon *}
		&\colon H^{2\kappa,\varepsilon}(\rho^{-1/2})
		\longrightarrow H_{\theta B}^{2\kappa-1,\varepsilon}(\rho^{-1/2+2\kappa}),\\
		(\mathcal{L}_{\theta B}^\varepsilon+1)^{-1}
		&\colon H_{\theta B}^{2\kappa-1,\varepsilon}(\rho^{-1/2+2\kappa})
		\longrightarrow H_{\theta B}^{1+2\kappa,\varepsilon}(\rho^{-1/2+2\kappa}),\\
		U_{-a\theta B,j}^\varepsilon D_{\theta B,j}^\varepsilon
		&\colon H_{\theta B}^{1+2\kappa,\varepsilon}(\rho^{-1/2+2\kappa})
		\longrightarrow H^{2\kappa,\varepsilon}(\rho^{-1/2+8\kappa}).
		\end{aligned}
	\end{equation}
	The first operator has norm
	$\lesssim\bigl(1+\lVert B\rVert_{H^{1,\varepsilon}(\rho)}\bigr)^{2\kappa}$
	by~\eqref{eq:weighted-covariant-negative-divergence} with $s=-2\kappa$.
	The resolvent bound is uniform by~\eqref{eq:weighted-covariant-resolvent}
	with $s=1+2\kappa$, whereas the third operator has norm
 $\lesssim\bigl(1+\lVert B\rVert_{H^{1,\varepsilon}(\rho)}\bigr)^{6\kappa}$
	by~\eqref{eq:weighted-covariant-derivative-free-output} with $s=2\kappa$
	and~\eqref{eq:covariant-derivative-translation-identities}.
	Since $\alpha>2\kappa$ and $-1/2+8\kappa\le4-\kappa$, combining
	\eqref{eq:appE-unified-comparison-R1-mapping},
	\eqref{eq:appE-unified-comparison-R2-mapping},
	\eqref{eq:appE-unified-comparison-projected-resolvent-neg-alpha},
	and~\eqref{eq:appE-unified-comparison-projected-resolvent-2kappa}
	with the weighted product estimate, we obtain
	\begin{equation}
		\label{eq:appE-unified-comparison-R3-estimate}
		\begin{aligned}
		\left\lVert\overline{\tau_j^a\phi_B^\varepsilon}\,
		R_{3,j}^{\theta,a,b}\mathcal{A}^{A,B}
		\right\rVert_{H^{-2\kappa,\varepsilon}(\rho^4)}
		&\lesssim
		\left\lVert\tau_j^a\phi_B^\varepsilon
		\right\rVert_{\mathcal{C}^{-\kappa,\varepsilon}(\rho^\kappa)}
		\Bigl[\bigl(1+\lVert B\rVert_{H^{1,\varepsilon}(\rho)}\bigr)^{1+13\alpha}\\[-2pt]
		&\hspace{72pt}+\bigl(1+\lVert B\rVert_{H^{1,\varepsilon}(\rho)}\bigr)^{2+2\alpha+12\kappa}\Bigr]
		\left\lVert\mathcal{A}^{A,B}
		\right\rVert_{\mathcal{C}^{-\alpha/2,\varepsilon}(\rho^{-4})}\\
		&\lesssim\lVert\Xi\rVert_\kappa.
		\end{aligned}
	\end{equation}

	\smallskip
	\noindent\textup{Step 4.} \textit{Estimate of $S_1$.}
	Recall from~\eqref{eq:appE-unified-comparison-S-terms} that
	\begin{equation}
		\label{eq:appE-unified-comparison-S1-recalled}
		S_1
		=
		\overline{\tau_j^a
		(\phi_B^\varepsilon-\phi^\varepsilon)}\,
		D_j^\varepsilon
		\bigl(\mathcal{L}_0^\varepsilon+1\bigr)^{-1}
		D^{\varepsilon *}\mathcal{A}^{A,B}.
	\end{equation}
	The free multiplier
	$D_j^\varepsilon(\mathcal{L}_0^\varepsilon+1)^{-1}
	D^{\varepsilon *}$ is bounded on
	$\mathcal{C}^{-\alpha/2,\varepsilon}(\rho^{-4})$. By the weighted product estimate
	and~\eqref{eq:appE-unified-comparison-source-bound}, we obtain
	\begin{equation}
		\label{eq:appE-unified-comparison-S1-estimate}
		\begin{aligned}
		\lVert S_1\rVert_{\mathcal{C}^{-\alpha/2,\varepsilon}(\rho^{-1})}
		&\lesssim
		\left\lVert
		\phi_B^\varepsilon-\phi^\varepsilon
		\right\rVert_{\mathcal{C}^{\alpha,\varepsilon}(\rho^{3})}
		\left\lVert
		D_j^\varepsilon
		\bigl(\mathcal{L}_0^\varepsilon+1\bigr)^{-1}
		D^{\varepsilon *}\mathcal{A}^{A,B}
		\right\rVert_{\mathcal{C}^{-\alpha/2,\varepsilon}(\rho^{-4})}\\
		&\lesssim
		\frac{
			\left\lVert
			\phi_B^\varepsilon-\phi^\varepsilon
			\right\rVert_{\mathcal{C}^{\alpha,\varepsilon}(\rho^{3})}
		}{
			\bigl(1+\lVert B\rVert_{H^{1,\varepsilon}(\rho)}\bigr)^{l_0}
		}\\
		&\lesssim
		\lVert\Xi\rVert_\kappa,
		\end{aligned}
	\end{equation}
	where the last inequality follows from
	\eqref{eq:stochastic-object-Xi-norm}.

	\smallskip
	\noindent\textup{Step 5.} \textit{Estimate of $S_2$.}
	By~\eqref{eq:appE-unified-comparison-S-terms},
	\begin{equation}
		\label{eq:appE-unified-comparison-S2-recalled}
		S_2
		=
		\overline{\tau_j^a\phi^\varepsilon}\,
		D_j^\varepsilon
		\bigl(\mathcal{L}_0^\varepsilon+1\bigr)^{-1}
		D^{\varepsilon *}
		\bigl(\mathcal{A}^{A,B}-\mathcal{A}^{A,0}\bigr).
	\end{equation}
	For each $\ell=1,2$, the definition
	\eqref{eq:appE-unified-comparison-sources} gives
	\begin{equation}
		\label{eq:appE-unified-comparison-source-difference}
		\mathcal{A}_\ell^{A,B}-\mathcal{A}_\ell^{A,0}
		=
		\int_0^1\mathbb{E}_t^\phi\left[
		\mathcal{U}_{>L_1}^\varepsilon
		\left(
		r_{A,\ell}^{\varepsilon,(b)}\tau_\ell^b
		(\phi_B^\varepsilon-\phi^\varepsilon)
		\right)
		\right]dt.
	\end{equation}
	For every $f\colon\Lambda_\varepsilon\to\mathbb{C}$, the discrete
	product rule and the paraproduct estimates give
	\begin{equation}
		\label{eq:appE-unified-comparison-divergence-product-bound}
		\begin{aligned}
		&\left\lVert
		D^{\varepsilon *}
		\bigl(
		(r_{A,\ell}^{\varepsilon,(b)}\tau_\ell^b f)_{\ell=1}^2
		\bigr)
		\right\rVert_{\mathcal{C}^{-2\kappa,\varepsilon}(\rho^{3+\kappa})}\\
		&\quad\lesssim
		\begin{cases}
		\left\lVert
		(D^{\varepsilon *}\overline{r_A^\varepsilon})\succ f
		\right\rVert_{\mathcal{C}^{-2\kappa,\varepsilon}(\rho^{3+\kappa})}
		+
		\displaystyle\sum_{\ell=1}^2
		\left\lVert
		\overline{r_{A,\ell}^\varepsilon}
		(\mathord\cdot-\varepsilon e_\ell)
		\succ D_\ell^{\varepsilon *}f
		\right\rVert_{\mathcal{C}^{-2\kappa,\varepsilon}(\rho^{3+\kappa})}
		+
		\left\lVert
		\overline{r_A^\varepsilon}\preccurlyeq f
		\right\rVert_{\mathcal{C}^{1-2\kappa,\varepsilon}(\rho^{3+\kappa})},
		& b=0,\\[6pt]
		\left\lVert
		(D^{\varepsilon *}r_A^\varepsilon)\succ f
		\right\rVert_{\mathcal{C}^{-2\kappa,\varepsilon}(\rho^{3+\kappa})}
		+
		\displaystyle\sum_{\ell=1}^2
		\left\lVert
		r_{A,\ell}^\varepsilon\succ D_\ell^\varepsilon f
		\right\rVert_{\mathcal{C}^{-2\kappa,\varepsilon}(\rho^{3+\kappa})}
		+
		\left\lVert
		(r_{A,\ell}^\varepsilon\preccurlyeq
		\tau_\ell f)_{\ell=1}^2
		\right\rVert_{\mathcal{C}^{1-2\kappa,\varepsilon}(\rho^{3+\kappa})},
		& b=1,
		\end{cases}\\
		&\quad\lesssim
		\left(
		\left\lVert
		D^{\varepsilon *}\overline{r_A^\varepsilon}
		\right\rVert_{\mathcal{C}^{-\kappa,\varepsilon}(\rho^{\kappa})}
		+
		\left\lVert r_A^\varepsilon
		\right\rVert_{\mathcal{C}^{-\kappa,\varepsilon}(\rho^{\kappa})}
		\right)
		\lVert f\rVert_{\mathcal{C}^{1-\kappa,\varepsilon}(\rho^{3})}\\
		&\quad\lesssim
		\lVert\Xi\rVert_\kappa
		\lVert f\rVert_{\mathcal{C}^{1-\kappa,\varepsilon}(\rho^{3})}.
		\end{aligned}
	\end{equation}
	The commutator is estimated by~
	\eqref{eq:high-frequency-truncation-commutator-estimate}, with
	$s=-2\kappa$, $\delta=\alpha-2\kappa$, and input weight $\rho^{3+\kappa}$.
	Since $(\alpha-2\kappa)/(\alpha/2-\kappa)=2$ and
	$3+\kappa-\sigma_1(\alpha-2\kappa)=-5-\kappa\le-4$,
	\eqref{eq:high-frequency-truncation-estimate},
	the cutoff choice in Table~\ref{tab:higgs-high-frequency-cutoff-choices},
	\eqref{eq:appE-unified-comparison-divergence-product-bound},
	and the weighted product estimate give
	\begin{equation}
		\label{eq:appE-unified-comparison-high-difference-bound}
		\begin{aligned}
		\left\lVert
		D^{\varepsilon *}\mathcal{U}_{>L_1}^\varepsilon
		\left(
		\bigl(
		r_{A,\ell}^{\varepsilon,(b)}\tau_\ell^b
		(\phi_B^\varepsilon-\phi^\varepsilon)
		\bigr)_{\ell=1}^2
		\right)
		\right\rVert_{\mathcal{C}^{-\alpha,\varepsilon}(\rho^{-4})}
		&\lesssim
		2^{-L_1(\alpha-2\kappa)}
		\left\lVert
		D^{\varepsilon *}
		\left(
		\bigl(
		r_{A,\ell}^{\varepsilon,(b)}\tau_\ell^b
		(\phi_B^\varepsilon-\phi^\varepsilon)
		\bigr)_{\ell=1}^2
		\right)
		\right\rVert_{\mathcal{C}^{-2\kappa,\varepsilon}(\rho^{3+\kappa})}\\
		&\quad+
		2^{-L_1(\alpha-2\kappa)}
		\left\lVert
		\bigl(
		r_{A,\ell}^{\varepsilon,(b)}\tau_\ell^b
		(\phi_B^\varepsilon-\phi^\varepsilon)
		\bigr)_{\ell=1}^2
		\right\rVert_{\mathcal{C}^{-2\kappa,\varepsilon}(\rho^{3+\kappa})}\\
		&\lesssim
		\frac{
			\lVert\Xi\rVert_\kappa
			\left\lVert
			\phi_B^\varepsilon-\phi^\varepsilon
			\right\rVert_{\mathcal{C}^{1-\kappa,\varepsilon}(\rho^{3})}
		}{
			\bigl(1+\lVert B\rVert_{H^{1,\varepsilon}(\rho)}\bigr)^{2l_0}
			\lVert\Xi\rVert_\kappa^2
		}
		\lesssim1.
		\end{aligned}
	\end{equation}
	The contraction property of conditional expectation and integration
	in time now yield
	\begin{equation}
		\label{eq:appE-unified-comparison-source-difference-bound}
		\left\lVert
		D^{\varepsilon *}
		\bigl(\mathcal{A}^{A,B}-\mathcal{A}^{A,0}\bigr)
		\right\rVert_{L^\infty(\Omega^\phi;
			\mathcal{C}^{-\alpha,\varepsilon}(\rho^{-4}))}
		\lesssim1.
	\end{equation}
	Finally,
	$D_j^\varepsilon(\mathcal{L}_0^\varepsilon+1)^{-1}$ maps
	$\mathcal{C}^{-\alpha,\varepsilon}(\rho^{-4})$ into
	$\mathcal{C}^{1-\alpha,\varepsilon}(\rho^{-4})$. Thus
	\begin{equation}
		\label{eq:appE-unified-comparison-S2-estimate}
		\begin{aligned}
		\lVert S_2\rVert_{\mathcal{C}^{-\kappa,\varepsilon}(\rho^{\kappa-4})}
		&\lesssim
		\lVert\phi^\varepsilon\rVert_{\mathcal{C}^{-\kappa,\varepsilon}(\rho^{\kappa})}
		\left\lVert
		D_j^\varepsilon
		\bigl(\mathcal{L}_0^\varepsilon+1\bigr)^{-1}
		D^{\varepsilon *}
		\bigl(\mathcal{A}^{A,B}-\mathcal{A}^{A,0}\bigr)
		\right\rVert_{\mathcal{C}^{1-\alpha,\varepsilon}(\rho^{-4})}\\
		&\lesssim
		\lVert\phi^\varepsilon\rVert_{\mathcal{C}^{-\kappa,\varepsilon}(\rho^{\kappa})}
		\left\lVert
		D^{\varepsilon *}
		\bigl(\mathcal{A}^{A,B}-\mathcal{A}^{A,0}\bigr)
		\right\rVert_{\mathcal{C}^{-\alpha,\varepsilon}(\rho^{-4})}.
		\end{aligned}
	\end{equation}
	Since $\rho^{8-\kappa}\in L^{2,\varepsilon}$, we have
	$\mathcal C^{-\kappa,\varepsilon}(\rho^{\kappa-4})
	\hookrightarrow H^{-\alpha,\varepsilon}(\rho^4)$.
	Together with~\eqref{eq:free-gff-wick-power-negative-holder-moment}
	with $(k_1,k_2)=(1,0)$ and
	\eqref{eq:appE-unified-comparison-source-difference-bound}, this gives
	\begin{equation}
		\label{eq:appE-unified-comparison-S2-moment}
		\mathbb{E}^\phi
		\lVert S_2\rVert_{H^{-\alpha,\varepsilon}(\rho^4)}
		\lesssim
		1
		\lesssim
		\mathbb{E}^\phi\left[
		\lVert\Xi\rVert_\kappa^{K_0}
		\right].
	\end{equation}

	Combining Steps 1--5, using
	$H^{-2\kappa,\varepsilon}(\rho^4)\hookrightarrow H^{-\alpha,\varepsilon}(\rho^4)$,
	$\mathcal{C}^{-\alpha/2,\varepsilon}(\rho^{-1})\hookrightarrow
	H^{-\alpha,\varepsilon}(\rho^4)$, and
	$\mathcal{C}^{-\kappa,\varepsilon}(\rho^{\kappa-4})\hookrightarrow
	H^{-\alpha,\varepsilon}(\rho^4)$, we obtain
	\begin{equation}
		\label{eq:appE-unified-comparison-conclusion}
		\begin{aligned}
		&\Bigg\lVert
		\mathbb{E}^\phi\Big[
		\overline{\tau_j^a\phi_B^\varepsilon}\,
		U_{-aB,j}^\varepsilon D_{B,j}^\varepsilon
		\bigl(\mathcal{L}_B^\varepsilon+1\bigr)^{-1}
		D_B^{\varepsilon *}U_{bB}^\varepsilon\mathcal{A}^{A,B}
		-
		\overline{\tau_j^a\phi^\varepsilon}\,
		D_j^\varepsilon
		\bigl(\mathcal{L}_0^\varepsilon+1\bigr)^{-1}
		D^{\varepsilon *}\mathcal{A}^{A,0}
		\Big]
		\Bigg\rVert_{H^{-\alpha,\varepsilon}(\rho^4)}\\
		&\quad\lesssim
		\int_0^1
		\mathbb{E}^\phi\Bigg[
		\left\lVert
		\overline{\tau_j^a\phi_B^\varepsilon}\,
		\bigl(
		R_{1,j}^{\theta,a,b}
		+R_{2,j}^{\theta,a,b}
		+R_{3,j}^{\theta,a,b}
		\bigr)\mathcal{A}^{A,B}
		\right\rVert_{H^{-\alpha,\varepsilon}(\rho^4)}
		\Bigg]d\theta
		+
		\mathbb{E}^\phi\left[
		\lVert S_1\rVert_{H^{-\alpha,\varepsilon}(\rho^4)}
		+
		\lVert S_2\rVert_{H^{-\alpha,\varepsilon}(\rho^4)}
		\right]\\
		&\quad\lesssim
		\mathbb{E}^\phi\left[
		\lVert\Xi\rVert_\kappa^{K_0}
		\right].
		\end{aligned}
	\end{equation}
	By~\eqref{eq:appE-unified-comparison-sources},
	\eqref{eq:appE-unified-comparison-conclusion} is precisely
	\eqref{eq:gauge-irregular-Zhi1-unified-free-comparison}.

	\bibliographystyle{amsalphaabbr}
	\bibliography{references}
\end{document}